\documentstyle[12pt,epic,eepic,epsf,here]{article}
\makeatletter
\if@twoside \def\ps@headings{\let\@mkboth\markboth
\def\@oddfoot{}\def\@evenfoot{}\def\@evenhead{\rm \thepage\hfil
\sl\leftmark}\def\@oddhead{\hbox{}{\sl \rightmark} \hfil
\rm\thepage}\def\sectionmark##1{\markboth {{\ifnum \c@secnumdepth
>\z@
 \thesection\hskip 1em\relax \fi ##1}}{}}\def\subsectionmark##1{\markright
{\ifnum \c@secnumdepth >\@ne
 \thesubsection\hskip 1em\relax \fi ##1}}}
\else \def\ps@headings{\let\@mkboth\markboth
\def\@oddfoot{}\def\@evenfoot{}\def\@oddhead{\hbox {}{\small \rightmark} \hfil
\rm\thepage}\def\sectionmark##1{\markright {{\ifnum \c@secnumdepth
>\z@
 \thesection\hskip 1em\relax \fi ##1}}}}
\fi
\def\ps@myheadings{\let\@mkboth\@gobbletwo
\def\@oddhead{\hbox{}{\small\rightmark} \hfil
\rm\thepage}\def\@oddfoot{}\def\@evenhead{\rm \thepage\hfil{\small\leftmark}\hbo\
x
{}}\def\@evenfoot{}\def\sectionmark##1{}\def\subsectionmark##1{}}
\makeatother

\newtheorem{lm}{Lemma}
\newtheorem{thm}{Theorem}
\newtheorem{claim}{Claim}
\newtheorem{cor}{Corollary}
\newtheorem{defi}{Definition}
\newtheorem{prop}{Proposition}

\begin{document}

\def\gb{{\beta^2}}
\def\gbf{{\beta^4}}
\def\gbs{{\beta}}
\def\os{{\rm OpString}}
\def\cho{\choose}
\def\ovl{\overline}
\def\sk{{*_{\cal K}}}
\def\sl{{*_{\cal L}}}
\def\sg{{*_{\cal G}}}
\def\lan{\langle}
\def\ran{\rangle}
\def\oX{{\ovl X}}
\def\begeq{\begin{eqnarray*}}
\def\endeq{\end{eqnarray*}}
\def\pha{\phantom}
\def\la{\lambda}
\def\lad{\lambda^2}
\def\R{R}
\def\Y{Y}
\def\Z{Z}
\def\sm{\small}
\def\s#1{$#1_{\sigma}$}
\def\ss#1{$#1_{\sigma^2}$}
\def\si#1{#1_{\sigma}}
\def\ssi#1{#1_{\sigma^2}}
\def\a{\alpha}
\def\ab{\tilde{\alpha}}
\def\dline(#1,#2){\multiput(#1,#2)(0,-3){8}{\line(0,-1){2}}}
\def\hl{\hline}
\def\mcol{\multicolumn}
\def\td#1{\tilde #1}
\def\bl{$\bullet$}
\def\qed{\phantom{XXXXXXXXXXXXXXXXXXXXXXXXXXXXXXXXXXXXXXXXXXXXX}q.e.d.}

\title{Exotic subfactors of finite depth with Jones indices 
$(5+\sqrt{13})/2$ and $(5+\sqrt{17})/2$}
 \author{\sm{M. Asaeda,  
         Graduate School of Mathematical Sciences, 
         University of Tokyo} \footnote{Komaba, Meguro-ku, Tokyo, 153-8924, 
         JAPAN, e-mail: asaeda@ms.u-tokyo.ac.jp}\\ 
         \and \sm{U. Haagerup, 
  Institut for Matematik og Datalogi, 
 Odense Universitet} \footnote{Campusvej 55, DK-5230 Odense M, DENMARK, e-mail: 
 haagerup@imada.ou.dk}}
\date{}
\maketitle
\date{}
\begin{abstract}
We prove existence of subfactors of finite depth of the hyperfinite 
 II$_1$ factor with indices ${(5+\sqrt{13})}/{2}=
4.302 \cdots$ and ${(5+\sqrt{17})}/{2}=4.561\cdots$. The existence of 
the former was announced by the second named author in 1993 and that 
of the latter has been conjectured since then. These are the only 
known subfactors with finite depth which do not arise from classical
 groups,  quantum 
groups or rational conformal field theory.
\end{abstract}

\section{Introduction}

In the theory of operator algebras, subfactor theory has been developing 
dynamically, involving various fields in mathematics and mathematical 
physics since its foundation by V. F. R. Jones in 1983 \cite{J}.  Above
all, the classification of subfactors is one of the most important
topics in the theory.  In the celebrated Jones index theory \cite{J},
Jones introduced the Jones index for subfactors of type II$_{1}$
 as an invariant.
Later, he also introduced a principal graph and a dual principal
graph as finer invariants of subfactors.
Since Jones proved in the middle of 1980's that subfactors with
index less than 4 have one of the Dynkin diagrams as their (dual)
principal graphs,
the classification of the
hyperfinite II$_1$ subfactors,
has been studied by A. Ocneanu and S. Popa, and also by M. Izumi,
Y. Kawahigashi, and a number of other mathematicians. \par
In this process, Ocneanu's {\em paragroup theory}
\cite{Oc1} has been quite effective.
He penetrated the algebraic, or rather combinatorial, nature of subfactors
and constructed a paragroup from a subfactor of type II$_1$.
A paragroup is a set of data consisting of four graphs made of (dual) principal
graph and assignment of complex numbers to ``cells'' arising from four graphs,
called a {\it connection}.  Thanks to the ``generating property'' for 
subfactors of finite depth proved by Popa in \cite{P1}, it has
turned out that the correspondence between paragroups and subfactors
of the hyperfinite II$_1$ factor with finite index and finite depth
is bijective, therefore the
classification of hyperfinite II$_1$ subfactors with finite index and
finite depth is reduced to that of
paragroups.  By checking the flatness condition for the connections on
the Dynkin diagrams,  Ocneanu has announced in \cite{Oc1}
that subfactors with index less
than 4 are completely classified by the Dynkin diagrams $A_n$,
$D_{2n}$, $E_6$, and $E_8$. (See also \cite{BN}, \cite{I1}, \cite{I2},
\cite{K}, \cite{SV}.)  After that, Popa (\cite{P}) extended the correspondence between 
paragroups and subfactors of the hyperfinite II$_{1}$ factor to the 
strongly amenable case, and gave a
classification of subfactors with indices equal to 4.  (In all the above 
mentioned cases, the dual principal graph of a subfactor is the same as the principal
graph.  See also \cite{IK}.) We refer readers to \cite{EK},
\cite{GHJ} for algebraic aspects of a general theory of subfactors.  \par

The second named author then tried to find subfactors with index {\it a
little bit beyond} 4.  Some subfactors with index larger than 4 had
been already constructed from other mathematical objects.  For example, we can
construct a subfactor from an arbitrary finite group by a crossed product
with an outer action,
and this subfactor has an index equal to the order of the original
finite group.  Trivially, the index is at least 5 if it is larger than 4. 
We also have subfactors constructed from quantum groups
$U_q(sl(n)),\; q=e^{2 \pi i/k}$ with index $\frac{{\rm sin}^2 (n \pi 
/k)}{{\rm sin}^2(\pi/k)}$ as in \cite{W} and these index values do not
fall in the interval $(4,5)$.  Unfortunately or
naturally, these subfactors do not contain more information about
the algebraic structure than the original mathematical objects themselves 
such as groups or quantum groups.
{\it ``Does there exist any subfactor not arising from (quantum) 
groups ?''}
If it is the case, we have a subfactor as a really new object producing
new mathematical structures.
 We expect that subfactors with index
{\it slightly} larger than 4 would be indeed those with {\it
exotic} nature and they do not to arise from other mathematical objects.
The second named author gave in 1993 a list of possible candidates of graphs
 which might be realized as
(dual) principal graphs of subfactors with index in 
$(4, 3+\sqrt{3})=(4, 4.732 \cdots)$ in \cite{H}.  We see four candidates of
pairs of graphs, including two pairs with parameters,
 in \S7 of \cite{H}.
At the same time, the second named author announced a proof of 
existence of 
the subfactor with index $(5+\sqrt{13})/2$ for
 the case $n=3$ of (2) in the list in \cite{H}, but the proof has not been 
 published until now.
 Ever since, nothing
had been known for the other cases for some years, until D. Bisch
recently proved that a subfactor with (dual) principal graph (4) in \S7 of \cite{H} does {\it not} exist \cite{B}.  \par
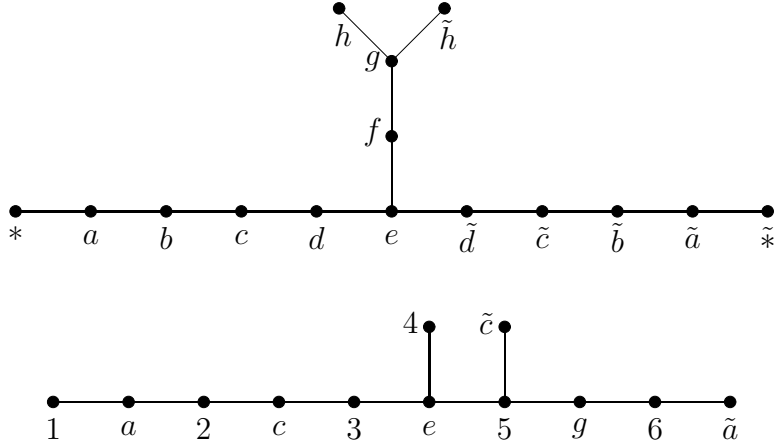
\begin{figure}[h]
\begin{center}
\thinlines
\unitlength 0.5mm
 \begin{picture}(220,80)
\multiput(10,20)(20,0){11}{\circle*{3}}
\multiput(110,40)(0,20){2}{\circle*{3}}
\multiput(96,74)(28,0){2}{\circle*{3}}
\multiput(10,20)(20,0){10}{\line(1,0){20}}
\multiput(110,20)(0,20){2}{\line(0,1){20}}
\put(110,60){\line(-1,1){14}}
\put(110,60){\line(1,1){14}}
\put(10,13){\makebox(0,0){*}}
\put(30,13){\makebox(0,0){$a$}}
\put(50,13){\makebox(0,0){$b$}}
\put(70,13){\makebox(0,0){$c$}}
\put(90,13){\makebox(0,0){$d$}}
\put(110,13){\makebox(0,0){$e$}}
\put(130,13){\makebox(0,0){$\tilde{d}$}}
\put(150,13){\makebox(0,0){$\tilde{c}$}}
\put(170,13){\makebox(0,0){$\tilde{b}$}}
\put(190,13){\makebox(0,0){$\tilde{a}$}}
\put(210,13){\makebox(0,0){$\tilde{*}$}}
\put(97,67){\makebox(0,0){$h$}}
\put(125,67){\makebox(0,0){$\tilde{h}$}}
\put(105,41){\makebox(0,0){$f$}}
\put(105,60){\makebox(0,0){$g$}}
\end{picture}
\thinlines
\unitlength 0.5mm
 \begin{picture}(200,50)
\multiput(10,20)(20,0){10}{\circle*{3}}
\multiput(110,40)(20,0){2}{\circle*{3}}
\multiput(10,20)(20,0){9}{\line(1,0){20}}
\multiput(110,20)(20,0){2}{\line(0,1){20}}
\put(10,13){\makebox(0,0){1}}
\put(30,13){\makebox(0,0){$a$}}
\put(50,13){\makebox(0,0){2}}
\put(70,13){\makebox(0,0){$c$}}
\put(90,13){\makebox(0,0){3}}
\put(110,13){\makebox(0,0){$e$}}
\put(130,13){\makebox(0,0){5}}
\put(150,13){\makebox(0,0){$g$}}
\put(170,13){\makebox(0,0){6}}
\put(190,13){\makebox(0,0){$\tilde{a}$}}
\put(105,41){\makebox(0,0){4}}
\put(125,41){\makebox(0,0){$\tilde{c}$}}
\end{picture}
\end{center}
\caption{The pair of graphs (3) in Haagerup's candidates list}
\label{17graph}
\end{figure}

About case (3) in \S7 of \cite{H} as in Figure \ref{17graph}, as well as the case $n=3$ of (2),
we can easily determine a biunitary connection uniquely
on the four graphs consisting of the graphs (3), and we thus have a
hyperfinite II$_1$ subfactor with index $(5+\sqrt{17})/2$ constructed 
from 
the connection by commuting square as in \cite{HS}.
The problem is whether this subfactor has (3) as
(dual) principal graphs or not.  This amounts to verifying the {\it flatness}
condition of the connection.
In 1996, K. Ikeda made a numerical check of 
the flatness of this connection by approximate computations on a computer in \cite{Ike} and 
showed
that the graphs (3) are very ``likely'' to exist
as (dual) principal graphs.
(He also made a numerical verification of the flatness for
the case $n=7$ of (2) in \S7 of \cite{H}.)  

In this paper, we will give the proof of the existence for the case of
 index $(5+\sqrt{13})/2$ previously announced by the second named author,
 and give the proof of the existence for the case of index
 $(5+\sqrt{17})/2$. The proof in the latter case was recently obtained
  by computations of
 the first named author based on a strategy of the second named author.
Our main result in
this paper is as follows.  \par

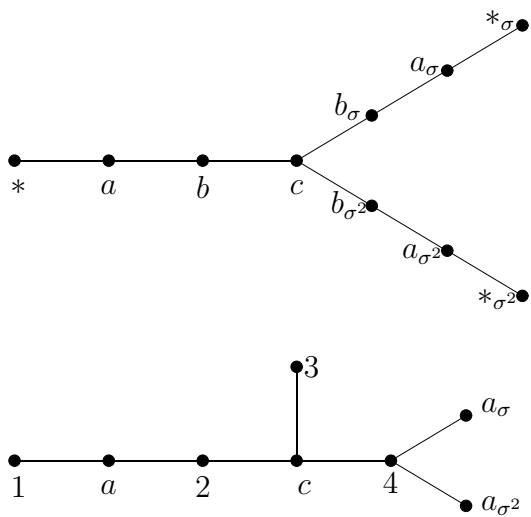
\begin{figure}[h]
\begin{center}
\thinlines
\unitlength 0.5mm
 \begin{picture}(170,80)(0,-20)
\multiput(10,20)(25,0){4}{\circle*{3}}
\multiput(10,20)(25,0){3}{\line(1,0){25}}
\multiput(85,20)(20,12){3}{\line(5,3){20}}
\multiput(85,20)(20,-12){3}{\line(5,-3){20}}
\multiput(105,32)(20,12){3}{\circle*{3}}
\multiput(105,8)(20,-12){3}{\circle*{3}}
\put(11,13){\makebox(0,0){$*$}}
\put(35,13){\makebox(0,0){$a$}}
\put(60,13){\makebox(0,0){$b$}}
\put(85,13){\makebox(0,0){$c$}}
\put(95,32){\s b}
\put(115,44){\s a}
\put(135,56){\s *}
\put(94,6){\ss b}
\put(113,-6){\ss a}
\put(133,-18){\ss *}
\end{picture}
\end{center}
\begin{center}
\thinlines
\unitlength 0.5mm
 \begin{picture}(170,50)(0,0)
\multiput(10,20)(25,0){5}{\circle*{3}}
\put(85,45){\circle*{3}}
\put(85,45){\line(0,-1){25}}
\multiput(10,20)(25,0){4}{\line(1,0){25}}
\put(110,20){\line(5,3){20}}
\put(110,20){\line(5,-3){20}}
\put(130,32){\circle*{3}}
\put(130,8){\circle*{3}}
\put(11,13){\makebox(0,0){$1$}}
\put(35,13){\makebox(0,0){$a$}}
\put(60,13){\makebox(0,0){$2$}}
\put(87,13){\makebox(0,0){$c$}}
\put(110,14){\makebox(0,0){$4$}}
\put(89,45){\makebox(0,0){$3$}}
\put(134,33){\s a}
\put(134,7){\ss a}
\end{picture}
\end{center}
\caption{The case $n=3$ of the pair of graphs (2) in the list of Haagerup}
\label{13graph}
\end{figure}

\begin{thm}{\rm (}$(5+\sqrt{13})/2$ case{\rm )} \\
The two graphs in Figure \ref{13graph} are realized as a pair of
(dual) principal graphs of a subfactor with index equal to 
$\frac{5+\sqrt{13}}{2}$ of the hyperfinite II$_1$ factor.
\end{thm}

\begin{thm}{\rm (}$(5+\sqrt{17})/2$ case{\rm )} \\
The two graphs in Figure \ref{17graph}
are realized as a pair of (dual) principal graphs of a
subfactor with index equal to $\frac{5+\sqrt{17}}{2}$
of the hyperfinite II$_1$ factor.
\end{thm}

 In Section 2, we will
give two key lemmas to prove our two main results respectively.  
In Section 3, we will give a
construction of {\it generalized open string bimodules} which
is a
generalization of Ocneanu's open string bimodules in \cite{Oc1}, 
\cite{SN}, and we will give a correspondence between bimodules and
 general biunitary connections on finite graphs.
In Sections 4 and 5, we will prove our two main theorems respectively. 

  The first named author acknowledges financial
supports and hospitality from Odense University and University of Copenhagen during her visit to Denmark in March/April and 
September, 1997. She also acknowledges a financial support 
from the Honda Heizaemon memorial fellowship.
She is much grateful to Y. Kawahigashi and
M. Izumi for constant advice and encouragement. 


\section{Key lemmas to the main results}

In this section, we will give the key lemmas which have been proved 
by the second named author. 

First of all, we will explain the motivation to the lemmas. Proofs of
our main  
theorems presented in Section 1 are reduced to verifying ``flatness'' of the 
biunitary connections which exist on the four graphs made of the pairs 
of the graphs in Figure \ref{13graph}, \ref{17graph} respectively. 
However, it is well known that, to verify {\it flatness} exactly is 
almost impossible, except for some easy cases such as the biunitary 
connections arising from the subfactors of crossed products of 
finite groups (\cite[10.6]{EK}). So far, in the history of 
classification of subfactors, several methods have been introduced to 
prove flatness/nonflatness of biunitary connections. Finding 
inconsistency of the fusion 
rule on the graph of a given biunitary connection has been sometimes effective 
to prove nonflatness, e.g. $D_{2n+1}, E_{7}$ (\cite{I1}, \cite{SV},
\ldots). On the other hand, since  
consistency of fusion rule does {\it never} mean flatness of given biunitary 
connection, several ideas have been introduced to prove flatness 
(\cite{EK1}, \cite{I2}, \cite{IK}, \cite{K} \ldots). 
The second named author, however, inspected the fusion rules of the
upper graphs in Figures \ref{13graph} and \ref{17graph} and noticed that 
if we can construct bimodules satisfying a part of the fusion rule, we
can conclude that there exists a subfactor having the desired principal
graph. 
These are the lemmas which we explain 
in this section.
Now we introduce the notation 
used in the lemmas.
\begin{defi}
Let $N$ and $M$ be II$_1$ factors and ${}_N X_M=X$ be an  $N$-$M$
bimodule (see \cite{EK}).
We denote by $R_X (M)$ 
the right action of $M$ on $X$, and by $L_X (N)$ the left 
action of $N$. We have the subfactor $R_X (M)' \supset
L_X (N)$. We denote its Jones index by $[X]$. We define the
 principal graph of  
the bimodule 
$X$ as that of the subfactor $R_X (M)' \supset L_X (N)$.
\end{defi}

\begin{defi}
For bimodules $X$ and $Y$ with common coefficient algebras, 
we define
$$ \langle X, Y \rangle = {\rm dim}\,{\rm Hom}(X, Y).$$ \\
A formal ${\bf Z}$-linear combination $Y$ of bimodules (of 
finite index) will be called {\em positive} if it is an actual 
bimodule, i.e. of $ \lan Y,Z \ran \geq 0$ for any irreducible bimodule $Z$ 
which appears in the direct sum decomposition of $X$. \\
When $Y \cong X \oplus Z$ for some positive bimodule $Z$, we write 
$X \prec Y$. 
\end{defi}

Hereafter we use the expression as follows, so far as it does not cause 
misunderstanding.

$$ {\bf 1}_{N} = {}_N N_{N}, $$
$$ 2X = X \oplus X, $$
$$ XY = X \otimes_{N} Y,$$
$$ X^{2}= X \otimes_{N} X,$$
where $N$ is a II$_{1}$ factor and $X$ and $Y$ are suitable bimodules.
\begin{lm}
Let $X={}_{N}X_{M}$ be a bimodule with finite Jones index larger than or 
equal to four. Then, \\
1) $X \oX - {\bf 1}_{N}$ and $(X \oX)^{2}-3 X 
 \oX \oplus {\bf 1}_{N}$ are positive $N$-$N$ bimodules. \\
2) $X \oX X -2X$ and $(X \oX)^{2} 
\otimes_{N} X - 4 X \oX X \oplus 3X$ are positive 
$N$-$M$ bimodules.
\end{lm}
{\it Proof} \\
Let ${\cal G}$ be the principal graph of $X$. We set
\begin{displaymath}
{\cal G}^{(0)}_{\rm even}:=\mbox{ set of all irreducible components 
of } \\
{\bf 1}_{N}, (X \oX)^{n}, \quad n=1,2,\ldots,
\end{displaymath}

\begin{displaymath}
{\cal G}^{(0)}_{\rm odd}:=\mbox{ set of all irreducible components 
of } \\
X, (X \oX)^{n}X, \quad n=1,2,\ldots,
\end{displaymath}
where ${\cal G}^{(0)}_{\rm even}$ (resp. ${\cal G}^{(0)}_{\rm odd}$) 
means the even (resp. odd) vertices of ${\cal G}$,
and $G=(G_{Y,Z})_{Y \in {\cal G}^{(0)}_{\rm even}, Z \in {\cal 
G}^{(0)}_{\rm odd}}$ to be the incidence matrix for ${\cal G}$ , i.e.,
$$G_{Y,Z}=\lan Y X, Z \ran, \quad Y \in {\cal G}^{(0)}_{\rm even}.$$
Since $4 \leq [X] < \infty$, We have $ 2 \leq ||G|| < \infty$. Put
$$ \Delta= \left( \matrix{
0 & G \cr
G^{t} & 0 } \right), $$
then $\Delta$ is the adjacency matrix of $G^{(0)}={\cal G}^{(0)}_{\rm 
even} \cup {\cal G}^{(0)}_{\rm odd}$ and $||\Delta|| \geq ||G|| \geq 
2$. Let $P_{0}$, $P_{1}$, $P_{2}$, \ldots be the sequence of the 
polynomials given by 
$$P_{0}(x)=1, \quad P_{1}(x)=x, \ldots 
, P_{n+1}(x)=P_{n}(x) x - P_{n-1}(x). $$
Then, by \cite{HW}, all of $P_{2}(\Delta)$, $P_{3}(\Delta)$, \ldots  
have non-negative entries. For $n=2, 3, 4, 5$, we get, in particular, that
$$GG^{t}-1, \, GG^{t}G-2G, \, (GG^{t})^{2}-3GG^{t}+1, \mbox{ and } 
(GG^{t})^{2}G-4GG^{t}G+3G $$
have non-negative entries. Hence, for $W \in {\cal G}^{(0)}_{\rm 
even}$, 
$$ \lan X  \oX - {\bf 1}_{N}, W \ran= (GG^{t}-1)_{{\bf 1}_{N},W} 
\geq 0,$$
$$ \lan (X \oX)^{2}-3 X \oX \oplus {\bf 1}_{N}, W 
\ran = 
((GG^{t})^{2}-3GG^{t}+1)_{{\bf 1}_{N}, W} \geq 0,$$
namely, $X \oX - {\bf 1}_{N}$ and $(X \oX)^{2}-3 
X \oX \oplus {\bf 1}_{N}$ are positive $N$-$N$ bimodules. The 
same argument (with $W \in {\cal G}^{(0)}_{\rm odd}$) shows that 
$X \oX X -2X$ and $(X \oX)^{2} 
 X - 4 X \oX X \oplus 3X$ are positive 
$N$-$M$ bimodules. \\
\qed

\subsection{Key lemma for the case of index $(5+ \sqrt{13})/2$}

In this subsection, we present the key lemma given by the second named
author to which the construction of the finite depth subfactor with
index $(5+ \sqrt{13})/2$ is reduced.
\begin{lm}
Let $M$ and $N$ be II$_1$ factors. Assume the following.  \\
1) We have an $N$-$M$ bimodule $X={}_N X_M$ of index 
$\frac{(5+\sqrt{13})}{2}$, \\
2) We have an $N$-$N$ bimodule $S = {}_N S_N$ of index $1$ satisfying
$ S^{3}  \cong {\bf 1}_N$ and 
$ S \not\cong {\bf 1}_{N}$, i.e.
$S$ is given by an automorphism of $N$ of outer period
$3$,  \\
3) The six bimodules \\
\begin{center}
${\bf 1}_{N}$, $S$, $S^{2}$, $X \oX - {\bf 1}_{N}$,\\
 $S(X{\ovl X}-{\bf 1}_{N})$, $S^{2}(X{\ovl X} - {\bf 1}_{N})$ 
 \end{center}
are irreducible and mutually inequivalent. \\
4) The four bimodules \\
\begin{center}
$X$, $SX$, $ S^{2}X$, $X{\ovl X}  X-2 X$  \\
\end{center}
are irreducible and mutually inequivalent,  \\
5)  (the most important assumption)
$$S (X{\ovl X} - {\bf 1}_{N})
\cong (X {\ovl X} -{\bf 1}_{N}) S^{2}.$$
Then the principal graph of $X$ and bimodules 
corresponding to the vertices on the graph are as follows;

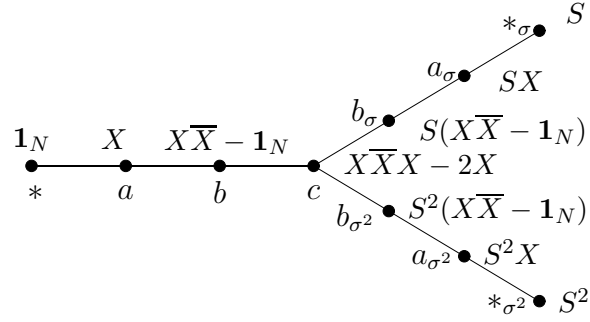
\begin{center}
\thinlines
\unitlength 0.5mm
 \begin{picture}(170,80)(0,-20)
\multiput(10,20)(25,0){4}{\circle*{3}}
\multiput(10,20)(25,0){3}{\line(1,0){25}}
\multiput(85,20)(20,12){3}{\line(5,3){20}}
\multiput(85,20)(20,-12){3}{\line(5,-3){20}}
\multiput(105,32)(20,12){3}{\circle*{3}}
\multiput(105,8)(20,-12){3}{\circle*{3}}
\put(11,13){\makebox(0,0){$*$}}
\put(10,27){\makebox(0,0){\sm{${\bf 1}_N$}}}
\put(35,13){\makebox(0,0){$a$}}
\put(32,27){\makebox(0,0){\sm{$X$}}}
\put(60,13){\makebox(0,0){$b$}}
\put(62,27){\makebox(0,0){\sm{$ X \oX- {\bf 1}_N$}}}
\put(85,13){\makebox(0,0){$c$}}
\put(93,18){\sm{$X \oX X-2X$}}
\put(95,32){\s b}
\put(115,44){\s a}
\put(135,56){\s *}
\put(91,4){\ss b}
\put(111,-6){\ss a}
\put(131,-18){\ss *}
\put(113,27){\sm{$S(X \oX-{\bf 1}_N)$}}
\put(134,40){\sm{$SX$}}
\put(152,58){\sm{$S$}}
\put(110,6){\sm{$S^2(X \oX-{\bf 1}_N)$}}
\put(130,-6){\sm{$S^2 X$}}
\put(150,-19){\sm{$S^2$}}
\end{picture}
\end{center}
\end{lm}
{\it Remark}  \\
By lemma 1, all the above bimodules are well-defined (i.e. positive 
in the sense of Def. 2).\\  \phantom{x} \\
{\it Proof} \\
We have 
\begeq
&&\lan (X \oX X - 2X)\oX, X \oX - {\bf 1}_N \ran   \\
&=& \lan X \oX X- 2X, (X \oX - {\bf 1}_N)X \ran \quad 
\mbox{(by Frobenius reciprocity)}\\
&=&\lan X \oX X- 2X, X \oX X- 2X \ran + \lan X \oX X- 2X, X
 \ran  \\
&=&1
\endeq
because $X$ and $X \oX X- 2X$ are irreducible and inequivalent.
Hence, by irreducibility of $X \oX-{\bf 1}_N$, we have  \\
(i) \phantom{XXXXXXXX} 
 $X \oX-{\bf 1}_N \prec  (X \oX X- 2X)\oX $. \\
We have \\
(ii) \phantom{XXXXXXXX}
   $(X \oX X- 2X)\oX  \cong (X \oX - {\bf 1}_N)^2-{\bf 1}_N$,  \\
and 5) says \\
(iii) \phantom{XXXXXXXX}
  $S(X \oX - {\bf 1}_N) \cong (X \oX - {\bf 1}_N)S^2.$   \\
Hence \\
(iv)        
\begeq
S^2(X \oX - {\bf 1}_N)& \cong & S(X \oX - {\bf 1}_N)S^2   \\
& \cong & (X \oX - {\bf 1}_N)S^4  \\
& \cong & (X \oX - {\bf 1}_N)S.
\endeq
Therefore, \\
$$S(X \oX - {\bf 1}_N)^2 \cong (X \oX - {\bf 1}_N)S^2(X \oX - {\bf 1}_N)
 \cong (X \oX - {\bf 1}_N)^2 S $$
and
$$S^2 (X \oX - {\bf 1}_N)^2 \cong (X \oX - {\bf 1}_N)^2 S^2. $$
Hence by (ii),
\begeq
S(X \oX X- 2X)\oX S^2 & \cong & (X \oX X- 2X)\oX S^3    \\
& \cong & (X \oX X- 2X)\oX,
\endeq
and similarly
$$ S^2 (X \oX X- 2X)\oX S \cong (X \oX X- 2X)\oX. $$
So, by (i), \\
(v) \pha{XXXXX}         
$S(X \oX X- 2X)\oX S^2 \prec (X \oX X- 2X)\oX, $ \\
(vi) \pha{XXXXX} 
$ S^2 (X \oX X- 2X)\oX S \prec (X \oX X-2X)\oX, $ \\
by (iii) and (iv). Hence, by (i), (v), (vi) and 3),
$$ (X \oX - {\bf 1}_N), \quad S(X \oX - {\bf 1}_N), \quad 
S^2 (X \oX - {\bf 1}_N)$$
are mutually inequivalent subbimodules of $(X \oX X- 2X)$,
 i.e.,
$$(X \oX X- 2X)\oX \cong (X \oX - {\bf 1}_N) \oplus S(X \oX - {\bf 1}_N)
\oplus S^2 (X \oX - {\bf 1}_N) \oplus Y $$
where $Y$ is an $N$-$N$ bimodule (possibly zero). \par
Since $X={}_N X_M$ is irreducible, the subfactor
$R_X (M)' \supset L_X (N)$ has the trivial relative commutant,
 hence extremal (see \cite{P}, p.176). Therefore, the square
 root of the Jones index of a bimodule $[\cdot]^{1/2}$ is 
additive and multiplicative on the bimodules expressed in 
terms of $X$ and $\oX$ (see \cite{P}).
Thus, we have 
$$ \sqrt{ [(X \oX X- 2X)\oX ] }= 3 \sqrt{ [ (X \oX - {\bf 1}_N) ] }
+ \sqrt{ [Y]}, $$
where the index of a zero bimodule is defined to be $0$. Hence 
with $\lambda = [X] = \sqrt{\frac{5+\sqrt{13}}{2}} $, we get

\begeq
\sqrt{ [Y]} &=& \la(\la^3-2\la)-3(\la^2-1) \\
&=& \la^4 -5\la^2 +3 \\
&=&  0.
\endeq
We must finally prove that \\
(vii) \pha{XXXXX} $ S (X \oX X- 2X) \cong X \oX X- 2X. $
\\
To see this, we compute 
\begeq
&&\lan S (X \oX X- 2X), X \oX X- 2X \ran  \\
&=& \lan S (X \oX X- X), X \oX X- 2X \ran - 
\lan SX, X \oX X- 2X \ran   \\
&=& \lan S (X \oX - {\bf 1}_N), (X \oX X- 2X) \oX \ran -
\lan SX, X \oX X- 2X \ran.
\endeq
The first bracket is $1$ because $S (X \oX - {\bf 1}_N)$ is contained
in $(X \oX X- 2X) \oX$ with multiplicity $1$, and the second
bracket is $0$ because $SX$ and $X \oX X- 2X$ are irreducible 
and inequivalent by 4), hence
$$ \lan S (X \oX X- 2X), X \oX X- 2X \ran =1, $$
thus, the equality of the irreducible bimodules (vii) holds. \par
From all the above, it follows easily that \\
(a) ${\bf 1}_N \in {\cal G}^{(0)}$, \\
(b) ${\cal G}$ is connected, \\
(c) Multiplication by $X$ (resp. $\oX$) from the right (resp. left)
on any bimodule $U$ 
 in ${\cal G}^{(0)}_{odd}$ (resp. ${\cal G}^{(0)}_{even}$) gives
a direct sum of the bimodules in ${\cal G}^{(0)}_{even}$ 
(resp. ${\cal G}^{(0)}_{even}$)  
connected to $U$ by edges, \\
namely, we find that ${\cal G}$ is the principal graph 
of $X$. \\
\qed

\subsection{Key lemma for the case of index $(5+\sqrt{17})/2$}

In this subsection we present the key lemma similar to the previous one
for the construction of the finite depth subfactor with index
$(5+\sqrt{17})/2$. 
\begin{lm}
Let $M$, $N$ be II$_1$ factors. Assume the following. \\
1) We have an $N$-$M$ bimodule $X$ of index $\frac{5+\sqrt{17}}{2}$. \\
2) We have an $N$-$N$ bimodule $S$ of index $1$ satisfying $S^2 
\cong {\bf 1}_N$ and $S \not\cong {\bf 1}_{N}$,
 i.e., $S$ is given by an automorphism of $N$ of outer period $2$, \\
3) The eight $N$-$N$ bimodules \\
\begin{center}
${\bf 1}_N, S, X \oX -{\bf 1}_N, S(X \oX -{\bf 1}_N), 
(X \oX -{\bf 1}_N)S$, \\
$S(X \oX -{\bf 1}_N)S, (X \oX)^2-3X \oX \oplus {\bf 1}_{N}, 
S((X \oX)^2-3X \oX \oplus {\bf 1}_N)$
\end{center}
are irreducible and mutually inequivalent. \\
4) The six $N$-$M$ bimodules \\
$$ X, SX, X \oX X -2X, S(X \oX X -2X), (X \oX)^2 X -
4 X \oX X \oplus 3X, (X \oX-{\bf 1}_N)SX $$
are irreducible and mutually inequivalent. \\
5) (The most important assumption)
$$ S(X \oX-{\bf 1}_N)SX \cong  (X \oX-{\bf 1}_N)SX. $$
Then the principal graph of $X$ and the bimodules corresponding
to the vertices on the graph are as follows: 
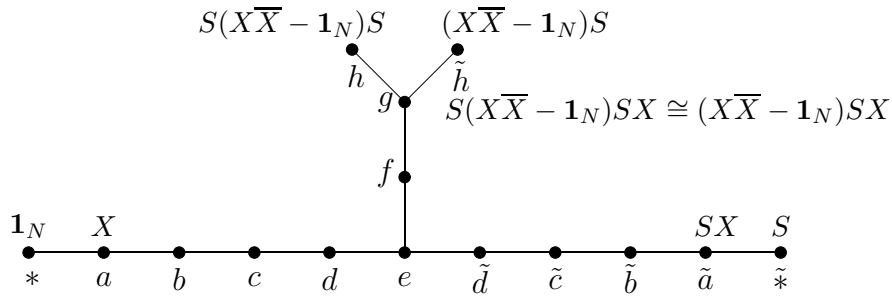
\begin{center}
\thinlines
\unitlength 0.5mm
 \begin{picture}(220,80)
\multiput(10,20)(20,0){11}{\circle*{3}}
\multiput(110,40)(0,20){2}{\circle*{3}}
\multiput(96,74)(28,0){2}{\circle*{3}}
\multiput(10,20)(20,0){10}{\line(1,0){20}}
\multiput(110,20)(0,20){2}{\line(0,1){20}}
\put(110,60){\line(-1,1){14}}
\put(110,60){\line(1,1){14}}
\put(11,13){\makebox(0,0){$*$}}
\put(10,27){\makebox(0,0){\sm{${\bf 1}_N$}}}
\put(30,13){\makebox(0,0){$a$}}
\put(30,27){\makebox(0,0){\sm{$X$}}}
\put(50,13){\makebox(0,0){$b$}}
\put(70,13){\makebox(0,0){$c$}}
\put(90,13){\makebox(0,0){$d$}}
\put(110,13){\makebox(0,0){$e$}}
\put(130,13){\makebox(0,0){$\tilde{d}$}}
\put(150,13){\makebox(0,0){$\tilde{c}$}}
\put(170,13){\makebox(0,0){$\tilde{b}$}}
\put(190,13){\makebox(0,0){$\tilde{a}$}}
\put(193,27){\makebox(0,0){\sm{$SX$}}}
\put(210,13){\makebox(0,0){$\tilde{*}$}}
\put(210,27){\makebox(0,0){\sm{$S$}}}
\put(97,67){\makebox(0,0){$h$}}
\put(80,81){\makebox(0,0){\sm{$ S(X \oX- {\bf 1}_N)S $}}}
\put(125,67){\makebox(0,0){$\tilde{h}$}}
\put(142,81){\makebox(0,0){\sm{$(X \oX- {\bf 1}_N)S$}}}
\put(105,41){\makebox(0,0){$f$}}
\put(105,60){\makebox(0,0){$g$}}
\put(180,58){\makebox(0,0){\sm{$S(X\oX-{\bf 1}_N)SX \cong 
(X\oX-{\bf 1}_N)SX$}}}
\end{picture}
\end{center}
where, 
\begeq
&& b \cdots X \oX- {\bf 1}_N \\
&& c \cdots X\oX X-2X, \\
&& d \cdots (X \oX)^2-3X \oX \oplus {\bf 1}_{N}, \\
&& e \cdots (X \oX)^2 X -4X \oX X \oplus 3X, \\
&& f \cdots (X \oX-{\bf 1}_{N})S(X \oX-{\bf 1}_{N})-S(X \oX-{\bf 1}_{N})S, \\
&& {\tilde d} \cdots S((X \oX)^2-3X \oX \oplus {\bf 1}_{N}), \\
&&{\tilde c} \cdots S(X\oX X-2X), \\
&& {\tilde{b}} \cdots S(X \oX- {\bf 1}_N)
\endeq
\end{lm}
{\it Remark} \par \noindent
By lemma 1, we know that all the bimodules above except for the 
``bimodule'' corresponding to $f$ are well-defined (i.e., positive in the 
sense of Def. 2). The well-definedness of the bimodule at $f$ will come 
out of the proof below. \\
{\it Proof} \\
In the proof we will sometimes use formal computations in the ${\bf 
Z}$-linear span of the $N$-$M$ bimodules or $N$-$N$ bimodules 
considered. The symbol $ \lan \phantom{x}, \phantom{x} \ran $ for computing 
the dimension of the space of intertwiners can be extended to ${\bf 
Z}$ bilinear maps and $ \lan Z,Z \ran =0$ implies $Z=0$ also for these generalized 
bimodule. \\
Since ${\bf 1}_N \prec X \oX$, we have by 5) that both 
$(X \oX -{\bf 1}_N)S$ and $S(X \oX -{\bf 1}_N)S$ are (equivalent to) 
subbimodules of $(X \oX -{\bf 1}_N)SX \oX$.  By 3), $(X \oX -{\bf 1}_N)S$ and
$S(X \oX -{\bf 1}_N)S$ are two non-equivalent irreducible bimodules.
Hence we can write
$$(X \oX -{\bf 1}_N)SX \oX \cong (X \oX -{\bf 1}_N)S \oplus 
S(X \oX -{\bf 1}_N)S \oplus \R $$
where $\R$ is an $N$-$N$ bimodule, and we have
$$\R \cong (X \oX -{\bf 1}_N)S(X \oX -{\bf 1}_N) - S(X \oX -{\bf 1}_N)S.$$
Note that $\R \not\cong 0$ because $\R \cong 0$ would imply
$[(X \oX -{\bf 1}_N)]=1$, which is impossible since $X$ is irreducible
and $[X] > 4$.\\
Next we will show the following. \\
(i) The bimodule $\R$ is irreducible, \\
(ii) $S((X \oX)^2-3X \oX \oplus {\bf 1}_N)S \cong 
(X \oX)^2-3X \oX \oplus {\bf 1}_N, $ \\
(iii) $S((X \oX)^2 X-4X \oX X \oplus 3X) \cong
 (X \oX)^2 X-4X \oX X \oplus 3X.$ \\
 \bigskip
We have
\begeq
\lan \R, \R \ran &=& \lan (X \oX -{\bf 1}_N)S(X \oX -{\bf 1}_N) ,
(X \oX -{\bf 1}_N)S(X \oX -{\bf 1}_N) \ran \\
&& - 2 \lan S(X \oX -{\bf 1}_N)S, (X \oX -{\bf 1}_N)S(X \oX -{\bf 1}_N) \ran \\
&& + \lan S(X \oX -{\bf 1}_N)S, S(X \oX -{\bf 1}_N)S \ran  \\
&=& t_1 +t_2 + t_3,
\endeq
where  \\
\pha{XXXXXXX} $ t_1=\lan S(X \oX -{\bf 1}_N)^2 S, (X \oX -{\bf 1}_N)^2 
\ran, $ \\
\pha{XXXXXXX} $t_2= \lan S(X \oX -{\bf 1}_N)S,
 (X \oX -{\bf 1}_N)S(X \oX -{\bf 1}_N) \ran, $ \\
\pha{XXXXXXX} $ t_3 = \lan (X \oX -{\bf 1}_N), (X \oX -{\bf 1}_N) \ran.$ \\
Note first that $t_3=1$ because $(X \oX -{\bf 1}_N)$ is irreducible.
Next 
\begeq
 t_2 =&  \lan S(X \oX -{\bf 1}_N)S, (X \oX -{\bf 1}_N)SX \oX \ran \\
  & -\lan S(X \oX -{\bf 1}_N)S, (X \oX -{\bf 1}_N)S \ran.
\endeq
The last term is $0$ because $S(X \oX -{\bf 1}_N)S$ and 
$(X \oX -{\bf 1}_N)S$ are irreducible and inequivalent by 3). 
Hence, using 4) and 5), we get
$$t_2 = \lan S(X \oX -{\bf 1}_N)SX, (X \oX -{\bf 1}_N)SX \ran =1. $$
To compute $t_1$, set irreducible bimodules $Y$, $Z$ as
$$ Y=X \oX-{\bf 1}_N, \quad Z=(X \oX)^2-3X \oX \oplus{\bf 1}_N. $$
Then
\begeq
t_1 &=& \lan S({\bf 1}_N \oplus \Y \oplus \Z)S, {\bf 1}_N \oplus
\Y \oplus \Z 
\ran \\
&=& \lan {\bf 1}_N, {\bf 1}_N \ran + \lan S \Y S, \Y \ran +\lan
S \Z S, \Z \ran + 2\lan {\bf 1}_N, \Y \ran+ \\
&& 2\lan {\bf 1}_N,  \Z \ran + 2 \lan S \Y S, \Z \ran.
\endeq
By 3), ${\bf 1}_N$, $\Y$, $S\Y S$, and $\Z$ are irreducible
and mutually inequivalent. Hence
$$ t_1= 1+ \lan S \Z S, \Z \ran. $$
Altogether, we have shown that
\begeq
 \lan \R, \R \ran &=& (1+\lan S\Z S, \Z \ran)-2+1 \\
&=& \lan S\Z S, \Z \ran.
\endeq
Since $\R \neq 0$, we have $\lan \R,R \ran \geq 1$. Moreover, 
since $\Z$ is irreducible, so is $S\Z S$. Hence 
$ \lan S \Z S, \Z \ran \leq 1$. Therefore
$$ \lan \R, R \ran = \lan S \Z S, \Z \ran =1, $$
which shows that $\R$ is irreducible, and using that $\Z$
and $S\Z S$ are irreducible, we also get that $S \Z S
\cong \Z.$ Hence we have verified (i) and (ii). \\
To prove (iii), put
$$ G=X \oX X-2X, \quad E=(X \oX)^2 X -4X \oX X  \oplus 3X. $$
Note that $E=((X \oX)^2  -3X \oX  \oplus {\bf 1}_N)X-(X \oX X \oplus 
2X)$,
then, by (ii)
\begeq
E & \cong & S((X \oX)^2  -3X \oX  \oplus {\bf 1}_N)SX-X \oX X \oplus
2X \\
 & \cong & S(X \oX-{\bf 1}_N)^2 SX - S(X \oX-{\bf 1}_N) SX-(X \oX-{\bf 1}_N)X.
\endeq
Using 5), we have
\begeq
E & \cong & S(X \oX-{\bf 1}_N) S(X \oX-{\bf 1}_N) SX-(X \oX-{\bf 1}_N) SX
-(X \oX-{\bf 1}_N)X \\
& \cong & S(X \oX-{\bf 1}_N) SX \oX SX-S(X \oX-{\bf 1}_N)X-(X \oX-{\bf 1}_N) SX
-(X \oX-{\bf 1}_N) X.
\endeq
Hence, again using 5) we get 
\begeq
E & \cong & (X \oX-{\bf 1}_N) SX \oX SX-S(X \oX-{\bf 1}_N) X-(X \oX-{\bf 1}_N) SX
-(X \oX-{\bf 1}_N) X  \\
& \cong & (X \oX-{\bf 1}_N) SX(\oX SX-{\bf 1}_N)-({\bf 1}_N \oplus S)(X \oX-{\bf 1}_N)X.
\endeq
From this expression of $E$ and 5), we clearly have $SE \cong E$,
 which proves (iii). \\
We next prove \\
(iv) \pha{XXXXXXX} $RX \cong (X \oX-{\bf 1}_N) SX \oplus E,$\\
where $R \cong (X \oX-{\bf 1}_N) S(X \oX-{\bf 1}_N) -S(X \oX-{\bf 1}_N)S $
is irreducible by (i). We put 
$$ E'=RX -  (X \oX-{\bf 1}_N) SX. $$
By 5), we have
\begeq
E' & \cong & (X \oX-{\bf 1}_N) S(X \oX-{\bf 1}_N)X-2(X \oX-{\bf 1}_N) SX \\
& \cong & (X \oX-{\bf 1}_N) S(X \oX-3 {\bf 1}_N)X.
\endeq
To prove (iv), we just have to show that $E \cong E'$, namely
$$ \lan E'-E, E'-E \ran =0. $$
Note that 
$$ \lan E'-E, E'-E \ran = s_1- 2s_2 +s_3, $$
where $s_1= \lan E', E' \ran$, $s_2=\lan E', E \ran$, and
$s_3 = \lan E,E \ran.$  \\
First, $s_3=1$ because $E$ is irreducible. Next,
\begeq
s_2 &=& \lan  (X \oX-{\bf 1}_N) S(X \oX-3 {\bf 1}_N)X,
 (X \oX-{\bf 1}_N) (X \oX-3 {\bf 1}_N)X \ran \\
&=& \lan S, (X \oX-{\bf 1}_N)^2 (X \oX-3 {\bf 1}_N)X \oX
(X \oX-3 {\bf 1}_N)(X \oX-{\bf 1}_N) \ran \\
&=& \lan S(X \oX-{\bf 1}_N) (X \oX-3 {\bf 1}_N)X, 
(X \oX-{\bf 1}_N) (X \oX-3 {\bf 1}_N)X \ran \\
&=& \lan SE,E \ran \\
&=& 1 \quad \mbox{(by (iii))}.
\endeq
Finally,
\begeq
s_1 &=& \lan E', E' \ran  \\
&=& \lan S(X \oX-{\bf 1}_N)^2 S, (X \oX-3 {\bf 1}_N)X \oX
(X \oX-3 {\bf 1}_N) \ran \\
&=& \lan S({\bf 1}_N \oplus Y \oplus Z )S, (X \oX-3 {\bf 1}_N)^2
X \oX \ran. 
\endeq
Here, using $SYS = S(X \oX-{\bf 1}_N)S$ and $SZS \cong Z$ 
by (ii), we have
\begeq
s_1 &=& \lan {\bf 1}_N \oplus S(X \oX-{\bf 1}_N)S \oplus 
Z, (X \oX-3 {\bf 1}_N)^2 X \oX \ran \\
&=& \lan X \oplus S(X \oX-{\bf 1}_N)SX \oplus ZX, 
(X \oX-3 {\bf 1}_N)^2 X \ran \\
&=& \lan X \oplus (X \oX-{\bf 1}_N)SX \oplus ZX,
 (X \oX-3 {\bf 1}_N)^2 X \ran,
\endeq
where we have used 5) again. \\
We expand $ZX$ and $(X \oX-3 {\bf 1}_N)^2 X$ in terms of the
irreducible bimodules $X$, $G=X \oX X-2X$, and $E=(X \oX)^2 X
-4X \oX X \oplus 3X$, and get
$$ZX \cong G \oplus E$$
and
$$(X \oX-3 {\bf 1}_N)^2 X \cong 2X - 2G \oplus E. $$
Hence 
$$s_1= \lan X \oplus (X \oX-{\bf 1}_N)SX \oplus G \oplus E,
2X - 2G \oplus E \ran. $$
By 4), $X$ $G$, $E$, and $(X \oX-{\bf 1}_N)SX$ are irreducible
and mutually inequivalent, hence
$$s_1 = 2 \lan X, X \ran-2\lan G,G \ran + \lan E,E \ran =1.$$
Altogether,
$$ \lan E'-E,E'-E \ran =s_1-2 s_2 + s_3 =1-2+1=0,$$
which proves (iv). \\
We need to prove one more relation \\
(v) \pha{XXXXXXXX} $E \oX \cong Z \oplus SZ \oplus R.$ \\
To prove (v), note first that $X$, $Z$, and $E$ all correspond
to the vertices in ${\cal G}^{(0)}_{odd}$, the set
 of the odd vertices of the principal graph of $X$. (We write 
 ${\cal G}^{(0)}_{even}$ for the even vertices.) Hence, $E \ovl{E} \in 
{\cal G}^{(0)}_{even}$, and since 
$$ \lan S, E \ovl{E} \ran = \lan SE, E \ran =1,$$
$S$ is an irreducible subbimodule of $E \ovl{E}$, so also 
$S \in {\cal G}^{(0)}_{even}$. Therefore, every irreducible $N$-
$N$ bimodule or $N$-$M$ bimodule that can be expressed in 
terms of $X$, $\oX$, and $S$, belong to the principal graph 
${\cal G}$ of $X$. Therefore, by the same argument in the 
proof of the previous lemma, the square root of the Jones index
is additive and multiplicative on the $N$-$N$ bimodules or 
$N$-$M$ bimodules which can be expressed in terms of $X$, 
$\oX$, and $S$, because it will occur as a submodule of 
$$(X \oX)^n, \quad n \geq 0, \quad {\rm or} \quad (X \oX)^n X,
\quad n \geq 1. $$
\\ Since $ZX \cong G \oplus E$, we have $E \prec ZX$ and 
therefore  \\
(vi) \pha{XXXXXXXXXXX} $Z \prec E \oX.$  \\
By (iii), also \\
(vii) \pha{XXXXXXXXXXX} $SZ \prec E \oX.$ \\
Moreover, in the same way, we have \\
(viii) \pha{XXXXXXXXXXX} $R \prec E \oX$ \quad
 by (iv). \\
We know that $Z$, $SZ$, and $S$ are irreducible and 
$Z \not\cong SZ$ by 3). Moreover, by a simple computation
using the additivity and multiplicativity of $[\cdot]^{1/2}$, we 
have
$$ [R]^{1/2} = [Z]^{1/2}-1=[SZ]^{1/2}-1.$$
Hence all of $R$, $Z$, $SZ$ are mutually inequivalent. Thus
$$E \oX \cong Z \oplus SZ \oplus R \oplus T,$$
where $T$ is an $N$-$N$ bimodule. By $[X]=(5+\sqrt{17})/2$,
we easily get 
$$[E \oX]^{1/2} = [Z]^{1/2} + [SZ]^{1/2} + [R]^{1/2},$$
hence, $T=0.$ \par
Putting everything together, we see that conditions
 (a), (b), (c) in the proof of the previous lemma hold,
namely, ${\cal G}$ is the principal graph of $X$. \\
\qed


\section{Generalized open string bimodules}

In section 2, we have reduced our construction problem to verification
of certain fusion rules, but we still have a problem of handling
bimodules in a concrete way.  For example, we do not know how to
represent $X$ or $S$, or how to verify equalities of infinite
dimensional bimodules.
In this section, we will introduce the item 
to make full use of the lemmas.
\begin{figure}[h]
\begin{center}
\thinlines
\unitlength 1.0mm
\begin{picture}(30,20)(0,-7)
\multiput(11,-2)(0,10){2}{\line(1,0){8}}
\multiput(10,7)(10,0){2}{\line(0,-1){8}}
\multiput(10,-2)(10,0){2}{\circle*{1}}
\multiput(10,8)(10,0){2}{\circle*{1}}
\put(6,12){\makebox(0,0){$V_0$}}
\put(24,12){\makebox(0,0){$V_1$}}
\put(6,-6){\makebox(0,0){$V_2$}}
\put(24,-6){\makebox(0,0){$V_3$}}
\put(15,12){\makebox(0,0){$\cal K$}}
\put(15,-6){\makebox(0,0){${\cal L}$}}
\put(15,3){\makebox(0,0){$\alpha$}}
\end{picture} 
\end{center}
\caption{the connection $\alpha$ with four graphs}
\label{alpha}
\end{figure}
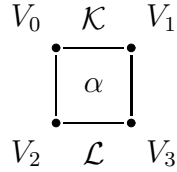
Consider a biunitary connection $\alpha$, as in Figure \ref{alpha}, on the four
graphs with upper graph 
${\cal K}$, lower graph ${\cal L}$ and the sets of vertices $V_0$, \dots ,
$V_3$. Note that, by the definition of biunitary connection, the graphs
${\cal K}$ and ${\cal L}$ should be connected, and the vertical graphs
are not nessesary to be connected.
 We fix the vertices $ \sk \in V_0$ and $*_{\cal L} \in
V_2$.  We will now construct the bimodule corresponding to $\a$. \par
First we construct AFD II$_1$ factors from the string algebras
$$ K =  \overline{\bigcup_{n=1}^{\infty}{{\rm String}_{*_{\cal K}}^{(n)} {\cal
K}}}^{\rm weak} $$ 
$$ L =  \overline{\bigcup_{n=1}^{\infty}{{\rm String}_{*_{\cal L}}^{(n)} {\cal
L}}}^{\rm weak} $$ 
by the GNS construction using the unique trace, where 
\begin{eqnarray*}  {\rm String }_{*_{\cal G} }^{(n)} {\cal G}= {\rm span} 
\{ (\xi, \eta) \vert
\mbox{ a pair of paths on the graph ${\cal G}$} \\
   \mbox{ $s(\xi)=s(\eta)=*_{\cal G}$, $r(\xi)=r(\eta)$, 
$|\xi|=|\eta|=n$}\}. 
\end{eqnarray*} 
Here for a path $\zeta$, we denote the initial vertex, the final vertex
and the length of the path by $s(\zeta)$, $r(\zeta)$ and $|\zeta|$
respectively. We define its $*$-algebra structure as 
$$ (\xi,\eta) \cdot (\xi',\eta') = \delta_{\eta, \xi'} (\xi,\eta'),$$
$$ (\xi,\eta)^* = (\eta,\xi).  $$
Now we have another AFD II$_1$ factor
$$ {\tilde L}=\overline{ \bigcup_{n=0}^{\infty}
{\rm span} \left\{ 
\left(
\thinlines
\unitlength 0.5mm
\begin{picture}(110,10)(0,5)
\multiput(7,7)(50,0){2}{$*_{\cal K}$}
\multiput(10,0)(50,0){2}{\line(1,0){40}}
\multiput(10,0)(50,0){2}{\circle*{1}}
\multiput(10,6)(50,0){2}{\line(0,-1){6}}
\multiput(50,0)(50,0){2}{\circle*{1}}
\multiput(48,-4.5)(50,0){2}{$x$}
\put(55,0){,}
\end{picture} 
\right)
\vert \begin{array}{c} 
\mbox{ a pair of paths, } x \in {\cal L}^{(0)}, \\
\mbox{horizontal paths are in ${\cal L}$, length $n$.}
\end{array} \right\} }^{\rm weak}, $$
where ${\cal L}^{(0)}$ denotes the set of vertices on ${\cal L}$.
We identify elements in $K$ with elements in $\tilde L$ by the embedding 
using
connection $\alpha$, and then have an AFD II$_1$ subfactor $ K \subset
{\tilde L}. $ (See \cite{EK}, Chapter 11). \par

Next we construct the $K$-$L$ bimodule corresponding to $\alpha$. 
Consider a pair of paths as follows:

$$
\left(
\thinlines
\unitlength 0.5mm
\begin{picture}(110,10)(0,5)
\put(10,7){$*_{\cal K}$}
\put(14,10){\line(1,0){36}}
\put(50,10){\circle*{1}}
\put(50,10){\line(0,-1){10}}
\put(50,0){\circle*{1}}
\put(55,0){,}
\put(60,-3){$*_{\cal L}$}
\put(64,0){\line(1,0){36}}
\put(100,0){\circle*{1}}
\end{picture}
\right),
$$
here the horizontal part of the left (resp. right) path consists of 
edges of the graph ${\cal K}$ (resp. ${\cal L}$), the vertical edge is
from one of the two vertical graphs of the four graphs of the connection $\alpha$, and the paths have a 
common final vertex.
In general, a pair of paths, as above, with a common final vertex, not
necessary with a common initial vertex, is called an {\it open
string}. It was first introduced by Ocneanu in \cite{Oc1} in more
restricted situations.
We embed an open string of length $k$ into the linear span of open
strings of length $k+1$ in a similar way to the embedding of string
algebras as follows:

\begin{eqnarray*}
&&\left(
\thinlines
\unitlength 0.5mm
\begin{picture}(110,10)(0,5)
\put(10,7){$*_{\cal K}$}
\put(14,10){\line(1,0){36}}
\put(50,10){\circle*{1}}
\put(50,10){\line(0,-1){10}}
\put(45,3){$\eta$}
\put(50,0){\circle*{1}}
\put(55,0){,}
\put(60,-3){$*_{\cal L}$}
\put(64,0){\line(1,0){36}}
\put(100,0){\circle*{1}}
\end{picture}
\right)  \\
&=&
\sum_{|\xi|=1} 
\left(
\thinlines
\unitlength 0.5mm
\begin{picture}(120,10)(5,5)
\put(10,7){$*_{\cal K}$}
\put(14,10){\line(1,0){36}}
\put(50,10){\circle*{1}}
\put(50,10){\line(0,-1){10}}
\put(50,0){\circle*{1}}
\put(50,0){\line(1,0){10}}
\put(60,0){\circle*{1}}
\put(53,-7){$\xi$}
\put(65,0){,}
\put(70,-3){$*_{\cal L}$}
\put(74,0){\line(1,0){36}}
\put(110,0){\circle*{1}}
\put(110,0){\line(1,0){10}}
\put(120,0){\circle*{1}}
\put(113,-7){$\xi$}
\end{picture}
\right)  \\
&=&
\sum_{\eta', \xi'} \sum_{|\xi|=1}
\thinlines
\unitlength 0.5mm
\begin{picture}(30,20)
\multiput(11,-2)(0,10){2}{\line(1,0){8}}
\multiput(10,7)(10,0){2}{\line(0,-1){8}}
\multiput(10,-2)(10,0){2}{\circle*{1}}
\multiput(10,8)(10,0){2}{\circle*{1}}
\put(15,13){\makebox(0,0){$\xi'$}}
\put(15,-7){\makebox(0,0){$\xi$}}
\put(5,2.5){$\eta$}
\put(22,2.5){$\eta'$}
\put(15,3){\makebox(0,0){$\alpha$}}
\end{picture}
\left(
\thinlines
\unitlength 0.5mm
\begin{picture}(120,10)(5,5)
\put(10,7){$*_{\cal K}$}
\put(14,10){\line(1,0){46}}
\put(50,10){\circle*{1}}
\put(53,12){$\xi'$}
\put(60,10){\line(0,-1){10}}
\put(61,4){$\eta'$}
\put(60,10){\circle*{1}}
\put(60,0){\circle*{1}}
\put(60,0){\circle*{1}}
\put(65,0){,}
\put(70,-3){$*_{\cal L}$}
\put(74,0){\line(1,0){36}}
\put(110,0){\circle*{1}}
\put(110,0){\line(1,0){10}}
\put(120,0){\circle*{1}}
\put(113,-7){$\xi$}
\end{picture}
\right),
\end{eqnarray*}
here the square marked with $\alpha$ means the value given
by the connection $\alpha$. \par

We define the vector space spanned by the above open strings with the
above embedding as follows.
\begin{eqnarray*}
 \stackrel{\circ}{X^{\alpha}} &=& \bigcup_{n}{\rm span} \{ (\xi, \eta) | s(\xi)=*_{\cal K},
s(\eta)=*_{\cal L}, r(\xi)=r(\eta)\} \\
&=& \bigcup {\rm span} \{ \left(
\thinlines
\unitlength 0.5mm
\begin{picture}(110,10)(0,5)
\put(10,7){$*_{\cal K}$}
\put(14,10){\line(1,0){36}}
\put(50,10){\circle*{1}}
\put(50,10){\line(0,-1){10}}
\put(50,0){\circle*{1}}
\put(55,0){,}
\put(60,-3){$*_{\cal L}$}
\put(64,0){\line(1,0){36}}
\put(100,0){\circle*{1}}
\end{picture}
\right)
\}.
\end{eqnarray*}

We define an inner product of $\stackrel{\circ}{X^{\alpha}}$ as 
the sesqui-linear extension of the following;

\begin{eqnarray*}
&&\langle (\xi \cdot \zeta, \eta), (\xi' \cdot \zeta', \eta') \rangle  \\
&=&
\langle \left(
\thinlines
\unitlength 0.5mm
\begin{picture}(110,10)(0,5)
\put(10,7){$*_{\cal K}$}
\put(14,10){\line(1,0){36}}
\put(27,3){$\xi$}
\put(50,10){\circle*{1}}
\put(50,10){\line(0,-1){10}}
\put(51,4){$\zeta$}
\put(50,0){\circle*{1}}
\put(55,0){,}
\put(60,-3){$*_{\cal L}$}
\put(64,0){\line(1,0){36}}
\put(81,4){$\eta$}
\put(100,0){\circle*{1}}
\end{picture}
\right),
\left(
\thinlines
\unitlength 0.5mm
\begin{picture}(110,10)(0,5)
\put(10,7){$*_{\cal K}$}
\put(14,10){\line(1,0){36}}
\put(27,3){$\xi'$}
\put(50,10){\circle*{1}}
\put(50,10){\line(0,-1){10}}
\put(51,4){$\zeta'$}
\put(50,0){\circle*{1}}
\put(55,0){,}
\put(60,-3){$*_{\cal L}$}
\put(64,0){\line(1,0){36}}
\put(81,4){$\eta'$}
\put(100,0){\circle*{1}}
\end{picture}
\right)
\rangle  \\
&=&
\frac{\mu_{\cal K} (s(\zeta))}{\mu_{\cal L} (r(\zeta))}
\delta_{\zeta, \zeta'}
{\rm tr}_{K} (\xi,\xi')
{\rm tr}_L (\eta', \eta),
\end{eqnarray*}
where $\xi \cdot \zeta$ denotes the concatenation of $\xi$ and $\zeta$,
$\mu_{\cal K}$ and $\mu_{\cal L}$ denotes the Perron-Frobenius eigen 
vector of the graphs,
tr$_K$ is the unique trace on $K$, and tr$_L$ is 
as well. We set $(\xi, \xi')$ and $(\eta',\eta)$ are the elements $0$ of 
${\tilde L}$ and $L$ respectively if the end 
points of each pair of paths do not coincide. \par
By this inner product, $\stackrel{\circ}{X^{\alpha}}$ is regarded as a
pre-Hilbert space, and then we complete it and denote the completion by
$X^{\alpha}$. \par
We have the natural left action of $K$ and the right action of $L$ as 
follows; for 

\begin{eqnarray*}
x &=& \left(
\thinlines
\unitlength 0.5mm
\begin{picture}(110,10)(0,5)
\put(10,7){$*_{\cal K}$}
\put(14,10){\line(1,0){36}}
\put(31,3){$\xi$}
\put(50,10){\circle*{1}}
\put(50,10){\line(0,-1){10}}
\put(50,0){\circle*{1}}
\put(55,0){,}
\put(60,-3){$*_{\cal L}$}
\put(64,0){\line(1,0){36}}
\put(81,4){$\eta$}
\put(100,0){\circle*{1}}
\end{picture}
\right)         \in       X^{\alpha}, \\
k &=& 
\left(
\thinlines
\unitlength 0.5mm
\begin{picture}(110,10)(0,5)
\put(10,7){$*_{\cal K}$}
\put(14,10){\line(1,0){36}}
\put(31,3){$\sigma$}
\put(50,10){\circle*{1}}
\put(55,0){,}
\put(60,7){$*_{\cal K}$}
\put(64,10){\line(1,0){36}}
\put(81,3){$\sigma'$}
\put(100,10){\circle*{1}}
\end{picture}
\right) 		\in 	K, \\
l &=& 
\left(
\thinlines
\unitlength 0.5mm
\begin{picture}(110,10)(0,5)
\put(10,-3){$*_{\cal L}$}
\put(14,0){\line(1,0){36}}
\put(31,4){$\rho$}
\put(50,0){\circle*{1}}
\put(55,0){,}
\put(60,-3){$*_{\cal L}$}
\put(64,0){\line(1,0){36}}
\put(81,4){$\rho'$}
\put(100,0){\circle*{1}}
\end{picture}
\right)        \in  L,  \\
\end{eqnarray*}
we have

\begin{eqnarray*}
k \cdot x &=& \sum_{ |\zeta|=1} 
\left(
\thinlines
\unitlength 0.5mm
\begin{picture}(110,10)(0,5)
\put(10,7){$*_{\cal K}$}
\put(14,10){\line(1,0){36}}
\put(31,3){$\sigma$}
\put(50,10){\circle*{1}}
\put(50,10){\line(0,-1){10}}
\put(51,4){$\zeta$}
\put(50,0){\circle*{1}}
\put(55,0){,}
\put(60,7){$*_{\cal K}$}
\put(64,10){\line(1,0){36}}
\put(81,3){$\sigma'$}
\put(100,10){\circle*{1}}
\put(100,10){\line(0,-1){10}}
\put(101,4){$\zeta$}
\put(100,0){\circle*{1}}
\end{picture}
\right) \cdot x  \\
&=&
\sum_{ |\zeta|=1} \delta_{\sigma' \cdot \zeta, \xi} 
\left(
\thinlines
\unitlength 0.5mm
\begin{picture}(110,10)(0,5)
\put(10,7){$*_{\cal K}$}
\put(14,10){\line(1,0){36}}
\put(31,3){$\sigma$}
\put(50,10){\circle*{1}}
\put(50,10){\line(0,-1){10}}
\put(51,4){$\zeta$}
\put(50,0){\circle*{1}}
\put(55,0){,}
\put(60,-3){$*_{\cal L}$}
\put(64,0){\line(1,0){36}}
\put(81,4){$\eta$}
\put(100,0){\circle*{1}}
\end{picture}
\right),  \\
x \cdot l &=& 
\delta_{\eta, \rho} 
\left(
\thinlines
\unitlength 0.5mm
\begin{picture}(110,10)(0,5)
\put(10,7){$*_{\cal K}$}
\put(14,10){\line(1,0){36}}
\put(31,3){$\xi$}
\put(50,10){\circle*{1}}
\put(50,10){\line(0,-1){10}}
\put(50,0){\circle*{1}}
\put(55,0){,}
\put(60,-3){$*_{\cal L}$}
\put(64,0){\line(1,0){36}}
\put(81,4){$\rho'$}
\put(100,0){\circle*{1}}
\end{picture}
\right). 
\end{eqnarray*}
By the extension of this action, the Hilbert space $X^{\alpha}$ is
considered as a Hilbert $K$-$L$ bimodule ${}_K X^{\alpha} _L$.
Then we have a $K$-$L$
bimodule ${}_K {X^{\a}}_L$ constructed from $\alpha$.  (We call this
bimodule made
of open strings an {\it open string bimodule\/}.  This is a
generalization of open string bimodules in \cite{Oc1} and \cite{SN}, 
which are the bimodules constructed from flat connections.) \par
We make the correspondence between direct sums, relative 
tensor products, and the contragredient map of bimodules and ``sums'',  
``products'', and the renormalization of connections, so that fusion 
rules on open string bimodules reduced to the operations of connections.
 (\cite{Oc3}) \par
First we introduce the sum of two connections.  Consider $\alpha$ and 
$\beta$ as connections on the four graphs with upper graph ${{\cal K}}$,  
lower graph ${{\cal L}}$ and sets of vertices $V_0$,\dots,$V_3$ as in
Figure \ref{alpha} (The side graphs of $\alpha$ and $\beta$ need not to be
identical), then they give rise two $K$-$L$ bimodules. We define the sum of 
the connections as follows:  

\begin{eqnarray*}
(\alpha + \beta) (
\thinlines
\unitlength 0.5mm
\begin{picture}(30,25)(0,0)
\multiput(11,-2)(0,10){2}{\line(1,0){8}}
\multiput(10,7)(10,0){2}{\line(0,-1){8}}
\multiput(10,-2)(10,0){2}{\circle*{1}}
\multiput(10,8)(10,0){2}{\circle*{1}}
\put(15,12){\makebox(0,0){$k$}}
\put(15,-6){\makebox(0,0){$l$}}
\put(6,3){\makebox(0,0){$m$}}
\put(24,3){\makebox(0,0){$n$}}
\end{picture}
) 
 = \quad 
\left\{ \begin{array}{c}
\alpha(\thinlines
\unitlength 0.5mm
\begin{picture}(30,25)(0,0)
\multiput(11,-2)(0,10){2}{\line(1,0){8}}
\multiput(10,7)(10,0){2}{\line(0,-1){8}}
\multiput(10,-2)(10,0){2}{\circle*{1}}
\multiput(10,8)(10,0){2}{\circle*{1}}
\put(15,12){\makebox(0,0){$k$}}
\put(15,-6){\makebox(0,0){$l$}}
\put(6,3){\makebox(0,0){$m$}}
\put(24,3){\makebox(0,0){$n$}}
\end{picture} 
), \qquad \mbox{ if both $m$, $n$ are edges appearing in $\alpha$,} \\ \\
\beta(\thinlines
\unitlength 0.5mm
\begin{picture}(30,20)(0,0)
\multiput(11,-2)(0,10){2}{\line(1,0){8}}
\multiput(10,7)(10,0){2}{\line(0,-1){8}}
\multiput(10,-2)(10,0){2}{\circle*{1}}
\multiput(10,8)(10,0){2}{\circle*{1}}
\put(15,12){\makebox(0,0){$k$}}
\put(15,-6){\makebox(0,0){$l$}}
\put(6,3){\makebox(0,0){$m$}}
\put(24,3){\makebox(0,0){$n$}}
\end{picture} 
), \qquad \mbox{ if both $m$, $n$ are edges appearing in $\beta$,} \\
\thinlines
\unitlength 0.5mm
\begin{picture}(30,20)(0,0)
\put(15,3){\makebox(0,0){0,}}
\end{picture} \qquad \quad \mbox{otherwise.} \phantom{XXXXXXXXXX}
\end{array}
\right.
\end{eqnarray*}

Obviously it satisfies the biunitarity.  We denote 
the bimodule constructed from a connection $\gamma$ by $X^{\gamma}$.
By considering the action of $K$ from the left, it is easy to see that
$$ {}_K X^{\alpha}_L \oplus {}_K X^{\beta}_L\quad =\quad 
{}_K X^{\alpha + \beta}_L, $$
thus, we can use the summation of connections instead of the direct
sum of bimodules. \par
Next we define the product of connections (\cite{Oc3}, \cite{SN}).
  Consider the connections
$\alpha$ and $\beta$, as in Figure \ref{algam}, which give rise to 
$K$-$L$ bimodule (resp. $L$-$M$ bimodule). 

\begin{figure}[h]
\begin{center}
\thinlines
\unitlength 1.0mm
\begin{picture}(30,20)(0,-7)
\multiput(11,-2)(0,10){2}{\line(1,0){8}}
\multiput(10,7)(10,0){2}{\line(0,-1){8}}
\multiput(10,-2)(10,0){2}{\circle*{1}}
\multiput(10,8)(10,0){2}{\circle*{1}}
\put(6,12){\makebox(0,0){$V_0$}}
\put(24,12){\makebox(0,0){$V_1$}}
\put(6,-6){\makebox(0,0){$V_2$}}
\put(24,-6){\makebox(0,0){$V_3$}}
\put(15,12){\makebox(0,0){${\cal K}$}}
\put(15,-6){\makebox(0,0){${\cal L}$}}
\put(15,3){\makebox(0,0){$\alpha$}}
\end{picture} 
\thinlines
\unitlength 1.0mm
\begin{picture}(30,20)(0,-7)
\multiput(11,-2)(0,10){2}{\line(1,0){8}}
\multiput(10,7)(10,0){2}{\line(0,-1){8}}
\multiput(10,-2)(10,0){2}{\circle*{1}}
\multiput(10,8)(10,0){2}{\circle*{1}}
\put(6,12){\makebox(0,0){$V_2$}}
\put(24,12){\makebox(0,0){$V_3$}}
\put(6,-6){\makebox(0,0){$V_4$}}
\put(24,-6){\makebox(0,0){$V_5$}}
\put(15,12){\makebox(0,0){${\cal L}$}}
\put(15,-6){\makebox(0,0){$\cal M$}}
\put(15,3){\makebox(0,0){$\beta$}}
\end{picture} 
\end{center}
\caption{}
\label{algam}
\end{figure}

Note that ${\cal L}$ is appears in the both four graphs.
Then we can define the product connection $\alpha\beta$ on the 
four graphs with upper graph ${\cal K}$, lower graph $\cal M$ and
the sets of vertices $V_0$, $V_1$, $V_4$, $V_5$.  The side graphs
 consist of the edges $\{p-q \; \vert \; p \in V_0({\rm resp.} 
\; V_1), \; q \in V_4({\rm resp.} \; V_5) \}$
with multiplicity $$ \sharp \{ p-x-q  \;\vert \mbox{ a path of length 2
from $p$ to $q$,  }  x \in V_2({\rm resp.} \; V_3) \}. $$  We have the
connection $\alpha 
\beta$ as follows: \par
\[
(\alpha  \beta) (
\thinlines
\unitlength 0.5mm
\begin{picture}(30,25)(0,0)
\multiput(11,-2)(0,10){2}{\line(1,0){8}}
\multiput(10,7)(10,0){2}{\line(0,-1){8}}
\multiput(10,-2)(10,0){2}{\circle*{1}}
\multiput(10,8)(10,0){2}{\circle*{1}}
\put(15,12){\makebox(0,0){$k$}}
\put(15,-6){\makebox(0,0){$m$}}
\put(6,3){\makebox(0,0){$n$}}
\put(24,3){\makebox(0,0){$o$}}
\end{picture} 
) 
=  (\alpha  \beta) \left(
\thinlines
\unitlength 0.5mm
\begin{picture}(30,25)(0,0)
\put(11,13){\line(1,0){8}}
\multiput(10,12)(10,0){2}{\line(0,-1){8}}
\multiput(10,3)(10,0){2}{\circle*{1}}
\multiput(10,13)(10,0){2}{\circle*{1}}
\put(15,17){\makebox(0,0){$k$}}
\put(6,8){\makebox(0,0){$n_1$}}
\put(24,8){\makebox(0,0){$o_1$}}
\put(11,-7){\line(1,0){8}}
\multiput(10,2)(10,0){2}{\line(0,-1){8}}
\multiput(10,-7)(10,0){2}{\circle*{1}}
\multiput(10,3)(10,0){2}{\circle*{1}}
\put(15,-11){\makebox(0,0){$m$}}
\put(6,-2){\makebox(0,0){$n_2$}}
\put(24,-2){\makebox(0,0){$o_2$}}
\multiput(11,3)(2,0){4}{\line(1,0){1}}
\end{picture}
\right)
 = 
\sum_{l}
\alpha(
\thinlines
\unitlength 0.5mm
\begin{picture}(30,25)(0,0)
\multiput(11,-2)(0,10){2}{\line(1,0){8}}
\multiput(10,7)(10,0){2}{\line(0,-1){8}}
\multiput(10,-2)(10,0){2}{\circle*{1}}
\multiput(10,8)(10,0){2}{\circle*{1}}
\put(15,12){\makebox(0,0){$k$}}
\put(15,-6){\makebox(0,0){$l$}}
\put(6,3){\makebox(0,0){$n_1$}}
\put(24,3){\makebox(0,0){$o_1$}}
\end{picture} 
)
\beta(
\thinlines
\unitlength 0.5mm
\begin{picture}(30,25)(0,0)
\multiput(11,-2)(0,10){2}{\line(1,0){8}}
\multiput(10,7)(10,0){2}{\line(0,-1){8}}
\multiput(10,-2)(10,0){2}{\circle*{1}}
\multiput(10,8)(10,0){2}{\circle*{1}}
\put(15,12){\makebox(0,0){$l$}}
\put(15,-6){\makebox(0,0){$m$}}
\put(6,3){\makebox(0,0){$n_2$}}
\put(24,3){\makebox(0,0){$o_2$}}
\end{picture} 
 ), 
\]
where $n_1$ and $n_2$ are edges such that their concatenation $n_1 \cdot
n_2 $ is $n$, and $o_1$ and $o_2$ are as well. 
We observe that this process corresponds to the following process of 
composing commuting squares of finite dimension.
$$
\begin{array}{ccccc}
\begin{array}{ccc}
A & \subset & B \\
\cap & & \cap \\
C & \subset & D
\end{array} & ,  &
\begin{array}{ccc}
C & \subset & D \\
\cap & &  \cap \\
E & \subset & F
\end{array} &
\Rightarrow &
\begin{array}{ccc}
A & \subset & B \\
\cap & & \cap \\
E & \subset & F
\end{array}
\end{array},
$$
where these three squares are finite dimensional commuting squares.
We will show that 
$_K X^{\alpha} {\otimes}_L X^{\beta}_M$ is isomorphic to 
$_K X^{\alpha\beta}_M$.  We define the map $\varphi$ from 
$_K X^{\alpha} {\otimes}_L X^{\beta}_M$ to $_K X^{\alpha\beta}_M$
as follows; For

\begin{eqnarray*}
x &=& \left(
\thinlines
\unitlength 0.5mm
\begin{picture}(110,10)(0,5)
\put(10,7){$*_{\cal K}$}
\put(14,10){\line(1,0){36}}
\put(31,3){$\xi$}
\put(50,10){\circle*{1}}
\put(50,10){\line(0,-1){10}}
\put(50,0){\circle*{1}}
\put(55,0){,}
\put(60,-3){$*_{\cal L}$}
\put(64,0){\line(1,0){36}}
\put(81,4){$\eta$}
\put(100,0){\circle*{1}}
\end{picture}
\right)  \\
&=& \sum_{|\zeta|=1}
 \left(
\thinlines
\unitlength 0.5mm
\begin{picture}(110,20)(0,10)
\put(10,17){$*_{\cal K}$}
\put(14,20){\line(1,0){36}}
\put(31,13){$\xi$}
\put(50,20){\circle*{1}}
\put(50,20){\line(0,-1){20}}
\put(50,10){\circle*{1}}
\put(50,0){\circle*{1}}
\put(51,4){$\zeta$}
\put(55,0){,}
\put(60,7){$*_{\cal L}$}
\put(64,10){\line(1,0){36}}
\put(81,14){$\eta$}
\put(100,10){\circle*{1}}
\put(100,10){\line(0,-1){10}}
\put(101,4){$\zeta$}
\put(100,0){\circle*{1}}
\end{picture}
\right)      \in {}_K \stackrel{\circ}{X^{\alpha}}_L,  \\
y &=& \left(
\thinlines
\unitlength 0.5mm
\begin{picture}(110,10)(0,5)
\put(10,7){$*_{\cal L}$}
\put(14,10){\line(1,0){36}}
\put(31,4){$\rho$}
\put(50,10){\circle*{1}}
\put(50,10){\line(0,-1){10}}
\put(50,0){\circle*{1}}
\put(55,0){,}
\put(60,-3){$*_{\cal M}$}
\put(64,0){\line(1,0){36}}
\put(81,4){$\sigma$}
\put(100,0){\circle*{1}}
\end{picture}
\right) \in {}_L \stackrel{\circ}{X^{\beta}}_M, 
\end{eqnarray*}
we define
\begin{eqnarray*}
\varphi (  x {\otimes}_L y ) &=& \; x \cdot y   \\
&=& \delta_{\eta \cdot \zeta, \rho} 
 \left(
\thinlines
\unitlength 0.5mm
\begin{picture}(110,20)(0,10)
\put(10,17){$*_{\cal K}$}
\put(14,20){\line(1,0){36}}
\put(31,13){$\xi$}
\put(50,20){\circle*{1}}
\put(50,20){\line(0,-1){20}}
\put(50,10){\circle*{1}}
\put(50,0){\circle*{1}}
\put(51,4){$\zeta$}
\put(55,0){,}
\put(60,-3){$*_{\cal M}$}
\put(64,0){\line(1,0){36}}
\put(81,4){$\sigma$}
\put(100,0){\circle*{1}}
\end{picture}
\right)    \in   {}_K X^{\alpha\beta}_M.
\end{eqnarray*}
Since
\begin{eqnarray*}
 (x {\otimes}_L y, \; x {\otimes}_L y) 
= (x (y, y)_L , x) &=& {\rm tr}_M( x^* \cdot x (y,y)_L),
 \\ 
 (x \cdot y, \; x \cdot y) 
&=& {\rm tr}_M(y^* \cdot x^* \cdot x \cdot y) \\
 &=& {\rm tr}_M(y^* \cdot (x^* \cdot x ) y) \quad ((x^* \cdot x )
 \in L) \\
&=&  ((x^* \cdot x)  y, y) \\
&=& {\rm tr}_M(x^* \cdot x (y, y)_L)=(x {\otimes}_L y, \; x 
{\otimes}_L y),
\end{eqnarray*}
where $x^*$ and $y^*$ means that we reverse the order of the pairs of
paths and also take the complex conjugate of their coefficients,
we see that $\varphi$ is an isometry, so it is well-defined as
 a linear map from
$_K X^{\alpha} {\otimes}_L X^{\beta}_M $ to  
$_K X^{\alpha\beta}_M $, and it is also injective.  Moreover,
it is surjective because, for an element 
$$
x =  
\left(
\thinlines
\unitlength 0.5mm
\begin{picture}(110,20)(0,10)
\put(10,17){$*_{\cal K}$}
\put(14,20){\line(1,0){36}}
\put(50,20){\circle*{1}}
\put(50,20){\line(0,-1){20}}
\put(50,10){\circle*{1}}
\put(50,0){\circle*{1}}
\put(51,4){$\rho$}
\put(55,0){,}
\put(60,-3){$*_{\cal M}$}
\put(64,0){\line(1,0){36}}
\put(100,0){\circle*{1}}
\end{picture}
\right)    \in   {}_K X^{\alpha\beta}_M,
$$
where we assume without loss of generality that $x$ is long 
enough that there is a path connecting $*_{\cal L}$ and 
$s(\rho)$, 

\begin{eqnarray*}
x &=& \sum_{\xi} \delta_{\rho,\xi}
\left(
\thinlines
\unitlength 0.5mm
\begin{picture}(120,20)(0,10)
\put(10,17){$*_{\cal K}$}
\put(14,20){\line(1,0){36}}
\put(50,20){\circle*{1}}
\put(50,20){\line(0,-1){20}}
\put(50,10){\circle*{1}}
\put(50,0){\circle*{1}}
\put(51,4){$\xi$}
\put(55,0){,}
\put(60,-3){$*_{\cal M}$}
\put(64,0){\line(1,0){36}}
\put(100,0){\circle*{1}}
\end{picture}
\right)   \\
&=& 
\sum_{\xi}
 \left(
\thinlines
\unitlength 0.5mm
\begin{picture}(110,20)(0,10)
\put(10,17){$*_{\cal K}$}
\put(14,20){\line(1,0){36}}
\put(50,20){\circle*{1}}
\put(50,20){\line(0,-1){20}}
\put(50,10){\circle*{1}}
\put(50,0){\circle*{1}}
\put(51,4){$\xi$}
\put(55,0){,}
\put(60,7){$*_{\cal L}$}
\put(64,10){\line(1,0){36}}
\put(81,14){$\eta$}
\put(100,10){\circle*{1}}
\put(100,10){\line(0,-1){10}}
\put(101,4){$\xi$}
\put(100,0){\circle*{1}}
\end{picture}
\right) 
\cdot
\left(
\thinlines
\unitlength 0.5mm
\begin{picture}(110,10)(0,5)
\put(10,7){$*_{\cal L}$}
\put(14,10){\line(1,0){36}}
\put(31,14){$\eta$}
\put(50,10){\circle*{1}}
\put(50,10){\line(0,-1){10}}
\put(51,4){$\rho$}
\put(50,0){\circle*{1}}
\put(55,0){,}
\put(60,-3){$*_{\cal M}$}
\put(64,0){\line(1,0){36}}
\put(100,0){\circle*{1}}
\end{picture}
\right)   \\
&=&
\varphi(
\left(
\thinlines
\unitlength 0.5mm
\begin{picture}(110,10)(0,5)
\put(10,7){$*_{\cal K}$}
\put(14,10){\line(1,0){36}}
\put(50,10){\circle*{1}}
\put(50,10){\line(0,-1){10}}
\put(50,0){\circle*{1}}
\put(55,0){,}
\put(60,-3){$*_{\cal L}$}
\put(64,0){\line(1,0){36}}
\put(81,4){$\eta$}
\put(100,0){\circle*{1}}
\end{picture}
\right) 
\otimes_L
\left(
\thinlines
\unitlength 0.5mm
\begin{picture}(110,10)(0,5)
\put(10,7){$*_{\cal L}$}
\put(14,10){\line(1,0){36}}
\put(31,14){$\eta$}
\put(50,10){\circle*{1}}
\put(50,10){\line(0,-1){10}}
\put(51,4){$\rho$}
\put(50,0){\circle*{1}}
\put(55,0){,}
\put(60,-3){$*_{\cal M}$}
\put(64,0){\line(1,0){36}}
\put(100,0){\circle*{1}}
\end{picture}
\right) 
)
\end{eqnarray*}
for some $\eta$ with $s(\eta)=*_{\cal L}$, $r(\eta)
=s(\rho)$. Therefore, we have the isomorphism
$$
{}_K X^{\alpha} \otimes_L X^{\beta}_M  \cong
{}_K X^{\alpha\beta}_M.
$$
Next we prove 
$$ \overline{{}_K X^{\alpha}_L} = {}_L X^{\tilde {\alpha}}_K,
 $$
here we denote the renormalization of the connection
$\alpha$ by ${\tilde \alpha}$. \par
Take an element 
$$ x=
\left(
\thinlines
\unitlength 0.5mm
\begin{picture}(110,10)(0,5)
\put(10,7){$*_{\cal K}$}
\put(14,10){\line(1,0){36}}
\put(31,14){$\zeta$}
\put(50,10){\circle*{1}}
\put(50,10){\line(0,-1){10}}
\put(51,4){$\xi$}
\put(50,0){\circle*{1}}
\put(55,0){,}
\put(60,-3){$*_{\cal L}$}
\put(64,0){\line(1,0){36}}
\put(81,4){$\eta$}
\put(100,0){\circle*{1}}
\end{picture}
\right)    \in   {}_K X^{\alpha}_L.
$$
For $x$, we easily see that its image by the contragredient map is given as
$$
{\bar x}=
\left(
\thinlines
\unitlength 0.5mm
\begin{picture}(110,10)(0,5)
\put(10,7){$*_{\cal L}$}
\put(14,10){\line(1,0){36}}
\put(30,14){$\eta$}
\put(50,10){\circle*{1}}
\put(50,10){\line(0,-1){10}}
\put(51,4){${\tilde \xi}$}
\put(50,0){\circle*{1}}
\put(55,0){,}
\put(60,-3){$*_{\cal K}$}
\put(64,0){\line(1,0){36}}
\put(81,4){$\zeta$}
\put(100,0){\circle*{1}}
\end{picture}
\right)    \in   \overline{{}_K X^{\alpha}_L},
$$
here ${\tilde \xi}$ means the upside down edge of $\xi$.
Since
\begin{eqnarray*}
{\bar x}
&=& 
\sum_{\sigma}
\overline{
\left(
\thinlines
\unitlength 0.5mm
\begin{picture}(130,14)(0,5)
\put(10,7){$*_{\cal K}$}
\put(14,10){\line(1,0){36}}
\put(30,14){$\zeta$}
\put(50,10){\circle*{1}}
\put(50,10){\line(0,-1){10}}
\put(51,4){$\xi$}
\put(50,0){\circle*{1}}
\put(50,0){\line(1,0){10}}
\put(55,1){$\sigma$}
\put(60,0){\circle*{1}}
\put(65,0){,}
\put(70,-3){$*_{\cal L}$}
\put(74,0){\line(1,0){36}}
\put(90,4){$\eta$}
\put(110,0){\circle*{1}}
\put(110,0){\line(1,0){10}}
\put(114,1){$\sigma$}
\put(120,0){\circle*{1}}
\end{picture}
\right)
} \\
&=&
\sum_{\sigma, \sigma', \xi'}
\thinlines
\unitlength 0.5mm
\begin{picture}(30,20)
\multiput(11,-2)(0,10){2}{\line(1,0){8}}
\multiput(10,7)(10,0){2}{\line(0,-1){8}}
\multiput(10,-2)(10,0){2}{\circle*{1}}
\multiput(10,8)(10,0){2}{\circle*{1}}
\put(15,13){\makebox(0,0){$\sigma'$}}
\put(15,-7){\makebox(0,0){$\sigma$}}
\put(5,2.5){$\xi$}
\put(22,2.5){$\xi'$}
\put(15,3){\makebox(0,0){$\alpha$}}
\end{picture}
\overline{
\left(
\thinlines
\unitlength 0.5mm
\begin{picture}(130,14)(0,5)
\put(10,7){$*_{\cal K}$}
\put(14,10){\line(1,0){36}}
\put(30,14){$\zeta$}
\put(50,10){\circle*{1}}
\put(60,10){\line(0,-1){10}}
\put(61,4){$\xi'$}
\put(60,0){\circle*{1}}
\put(50,10){\line(1,0){10}}
\put(54,11){$\sigma'$}
\put(65,0){,}
\put(70,-3){$*_{\cal L}$}
\put(74,0){\line(1,0){36}}
\put(90,4){$\eta$}
\put(110,0){\circle*{1}}
\put(110,0){\line(1,0){10}}
\put(114,1){$\sigma$}
\put(120,0){\circle*{1}}
\end{picture}
\right)
}   \\
&=&
\sum_{\sigma, \sigma', \xi'}
\thinlines
\unitlength 0.5mm
\begin{picture}(30,20)
\multiput(11,-2)(0,10){2}{\line(1,0){8}}
\multiput(10,7)(10,0){2}{\line(0,-1){8}}
\multiput(10,-2)(10,0){2}{\circle*{1}}
\multiput(10,8)(10,0){2}{\circle*{1}}
\put(15,13){\makebox(0,0){$\sigma$}}
\put(15,-7){\makebox(0,0){$\sigma'$}}
\put(5,2.5){${\tilde \xi}$}
\put(22,2.5){${\tilde \xi'}$}
\put(15,3){\makebox(0,0){${\tilde \alpha}$}}
\end{picture}
\left(
\thinlines
\unitlength 0.5mm
\begin{picture}(130,10)(0,5)
\put(10,7){$*_{\cal K}$}
\put(14,10){\line(1,0){36}}
\put(30,14){$\eta$}
\put(50,10){\circle*{1}}
\put(60,10){\line(0,-1){10}}
\put(61,4){${\tilde \xi}$}
\put(60,0){\circle*{1}}
\put(50,10){\line(1,0){10}}
\put(54,11){$\sigma$}
\put(65,0){,}
\put(70,-3){$*_{\cal L}$}
\put(74,0){\line(1,0){36}}
\put(90,4){$\zeta$}
\put(110,0){\circle*{1}}
\put(110,0){\line(1,0){10}}
\put(114,1){$\sigma'$}
\put(120,0){\circle*{1}}
\end{picture}
\right),
\end{eqnarray*}
${\bar x}$ is regarded as the element of ${}_L X^{\tilde
 \alpha}_K$, thus we have
$$ \overline{{}_K X^{\alpha}_L}
 \cong {}_L X^{\tilde \alpha}_K. $$
Now we have a good correspondence between the operations 
of certain bimodules and those of connections.  To complete 
it, we should check that the construction of bimodules from 
connections is a {\it one to one\/} correspondence of the 
equivalent classes.

\begin{thm}
Let $\alpha$ and $\beta$ be two connections as below;

\begin{center}
\thinlines
\unitlength 1.0mm
\begin{picture}(30,20)(0,-7)
\multiput(11,-2)(0,10){2}{\line(1,0){8}}
\multiput(10,7)(10,0){2}{\line(0,-1){8}}
\multiput(10,-2)(10,0){2}{\circle*{1}}
\multiput(10,8)(10,0){2}{\circle*{1}}
\put(6,12){\makebox(0,0){$V_0$}}
\put(24,12){\makebox(0,0){$V_1$}}
\put(6,-6){\makebox(0,0){$V_2$}}
\put(24,-6){\makebox(0,0){$V_3$}}
\put(15,12){\makebox(0,0){${\cal K}$}}
\put(15,-6){\makebox(0,0){${\cal L}$}}
\put(6,3){\makebox(0,0){${\cal S}_1$}}
\put(24,3){\makebox(0,0){${\cal T}_1$}}
\put(15,3){\makebox(0,0){$\alpha$}}
\end{picture} 
\thinlines
\unitlength 1.0mm
\begin{picture}(30,20)(0,-7)
\multiput(11,-2)(0,10){2}{\line(1,0){8}}
\multiput(10,7)(10,0){2}{\line(0,-1){8}}
\multiput(10,-2)(10,0){2}{\circle*{1}}
\multiput(10,8)(10,0){2}{\circle*{1}}
\put(6,12){\makebox(0,0){$V_0$}}
\put(24,12){\makebox(0,0){$V_1$}}
\put(6,-6){\makebox(0,0){$V_2$}}
\put(24,-6){\makebox(0,0){$V_3$}}
\put(15,12){\makebox(0,0){${\cal K}$}}
\put(15,-6){\makebox(0,0){$\cal L$}}
\put(6,3){\makebox(0,0){${\cal S}_2$}}
\put(24,3){\makebox(0,0){${\cal T}_2$}}
\put(15,3){\makebox(0,0){$\beta$}}
\end{picture} 
\end{center}
then the $K$-$L$ bimodules ${}_K X^{\alpha}_L$ and 
${}_K X^{\beta}_L$ are isomorphic if and only if $\alpha$ 
and $\beta$ are equivalent to each other up to gauge choice
for the vertical edges, in particular the pairs $({\cal S}_{1},{\cal 
T}_{1})$ and $({\cal S}_{2},{\cal 
T}_{2})$ of the vertical graphs must coincide. 
\end{thm}
{\it Remark} \par \noindent
In \cite{Oc3}, the same correspondence of bimodules and equivalent 
classes of connections has been introduced for limited objects, and there 
an equivalent class of connections is defined as that of a gauge 
transform not only by vertical gauges but also horizontal ones. If the 
horizontal graphs are ``trees'', the equivalent class by total gauges is 
the same as that by vertical gauges, however, for general biunitary 
connections, we should limit the gauge choices only to vertical
ones.\par \pha{x} \par

\noindent
{\it Proof.} \par
First assume that $\alpha$ and $\beta$ are equivalent up to
 gauge choice for the vertical edges. Now $\alpha$ and
$\beta$ are on the common four graphs, namely 
${\cal S}_1={\cal S}_2={\cal S}$, 
${\cal T}_1={\cal T}_2={\cal T}$. 
From the assumption, we have two unitary matrices 
$u_{\cal S}$, $u_{\cal T}$ corresponding to the graphs 
${\cal S}$, ${\cal T}$ respectively, such that
$$ u_{\cal S}^* \alpha u_{\cal T} = \beta, $$
where $\alpha$ and $\beta$ represent the matrices 
corresponding to the connections. Now we define the 
isomorphism $\Phi$ from 
${}_K X^{\alpha}_L$ to ${}_K X^{\beta}_L$ as follows;

\begin{eqnarray*}
 &x&=
\left(
\thinlines
\unitlength 0.5mm
\begin{picture}(110,10)(0,5)
\put(10,7){$*_{\cal K}$}
\put(14,10){\line(1,0){36}}
\put(50,10){\circle*{1}}
\put(50,10){\line(0,-1){10}}
\put(50,0){\circle*{1}}
\put(55,0){,}
\put(60,-3){$*_{\cal L}$}
\put(64,0){\line(1,0){36}}
\put(100,0){\circle*{1}}
\end{picture}
\right)    \in   {}_K X^{\alpha}_L, \quad |x|=n \\
&{\downarrow}&   \\
& \Phi(x) &=
\left\{  \begin{array}{c}
({\rm id}^{(n)} \cdot u_{\cal S})
\left(
\thinlines
\unitlength 0.5mm
\begin{picture}(110,10)(0,5)
\put(10,7){$*_{\cal K}$}
\put(14,10){\line(1,0){36}}
\put(50,10){\circle*{1}}
\put(50,10){\line(0,-1){10}}
\put(51,4){$\xi$}
\put(50,0){\circle*{1}}
\put(55,0){,}
\put(60,-3){$*_{\cal L}$}
\put(64,0){\line(1,0){36}}
\put(100,0){\circle*{1}}
\end{picture}
\right),    \quad \mbox{if $n$ is even},  \\ \\
({\rm id}^{(n)} \cdot u_{\cal T})
\left(
\thinlines
\unitlength 0.5mm
\begin{picture}(110,10)(0,5)
\put(10,7){$*_{\cal K}$}
\put(14,10){\line(1,0){36}}
\put(50,10){\circle*{1}}
\put(50,10){\line(0,-1){10}}
\put(51,4){$\xi$}
\put(50,0){\circle*{1}}
\put(55,0){,}
\put(60,-3){$*_{\cal L}$}
\put(64,0){\line(1,0){36}}
\put(100,0){\circle*{1}}
\end{picture}
\right),    \quad \mbox{if $n$ is odd},
\end{array}  \right.   \\
&&
\in {}_K X^{\beta}_L,
\end{eqnarray*}
where id$^{(n)}$ represents the identity of 
String$^{(n)}_{*_{\cal K}}{\cal K}$, and id$^{(n)} \cdot
u_{\cal S}$ is the concatenation, regarding $u_{\cal S}$
 as an element of 
$\bigoplus_{p \in V_0}{\rm String}^{(1)}_{p}{\cal S}$,
and $u_{\cal T}$ is as well. Note that this map changes
only the vertical part of the elements of bimodule.  Now we
check that $\Phi$ is a well-defined linear map, i.e., does not
depend on the length of the expression of $x$. Here we 
assume $n$ is even. We have
\begin{eqnarray*}
\Phi(x)&=&
\sum_{\eta}
u^{\eta,\xi}_{\cal S}
\left(
\thinlines
\unitlength 0.5mm
\begin{picture}(110,10)(0,5)
\put(10,7){$*_{\cal K}$}
\put(14,10){\line(1,0){36}}
\put(50,10){\circle*{1}}
\put(50,10){\line(0,-1){10}}
\put(51,4){$\eta$}
\put(50,0){\circle*{1}}
\put(55,0){,}
\put(60,-3){$*_{\cal L}$}
\put(64,0){\line(1,0){36}}
\put(100,0){\circle*{1}}
\end{picture}
\right)   \\
&=&
\sum_{\eta}
u^{\eta,\xi}_{\cal S}
\sum_{\sigma}
\left(
\thinlines
\unitlength 0.5mm
\begin{picture}(130,10)(0,5)
\put(10,7){$*_{\cal K}$}
\put(14,10){\line(1,0){36}}
\put(50,10){\circle*{1}}
\put(50,10){\line(0,-1){10}}
\put(45,4){$\eta$}
\put(50,0){\circle*{1}}
\put(50,0){\line(1,0){10}}
\put(54,1){$\sigma$}
\put(60,0){\circle*{1}}
\put(65,0){,}
\put(70,-3){$*_{\cal L}$}
\put(74,0){\line(1,0){36}}
\put(110,0){\circle*{1}}
\put(110,0){\line(1,0){10}}
\put(114,1){$\sigma$}
\put(120,0){\circle*{1}}
\end{picture}
\right)   \\
&=&
\sum_{\eta}
u^{\eta,\xi}_{\cal S}
\sum_{\sigma, \sigma', \eta'}
\thinlines
\unitlength 0.5mm
\begin{picture}(30,20)
\multiput(11,-2)(0,10){2}{\line(1,0){8}}
\multiput(10,7)(10,0){2}{\line(0,-1){8}}
\multiput(10,-2)(10,0){2}{\circle*{1}}
\multiput(10,8)(10,0){2}{\circle*{1}}
\put(15,13){\makebox(0,0){$\sigma'$}}
\put(15,-7){\makebox(0,0){$\sigma$}}
\put(5,2.5){$\eta$}
\put(22,2.5){$\eta'$}
\put(15,3){\makebox(0,0){$\beta$}}
\end{picture}
\left(
\thinlines
\unitlength 0.5mm
\begin{picture}(130,10)(0,5)
\put(10,7){$*_{\cal K}$}
\put(14,10){\line(1,0){36}}
\put(50,10){\circle*{1}}
\put(60,10){\line(0,-1){10}}
\put(60,10){\circle*{1}}
\put(61,4){$\eta'$}
\put(60,0){\circle*{1}}
\put(50,10){\line(1,0){10}}
\put(54,11){$\sigma'$}
\put(65,0){,}
\put(70,-3){$*_{\cal L}$}
\put(74,0){\line(1,0){36}}
\put(110,0){\circle*{1}}
\put(110,0){\line(1,0){10}}
\put(114,1){$\sigma$}
\put(120,0){\circle*{1}}
\end{picture}
\right),
\end{eqnarray*}
where $u^{\eta,\xi}_{\cal S}$ denotes the $\eta$-$\xi$ entry of the
matrix $u_{\cal S}$.
On the other hand, we have
\begin{eqnarray*}
\Phi(x)
&=&
 \Phi(
\sum_{\sigma}
\left(
\thinlines
\unitlength 0.5mm
\begin{picture}(130,10)(0,5)
\put(10,7){$*_{\cal K}$}
\put(14,10){\line(1,0){36}}
\put(50,10){\circle*{1}}
\put(50,10){\line(0,-1){10}}
\put(45,4){$\xi$}
\put(50,0){\circle*{1}}
\put(50,0){\line(1,0){10}}
\put(54,1){$\sigma$}
\put(60,0){\circle*{1}}
\put(65,0){,}
\put(70,-3){$*_{\cal L}$}
\put(74,0){\line(1,0){36}}
\put(110,0){\circle*{1}}
\put(110,0){\line(1,0){10}}
\put(114,1){$\sigma$}
\put(120,0){\circle*{1}}
\end{picture}
\right)
)   \\
&=&
{\rm id}^{(n+1)} \cdot u_{\cal T}
\sum_{\sigma, \sigma', \xi'}
\thinlines
\unitlength 0.5mm
\begin{picture}(30,20)
\multiput(11,-2)(0,10){2}{\line(1,0){8}}
\multiput(10,7)(10,0){2}{\line(0,-1){8}}
\multiput(10,-2)(10,0){2}{\circle*{1}}
\multiput(10,8)(10,0){2}{\circle*{1}}
\put(15,13){\makebox(0,0){$\sigma'$}}
\put(15,-7){\makebox(0,0){$\sigma$}}
\put(5,2.5){$\xi$}
\put(22,2.5){$\xi'$}
\put(15,3){\makebox(0,0){$\alpha$}}
\end{picture}
\left(
\thinlines
\unitlength 0.5mm
\begin{picture}(130,10)(0,5)
\put(10,7){$*_{\cal K}$}
\put(14,10){\line(1,0){36}}
\put(50,10){\circle*{1}}
\put(60,10){\circle*{1}}
\put(60,10){\line(0,-1){10}}
\put(61,4){$\eta'$}
\put(60,0){\circle*{1}}
\put(50,10){\line(1,0){10}}
\put(54,11){$\sigma'$}
\put(65,0){,}
\put(70,-3){$*_{\cal L}$}
\put(74,0){\line(1,0){36}}
\put(110,0){\circle*{1}}
\put(110,0){\line(1,0){10}}
\put(114,1){$\sigma$}
\put(120,0){\circle*{1}}
\end{picture}
\right)   \\
&=&
\sum_{\sigma, \sigma', \xi'}
\thinlines
\unitlength 0.5mm
\begin{picture}(30,20)
\multiput(11,-2)(0,10){2}{\line(1,0){8}}
\multiput(10,7)(10,0){2}{\line(0,-1){8}}
\multiput(10,-2)(10,0){2}{\circle*{1}}
\multiput(10,8)(10,0){2}{\circle*{1}}
\put(15,13){\makebox(0,0){$\sigma'$}}
\put(15,-7){\makebox(0,0){$\sigma$}}
\put(5,2.5){$\xi$}
\put(22,2.5){$\xi'$}
\put(15,3){\makebox(0,0){$\alpha$}}
\end{picture}
\sum_{\eta'}
u^{\eta',\xi'}_{\cal T}
\left(
\thinlines
\unitlength 0.5mm
\begin{picture}(130,10)(0,5)
\put(10,7){$*_{\cal K}$}
\put(14,10){\line(1,0){36}}
\put(50,10){\circle*{1}}
\put(60,10){\line(0,-1){10}}
\put(60,10){\circle*{1}}
\put(61,4){$\xi'$}
\put(60,0){\circle*{1}}
\put(50,10){\line(1,0){10}}
\put(54,11){$\sigma'$}
\put(65,0){,}
\put(70,-3){$*_{\cal L}$}
\put(74,0){\line(1,0){36}}
\put(110,0){\circle*{1}}
\put(110,0){\line(1,0){10}}
\put(114,1){$\sigma$}
\put(120,0){\circle*{1}}
\end{picture}
\right).
\end{eqnarray*}
By $u_{\cal S}^* \alpha u_{\cal T} = \beta $, the above two 
expressions of $\Phi(x)$ coincide. When $n$ is odd, it follows
from the same argument. Therefore, $\Phi$ is a
well-defined linear map. \par
Here, $\Phi$ is obviously a right $L$-homomorphism,
and, since id$\cdot u_{\cal S}$ (resp. id$\cdot u_{\cal T}$)
of any length commutes with the element of $K$ of the
same length, $\Phi$ is a left $K$-homomorphism, too. 
 Since $u_{\cal S}$ and $u_{\cal T}$ are unitaries, $\Phi$ is an
isomorphism. Then, we have
$ {}_K X^{\alpha}_L \cong {}_K X^{\beta}_L. $

Next we prove the converse. Assume 
$ {}_K X^{\alpha}_L \cong {}_K X^{\beta}_L. $
Then we have a partial isometry
$$ u \in {\rm End}({}_K X^{\alpha}_L \oplus {}_K X^{\beta}_L )
  = {\rm End}({}_K X^{\alpha + \beta}_L ) $$
such that 
$$ u: {}_K X^{\alpha}_L \stackrel{\sim}{\longrightarrow}
 {}_K X^{\beta}_L,   \quad
 uu^*+u^*u={\rm id}.$$

Our aim is to prove ${\cal S}_1={\cal S}_2$, ${\cal T}_1
={\cal T}_2$ and construct a gauge transform between
$\alpha$ and $\beta$ from $u$.

\begin{claim}
Consider a connection $\gamma$ with four graphs
as below.
\begin{center}
\thinlines
\unitlength 1.0mm
\begin{picture}(30,20)(0,-7)
\multiput(11,-2)(0,10){2}{\line(1,0){8}}
\multiput(10,7)(10,0){2}{\line(0,-1){8}}
\multiput(10,-2)(10,0){2}{\circle*{1}}
\multiput(10,8)(10,0){2}{\circle*{1}}
\put(6,12){\makebox(0,0){$V_0$}}
\put(24,12){\makebox(0,0){$V_1$}}
\put(6,-6){\makebox(0,0){$V_2$}}
\put(24,-6){\makebox(0,0){$V_3$}}
\put(15,12){\makebox(0,0){${\cal K}$}}
\put(15,-6){\makebox(0,0){$\cal L$}}
\put(6,3){\makebox(0,0){${\cal S}$}}
\put(24,3){\makebox(0,0){${\cal T}$}}
\put(15,3){\makebox(0,0){$\gamma$}}
\end{picture} ,
\end{center}
and three AFD II$_1$ factors as in the beginning of this 
section. Then we have  
$$ {\rm End}({}_K X^{\gamma}_L)=K' \cap {\tilde L}, $$
where the embedding of $K \subset {\td L}$ is given by $\gamma$.
\end{claim}
{\it Proof}  \par
First we have 
\begin{eqnarray*}
 {\rm End}({}_K X^{\gamma}_L)
&=&(\mbox{the left action of $K$ on $X^{\gamma}$})' \\
&& \cap (\mbox{the right action of $L$ on 
$X^{\gamma}$})'.
\end{eqnarray*}
We have a natural left action of ${\tilde L}$ on
 ${}_K X^{\gamma}_L$. Now we prove 
\begin{eqnarray*}
 (\mbox{the right action of $L$ on $X^{\gamma}$})' 
 \phantom{XXXXXXXXXX} \\
\phantom{XXXXXXX} 
=(\mbox{the left action of ${\tilde L}$ on 
$X^{\gamma}$}). \phantom{XXXXXXX} (\spadesuit)
\end{eqnarray*}
Obviously we have the inclusion $\subset$, so we prove
the equality by comparing dimensions of $X^{\gamma}$
 as modules of both algebras.
Take a vertex $x$ on ${\cal L}$ and consider projections as 
below;
\begin{eqnarray*}
p &=& 
\left(
\thinlines
\unitlength 0.5mm
\begin{picture}(110,10)(0,5)
\multiput(7,7)(50,0){2}{$*_{\cal K}$}
\multiput(10,0)(50,0){2}{\line(1,0){40}}
\multiput(10,0)(50,0){2}{\circle*{1}}
\multiput(10,6)(50,0){2}{\line(0,-1){6}}
\multiput(50,0)(50,0){2}{\circle*{1}}
\multiput(48,-4.5)(50,0){2}{$x$}
\put(55,0){,}
\end{picture} 
\right)  \in {\tilde L}   \\
q &=& 
\left(
\thinlines
\unitlength 0.5mm
\begin{picture}(110,10)(0,5)
\multiput(10,-3)(50,0){2}{$*_{\cal L}$}
\multiput(14,0)(50,0){2}{\line(1,0){36}}
\multiput(50,0)(50,0){2}{\circle*{1}}
\multiput(48,-4.5)(50,0){2}{$x$}
\put(55,0){,}
\end{picture}  
\right)    \in L.
\end{eqnarray*}
We see that $p {\tilde L}p$ consists of the strings such as
$$p \cdot
\left(
\thinlines
\unitlength 0.5mm
\begin{picture}(110,10)(0,5)
\multiput(10,0)(50,0){2}{\circle*{1}}
\multiput(8,-5)(50,0){2}{$x$}
\multiput(10,0)(50,0){2}{\line(1,0){40}}
\multiput(50,0)(50,0){2}{\circle*{1}}
\put(55,0){,}
\end{picture}  
\right),
$$
where $\cdot$ means the concatenation.  Namely, $p {\tilde L}p$
 essentially consists of the strings of ${\cal L}$ with the initial
 vertex $x$.
It is the case for $p X^{\gamma} q$ and $qLq$ by similar
argument, thus we have
$$ {\rm dim}_{p {\tilde L}p}(p X^{\gamma} q) =1. $$
On the other hand, we have 
$${\rm dim}_{p {\tilde L}p}(p X^{\gamma} q)
= \frac{{\rm tr}_L q}{{\rm tr}_{\tilde L} p}
 {\rm dim}_{\tilde L} X^{\gamma},$$
then we have 
$$ {\rm dim}_{\tilde L} X^{\gamma}
=\frac{{\rm tr}_{\tilde L} p}{{\rm tr}_L q}. $$
By the same argument, we have 
$$ {\rm dim}(p X^{\gamma} q)_{qLq}=1 $$
and
$$ {\rm dim}X^{\gamma}_L 
= \frac{{\rm tr}_L q}{{\rm tr}_{\tilde L} p}. $$
Thus, we have 
$$ {\rm dim}_{\tilde L} X^{\gamma}
=\frac{1}{{\rm dim}X^{\gamma}_L }
={\rm dim}X^{\gamma}_{L'}  $$
and the equality in $(\spadesuit)$ holds. \\
\qed

By applying this claim to $\alpha + \gbs$,
we see that the partial isometry $u$ is in $K' \cap {\tilde L}$
and the map $X^{\alpha} \longrightarrow X^{\beta}$ is
given by the natural left action of ${\tilde L}$ on 
$X^{\alpha}$. To construct the gauge matrices which 
transfer $\alpha$ to $\beta$, we use the {\it compactness
argument of Ocneanu.} (\cite{Oc2}, \cite[Section 11.4]{EK}) 
We introduce some necessary notions and facts. \\ \phantom{x} \\
{\bf Definition 3} (Flat element, Flat field, Ocneanu \cite{Oc2}, \cite{EK}) \\
 \bigskip
Consider a connection $\gamma$ on the four graphs as in the previous claim, 
and three AFD II$_1$ factors $K$, $L$, and ${\td L}$  
as at the beginning of this section.
Take an element 
 $\xi \in \bigoplus_{p \in V_0}
{\rm String}^{(1)}_p {\cal S}$.
It is called a {\it flat element} if
$$
\thinlines
\unitlength 0.5mm
\begin{picture}(35,10)(0,5)
\put(10,10){\circle*{1}}
\put(10,10){\line(0,-1){10}}
\put(5,4){$\xi$}
\put(10,0){\circle*{1}}
\put(10,0){\line(1,0){15}}
\put(25,0){\circle*{1}}
\put(15,1){id$^{(2l)}$}
\end{picture}
=
\thinlines
\unitlength 0.5mm
\begin{picture}(35,10)(0,5)
\put(10,10){\circle*{1}}
\put(25,10){\line(0,-1){10}}
\put(26,4){$\xi$}
\put(25,0){\circle*{1}}
\put(10,10){\line(1,0){15}}
\put(25,10){\circle*{1}}
\put(15,11){id$^{(2l)}$}
\end{picture} , \quad  l \in {\bf N},
$$
under the identification by the connection $\gamma$,
where id$^{(2l)}$ denotes the string 
$\sum_{|\sigma|=2l}(\sigma, \sigma)$ on the graph
${\cal K}$ (resp. ${\cal L}$). We use
this notation often hereafter under similar conditions. \par
It is known that, for a flat element $\xi$, 
there is the element 
$\eta \in \bigoplus_{p \in V_1} {\rm String}_p^{(1)}{\cal T}$
such that 
$$
\thinlines
\unitlength 0.5mm
\begin{picture}(30,10)(0,5)
\put(10,10){\circle*{1}}
\put(10,10){\line(0,-1){10}}
\put(5,4){$\xi$}
\put(10,0){\circle*{1}}
\put(10,0){\line(1,0){10}}
\put(20,0){\circle*{1}}
\put(15,1){id$^{(1)}$}
\end{picture}
=
\thinlines
\unitlength 0.5mm
\begin{picture}(30,10)(0,5)
\put(10,10){\circle*{1}}
\put(20,10){\line(0,-1){10}}
\put(21,4){$\eta$}
\put(20,0){\circle*{1}}
\put(10,10){\line(1,0){10}}
\put(20,10){\circle*{1}}
\put(15,11){id$^{(1)}$}
\end{picture},
$$
and
$$ 
\thinlines
\unitlength 0.5mm
\begin{picture}(35,10)(0,5)
\put(10,10){\circle*{1}}
\put(10,10){\line(0,-1){10}}
\put(5,4){$\eta$}
\put(10,0){\circle*{1}}
\put(10,0){\line(1,0){15}}
\put(25,0){\circle*{1}}
\put(15,1){id$^{(2l)}$}
\end{picture}
=
\thinlines
\unitlength 0.5mm
\begin{picture}(35,10)(0,5)
\put(10,10){\circle*{1}}
\put(25,10){\line(0,-1){10}}
\put(26,4){$\eta$}
\put(25,0){\circle*{1}}
\put(10,10){\line(1,0){15}}
\put(25,10){\circle*{1}}
\put(15,11){id$^{(2l)}$}
\end{picture}.
$$
by the connection $\gamma$.  We call $\eta$ a flat element,
too. This ``couple'' of flat elements represents an element of the string 
algebra with identification by $\gamma$, namely,
$$
\thinlines
\unitlength 0.5mm
\begin{picture}(55,10)(0,5)
\put(10,7){$*_{\cal K}$}
\put(14,10){\line(1,0){36}}
\put(24,11){$id^{(k)}$}
\put(50,10){\circle*{1}}
\put(50,10){\line(0,-1){10}}
\put(51,4){$\xi$}
\put(50,0){\circle*{1}}
\end{picture}
=
\thinlines
\unitlength 0.5mm
\begin{picture}(55,10)(0,5)
\put(10,7){$*_{\cal K}$}
\put(14,10){\line(1,0){36}}
\put(24,11){$id^{(l)}$}
\put(50,10){\circle*{1}}
\put(50,10){\line(0,-1){10}}
\put(51,4){$\eta$}
\put(50,0){\circle*{1}}
\end{picture}
$$
for any sufficiently large $k$:even and $l$:odd, that is, large 
enough that the set of the end points of $id^{(k)}$ (resp. $id^{(l)}$) 
coincides with $V_{0}$ (resp. $V_{1}$). Now we define $z$ to be a function  
on $V_{0} \cup V_{1}$ such that 
$z(p) \in  {\rm String}_p^{(1)} {\cal S}$ (resp. ${\rm 
String}_p^{(1)} {\cal T}$) for $p \in  V_{0}$ 
(resp. $V_{1}$) and $\oplus_{p \in V_{0}}z(p)=\xi$ (resp. $\oplus_{p \in V_{1}}z(p)
=\eta$), and call it {\it flat field}.
 Let $V_{n}'$ to be a proper subset of $V_{n}$, where $n=0, 1$, 
 and put $\xi_{0} = \oplus_{p \in V_{0}'}z(p)$ (resp. 
 $\eta_{0} = \oplus_{o \in V_{1}'}z(p)$). It is known that
 $$
 \thinlines
\unitlength 0.5mm
\begin{picture}(55,10)(0,5)
\put(10,10){\circle*{1}}
\put(10,10){\line(0,-1){10}}
\put(11,4){$\xi_0$(resp. $\eta_0$)}
\put(10,0){\circle*{1}}
\put(10,0){\line(1,0){35}}
\put(45,0){\circle*{1}}
\put(27,-8){id$^{(2j)}$}
\end{picture}
= 
\thinlines
\unitlength 0.5mm
\begin{picture}(65,10)(0,5)
\put(10,10){\circle*{1}}
\put(45,10){\line(0,-1){10}}
\put(46,4){$\xi$(resp. $\eta$),}
\put(45,0){\circle*{1}}
\put(10,10){\line(1,0){35}}
\put(45,10){\circle*{1}}
\put(27,11){id$^{(2j)}$}
\put(45,0){\circle*{1}}
\end{picture}
$$
\pha{x} \par \noindent
for sufficiently large $j$. We call such elements as $\xi_{0}$ and 
$\eta_{0}$ {\it flat}, too, though they are not flat elements by the 
definition above.
\begin{thm}{\rm (Ocneanu \cite{Oc2},\cite{EK})}  \\
Let $K \subset {\tilde L}$ be the AFD II$_1$ subfactor 
constructed from the connection $\gamma$. Then,
$$ K' \cap {\tilde L} = \{ \mbox{flat field} \}.  $$
The correspondence of elements as follows; \\
Take a flat field $z$ and let 
$\xi=\bigoplus_{p \in V_0} z(p)$, then
$$
\thinlines
\unitlength 0.5mm
\begin{picture}(35,10)(0,5)
\put(10,7){$*_{\cal K}$}
\put(25,10){\line(0,-1){10}}
\put(26,4){$\xi$}
\put(25,0){\circle*{1}}
\put(14,10){\line(1,0){11}}
\put(25,10){\circle*{1}}
\put(15,11){id$^{(2k)}$}
\end{picture} 
\in K' \cap {\tilde L},
$$
and conversely, for $x \in K' \cap {\tilde L}$, it turns out
that $x$ is written as
$$
x= \thinlines
\unitlength 0.7mm
\begin{picture}(20,10)(0,5)
\put(8,7){$*_{\cal K}$}
\put(10,7){\line(0,-1){7}}
\put(11,1){$z(*)$}
\put(10,0){\circle*{1}}
\end{picture} 
\in \mbox{String}_{*}^{(1)}{\cal S} \subset {\tilde L}
$$
for some flat field $z$.
\end{thm}

This theorem is proved by {\it the compactness argument
of Ocneanu}, see \cite{Oc2} and \cite{EK}. (Generally, the
length of flat field/element can be arbitrary.) \bigskip

Now we continue the proof of Theorem 3, using the above notions.
Let $\gamma= \alpha + \beta$ and ${\cal S}= {\cal S}_1 \cup
{\cal S}_2$, ${\cal T}= {\cal T}_1 \cup {\cal T}_2$. By the
above theorem, we consider the partial isometry $u \in K' \cap {\td 
L}$ which gives the isometry $X^{\a} \longrightarrow X^{\gbs}$ as a flat 
field for the connection $\gamma$.
Take $p \in V_0$ and $q \in V_2$ so that they are connected 
in ${\cal S}$. Since $uu^*+u^*u=1$, we have
$$u(p,q)u(p,q)^*+u(p,q)^*u(p,q)=1$$
in the algebra String$_{(p,q)}^{(1)}{\cal S}$ =span$\{(\sigma, \rho) | 
\quad |\sigma|=|\rho|=1, \quad s(\sigma)=s(\rho)=p,
r(\sigma)=r(\rho)=q \} $, where $u(p,q) \in 
\mbox{String}_{(p,q)}^{(1)}{\cal S}$ such that $\oplus_{q \in 
V_{2}}u(p,q) = u(p)$.  Take an element of $X^{\alpha}$
$$x=
\left(
\thinlines
\unitlength 0.5mm
\begin{picture}(110,10)(0,5)
\put(10,7){$*_{\cal K}$}
\put(14,10){\line(1,0){36}}
\put(31,13){$\zeta$}
\put(50,10){\circle*{1}}
\put(48,13){$p$}
\put(50,10){\line(0,-1){10}}
\put(51,4){$\xi$}
\put(50,0){\circle*{1}}
\put(48,-7){$q$}
\put(55,0){,}
\put(60,-3){$*_{\cal L}$}
\put(64,0){\line(1,0){36}}
\put(81,1){$\varepsilon$}
\put(100,0){\circle*{1}}
\put(98,-7){$q$}
\end{picture}
\right), \qquad \xi \in {\cal S}_1.
 $$ 
From the definition of $u$, we have $ux \in X^{\beta}$, then
\begin{eqnarray*}
(id \cdot u(p,q)) \cdot x =  
\thinlines
\unitlength 0.5mm
\begin{picture}(65,10)(10,5)
\put(10,7){$*_{\cal K}$}
\put(14,10){\line(1,0){36}}
\put(31,3){id}
\put(50,10){\circle*{1}}
\put(50,10){\line(0,-1){10}}
\put(51,4){$u(p,q)$}
\put(50,0){\circle*{1}} 
\end{picture} 
\cdot
\left(
\thinlines
\unitlength 0.5mm
\begin{picture}(110,10)(0,5)
\put(10,7){$*_{\cal K}$}
\put(14,10){\line(1,0){36}}
\put(31,13){$\zeta$}
\put(50,10){\circle*{1}}
\put(48,13){$p$}
\put(50,10){\line(0,-1){10}}
\put(51,4){$\xi$}
\put(50,0){\circle*{1}}
\put(48,-7){$q$}
\put(55,0){,}
\put(60,-3){$*_{\cal L}$}
\put(64,0){\line(1,0){36}}
\put(81,1){$\varepsilon$}
\put(100,0){\circle*{1}}
\put(98,-7){$q$}
\end{picture}
\right)   \\  \phantom{X} \\
= 
\sum_{\eta} u(p,q)^{\eta,\xi} 
\left(
\thinlines
\unitlength 0.5mm
\begin{picture}(110,10)(0,5)
\put(10,7){$*_{\cal K}$}
\put(14,10){\line(1,0){36}}
\put(31,13){$\zeta$}
\put(50,10){\circle*{1}}
\put(48,13){$p$}
\put(50,10){\line(0,-1){10}}
\put(51,4){$\eta$}
\put(50,0){\circle*{1}}
\put(48,-7){$q$}
\put(55,0){,}
\put(60,-3){$*_{\cal L}$}
\put(64,0){\line(1,0){36}}
\put(81,1){$\varepsilon$}
\put(100,0){\circle*{1}}
\put(98,-7){$q$}
\end{picture}
\right) 
\in X^{\beta},  \quad  u(p,q)^{\eta,\xi} \in {\bf C}.  
\end{eqnarray*}
Note that $\eta \in {\cal S}_2$ if $u(p,q)^{\eta,\xi} \not=0$. 
Since $u$ gives an isometry of $X^{\alpha}$ and $X^{\beta}$,
\begin{eqnarray*}
{\rm dim}\;{\rm span} \{
\left(
\thinlines
\unitlength 0.5mm
\begin{picture}(110,10)(0,5)
\put(10,7){$*_{\cal K}$}
\put(14,10){\line(1,0){36}}
\put(31,13){$\zeta$}
\put(50,10){\circle*{1}}
\put(48,13){$p$}
\put(50,10){\line(0,-1){10}}
\put(51,4){$\xi$}
\put(50,0){\circle*{1}}
\put(48,-7){$q$}
\put(55,0){,}
\put(60,-3){$*_{\cal L}$}
\put(64,0){\line(1,0){36}}
\put(81,1){$\varepsilon$}
\put(100,0){\circle*{1}}
\put(98,-7){$q$}
\end{picture}
\right), \quad \xi \in {\cal S}_1 \} \phantom{XXXXXX} \\
\phantom{X} \\
=  
{\rm dim}\;{\rm span} \{
\left(
\thinlines
\unitlength 0.5mm
\begin{picture}(110,10)(0,5)
\put(10,7){$*_{\cal K}$}
\put(14,10){\line(1,0){36}}
\put(31,13){$\zeta$}
\put(50,10){\circle*{1}}
\put(48,13){$p$}
\put(50,10){\line(0,-1){10}}
\put(51,4){$\eta$}
\put(50,0){\circle*{1}}
\put(48,-7){$q$}
\put(55,0){,}
\put(60,-3){$*_{\cal L}$}
\put(64,0){\line(1,0){36}}
\put(81,1){$\varepsilon$}
\put(100,0){\circle*{1}}
\put(98,-7){$q$}
\end{picture}
\right), \quad \eta \in {\cal S}_2 \}
\end{eqnarray*}
for each $\zeta$ and $\varepsilon$. This means
$$
{\sharp} \left\{
\thinlines
\unitlength 0.5mm
\begin{picture}(20,10)(0,5)
\put(10,10){\circle*{1}}
\put(8,13){$p$}
\put(10,10){\line(0,-1){10}}
\put(11,4){$\xi$}
\put(10,0){\circle*{1}} 
\put(8,-7){$q$} 
\end{picture}
\in {\cal S}_1
\right\}
=
{\sharp}  \left\{
\thinlines
\unitlength 0.5mm
\begin{picture}(20,10)(0,5)
\put(10,10){\circle*{1}}
\put(8,13){$p$}
\put(10,10){\line(0,-1){10}}
\put(11,4){$\eta$}
\put(10,0){\circle*{1}} 
\put(8,-7){$q$} 
\end{picture}
\in {\cal S}_2
\right\}. $$
By seeing all the possible pairs of vertices $p$ and $q$, we have
$${\cal S}_1={\cal S}_2.$$
By the same discussion, we have also
$${\cal T}_1={\cal T}_2,$$ 
then we know that $\alpha$ and $\beta$ are on the same
four graphs. \par
Now we see that $u(p,q)$ gives the gauge matrix for the 
edges which connect $p$ and $q$. Let $u_{\cal S}$ and 
$u_{\cal T}$ be ``stable" flat elements on ${\cal S}$
and {\cal T} corresponding to the flat field $u$. Since the 
isomorphism $x \in X^{\alpha} \rightarrow u \cdot x \in
X^{\beta}$ is well-defined, from the same deformation as
we proved the well-definedness of $\Phi$ in the first half
proof of our main statement here, 
$$ u_{\cal S}^* \alpha u_{\cal T} = \beta $$
follows. Under the identification of ${\cal S}_1={\cal S}_2$
 and ${\cal T}_1={\cal T}_2$, $u_{\cal S}$ and $ u_{\cal T}$
are considered as unitary matrices corresponding to the
gauge transform action of $\alpha$ and $\beta$. Thus we have
$$ \alpha \cong \beta \qquad \mbox{up to vertical gauge choice}. $$
\begin{cor}  \par
$ {}_K X^{\alpha}_L$ is irreducible if and only if $\alpha$
is indecomposable.
\end{cor}
{\it Proof} \\
Assume $\alpha$ is decomposable, i.e. there exist gauge unitaries
$u_{\cal S}$, $u_{\cal T}$ and connections $\beta$, $\gamma$
such that 
$$ u_{\cal S}^* \alpha u_{\cal T}=\beta+\gamma. $$
Then we have
$$ {}_K X^{u_{\cal S}^* \alpha u_{\cal T}}_L \cong
{}_K X^{\beta}_L \oplus {}_K X^{\gamma}_L 
\cong {}_K X^{\alpha}_L, $$
namely, ${}_K X^{\alpha}_L$ is reducible. \\
Conversely, assume ${}_K X^{\alpha}_L$ is reducible. 
then we have bimodules ${}_K Y_L$ and ${}_K Z_L$
such that 
$$ {}_K X^{\alpha}_L \cong {}_K Y_L \oplus {}_K Z_L $$ 
and a projection
\begeq
 p \in {\rm End}({}_K X^{\alpha}_L)=K' \cap {\tilde L} \quad
\mbox{with  }  p: {}_K X^{\alpha}_L \longrightarrow {}_K Y_L. 
 \endeq
Along the same argument as in the proof of the previous 
theorem, we consider $p$ as a flat field and make the
projections $p_{\cal S}$ and $p_{\cal T}$ which project
 elements of ${}_K X^{\alpha}_L$ to ${}_K Y_L$ at the 
finite level, and they act as the ``projections" of the
connection matrix, and we have a ``sub connection" of 
$\alpha$
$$ \beta= p_{\cal S}^* \alpha p_{\cal T} $$
so that
$$    {}_K X^{\beta}_L \cong {}_K Y_L,$$
thus, $\alpha$ is decomposable.      \\      
\qed

\begin{cor} \par
\label{irred}
Let $\gamma$ be a connection as in Claim 1, i.e.,
\begin{center}
\thinlines
\unitlength 1.0mm
\begin{picture}(30,20)(0,-7)
\multiput(11,-2)(0,10){2}{\line(1,0){8}}
\multiput(10,7)(10,0){2}{\line(0,-1){8}}
\multiput(10,-2)(10,0){2}{\circle*{1}}
\multiput(10,8)(10,0){2}{\circle*{1}}
\put(6,12){\makebox(0,0){$V_0$}}
\put(24,12){\makebox(0,0){$V_1$}}
\put(6,-6){\makebox(0,0){$V_2$}}
\put(24,-6){\makebox(0,0){$V_3$}}
\put(15,12){\makebox(0,0){${\cal K}$}}
\put(15,-6){\makebox(0,0){$\cal L$}}
\put(6,3){\makebox(0,0){${\cal S}$}}
\put(24,3){\makebox(0,0){${\cal T}$}}
\put(15,3){\makebox(0,0){$\gamma$}}
\end{picture} .
\end{center}
If there exists a vertex $p$ to which only one vertical edge
is connected. Then the bimodule ${}_{K} X^{\gamma}_{L}$ is 
irreducible.
\end{cor} 
{\it Proof} \\
Assume $p \in V_{0}$ without missing generality. Let $\xi$ be the only one 
vertical edge in ${\cal S}$ connected to $p$. 
Assume ${}_{K} X^{\gamma}_{L}$ is not irreducible. Then, by the 
above argument, we have connections $\gamma_{1}$ and $\gamma_{2}$ 
with the four graphs
\begin{center}
\thinlines
\unitlength 1.0mm
\begin{picture}(30,20)(0,-7)
\multiput(11,-2)(0,10){2}{\line(1,0){8}}
\multiput(10,7)(10,0){2}{\line(0,-1){8}}
\multiput(10,-2)(10,0){2}{\circle*{1}}
\multiput(10,8)(10,0){2}{\circle*{1}}
\put(6,12){\makebox(0,0){$V_0$}}
\put(24,12){\makebox(0,0){$V_1$}}
\put(6,-6){\makebox(0,0){$V_2$}}
\put(24,-6){\makebox(0,0){$V_3$}}
\put(15,12){\makebox(0,0){${\cal K}$}}
\put(15,-6){\makebox(0,0){$\cal L$}}
\put(6,3){\makebox(0,0){${\cal S}_{1}$}}
\put(24,3){\makebox(0,0){${\cal T}_{1}$}}
\put(15,3){\makebox(0,0){$\gamma_{1}$}}
\end{picture} ,
\thinlines
\unitlength 1.0mm
\begin{picture}(30,20)(0,-7)
\multiput(11,-2)(0,10){2}{\line(1,0){8}}
\multiput(10,7)(10,0){2}{\line(0,-1){8}}
\multiput(10,-2)(10,0){2}{\circle*{1}}
\multiput(10,8)(10,0){2}{\circle*{1}}
\put(6,12){\makebox(0,0){$V_0$}}
\put(24,12){\makebox(0,0){$V_1$}}
\put(6,-6){\makebox(0,0){$V_2$}}
\put(24,-6){\makebox(0,0){$V_3$}}
\put(15,12){\makebox(0,0){${\cal K}$}}
\put(15,-6){\makebox(0,0){$\cal L$}}
\put(6,3){\makebox(0,0){${\cal S}_{2}$}}
\put(24,3){\makebox(0,0){${\cal T}_{2}$}}
\put(15,3){\makebox(0,0){$\gamma_{2}$}}
\end{picture} 
\end{center}
respectively, so that 
$$\gamma \cong \gamma_{1} + \gamma_{2}, \quad 
{\cal S}={\cal S}_{1} \cup {\cal S}_{2}, \quad 
{\cal T}={\cal T}_{1} \cup {\cal T}_{2}. $$
$\xi$ should be contained either in ${\cal S}_{1}$ or ${\cal S}_{2}$.
Assume $\xi \in {\cal S}_{2}$, then no edge in ${\cal S}_{1}$ 
connects to $p$. This contradicts to the unitarity of $\gamma_{1}$. \\
\qed \\
{\it Remark} \\
Corollary 2 is a generalization of Wenzl's Criterion for irreducibility 
of subfactors obtained from a periodic sequence of commuting squares 
(c.f. \cite{W}).


\section{Main theorem for the case of $(5+\sqrt{13})/2$}

In this section, we give a proof for our main theorem 
for the case of index $(5+\sqrt{13})/2$ due to the second named author.

\begin{thm}
A subfactor with principal graph and dual principal graph as in Figure
\ref{13graph} exists.
\end{thm}
From the key lemma, we know that the above theorem follows from the 
next proposition. We define the connection $\sigma$ as
\begin{displaymath}
\sigma
\left(
\thinlines
\unitlength 0.5mm
\begin{picture}(30,20)
\multiput(11,-2)(0,10){2}{\line(1,0){8}}
\multiput(10,7)(10,0){2}{\line(0,-1){8}}
\multiput(10,-2)(10,0){2}{\circle*{1}}
\multiput(10,8)(10,0){2}{\circle*{1}}
\put(6,12){\makebox(0,0){$p$}}
\put(24,12){\makebox(0,0){$q$}}
\put(6,-6){\makebox(0,0){$r$}}
\put(24,-6){\makebox(0,0){$s$}}
\end{picture}
\right)
= \delta_{\sigma (p),r}\delta_{\sigma (q),s},
\end{displaymath}
where $p, q, r, s$ are the vertices on the upper graph in 
Figure \ref{13graph}, and we define $\sigma(\cdot)$ as 
$\sigma(x)=x_\sigma$, $\sigma(x_\sigma)=x_{\sigma^2}$, 
and $x_{\sigma^3}=x$.
Note that, for the vertex $c$ we put $\sigma(c)=c$. 

\begin{prop}
Let $\alpha$ be the unique connection on the four graphs 
consisting of the pair of the graphs appearing in Figure 
\ref{13graph},
and $\sigma$ be the connection defined above.
Then, the following hold. \par \noindent
1) The six connections
$$
{\bf 1}, \; \sigma, \;  \sigma^2, \;  (\a \ab -{\bf 1}), \; 
\sigma(\a \ab -{\bf 1}),
\; \sigma^2(\a \ab -{\bf 1})
$$
are indecomposable and mutually inequivalent. \\
2) The four connections
$$ \a, \; \sigma \a, \; \sigma^2 \a, \; \a \ab \a-2\a
$$
are irreducible and mutually inequivalent. \\
3) 
$$ \sigma(\a \ab -{\bf 1}) \cong (\a \ab -{\bf 1})\sigma^2. $$
\end{prop}
{\it Proof}

 The four graphs of the connection $\a$ are as in Figure \ref{13fourgraphs}. 

\begin{figure}[H]
\begin{center}
\thinlines
\unitlength 1.0mm
\begin{picture}(30,20)(0,-5)
\multiput(11,-2)(0,10){2}{\line(1,0){8}}
\multiput(10,7)(10,0){2}{\line(0,-1){8}}
\multiput(10,-2)(10,0){2}{\circle*{1}}
\multiput(10,8)(10,0){2}{\circle*{1}}
\put(6,12){\makebox(0,0){$V_0$}}
\put(24,12){\makebox(0,0){$V_1$}}
\put(6,-6){\makebox(0,0){$V_3$}}
\put(24,-6){\makebox(0,0){$V_2$}}
\put(15,12){\makebox(0,0){${\cal G}_0$}}
\put(15,-6){\makebox(0,0){${\cal G}_2$}}
\put(6,3){\makebox(0,0){${\cal G}_3$}}
\put(24,3){\makebox(0,0){${\cal G}_1$}}
\put(15,3){\makebox(0,0){$\a$}}
\end{picture} 
\end{center}
\begin{center}
	\setlength{\unitlength}{0.5mm}
	\begin{picture}(200,120)(-20,-10)
	\multiput(0,100)(30,0){6}{\circle*{2}}
	\multiput(15,75)(30,0){3}{\circle*{2}}
	\put(135,75){\circle*{2}}
	\multiput(0,50)(30,0){4}{\circle*{2}}
	\multiput(15,25)(30,0){3}{\circle*{2}}
	\put(135,25){\circle*{2}}
	\multiput(0,0)(30,0){6}{\circle*{2}}
	\multiput(0,100)(30,0){3}{\line(3,-5){30}}
	\multiput(30,100)(30,0){2}{\line(-3,-5){30}}
	\put(150,100){\line(-3,-5){15}}
   \put(90,100){\line(9,-5){45}}
    \multiput(90,100)(30,0){2}{\line(-9,-5){45}}
   \put(45,75){\line(9,-5){90}}
	\put(135,75){\line(-9,-5){90}}
	\multiput(0,50)(30,0){2}{\line(3,-5){30}}
	\multiput(30,50)(30,0){3}{\line(-3,-5){30}}
	
   	\put(45,25){\line(9,-5){45}}
   	\put(75,25){\line(15,-5){75}}
   	\put(135,25){\line(-9,-5){45}}
	\put(135,25){\line(-3,-5){15}}

    \put(5,100){\makebox(1,1){$*$}}
	\put(35,100){\makebox(1,1){$b$}}
	\put(68,100){\makebox(1,1){\s{b}}}
	\put(98,102){\makebox(1,1){\ss{b}}}
	\put(128,100){\makebox(1,1){\s{*}}}
	\put(158,100){\makebox(1,1){\ss{*}}}
	\put(21,75){\makebox(1,1){$a$}}
	\put(52,75){\makebox(1,1){$c$}}
	\put(83,70){\makebox(1,1){\s{a}}}
	\put(142,75){\makebox(1,1){\ss{a}}}
	\put(5,50){\makebox(1,1){1}}
	\put(35,50){\makebox(1,1){2}}
	\put(65,50){\makebox(1,1){3}}
	\put(100,50){\makebox(1,1){4}}
	\put(21,25){\makebox(1,1){$a$}}
	\put(54,25){\makebox(1,1){$c$}}
	\put(84,25){\makebox(1,1){\s{a}}}
	\put(144,25){\makebox(1,1){\ss{a}}}
	\put(5,0){\makebox(1,1){$*$}}
	\put(35,0){\makebox(1,1){$b$}}
	\put(68,0){\makebox(1,1){\s{b}}}
	\put(99,-5){\makebox(1,1){\ss{b}}}
	\put(128,0){\makebox(1,1){\s{*}}}
	\put(158,0){\makebox(1,1){\ss{*}}}
	\put(-10,87){\makebox(2,2){${\cal G}_0$}}
	\put(-10,62){\makebox(2,2){ ${\cal G}_1$}}
	\put(-10,37){\makebox(2,2){${\cal G}_2$}}
	\put(-10,12){\makebox(2,2){${\cal G}_3$}}
	\put(-15,100){\makebox(2,2){ $V_0$}}
	\put(-15,75){\makebox(2,2){$V_1$}}
	\put(-15,50){\makebox(2,2){$V_2$}}
	\put(-15,25){\makebox(2,2){ $V_3$}}
	\put(-15,0){\makebox(2,2){ $V_0$}}
	\end{picture}
\caption{Four graphs of the connection $\alpha$}
\label{13fourgraphs}
\end{center}
\end{figure}

The Perron-Frobenius weights of the vertices can easily be computed 
as follows,
\begeq
\mu(*)=1, \quad \mu(a)=\mu(\si{a})=\mu(\ssi{a})=\la, \\
\mu(b)=\mu(\si{b})=\mu(\ssi{b})=\lad-1, \quad \mu(c)=\la^{3}-2\la,\\
\mu(1)=1, \quad \mu(2)=\lad-1, \quad \mu(3)= \lad-2, \quad 
\mu(4)=\lad,
\endeq
where $\la=\sqrt{\frac{5+\sqrt{13}}{2}}$.
One can check that the Table \ref{a13} defines a connection $\a$ on 
the four graphs (Figure \ref{13fourgraphs}) which satisfies Ocneanu's 
biunitary conditions, i.e., 
$$
\label{biunit}
\left( \alpha \left( 
\thinlines
\unitlength 0.5mm
\begin{picture}(30,20)
\multiput(11,-2)(0,10){2}{\line(1,0){8}}
\multiput(10,7)(10,0){2}{\line(0,-1){8}}
\multiput(10,-2)(10,0){2}{\circle*{1}}
\multiput(10,8)(10,0){2}{\circle*{1}}
\put(6,12){\makebox(0,0){$p$}}
\put(24,-6){\makebox(0,0){$s$}}
\put(6,3){\makebox(0,0){$\xi$}}
\put(24,3){\makebox(0,0){$\xi'$}}
\put(15,12){\makebox(0,0){$\eta'$}}
\put(15,-6){\makebox(0,0){$\eta$}}
\end{picture} \right) \right)_{\xi \cdot \eta, \xi \cdot \xi'}
 \quad \mbox{  is a unitary matrix for each fixed $p$, $s$, (unitarity)}
$$
and
$$
\thinlines
\unitlength 0.5mm
\begin{picture}(30,20)
\multiput(11,-2)(0,10){2}{\line(1,0){8}}
\multiput(10,7)(10,0){2}{\line(0,-1){8}}
\multiput(10,-2)(10,0){2}{\circle*{1}}
\multiput(10,8)(10,0){2}{\circle*{1}}
\put(6,12){\makebox(0,0){$x$}}
\put(24,12){\makebox(0,0){$z$}}
\put(6,-6){\makebox(0,0){$y$}}
\put(24,-6){\makebox(0,0){$w$}}
\put(6,3){\makebox(0,0){$\xi$}}
\put(24,3){\makebox(0,0){$\xi'$}}
\put(15,12){\makebox(0,0){$\eta'$}}
\put(15,-6){\makebox(0,0){$\eta$}}
\end{picture}
=
\sqrt{\frac{\mu(y) \mu(z)}{\mu(x) \mu(w)}} 
\cdot
\overline{
\thinlines
\unitlength 0.5mm
\begin{picture}(30,20)
\multiput(11,-2)(0,10){2}{\line(1,0){8}}
\multiput(10,7)(10,0){2}{\line(0,-1){8}}
\multiput(10,-2)(10,0){2}{\circle*{1}}
\multiput(10,8)(10,0){2}{\circle*{1}}
\put(6,12){\makebox(0,0){$y$}}
\put(24,12){\makebox(0,0){$w$}}
\put(6,-6){\makebox(0,0){$x$}}
\put(24,-6){\makebox(0,0){$z$}}
\put(6,3){\makebox(0,0){$\td{\xi}$}}
\put(24,3){\makebox(0,0){$\td{\xi'}$}}
\put(15,12){\makebox(0,0){$\eta$}}
\put(15,-6){\makebox(0,0){$\eta'$}}
\end{picture}
}, \qquad \mbox{  (renormalization)}
$$
see \cite{Oc1} and \cite[chapter 10]{EK}.
We see that such a biunitary connection $\alpha$ on these four graphs is 
determined uniquely up to complex conjugate arising from 
the symmetricity 
of the graphs, namely it is essentially unique. The connection 
$\alpha$ is as in Table \ref{a13}. Note
$$
\mbox{$(xy, zw)$-entry in the table}=
\alpha \left(
\thinlines
\unitlength 0.5mm
\begin{picture}(30,20)
\multiput(11,-2)(0,10){2}{\line(1,0){8}}
\multiput(10,7)(10,0){2}{\line(0,-1){8}}
\multiput(10,-2)(10,0){2}{\circle*{1}}
\multiput(10,8)(10,0){2}{\circle*{1}}
\put(6,12){\makebox(0,0){$x$}}
\put(24,12){\makebox(0,0){$z$}}
\put(6,-6){\makebox(0,0){$y$}}
\put(24,-6){\makebox(0,0){$w$}}
\end{picture}
\right),
$$
where we note that, since all the graphs which consist the four 
graphs in Figure \ref{13fourgraphs} are ``tree'', all the edges are 
expressed by the both ends. For example, in the Table \ref{a13} one 
can find 
$$
\mbox{$(*a, a2)$-entry}=
\alpha \left(
\thinlines
\unitlength 0.5mm
\begin{picture}(30,20)
\multiput(11,-2)(0,10){2}{\line(1,0){8}}
\multiput(10,7)(10,0){2}{\line(0,-1){8}}
\multiput(10,-2)(10,0){2}{\circle*{1}}
\multiput(10,8)(10,0){2}{\circle*{1}}
\put(6,12){\makebox(0,0){$*$}}
\put(24,12){\makebox(0,0){$a$}}
\put(6,-6){\makebox(0,0){$a$}}
\put(24,-6){\makebox(0,0){$2$}}
\end{picture}
\right)=1. $$
We also note that blank entries are all 0's, and   
$\rho = \frac{1}{2}(-\sqrt{\lad-4}+i \sqrt{8-\lad})$, 
   $\tau = \frac{1}{2}(-\sqrt{\lad-1}-i\sqrt{5-\lad})$, 
 $|\rho|=|\tau|=1. \quad ({\bar \tau}^3=\rho).$ \\
\begin{table}[h]
\begin{center}
\begin{tabular}[pos]{|c||c||c|c||c||c|c|c|}
\hl
& $a1$ & $a2$ & $c2$ & $c3$ & $c4$ & {\ss a}4 & {\s a}4
 \\ \hline\hline
$*a$ & 1 & 1 & &&&& \\ \hline\hline
$ba$ & 1 & $\frac{-1}{\lad-1}$ & $\frac{\la \sqrt{\lad-2}}{\lad-1}$
 & &&& \\ \hl
$bc$ & & $\frac{\la \sqrt{\lad-2}}{\lad-1}$ & $\frac{1}{\lad-1}$
 & 1 & 1 & & \\ \hline\hline
{\s b}$c$ & && ${\bar \rho}$ & ${\bar \tau}$ & 
$\frac{1}{\sqrt{\lad-1}}$ & $\sqrt{\frac{\lad-2}{\lad-1}}$
 &  \\ \hl
{\s b}{\s a} & &&&& $\sqrt{\frac{\lad-2}{\lad-1}}$ &
 $\frac{-1}{\sqrt{\lad-1}}$ &  \\ \hline\hline
{\ss b}$c$ & && $\rho$ & $\tau$ & $\frac{1}{\sqrt{\lad-1}}$ & & 
$\sqrt{\frac{\lad-2}{\lad-1}}$ \\ \hl
{\ss b}{\s a} & &&&&  $\sqrt{\frac{\lad-2}{\lad-1}}$ && 
$\frac{-1}{\sqrt{\lad-1}}$ \\ \hline\hline
{\s *}{\s a} &&&&&& 1 & \\ \hline\hline
{\ss *}{\ss a} &&&&&&& 1 \\ \hl
\end{tabular}
\end{center}
\caption{Connection $\alpha$}
\label{a13}
\end{table}
Now we display the table of the connection $\tilde{\alpha}$ computed 
by ``renormalization'' in \framebox{Table \ref{ab13}}
for use of the later computations. \\

\begin{table}[h]
\begin{center}
\begin{tabular}[pos]{|c||c|c||c|c|c||c||c|}
\hl
& 1$a$ & 2$a$ & 2$c$ & 3$c$ & 4$c$ & 4\ss{a} & 4\s{a} \\ \hline\hline
$a*$ & $\frac{1}{\la}$ & $\frac{\sqrt{\lad-1}}{\la}$ & &&&& \\ \hl
$ab$ & $\frac{\sqrt{\lad-1}}{\la}$ & $\frac{-1}{\la}$ & 1 & 
&&& \\ \hl \hl
$cb$ & & 1 & $\frac{1}{\la(\lad-2)}$ & $\frac{1}{\sqrt{3}}$ & 
$\frac{\sqrt{\lad-1}}{\lad-2}$ & & \\ \hl
$c$\s{b} & &&  $\frac{\lad-1}{\la(\lad-2)}\rho$  &
 $\frac{1}{\sqrt{3}}\tau$ & $\frac{1}{\lad-2}$ & 1 & \\ \hl
$c$\ss{b} &&& $\frac{\lad-1}{\la(\lad-2)}\bar{\rho}$ & 
$\frac{1}{\sqrt{3}}\bar{\tau}$ & $\frac{1}{\lad-2}$ && 1 \\ 
\hline\hline
\s{a}\s{b} &&&&& 1 & $-1$ & \\ \hl
\s{a}\ss{*} &&&&&&& 1 \\ \hline\hline
\ss{a}\ss{b} &&&&& 1 && $-1$ \\ \hl
\ss{a}\s{*} &&&&&& 1 & \\ \hl
\end{tabular}
\caption{Connection $\tilde{\alpha}$}
\label{ab13}
\end{center}
\end{table}

First we check condition 3), namely we prove
$$ \sigma(\a \ab-{\bf 1}) \cong (\a \ab-{\bf 1}) \sigma^2$$
up to vertical gauge choice. It is enough to show
$$ (\a \ab-{\bf 1}) \cong \sigma(\a \ab-{\bf 1}) \sigma,$$
so now we will prove this equivalence. \\ 
First we compute the connection $\a \ab$. The four graphs on 
which the connection $\a \ab$ exists are as in Figure \ref{aab13graph}. 
The vertical graphs are constructed as in Figure \ref{g1g1t}, where we 
explain it only by ${\cal G}{\cal G}^{t}$.
\pagebreak
\begin{figure}[H]
\begin{center}
	\setlength{\unitlength}{0.5mm}
	\begin{picture}(200,60)(-20,35)
	\multiput(15,100)(30,0){3}{\circle*{2}}
	\put(135,100){\circle*{2}}
	\multiput(15,50)(30,0){3}{\circle*{2}}
	\put(135,50){\circle*{2}}
	\multiput(15,100)(30,0){3}{\line(3,-5){15}}
	\multiput(15,100)(30,0){2}{\line(-3,-5){15}}
	 \put(45,100){\line(9,-5){90}}
	\put(135,100){\line(-9,-5){90}}
	\multiput(0,75)(30,0){2}{\line(3,-5){15}}
	\multiput(30,75)(30,0){3}{\line(-3,-5){15}}	
	\put(21,100){\makebox(1,1){$a$}}
	\put(52,100){\makebox(1,1){$c$}}
	\put(83,100){\makebox(1,1){\s{a}}}
	\put(142,100){\makebox(1,1){\ss{a}}}
	\put(5,75){\makebox(1,1){1}}
	\put(35,75){\makebox(1,1){2}}
	\put(65,75){\makebox(1,1){3}}
	\put(100,75){\makebox(1,1){4}}
	\put(21,50){\makebox(1,1){$a$}}
	\put(54,50){\makebox(1,1){$c$}}
	\put(84,50){\makebox(1,1){\s{a}}}
	\put(144,50){\makebox(1,1){\ss{a}}}
	\put(-10,87){\makebox(2,2){ ${\cal G}_1$}}
	\put(-10,62){\makebox(2,2){${{\cal G}_1}^{t}$}}
	\put(-15,100){\makebox(2,2){$V_1$}}
	\put(-15,75){\makebox(2,2){$V_2$}}
	\put(-15,50){\makebox(2,2){ $V_1$}}
	\put(75,30){\makebox(1,2){$\Downarrow$}}
	\end{picture}
\begin{picture}(200,55)(-20,40)
 		\multiput(15,75)(30,0){3}{\circle*{2}}
	\put(135,75){\circle*{2}}
\multiput(15,50)(30,0){3}{\circle*{2}}
	\put(135,50){\circle*{2}}
	  \dline(15,75) \dline(45,75) \dline(75,75)  \dline(135,75) 
   \put(14,75){\line(0,-1){25}}
  \multiput(15,75)(30,0){2}{\line(6,-5){30}}
   \multiput(44,75)(2,0){2}{\line(0,-1){25}}
  \multiput(45,75)(30,0){2}{\line(-6,-5){30}}
 \put(45,75){\line(18,-5){90}}
  \put(75,75){\line(12,-5){60}}
  \put(135,75){\line(-18,-5){90}}
    \put(135,75){\line(-12,-5){60}}
		\put(21,75){\makebox(1,1){$a$}}
	\put(52,75){\makebox(1,1){$c$}}
	\put(84,75){\makebox(1,1){\s{a}}}
	\put(143,75){\makebox(1,1){\ss{a}}}
\put(21,50){\makebox(1,1){$a$}}
	\put(52,49){\makebox(1,1){$c$}}
	\put(83,49){\makebox(1,1){\s{a}}}
	\put(143,50){\makebox(1,1){\ss{a}}}
		\put(-10,62){\makebox(2,2){${\cal G}_1 {{\cal G}_1}^t$}}
			\put(-15,75){\makebox(2,2){$V_1$}}
	\put(-15,50){\makebox(2,2){$V_1$}}
\end{picture}
\end{center}
\caption{construction of the vertical graphs of $\a \ab$.}
\label{g1g1t}
\end{figure}
\pagebreak
\begin{figure}[H]
\begin{center}
$$
\thinlines
\unitlength 1.0mm
\begin{picture}(25,20)(5,0)
\put(11,13){\line(1,0){8}}
\multiput(11,3)(2,0){4}{\line(1,0){1}}
\multiput(10,12)(10,0){2}{\line(0,-1){8}}
\multiput(10,3)(10,0){2}{\circle*{1}}
\multiput(10,13)(10,0){2}{\circle*{1}}
\put(6,17){\makebox(0,0){$V_0$}}
\put(24,17){\makebox(0,0){$V_1$}}
\put(6,3){\makebox(0,0){$V_2$}}
\put(24,3){\makebox(0,0){$V_3$}}
\put(15,17){\makebox(0,0){${\cal G}_0$}}
\put(15,3){\makebox(0,0){\sm{${\cal G}_2$}}}
\put(6,8){\makebox(0,0){\sm{${\cal G}_3$}}}
\put(24,8){\makebox(0,0){\sm{${\cal G}_1$}}}
\put(15,8){\makebox(0,0){$\a$}}
\put(11,-7){\line(1,0){8}}
\multiput(10,2)(10,0){2}{\line(0,-1){8}}
\multiput(10,-7)(10,0){2}{\circle*{1}}
\put(6,-11){\makebox(0,0){$V_0$}}
\put(24,-11){\makebox(0,0){$V_1$}}
\put(15,-11){\makebox(0,0){${\cal G}_0$}}
\put(6,-2){\makebox(0,0){\sm{${{\cal G}_3}^t$}}}
\put(24,-2){\makebox(0,0){\sm{${{\cal G}_1}^t$}}}
\put(15,-2){\makebox(0,0){$\ab$}}
\end{picture}
=
\thinlines
\unitlength 1.0mm
\begin{picture}(30,20)(-5,0)
\multiput(11,-2)(0,10){2}{\line(1,0){8}}
\multiput(10,7)(10,0){2}{\line(0,-1){8}}
\multiput(10,-2)(10,0){2}{\circle*{1}}
\multiput(10,8)(10,0){2}{\circle*{1}}
\put(6,12){\makebox(0,0){$V_0$}}
\put(24,12){\makebox(0,0){$V_1$}}
\put(6,-6){\makebox(0,0){$V_0$}}
\put(24,-6){\makebox(0,0){$V_1$}}
\put(15,12){\makebox(0,0){${\cal G}_0$}}
\put(15,-6){\makebox(0,0){${\cal G}_0$}}
\put(4,3){\makebox(0,0){${\cal G}_3 {{\cal G}_3}^t$}}
\put(27,3){\makebox(0,0){${\cal G}_1 {{\cal G}_1}^t$}}
\put(15,3){\makebox(0,0){$\a \ab$}}
\end{picture}
$$ 
\end{center}
\begin{center}
	\setlength{\unitlength}{0.5mm}
	\begin{picture}(200,145)(-20,-5)
	\multiput(0,100)(30,0){6}{\circle*{2}}
	\multiput(15,75)(30,0){3}{\circle*{2}}
	\put(135,75){\circle*{2}}
\multiput(15,50)(30,0){3}{\circle*{2}}
	\put(135,50){\circle*{2}}
	\multiput(0,25)(30,0){6}{\circle*{2}}
	\multiput(0,0)(30,0){6}{\circle*{2}}

    \multiput(0,100)(30,0){2}{\line(3,-5){15}}
 \put(60,100){\line(3,-5){15}}
	\multiput(30,100)(30,0){2}{\line(-3,-5){15}}
	\put(150,100){\line(-3,-5){15}}
   \put(90,100){\line(9,-5){45}}
    \multiput(90,100)(30,0){2}{\line(-9,-5){45}}
  \dline(15,75) \dline(45,75) \dline(75,75)  \dline(135,75) 
   \put(14,75){\line(0,-1){25}}
  \multiput(15,75)(30,0){2}{\line(6,-5){30}}
   \multiput(44,75)(2,0){2}{\line(0,-1){25}}
  \multiput(45,75)(30,0){2}{\line(-6,-5){30}}
 \put(45,75){\line(18,-5){90}}
  \put(75,75){\line(12,-5){60}}
  \put(135,75){\line(-18,-5){90}}
    \put(135,75){\line(-12,-5){60}}
\multiput(15,50)(30,0){3}{\line(-3,-5){15}}
\multiput(15,50)(30,0){2}{\line(3,-5){15}}
\multiput(45,50)(30,0){2}{\line(9,-5){45}}
\put(135,50){\line(-9,-5){45}}
\put(135,50){\line(3,-5){15}}

\dline(0,25) \dline(30,25) \dline(60,25) \dline(90,25)
\dline(120,25) \dline(150,25) 
\multiput(0,25)(30,0){4}{\line(6,-5){30}}
\multiput(30,25)(30,0){4}{\line(-6,-5){30}}
 \put(30,25){\line(12,-5){60}}
 \put(60,25){\line(18,-5){90}}
 \put(90,25){\line(-12,-5){60}}
 \put(150,25){\line(-18,-5){90}}
\multiput(29,25)(30,0){3}{\line(0,-1){25}}

	\put(5,100){\makebox(1,1){$*$}}
	\put(35,100){\makebox(1,1){$b$}}
	\put(68,100){\makebox(1,1){\s{b}}}
	\put(98,102){\makebox(1,1){\ss{b}}}
	\put(127,100){\makebox(1,1){\s{*}}}
	\put(157,100){\makebox(1,1){\ss{*}}}
	\put(21,75){\makebox(1,1){$a$}}
	\put(52,75){\makebox(1,1){$c$}}
	\put(84,75){\makebox(1,1){\s{a}}}
	\put(143,75){\makebox(1,1){\ss{a}}}
\put(21,50){\makebox(1,1){$a$}}
	\put(52,49){\makebox(1,1){$c$}}
	\put(83,49){\makebox(1,1){\s{a}}}
	\put(143,50){\makebox(1,1){\ss{a}}}
\put(5,25){\makebox(1,1){$*$}}
	\put(35,28){\makebox(1,1){$b$}}
	\put(68,28){\makebox(1,1){\s{b}}}
	\put(100,23){\makebox(1,1){\ss{b}}}
	\put(126,25){\makebox(1,1){\s{*}}}
	\put(157,25){\makebox(1,1){\ss{*}}}
	
	\put(4,-2){\makebox(1,1){$*$}}
	\put(33,-4){\makebox(1,1){$b$}}
	\put(65,-4){\makebox(1,1){\s{b}}}
	\put(96,-4){\makebox(1,1){\ss{b}}}
	\put(126,-2){\makebox(1,1){\s{*}}}
	\put(157,-2){\makebox(1,1){\ss{*}}}

	\put(-10,87){\makebox(2,2){${\cal G}_0$}}
	\put(-10,62){\makebox(2,2){${\cal G}_1 {{\cal G}_1}^t$}}
	\put(-10,37){\makebox(2,2){${\cal G}_0$}}
	\put(-14,12){\makebox(2,2){${\cal G}_3 {{\cal G}_3}^t$}}
	\put(-15,100){\makebox(2,2){$V_0$}}
	\put(-15,75){\makebox(2,2){$V_1$}}
	\put(-15,50){\makebox(2,2){$V_1$}}
	\put(-15,25){\makebox(2,2){$V_0$}}
	\put(-15,0){\makebox(2,2){$V_0$}}
	\end{picture}
\caption{the four graphs of $\a \ab$}
\label{aab13graph}
\end{center}
\end{figure}
To obtain the connection $\a \ab-{\bf 1}$, we multiply the entries of
the connections $\a$ and $\ab$ properly (We call this sort of 
 computations of the multiplication of the connections 
``actual" multiplication.), transform it by vertical gauge so that
the entries corresponding to the trivial connection ${\bf 1}$ are 1, 
and subtract ${\bf 1}$.
(In Figure \ref{aab13graph}, the broken lines corresponds to 
this trivial summand.) Here, in the \framebox{Table \ref{land13.tex}},
we show the landscape of $\a \ab$ with 1's in the entries 
corresponding to ${\bf 1}$.

\def\odo{$\odot$}
\def\ci{$\circ$}
\def\dia{$\diamond$}
\begin{table}[h]
\begin{center}
\begin{tabular}[pos]{|@{}c@{}||@{}c@{}|@{}c@{}|@{}c@{}|| @{}c@{}|@{}c@{}
|@{}c@{}| @{}c@{}|@{}c@{}|@{}c@{}|| @{}c@{}|@{}c@{}|@{}c@{}|| @{}c@{}
|@{}c@{}|@{}c@{}|@{}c@{} }
\hl
& $aa^1$ & $aa^2$ & $ca$ & $ac$ & $cc^1$ & $cc^2$ & $cc^3$ & 
\s{a}$c$ &
\ss{a}$c$ & $c$\s{a} & \s{a}\s{a} & \ss{a}\s{a} & $c$\ss{a} & 
\s{a}\ss{a} & \ss{a}\ss{a} \\ \hline\hline
$**$ & 1 & 0 &&&&&&&&&&&&&  \\ \hl
$*b$ & 0 & \ci & & \odo &&&&&&&&&&&  \\ \hline\hline
$b*$ & 0 & \ci & \odo &&&&&&&&&&&& \\ \hl
$bb^1$ & 1 & 0 & 0 & 0 & 1 & 0 & 0 &&&&&&&& \\ \hl
$bb^2$ & 0 & \bl & \ci & \ci & 0 & \bl& \bl &&&&&&&& \\ \hl
$b$\s{b} &&&& \odo & 0 & \dia & \dia &&& \odo &&&&& \\ \hl
$b$\ss{b} &&&&\odo  & 0 & \dia & \dia & && &&& \odo && \\ \hline\hline
\s{b}$b$ & && \odo & & 0 & \dia & \dia & \odo & & &&&&& \\ \hl
\s{b}\s{b}$^1$ & &&&& 1 & 0 & 0 & 0 & & 0 & 1 &&&& \\ \hl
\s{b}\s{b}$^2$ & &&&& 0 & \bl & \bl & \ci & & \ci & 0 &&&& \\ \hl
\s{b}\ss{b} & &&&& 0 & \dia & \dia & \odo & &&&& \odo & \odo & \\ \hl
\s{b}\ss{*} &&&&& &&&&& &&& \odo & \odo & \\ \hline\hline
\ss{b}$b$ &&& \odo & & 0 & \dia & \dia & & \odo & &&&&& \\ \hl
\ss{b}\s{b} &&&&& 0 & \dia & \dia && \odo & \odo &  & \odo &&& \\ \hl
\ss{b}\ss{b}$^1$ &&&&& 1 & \bl & \bl & & \ci & & & & \ci & & 1 \\ \hl
\ss{b}\ss{b}$^2$ &&&&& 0 & \bl & \bl & & \ci & & & & \ci & & 0 \\ \hl
\ss{b}\s{*} & &&&& &&&&& \odo &  & \odo & && \\ \hline\hline
\s{*}\ss{b} &&&& &&&& \odo & & &&&& \odo & \\ \hl
\s{*}\s{*} &&&&& &&&&& & 1 &&&& \\ \hline\hline
\ss{*}\s{b} &&&& &&&& & \odo &&& \odo &&& \\ \hl
\ss{*}\ss{*} &&&&& &&&&& &&&&& 1 \\ \hl
\end{tabular}
\end{center}
\caption{Landscape of the connection $\alpha{\tilde \alpha}$ after a 
gauge transform}
\label{land13.tex}
\end{table}
First we will compute the entries marked \odo \ in the Table 
\ref{land13.tex}.
 We assume that $1 \times 1$ gauge transform unitaries 
corresponding to single vertical edges which connect different 
vertices in the graph ${\cal G}$ to be 1 without losing 
generality, because they are not involved in the trivial 
connection ${\bf 1}$. We compute such entries by ``actual'' 
multiplication. \\

$$
\thinlines
\unitlength 0.5mm
\begin{picture}(30,20)(0,0)
\multiput(11,-2)(0,10){2}{\line(1,0){8}}
\multiput(10,7)(10,0){2}{\line(0,-1){8}}
\multiput(10,-2)(10,0){2}{\circle*{1}}
\multiput(10,8)(10,0){2}{\circle*{1}}
\put(15,3){\makebox(0,0){{\sm $\a \ab$}}}
\put(6,12){\makebox(0,0){$*$}}
\put(24,12){\makebox(0,0){$a$}}
\put(6,-6){\makebox(0,0){$b$}}
\put(24,-6){\makebox(0,0){$c$}}
\end{picture}
=
\thinlines
\unitlength 0.5mm
\begin{picture}(30,20)(0,0)
\multiput(11,-2)(0,10){2}{\line(1,0){8}}
\multiput(10,7)(10,0){2}{\line(0,-1){8}}
\multiput(10,-2)(10,0){2}{\circle*{1}}
\multiput(10,8)(10,0){2}{\circle*{1}}
\put(15,3){\makebox(0,0){$\a$}}
\put(6,12){\makebox(0,0){$*$}}
\put(24,12){\makebox(0,0){$a$}}
\put(6,-6){\makebox(0,0){$a$}}
\put(24,-6){\makebox(0,0){$2$}}
\end{picture}
\cdot
\thinlines
\unitlength 0.5mm
\begin{picture}(30,20)(0,0)
\multiput(11,-2)(0,10){2}{\line(1,0){8}}
\multiput(10,7)(10,0){2}{\line(0,-1){8}}
\multiput(10,-2)(10,0){2}{\circle*{1}}
\multiput(10,8)(10,0){2}{\circle*{1}}
\put(15,3){\makebox(0,0){$\ab$}}
\put(6,12){\makebox(0,0){$a$}}
\put(24,12){\makebox(0,0){$2$}}
\put(6,-6){\makebox(0,0){$b$}}
\put(24,-6){\makebox(0,0){$c$}}
\end{picture}
= 1\cdot1 = 1,
$$

$$
\thinlines
\unitlength 0.5mm
\begin{picture}(30,20)(0,0)
\multiput(11,-2)(0,10){2}{\line(1,0){8}}
\multiput(10,7)(10,0){2}{\line(0,-1){8}}
\multiput(10,-2)(10,0){2}{\circle*{1}}
\multiput(10,8)(10,0){2}{\circle*{1}}
\put(15,3){\makebox(0,0){{\sm $\a \ab$}}}
\put(6,12){\makebox(0,0){$b$}}
\put(24,12){\makebox(0,0){$c$}}
\put(6,-6){\makebox(0,0){$*$}}
\put(24,-6){\makebox(0,0){$a$}}
\end{picture}
=
\thinlines
\unitlength 0.5mm
\begin{picture}(30,20)(0,0)
\multiput(11,-2)(0,10){2}{\line(1,0){8}}
\multiput(10,7)(10,0){2}{\line(0,-1){8}}
\multiput(10,-2)(10,0){2}{\circle*{1}}
\multiput(10,8)(10,0){2}{\circle*{1}}
\put(15,3){\makebox(0,0){$\a$}}
\put(6,12){\makebox(0,0){$b$}}
\put(24,12){\makebox(0,0){$c$}}
\put(6,-6){\makebox(0,0){$a$}}
\put(24,-6){\makebox(0,0){$2$}}
\end{picture}
\cdot
\thinlines
\unitlength 0.5mm
\begin{picture}(30,20)(0,0)
\multiput(11,-2)(0,10){2}{\line(1,0){8}}
\multiput(10,7)(10,0){2}{\line(0,-1){8}}
\multiput(10,-2)(10,0){2}{\circle*{1}}
\multiput(10,8)(10,0){2}{\circle*{1}}
\put(15,3){\makebox(0,0){$\ab$}}
\put(6,12){\makebox(0,0){$a$}}
\put(24,12){\makebox(0,0){$2$}}
\put(6,-6){\makebox(0,0){$*$}}
\put(24,-6){\makebox(0,0){$a$}}
\end{picture}
= \frac{\la \sqrt{\lad-2}}{\lad-1} \cdot \frac{\sqrt{\lad-1}}{\la}
=\frac{\sqrt{\lad-2}}{\sqrt{\lad-1}}.
$$
From here, we only write the result of multiplication.

$$
\thinlines
\unitlength 0.5mm
\begin{picture}(30,20)(0,0)
\multiput(11,-2)(0,10){2}{\line(1,0){8}}
\multiput(10,7)(10,0){2}{\line(0,-1){8}}
\multiput(10,-2)(10,0){2}{\circle*{1}}
\multiput(10,8)(10,0){2}{\circle*{1}}
\put(15,3){\makebox(0,0){{\sm $\a \ab$}}}
\put(6,12){\makebox(0,0){$b$}}
\put(24,12){\makebox(0,0){$a$}}
\put(6,-6){\makebox(0,0){\s{b}}}
\put(24,-6){\makebox(0,0){$c$}}
\end{picture}
=
\thinlines
\unitlength 0.5mm
\begin{picture}(30,20)(0,0)
\multiput(11,-2)(0,10){2}{\line(1,0){8}}
\multiput(10,7)(10,0){2}{\line(0,-1){8}}
\multiput(10,-2)(10,0){2}{\circle*{1}}
\multiput(10,8)(10,0){2}{\circle*{1}}
\put(15,3){\makebox(0,0){$\a$}}
\put(6,12){\makebox(0,0){$b$}}
\put(24,12){\makebox(0,0){$a$}}
\put(6,-6){\makebox(0,0){$c$}}
\put(24,-6){\makebox(0,0){$2$}}
\end{picture}
\cdot
\thinlines
\unitlength 0.5mm
\begin{picture}(30,20)(0,0)
\multiput(11,-2)(0,10){2}{\line(1,0){8}}
\multiput(10,7)(10,0){2}{\line(0,-1){8}}
\multiput(10,-2)(10,0){2}{\circle*{1}}
\multiput(10,8)(10,0){2}{\circle*{1}}
\put(15,3){\makebox(0,0){$\ab$}}
\put(6,12){\makebox(0,0){$c$}}
\put(24,12){\makebox(0,0){$2$}}
\put(6,-6){\makebox(0,0){\s{b}}}
\put(24,-6){\makebox(0,0){$c$}}
\end{picture}
= \frac{\rho}{\sqrt{\lad-2}},
$$

$$
\thinlines
\unitlength 0.5mm
\begin{picture}(30,20)(0,0)
\multiput(11,-2)(0,10){2}{\line(1,0){8}}
\multiput(10,7)(10,0){2}{\line(0,-1){8}}
\multiput(10,-2)(10,0){2}{\circle*{1}}
\multiput(10,8)(10,0){2}{\circle*{1}}
\put(15,3){\makebox(0,0){{\sm $\a \ab$}}}
\put(6,12){\makebox(0,0){$b$}}
\put(24,12){\makebox(0,0){$c$}}
\put(6,-6){\makebox(0,0){\s{b}}}
\put(24,-6){\makebox(0,0){\s{a}}}
\end{picture}
=
\thinlines
\unitlength 0.5mm
\begin{picture}(30,20)(0,0)
\multiput(11,-2)(0,10){2}{\line(1,0){8}}
\multiput(10,7)(10,0){2}{\line(0,-1){8}}
\multiput(10,-2)(10,0){2}{\circle*{1}}
\multiput(10,8)(10,0){2}{\circle*{1}}
\put(15,3){\makebox(0,0){$\a$}}
\put(6,12){\makebox(0,0){$b$}}
\put(24,12){\makebox(0,0){$c$}}
\put(6,-6){\makebox(0,0){$c$}}
\put(24,-6){\makebox(0,0){$4$}}
\end{picture}
\cdot
\thinlines
\unitlength 0.5mm
\begin{picture}(30,20)(0,0)
\multiput(11,-2)(0,10){2}{\line(1,0){8}}
\multiput(10,7)(10,0){2}{\line(0,-1){8}}
\multiput(10,-2)(10,0){2}{\circle*{1}}
\multiput(10,8)(10,0){2}{\circle*{1}}
\put(15,3){\makebox(0,0){$\ab$}}
\put(6,12){\makebox(0,0){$c$}}
\put(24,12){\makebox(0,0){$4$}}
\put(6,-6){\makebox(0,0){\s{b}}}
\put(24,-6){\makebox(0,0){\s{a}}}
\end{picture}
= 1,
$$

$$
\thinlines
\unitlength 0.5mm
\begin{picture}(30,20)(0,0)
\multiput(11,-2)(0,10){2}{\line(1,0){8}}
\multiput(10,7)(10,0){2}{\line(0,-1){8}}
\multiput(10,-2)(10,0){2}{\circle*{1}}
\multiput(10,8)(10,0){2}{\circle*{1}}
\put(15,3){\makebox(0,0){{\sm $\a \ab$}}}
\put(6,12){\makebox(0,0){$b$}}
\put(24,12){\makebox(0,0){$a$}}
\put(6,-6){\makebox(0,0){\ss{b}}}
\put(24,-6){\makebox(0,0){$c$}}
\end{picture}
=
\thinlines
\unitlength 0.5mm
\begin{picture}(30,20)(0,0)
\multiput(11,-2)(0,10){2}{\line(1,0){8}}
\multiput(10,7)(10,0){2}{\line(0,-1){8}}
\multiput(10,-2)(10,0){2}{\circle*{1}}
\multiput(10,8)(10,0){2}{\circle*{1}}
\put(15,3){\makebox(0,0){$\a$}}
\put(6,12){\makebox(0,0){$b$}}
\put(24,12){\makebox(0,0){$a$}}
\put(6,-6){\makebox(0,0){$c$}}
\put(24,-6){\makebox(0,0){$2$}}
\end{picture}
\cdot
\thinlines
\unitlength 0.5mm
\begin{picture}(30,20)(0,0)
\multiput(11,-2)(0,10){2}{\line(1,0){8}}
\multiput(10,7)(10,0){2}{\line(0,-1){8}}
\multiput(10,-2)(10,0){2}{\circle*{1}}
\multiput(10,8)(10,0){2}{\circle*{1}}
\put(15,3){\makebox(0,0){$\ab$}}
\put(6,12){\makebox(0,0){$c$}}
\put(24,12){\makebox(0,0){$2$}}
\put(6,-6){\makebox(0,0){\ss{b}}}
\put(24,-6){\makebox(0,0){$c$}}
\end{picture}
= \frac{\bar \rho}{\sqrt{\lad-2}},
$$

$$
\thinlines
\unitlength 0.5mm
\begin{picture}(30,20)(0,0)
\multiput(11,-2)(0,10){2}{\line(1,0){8}}
\multiput(10,7)(10,0){2}{\line(0,-1){8}}
\multiput(10,-2)(10,0){2}{\circle*{1}}
\multiput(10,8)(10,0){2}{\circle*{1}}
\put(15,3){\makebox(0,0){{\sm $\a \ab$}}}
\put(6,12){\makebox(0,0){$b$}}
\put(24,12){\makebox(0,0){$c$}}
\put(6,-6){\makebox(0,0){\ss{b}}}
\put(24,-6){\makebox(0,0){\ss{a}}}
\end{picture}
=
\thinlines
\unitlength 0.5mm
\begin{picture}(30,20)(0,0)
\multiput(11,-2)(0,10){2}{\line(1,0){8}}
\multiput(10,7)(10,0){2}{\line(0,-1){8}}
\multiput(10,-2)(10,0){2}{\circle*{1}}
\multiput(10,8)(10,0){2}{\circle*{1}}
\put(15,3){\makebox(0,0){$\a$}}
\put(6,12){\makebox(0,0){$b$}}
\put(24,12){\makebox(0,0){$c$}}
\put(6,-6){\makebox(0,0){$c$}}
\put(24,-6){\makebox(0,0){$4$}}
\end{picture}
\cdot
\thinlines
\unitlength 0.5mm
\begin{picture}(30,20)(0,0)
\multiput(11,-2)(0,10){2}{\line(1,0){8}}
\multiput(10,7)(10,0){2}{\line(0,-1){8}}
\multiput(10,-2)(10,0){2}{\circle*{1}}
\multiput(10,8)(10,0){2}{\circle*{1}}
\put(15,3){\makebox(0,0){$\ab$}}
\put(6,12){\makebox(0,0){$c$}}
\put(24,12){\makebox(0,0){$4$}}
\put(6,-6){\makebox(0,0){\ss{b}}}
\put(24,-6){\makebox(0,0){\ss{a}}}
\end{picture}
= 1,
$$
$$
\thinlines
\unitlength 0.5mm
\begin{picture}(30,20)(0,0)
\multiput(11,-2)(0,10){2}{\line(1,0){8}}
\multiput(10,7)(10,0){2}{\line(0,-1){8}}
\multiput(10,-2)(10,0){2}{\circle*{1}}
\multiput(10,8)(10,0){2}{\circle*{1}}
\put(15,3){\makebox(0,0){{\sm $\a \ab$}}}
\put(6,12){\makebox(0,0){\s{b}}}
\put(24,12){\makebox(0,0){$c$}}
\put(6,-6){\makebox(0,0){$b$}}
\put(24,-6){\makebox(0,0){$a$}}
\end{picture}
=
\thinlines
\unitlength 0.5mm
\begin{picture}(30,20)(0,0)
\multiput(11,-2)(0,10){2}{\line(1,0){8}}
\multiput(10,7)(10,0){2}{\line(0,-1){8}}
\multiput(10,-2)(10,0){2}{\circle*{1}}
\multiput(10,8)(10,0){2}{\circle*{1}}
\put(15,3){\makebox(0,0){$\a$}}
\put(6,12){\makebox(0,0){\s{b}}}
\put(24,12){\makebox(0,0){$c$}}
\put(6,-6){\makebox(0,0){$c$}}
\put(24,-6){\makebox(0,0){$2$}}
\end{picture}
\cdot
\thinlines
\unitlength 0.5mm
\begin{picture}(30,20)(0,0)
\multiput(11,-2)(0,10){2}{\line(1,0){8}}
\multiput(10,7)(10,0){2}{\line(0,-1){8}}
\multiput(10,-2)(10,0){2}{\circle*{1}}
\multiput(10,8)(10,0){2}{\circle*{1}}
\put(15,3){\makebox(0,0){$\ab$}}
\put(6,12){\makebox(0,0){$c$}}
\put(24,12){\makebox(0,0){$2$}}
\put(6,-6){\makebox(0,0){$b$}}
\put(24,-6){\makebox(0,0){$a$}}
\end{picture}
= {\bar \rho},
$$

$$
\thinlines
\unitlength 0.5mm
\begin{picture}(30,20)(0,0)
\multiput(11,-2)(0,10){2}{\line(1,0){8}}
\multiput(10,7)(10,0){2}{\line(0,-1){8}}
\multiput(10,-2)(10,0){2}{\circle*{1}}
\multiput(10,8)(10,0){2}{\circle*{1}}
\put(15,3){\makebox(0,0){{\sm $\a \ab$}}}
\put(6,12){\makebox(0,0){\s{b}}}
\put(24,12){\makebox(0,0){\s{a}}}
\put(6,-6){\makebox(0,0){$b$}}
\put(24,-6){\makebox(0,0){$c$}}
\end{picture}
=
\thinlines
\unitlength 0.5mm
\begin{picture}(30,20)(0,0)
\multiput(11,-2)(0,10){2}{\line(1,0){8}}
\multiput(10,7)(10,0){2}{\line(0,-1){8}}
\multiput(10,-2)(10,0){2}{\circle*{1}}
\multiput(10,8)(10,0){2}{\circle*{1}}
\put(15,3){\makebox(0,0){$\a$}}
\put(6,12){\makebox(0,0){\s{b}}}
\put(24,12){\makebox(0,0){\s{a}}}
\put(6,-6){\makebox(0,0){$c$}}
\put(24,-6){\makebox(0,0){$4$}}
\end{picture}
\cdot
\thinlines
\unitlength 0.5mm
\begin{picture}(30,20)(0,0)
\multiput(11,-2)(0,10){2}{\line(1,0){8}}
\multiput(10,7)(10,0){2}{\line(0,-1){8}}
\multiput(10,-2)(10,0){2}{\circle*{1}}
\multiput(10,8)(10,0){2}{\circle*{1}}
\put(15,3){\makebox(0,0){$\ab$}}
\put(6,12){\makebox(0,0){$c$}}
\put(24,12){\makebox(0,0){$4$}}
\put(6,-6){\makebox(0,0){$b$}}
\put(24,-6){\makebox(0,0){$c$}}
\end{picture}
= \frac{1}{\sqrt{\lad-2}},
$$

$$
\thinlines
\unitlength 0.5mm
\begin{picture}(30,20)(0,0)
\multiput(11,-2)(0,10){2}{\line(1,0){8}}
\multiput(10,7)(10,0){2}{\line(0,-1){8}}
\multiput(10,-2)(10,0){2}{\circle*{1}}
\multiput(10,8)(10,0){2}{\circle*{1}}
\put(15,3){\makebox(0,0){{\sm $\a \ab$}}}
\put(6,12){\makebox(0,0){\s{b}}}
\put(24,12){\makebox(0,0){\s{a}}}
\put(6,-6){\makebox(0,0){\ss{b}}}
\put(24,-6){\makebox(0,0){$c$}}
\end{picture}
=
\thinlines
\unitlength 0.5mm
\begin{picture}(30,20)(0,0)
\multiput(11,-2)(0,10){2}{\line(1,0){8}}
\multiput(10,7)(10,0){2}{\line(0,-1){8}}
\multiput(10,-2)(10,0){2}{\circle*{1}}
\multiput(10,8)(10,0){2}{\circle*{1}}
\put(15,3){\makebox(0,0){$\a$}}
\put(6,12){\makebox(0,0){\s{b}}}
\put(24,12){\makebox(0,0){\s{a}}}
\put(6,-6){\makebox(0,0){$c$}}
\put(24,-6){\makebox(0,0){$4$}}
\end{picture}
\cdot
\thinlines
\unitlength 0.5mm
\begin{picture}(30,20)(0,0)
\multiput(11,-2)(0,10){2}{\line(1,0){8}}
\multiput(10,7)(10,0){2}{\line(0,-1){8}}
\multiput(10,-2)(10,0){2}{\circle*{1}}
\multiput(10,8)(10,0){2}{\circle*{1}}
\put(15,3){\makebox(0,0){$\ab$}}
\put(6,12){\makebox(0,0){$c$}}
\put(24,12){\makebox(0,0){$4$}}
\put(6,-6){\makebox(0,0){\ss{b}}}
\put(24,-6){\makebox(0,0){$c$}}
\end{picture}
= \frac{\sqrt{\lad-4}}{\sqrt{\lad-2}},
$$

$$
\thinlines
\unitlength 0.5mm
\begin{picture}(30,20)(0,0)
\multiput(11,-2)(0,10){2}{\line(1,0){8}}
\multiput(10,7)(10,0){2}{\line(0,-1){8}}
\multiput(10,-2)(10,0){2}{\circle*{1}}
\multiput(10,8)(10,0){2}{\circle*{1}}
\put(15,3){\makebox(0,0){{\sm $\a \ab$}}}
\put(6,12){\makebox(0,0){\s{b}}}
\put(24,12){\makebox(0,0){$c$}}
\put(6,-6){\makebox(0,0){\ss{b}}}
\put(24,-6){\makebox(0,0){\ss{a}}}
\end{picture}
=
\thinlines
\unitlength 0.5mm
\begin{picture}(30,20)(0,0)
\multiput(11,-2)(0,10){2}{\line(1,0){8}}
\multiput(10,7)(10,0){2}{\line(0,-1){8}}
\multiput(10,-2)(10,0){2}{\circle*{1}}
\multiput(10,8)(10,0){2}{\circle*{1}}
\put(15,3){\makebox(0,0){$\a$}}
\put(6,12){\makebox(0,0){\s{b}}}
\put(24,12){\makebox(0,0){$c$}}
\put(6,-6){\makebox(0,0){$c$}}
\put(24,-6){\makebox(0,0){$4$}}
\end{picture}
\cdot
\thinlines
\unitlength 0.5mm
\begin{picture}(30,20)(0,0)
\multiput(11,-2)(0,10){2}{\line(1,0){8}}
\multiput(10,7)(10,0){2}{\line(0,-1){8}}
\multiput(10,-2)(10,0){2}{\circle*{1}}
\multiput(10,8)(10,0){2}{\circle*{1}}
\put(15,3){\makebox(0,0){$\ab$}}
\put(6,12){\makebox(0,0){$c$}}
\put(24,12){\makebox(0,0){$4$}}
\put(6,-6){\makebox(0,0){\ss{b}}}
\put(24,-6){\makebox(0,0){\ss{a}}}
\end{picture}
= \frac{1}{\sqrt{\lad-1}},
$$

$$
\thinlines
\unitlength 0.5mm
\begin{picture}(30,20)(0,0)
\multiput(11,-2)(0,10){2}{\line(1,0){8}}
\multiput(10,7)(10,0){2}{\line(0,-1){8}}
\multiput(10,-2)(10,0){2}{\circle*{1}}
\multiput(10,8)(10,0){2}{\circle*{1}}
\put(15,3){\makebox(0,0){{\sm $\a \ab$}}}
\put(6,12){\makebox(0,0){\s{b}}}
\put(24,12){\makebox(0,0){\s{a}}}
\put(6,-6){\makebox(0,0){\ss{b}}}
\put(24,-6){\makebox(0,0){\ss{a}}}
\end{picture}
=
\thinlines
\unitlength 0.5mm
\begin{picture}(30,20)(0,0)
\multiput(11,-2)(0,10){2}{\line(1,0){8}}
\multiput(10,7)(10,0){2}{\line(0,-1){8}}
\multiput(10,-2)(10,0){2}{\circle*{1}}
\multiput(10,8)(10,0){2}{\circle*{1}}
\put(15,3){\makebox(0,0){$\a$}}
\put(6,12){\makebox(0,0){\s{b}}}
\put(24,12){\makebox(0,0){\s{a}}}
\put(6,-6){\makebox(0,0){$c$}}
\put(24,-6){\makebox(0,0){$4$}}
\end{picture}
\cdot
\thinlines
\unitlength 0.5mm
\begin{picture}(30,20)(0,0)
\multiput(11,-2)(0,10){2}{\line(1,0){8}}
\multiput(10,7)(10,0){2}{\line(0,-1){8}}
\multiput(10,-2)(10,0){2}{\circle*{1}}
\multiput(10,8)(10,0){2}{\circle*{1}}
\put(15,3){\makebox(0,0){$\ab$}}
\put(6,12){\makebox(0,0){$c$}}
\put(24,12){\makebox(0,0){$4$}}
\put(6,-6){\makebox(0,0){\ss{b}}}
\put(24,-6){\makebox(0,0){\ss{a}}}
\end{picture}
= \frac{\sqrt{\lad-2}}{\sqrt{\lad-1}},
$$

$$
\thinlines
\unitlength 0.5mm
\begin{picture}(30,20)(0,0)
\multiput(11,-2)(0,10){2}{\line(1,0){8}}
\multiput(10,7)(10,0){2}{\line(0,-1){8}}
\multiput(10,-2)(10,0){2}{\circle*{1}}
\multiput(10,8)(10,0){2}{\circle*{1}}
\put(15,3){\makebox(0,0){{\sm $\a \ab$}}}
\put(6,12){\makebox(0,0){\s{b}}}
\put(24,12){\makebox(0,0){$c$}}
\put(6,-6){\makebox(0,0){\ss{*}}}
\put(24,-6){\makebox(0,0){\ss{a}}}
\end{picture}
=
\thinlines
\unitlength 0.5mm
\begin{picture}(30,20)(0,0)
\multiput(11,-2)(0,10){2}{\line(1,0){8}}
\multiput(10,7)(10,0){2}{\line(0,-1){8}}
\multiput(10,-2)(10,0){2}{\circle*{1}}
\multiput(10,8)(10,0){2}{\circle*{1}}
\put(15,3){\makebox(0,0){$\a$}}
\put(6,12){\makebox(0,0){\s{b}}}
\put(24,12){\makebox(0,0){$c$}}
\put(6,-6){\makebox(0,0){\s{a}}}
\put(24,-6){\makebox(0,0){$4$}}
\end{picture}
\cdot
\thinlines
\unitlength 0.5mm
\begin{picture}(30,20)(0,0)
\multiput(11,-2)(0,10){2}{\line(1,0){8}}
\multiput(10,7)(10,0){2}{\line(0,-1){8}}
\multiput(10,-2)(10,0){2}{\circle*{1}}
\multiput(10,8)(10,0){2}{\circle*{1}}
\put(15,3){\makebox(0,0){$\ab$}}
\put(6,12){\makebox(0,0){\s{a}}}
\put(24,12){\makebox(0,0){$4$}}
\put(6,-6){\makebox(0,0){\ss{*}}}
\put(24,-6){\makebox(0,0){\ss{a}}}
\end{picture}
= \frac{\sqrt{\lad-2}}{\sqrt{\lad-1}},
$$
$$
\thinlines
\unitlength 0.5mm
\begin{picture}(30,20)(0,0)
\multiput(11,-2)(0,10){2}{\line(1,0){8}}
\multiput(10,7)(10,0){2}{\line(0,-1){8}}
\multiput(10,-2)(10,0){2}{\circle*{1}}
\multiput(10,8)(10,0){2}{\circle*{1}}
\put(15,3){\makebox(0,0){{\sm $\a \ab$}}}
\put(6,12){\makebox(0,0){\s{b}}}
\put(24,12){\makebox(0,0){\s{a}}}
\put(6,-6){\makebox(0,0){\ss{*}}}
\put(24,-6){\makebox(0,0){\ss{a}}}
\end{picture}
=
\thinlines
\unitlength 0.5mm
\begin{picture}(30,20)(0,0)
\multiput(11,-2)(0,10){2}{\line(1,0){8}}
\multiput(10,7)(10,0){2}{\line(0,-1){8}}
\multiput(10,-2)(10,0){2}{\circle*{1}}
\multiput(10,8)(10,0){2}{\circle*{1}}
\put(15,3){\makebox(0,0){$\a$}}
\put(6,12){\makebox(0,0){\s{b}}}
\put(24,12){\makebox(0,0){\s{a}}}
\put(6,-6){\makebox(0,0){\s{a}}}
\put(24,-6){\makebox(0,0){$4$}}
\end{picture}
\cdot
\thinlines
\unitlength 0.5mm
\begin{picture}(30,20)(0,0)
\multiput(11,-2)(0,10){2}{\line(1,0){8}}
\multiput(10,7)(10,0){2}{\line(0,-1){8}}
\multiput(10,-2)(10,0){2}{\circle*{1}}
\multiput(10,8)(10,0){2}{\circle*{1}}
\put(15,3){\makebox(0,0){$\ab$}}
\put(6,12){\makebox(0,0){\s{a}}}
\put(24,12){\makebox(0,0){$4$}}
\put(6,-6){\makebox(0,0){\ss{*}}}
\put(24,-6){\makebox(0,0){\ss{a}}}
\end{picture}
= \frac{-1}{\sqrt{\lad-1}},
$$
$$
\thinlines
\unitlength 0.5mm
\begin{picture}(30,20)(0,0)
\multiput(11,-2)(0,10){2}{\line(1,0){8}}
\multiput(10,7)(10,0){2}{\line(0,-1){8}}
\multiput(10,-2)(10,0){2}{\circle*{1}}
\multiput(10,8)(10,0){2}{\circle*{1}}
\put(15,3){\makebox(0,0){{\sm $\a \ab$}}}
\put(6,12){\makebox(0,0){\ss{b}}}
\put(24,12){\makebox(0,0){$c$}}
\put(6,-6){\makebox(0,0){$b$}}
\put(24,-6){\makebox(0,0){$a$}}
\end{picture}
=
\thinlines
\unitlength 0.5mm
\begin{picture}(30,20)(0,0)
\multiput(11,-2)(0,10){2}{\line(1,0){8}}
\multiput(10,7)(10,0){2}{\line(0,-1){8}}
\multiput(10,-2)(10,0){2}{\circle*{1}}
\multiput(10,8)(10,0){2}{\circle*{1}}
\put(15,3){\makebox(0,0){$\a$}}
\put(6,12){\makebox(0,0){\ss{b}}}
\put(24,12){\makebox(0,0){$c$}}
\put(6,-6){\makebox(0,0){$c$}}
\put(24,-6){\makebox(0,0){$2$}}
\end{picture}
\cdot
\thinlines
\unitlength 0.5mm
\begin{picture}(30,20)(0,0)
\multiput(11,-2)(0,10){2}{\line(1,0){8}}
\multiput(10,7)(10,0){2}{\line(0,-1){8}}
\multiput(10,-2)(10,0){2}{\circle*{1}}
\multiput(10,8)(10,0){2}{\circle*{1}}
\put(15,3){\makebox(0,0){$\ab$}}
\put(6,12){\makebox(0,0){$c$}}
\put(24,12){\makebox(0,0){$2$}}
\put(6,-6){\makebox(0,0){$b$}}
\put(24,-6){\makebox(0,0){$a$}}
\end{picture}
= \rho,
$$
$$
\thinlines
\unitlength 0.5mm
\begin{picture}(30,20)(0,0)
\multiput(11,-2)(0,10){2}{\line(1,0){8}}
\multiput(10,7)(10,0){2}{\line(0,-1){8}}
\multiput(10,-2)(10,0){2}{\circle*{1}}
\multiput(10,8)(10,0){2}{\circle*{1}}
\put(15,3){\makebox(0,0){{\sm $\a \ab$}}}
\put(6,12){\makebox(0,0){\ss{b}}}
\put(24,12){\makebox(0,0){\ss{a}}}
\put(6,-6){\makebox(0,0){$b$}}
\put(24,-6){\makebox(0,0){$c$}}
\end{picture}
=
\thinlines
\unitlength 0.5mm
\begin{picture}(30,20)(0,0)
\multiput(11,-2)(0,10){2}{\line(1,0){8}}
\multiput(10,7)(10,0){2}{\line(0,-1){8}}
\multiput(10,-2)(10,0){2}{\circle*{1}}
\multiput(10,8)(10,0){2}{\circle*{1}}
\put(15,3){\makebox(0,0){\ss{b}}}
\put(6,12){\makebox(0,0){\ss{a}}}
\put(24,12){\makebox(0,0){$c$}}
\put(6,-6){\makebox(0,0){$c$}}
\put(24,-6){\makebox(0,0){$4$}}
\end{picture}
\cdot
\thinlines
\unitlength 0.5mm
\begin{picture}(30,20)(0,0)
\multiput(11,-2)(0,10){2}{\line(1,0){8}}
\multiput(10,7)(10,0){2}{\line(0,-1){8}}
\multiput(10,-2)(10,0){2}{\circle*{1}}
\multiput(10,8)(10,0){2}{\circle*{1}}
\put(15,3){\makebox(0,0){$\ab$}}
\put(6,12){\makebox(0,0){$c$}}
\put(24,12){\makebox(0,0){$4$}}
\put(6,-6){\makebox(0,0){$b$}}
\put(24,-6){\makebox(0,0){$c$}}
\end{picture}
= \frac{1}{\sqrt{\lad-2}},
$$
$$
\thinlines
\unitlength 0.5mm
\begin{picture}(30,20)(0,0)
\multiput(11,-2)(0,10){2}{\line(1,0){8}}
\multiput(10,7)(10,0){2}{\line(0,-1){8}}
\multiput(10,-2)(10,0){2}{\circle*{1}}
\multiput(10,8)(10,0){2}{\circle*{1}}
\put(15,3){\makebox(0,0){{\sm $\a \ab$}}}
\put(6,12){\makebox(0,0){\ss{b}}}
\put(24,12){\makebox(0,0){\ss{a}}}
\put(6,-6){\makebox(0,0){\ss{b}}}
\put(24,-6){\makebox(0,0){$c$}}
\end{picture}
=
\thinlines
\unitlength 0.5mm
\begin{picture}(30,20)(0,0)
\multiput(11,-2)(0,10){2}{\line(1,0){8}}
\multiput(10,7)(10,0){2}{\line(0,-1){8}}
\multiput(10,-2)(10,0){2}{\circle*{1}}
\multiput(10,8)(10,0){2}{\circle*{1}}
\put(15,3){\makebox(0,0){$\a$}}
\put(6,12){\makebox(0,0){\ss{b}}}
\put(24,12){\makebox(0,0){\ss{a}}}
\put(6,-6){\makebox(0,0){$c$}}
\put(24,-6){\makebox(0,0){$4$}}
\end{picture}
\cdot
\thinlines
\unitlength 0.5mm
\begin{picture}(30,20)(0,0)
\multiput(11,-2)(0,10){2}{\line(1,0){8}}
\multiput(10,7)(10,0){2}{\line(0,-1){8}}
\multiput(10,-2)(10,0){2}{\circle*{1}}
\multiput(10,8)(10,0){2}{\circle*{1}}
\put(15,3){\makebox(0,0){$\ab$}}
\put(6,12){\makebox(0,0){$c$}}
\put(24,12){\makebox(0,0){$4$}}
\put(6,-6){\makebox(0,0){\s{b}}}
\put(24,-6){\makebox(0,0){$c$}}
\end{picture}
= \frac{\sqrt{\lad-4}}{\sqrt{\lad-2}},
$$
$$
\thinlines
\unitlength 0.5mm
\begin{picture}(30,20)(0,0)
\multiput(11,-2)(0,10){2}{\line(1,0){8}}
\multiput(10,7)(10,0){2}{\line(0,-1){8}}
\multiput(10,-2)(10,0){2}{\circle*{1}}
\multiput(10,8)(10,0){2}{\circle*{1}}
\put(15,3){\makebox(0,0){{\sm $\a \ab$}}}
\put(6,12){\makebox(0,0){\ss{b}}}
\put(24,12){\makebox(0,0){$c$}}
\put(6,-6){\makebox(0,0){\s{b}}}
\put(24,-6){\makebox(0,0){\s{a}}}
\end{picture}
=
\thinlines
\unitlength 0.5mm
\begin{picture}(30,20)(0,0)
\multiput(11,-2)(0,10){2}{\line(1,0){8}}
\multiput(10,7)(10,0){2}{\line(0,-1){8}}
\multiput(10,-2)(10,0){2}{\circle*{1}}
\multiput(10,8)(10,0){2}{\circle*{1}}
\put(15,3){\makebox(0,0){$\a$}}
\put(6,12){\makebox(0,0){\ss{b}}}
\put(24,12){\makebox(0,0){$c$}}
\put(6,-6){\makebox(0,0){$c$}}
\put(24,-6){\makebox(0,0){$4$}}
\end{picture}
\cdot
\thinlines
\unitlength 0.5mm
\begin{picture}(30,20)(0,0)
\multiput(11,-2)(0,10){2}{\line(1,0){8}}
\multiput(10,7)(10,0){2}{\line(0,-1){8}}
\multiput(10,-2)(10,0){2}{\circle*{1}}
\multiput(10,8)(10,0){2}{\circle*{1}}
\put(15,3){\makebox(0,0){$\ab$}}
\put(6,12){\makebox(0,0){$c$}}
\put(24,12){\makebox(0,0){$4$}}
\put(6,-6){\makebox(0,0){\s{b}}}
\put(24,-6){\makebox(0,0){\s{a}}}
\end{picture}
= \frac{1}{\sqrt{\lad-1}},
$$
$$
\thinlines
\unitlength 0.5mm
\begin{picture}(30,20)(0,0)
\multiput(11,-2)(0,10){2}{\line(1,0){8}}
\multiput(10,7)(10,0){2}{\line(0,-1){8}}
\multiput(10,-2)(10,0){2}{\circle*{1}}
\multiput(10,8)(10,0){2}{\circle*{1}}
\put(15,3){\makebox(0,0){{\sm $\a \ab$}}}
\put(6,12){\makebox(0,0){\ss{b}}}
\put(24,12){\makebox(0,0){\ss{a}}}
\put(6,-6){\makebox(0,0){\s{b}}}
\put(24,-6){\makebox(0,0){\s{a}}}
\end{picture}
=
\thinlines
\unitlength 0.5mm
\begin{picture}(30,20)(0,0)
\multiput(11,-2)(0,10){2}{\line(1,0){8}}
\multiput(10,7)(10,0){2}{\line(0,-1){8}}
\multiput(10,-2)(10,0){2}{\circle*{1}}
\multiput(10,8)(10,0){2}{\circle*{1}}
\put(15,3){\makebox(0,0){$\a$}}
\put(6,12){\makebox(0,0){\ss{b}}}
\put(24,12){\makebox(0,0){\ss{a}}}
\put(6,-6){\makebox(0,0){$c$}}
\put(24,-6){\makebox(0,0){$4$}}
\end{picture}
\cdot
\thinlines
\unitlength 0.5mm
\begin{picture}(30,20)(0,0)
\multiput(11,-2)(0,10){2}{\line(1,0){8}}
\multiput(10,7)(10,0){2}{\line(0,-1){8}}
\multiput(10,-2)(10,0){2}{\circle*{1}}
\multiput(10,8)(10,0){2}{\circle*{1}}
\put(15,3){\makebox(0,0){$\ab$}}
\put(6,12){\makebox(0,0){$c$}}
\put(24,12){\makebox(0,0){$4$}}
\put(6,-6){\makebox(0,0){\s{b}}}
\put(24,-6){\makebox(0,0){\s{a}}}
\end{picture}
= \frac{\sqrt{\lad-2}}{\sqrt{\lad-1}},
$$
$$
\thinlines
\unitlength 0.5mm
\begin{picture}(30,20)(0,0)
\multiput(11,-2)(0,10){2}{\line(1,0){8}}
\multiput(10,7)(10,0){2}{\line(0,-1){8}}
\multiput(10,-2)(10,0){2}{\circle*{1}}
\multiput(10,8)(10,0){2}{\circle*{1}}
\put(15,3){\makebox(0,0){{\sm $\a \ab$}}}
\put(6,12){\makebox(0,0){\ss{b}}}
\put(24,12){\makebox(0,0){$c$}}
\put(6,-6){\makebox(0,0){\s{*}}}
\put(24,-6){\makebox(0,0){\s{a}}}
\end{picture}
=
\thinlines
\unitlength 0.5mm
\begin{picture}(30,20)(0,0)
\multiput(11,-2)(0,10){2}{\line(1,0){8}}
\multiput(10,7)(10,0){2}{\line(0,-1){8}}
\multiput(10,-2)(10,0){2}{\circle*{1}}
\multiput(10,8)(10,0){2}{\circle*{1}}
\put(15,3){\makebox(0,0){$\a$}}
\put(6,12){\makebox(0,0){\ss{b}}}
\put(24,12){\makebox(0,0){$c$}}
\put(6,-6){\makebox(0,0){\ss{a}}}
\put(24,-6){\makebox(0,0){$4$}}
\end{picture}
\cdot
\thinlines
\unitlength 0.5mm
\begin{picture}(30,20)(0,0)
\multiput(11,-2)(0,10){2}{\line(1,0){8}}
\multiput(10,7)(10,0){2}{\line(0,-1){8}}
\multiput(10,-2)(10,0){2}{\circle*{1}}
\multiput(10,8)(10,0){2}{\circle*{1}}
\put(15,3){\makebox(0,0){$\ab$}}
\put(6,12){\makebox(0,0){\ss{a}}}
\put(24,12){\makebox(0,0){$4$}}
\put(6,-6){\makebox(0,0){\s{*}}}
\put(24,-6){\makebox(0,0){\s{a}}}
\end{picture}
= \frac{\sqrt{\lad-2}}{\sqrt{\lad-1}},
$$
$$
\thinlines
\unitlength 0.5mm
\begin{picture}(30,20)(0,0)
\multiput(11,-2)(0,10){2}{\line(1,0){8}}
\multiput(10,7)(10,0){2}{\line(0,-1){8}}
\multiput(10,-2)(10,0){2}{\circle*{1}}
\multiput(10,8)(10,0){2}{\circle*{1}}
\put(15,3){\makebox(0,0){{\sm $\a \ab$}}}
\put(6,12){\makebox(0,0){\ss{b}}}
\put(24,12){\makebox(0,0){\ss{a}}}
\put(6,-6){\makebox(0,0){\s{*}}}
\put(24,-6){\makebox(0,0){\s{a}}}
\end{picture}
=
\thinlines
\unitlength 0.5mm
\begin{picture}(30,20)(0,0)
\multiput(11,-2)(0,10){2}{\line(1,0){8}}
\multiput(10,7)(10,0){2}{\line(0,-1){8}}
\multiput(10,-2)(10,0){2}{\circle*{1}}
\multiput(10,8)(10,0){2}{\circle*{1}}
\put(15,3){\makebox(0,0){$\a$}}
\put(6,12){\makebox(0,0){\ss{b}}}
\put(24,12){\makebox(0,0){\ss{a}}}
\put(6,-6){\makebox(0,0){\ss{a}}}
\put(24,-6){\makebox(0,0){$4$}}
\end{picture}
\cdot
\thinlines
\unitlength 0.5mm
\begin{picture}(30,20)(0,0)
\multiput(11,-2)(0,10){2}{\line(1,0){8}}
\multiput(10,7)(10,0){2}{\line(0,-1){8}}
\multiput(10,-2)(10,0){2}{\circle*{1}}
\multiput(10,8)(10,0){2}{\circle*{1}}
\put(15,3){\makebox(0,0){$\ab$}}
\put(6,12){\makebox(0,0){\ss{a}}}
\put(24,12){\makebox(0,0){$4$}}
\put(6,-6){\makebox(0,0){\s{*}}}
\put(24,-6){\makebox(0,0){\s{a}}}
\end{picture}
= \frac{-1}{\sqrt{\lad-1}},
$$
$$
\thinlines
\unitlength 0.5mm
\begin{picture}(30,20)(0,0)
\multiput(11,-2)(0,10){2}{\line(1,0){8}}
\multiput(10,7)(10,0){2}{\line(0,-1){8}}
\multiput(10,-2)(10,0){2}{\circle*{1}}
\multiput(10,8)(10,0){2}{\circle*{1}}
\put(15,3){\makebox(0,0){{\sm $\a \ab$}}}
\put(6,12){\makebox(0,0){\s{*}}}
\put(24,12){\makebox(0,0){\s{a}}}
\put(6,-6){\makebox(0,0){\ss{b}}}
\put(24,-6){\makebox(0,0){$c$}}
\end{picture}
=
\thinlines
\unitlength 0.5mm
\begin{picture}(30,20)(0,0)
\multiput(11,-2)(0,10){2}{\line(1,0){8}}
\multiput(10,7)(10,0){2}{\line(0,-1){8}}
\multiput(10,-2)(10,0){2}{\circle*{1}}
\multiput(10,8)(10,0){2}{\circle*{1}}
\put(15,3){\makebox(0,0){$\a$}}
\put(6,12){\makebox(0,0){\s{*}}}
\put(24,12){\makebox(0,0){\s{a}}}
\put(6,-6){\makebox(0,0){\ss{a}}}
\put(24,-6){\makebox(0,0){$4$}}
\end{picture}
\cdot
\thinlines
\unitlength 0.5mm
\begin{picture}(30,20)(0,0)
\multiput(11,-2)(0,10){2}{\line(1,0){8}}
\multiput(10,7)(10,0){2}{\line(0,-1){8}}
\multiput(10,-2)(10,0){2}{\circle*{1}}
\multiput(10,8)(10,0){2}{\circle*{1}}
\put(15,3){\makebox(0,0){$\ab$}}
\put(6,12){\makebox(0,0){\ss{a}}}
\put(24,12){\makebox(0,0){$4$}}
\put(6,-6){\makebox(0,0){\ss{b}}}
\put(24,-6){\makebox(0,0){$c$}}
\end{picture}
= 1,
$$
$$
\thinlines
\unitlength 0.5mm
\begin{picture}(30,20)(0,0)
\multiput(11,-2)(0,10){2}{\line(1,0){8}}
\multiput(10,7)(10,0){2}{\line(0,-1){8}}
\multiput(10,-2)(10,0){2}{\circle*{1}}
\multiput(10,8)(10,0){2}{\circle*{1}}
\put(15,3){\makebox(0,0){{\sm $\a \ab$}}}
\put(6,12){\makebox(0,0){\s{*}}}
\put(24,12){\makebox(0,0){\s{a}}}
\put(6,-6){\makebox(0,0){\ss{b}}}
\put(24,-6){\makebox(0,0){\ss{a}}}
\end{picture}
=
\thinlines
\unitlength 0.5mm
\begin{picture}(30,20)(0,0)
\multiput(11,-2)(0,10){2}{\line(1,0){8}}
\multiput(10,7)(10,0){2}{\line(0,-1){8}}
\multiput(10,-2)(10,0){2}{\circle*{1}}
\multiput(10,8)(10,0){2}{\circle*{1}}
\put(15,3){\makebox(0,0){$\a$}}
\put(6,12){\makebox(0,0){\s{*}}}
\put(24,12){\makebox(0,0){\s{a}}}
\put(6,-6){\makebox(0,0){\ss{a}}}
\put(24,-6){\makebox(0,0){$4$}}
\end{picture}
\cdot
\thinlines
\unitlength 0.5mm
\begin{picture}(30,20)(0,0)
\multiput(11,-2)(0,10){2}{\line(1,0){8}}
\multiput(10,7)(10,0){2}{\line(0,-1){8}}
\multiput(10,-2)(10,0){2}{\circle*{1}}
\multiput(10,8)(10,0){2}{\circle*{1}}
\put(15,3){\makebox(0,0){$\ab$}}
\put(6,12){\makebox(0,0){\ss{a}}}
\put(24,12){\makebox(0,0){$4$}}
\put(6,-6){\makebox(0,0){\ss{b}}}
\put(24,-6){\makebox(0,0){\ss{a}}}
\end{picture}
= -1,
$$
$$
\thinlines
\unitlength 0.5mm
\begin{picture}(30,20)(0,0)
\multiput(11,-2)(0,10){2}{\line(1,0){8}}
\multiput(10,7)(10,0){2}{\line(0,-1){8}}
\multiput(10,-2)(10,0){2}{\circle*{1}}
\multiput(10,8)(10,0){2}{\circle*{1}}
\put(15,3){\makebox(0,0){{\sm $\a \ab$}}}
\put(6,12){\makebox(0,0){\ss{*}}}
\put(24,12){\makebox(0,0){\ss{a}}}
\put(6,-6){\makebox(0,0){\s{b}}}
\put(24,-6){\makebox(0,0){$c$}}
\end{picture}
=
\thinlines
\unitlength 0.5mm
\begin{picture}(30,20)(0,0)
\multiput(11,-2)(0,10){2}{\line(1,0){8}}
\multiput(10,7)(10,0){2}{\line(0,-1){8}}
\multiput(10,-2)(10,0){2}{\circle*{1}}
\multiput(10,8)(10,0){2}{\circle*{1}}
\put(15,3){\makebox(0,0){$\a$}}
\put(6,12){\makebox(0,0){\ss{*}}}
\put(24,12){\makebox(0,0){\ss{a}}}
\put(6,-6){\makebox(0,0){\s{a}}}
\put(24,-6){\makebox(0,0){$4$}}
\end{picture}
\cdot
\thinlines
\unitlength 0.5mm
\begin{picture}(30,20)(0,0)
\multiput(11,-2)(0,10){2}{\line(1,0){8}}
\multiput(10,7)(10,0){2}{\line(0,-1){8}}
\multiput(10,-2)(10,0){2}{\circle*{1}}
\multiput(10,8)(10,0){2}{\circle*{1}}
\put(15,3){\makebox(0,0){$\ab$}}
\put(6,12){\makebox(0,0){\s{a}}}
\put(24,12){\makebox(0,0){$4$}}
\put(6,-6){\makebox(0,0){\s{b}}}
\put(24,-6){\makebox(0,0){$c$}}
\end{picture}
= 1,
$$
$$
\thinlines
\unitlength 0.5mm
\begin{picture}(30,20)(0,0)
\multiput(11,-2)(0,10){2}{\line(1,0){8}}
\multiput(10,7)(10,0){2}{\line(0,-1){8}}
\multiput(10,-2)(10,0){2}{\circle*{1}}
\multiput(10,8)(10,0){2}{\circle*{1}}
\put(15,3){\makebox(0,0){{\sm $\a \ab$}}}
\put(6,12){\makebox(0,0){\ss{*}}}
\put(24,12){\makebox(0,0){\ss{a}}}
\put(6,-6){\makebox(0,0){\s{b}}}
\put(24,-6){\makebox(0,0){\s{a}}}
\end{picture}
=
\thinlines
\unitlength 0.5mm
\begin{picture}(30,20)(0,0)
\multiput(11,-2)(0,10){2}{\line(1,0){8}}
\multiput(10,7)(10,0){2}{\line(0,-1){8}}
\multiput(10,-2)(10,0){2}{\circle*{1}}
\multiput(10,8)(10,0){2}{\circle*{1}}
\put(15,3){\makebox(0,0){$\a$}}
\put(6,12){\makebox(0,0){\ss{*}}}
\put(24,12){\makebox(0,0){\ss{a}}}
\put(6,-6){\makebox(0,0){\ss{a}}}
\put(24,-6){\makebox(0,0){$4$}}
\end{picture}
\cdot
\thinlines
\unitlength 0.5mm
\begin{picture}(30,20)(0,0)
\multiput(11,-2)(0,10){2}{\line(1,0){8}}
\multiput(10,7)(10,0){2}{\line(0,-1){8}}
\multiput(10,-2)(10,0){2}{\circle*{1}}
\multiput(10,8)(10,0){2}{\circle*{1}}
\put(15,3){\makebox(0,0){$\ab$}}
\put(6,12){\makebox(0,0){\ss{a}}}
\put(24,12){\makebox(0,0){$4$}}
\put(6,-6){\makebox(0,0){\s{b}}}
\put(24,-6){\makebox(0,0){\s{a}}}
\end{picture}
= -1.
$$
\phantom{x} \\
Next, we will obtain the entries marked \ci. We have two 
vectors of connection $\a \ab$ concerning to $b$-$b$ double edges by 
``actual'' multiplication as follows:
$$
\thinlines
\unitlength 0.5mm
\begin{picture}(30,20)(0,0)
\multiput(11,-2)(0,10){2}{\line(1,0){8}}
\multiput(10,7)(10,0){2}{\line(0,-1){8}}
\multiput(10,-2)(10,0){2}{\circle*{1}}
\multiput(10,8)(10,0){2}{\circle*{1}}
\put(15,3){\makebox(0,0){{\sm $\a \ab$}}}
\put(6,12){\makebox(0,0){$b$}}
\put(24,12){\makebox(0,0){$c$}}
\put(6,-6){\makebox(0,0){$b$}}
\put(24,-6){\makebox(0,0){$a$}}
\end{picture}
=
\left(
\begin{array}{c}
\thinlines
\unitlength 0.5mm
\begin{picture}(30,15)(0,0)
\multiput(11,-2)(0,10){2}{\line(1,0){8}}
\multiput(10,7)(10,0){2}{\line(0,-1){8}}
\multiput(10,-2)(10,0){2}{\circle*{1}}
\multiput(10,8)(10,0){2}{\circle*{1}}
\put(15,3){\makebox(0,0){$\a$}}
\put(6,12){\makebox(0,0){$b$}}
\put(24,12){\makebox(0,0){$c$}}
\put(6,-6){\makebox(0,0){$a$}}
\put(24,-6){\makebox(0,0){$2$}}
\end{picture}
\cdot
\thinlines
\unitlength 0.5mm
\begin{picture}(30,15)(0,0)
\multiput(11,-2)(0,10){2}{\line(1,0){8}}
\multiput(10,7)(10,0){2}{\line(0,-1){8}}
\multiput(10,-2)(10,0){2}{\circle*{1}}
\multiput(10,8)(10,0){2}{\circle*{1}}
\put(15,3){\makebox(0,0){$\ab$}}
\put(6,12){\makebox(0,0){$a$}}
\put(24,12){\makebox(0,0){$2$}}
\put(6,-6){\makebox(0,0){$b$}}
\put(24,-6){\makebox(0,0){$a$}}
\end{picture}
\\ \\
\thinlines
\unitlength 0.5mm
\begin{picture}(30,15)(0,0)
\multiput(11,-2)(0,10){2}{\line(1,0){8}}
\multiput(10,7)(10,0){2}{\line(0,-1){8}}
\multiput(10,-2)(10,0){2}{\circle*{1}}
\multiput(10,8)(10,0){2}{\circle*{1}}
\put(15,3){\makebox(0,0){$\a$}}
\put(6,12){\makebox(0,0){$b$}}
\put(24,12){\makebox(0,0){$c$}}
\put(6,-6){\makebox(0,0){$c$}}
\put(24,-6){\makebox(0,0){$2$}}
\end{picture}
\cdot
\thinlines
\unitlength 0.5mm
\begin{picture}(30,15)(0,0)
\multiput(11,-2)(0,10){2}{\line(1,0){8}}
\multiput(10,7)(10,0){2}{\line(0,-1){8}}
\multiput(10,-2)(10,0){2}{\circle*{1}}
\multiput(10,8)(10,0){2}{\circle*{1}}
\put(15,3){\makebox(0,0){$\ab$}}
\put(6,12){\makebox(0,0){$c$}}
\put(24,12){\makebox(0,0){$2$}}
\put(6,-6){\makebox(0,0){$b$}}
\put(24,-6){\makebox(0,0){$a$}}
\end{picture}
\end{array}
\right)
= 
\left(
\begin{array}{c}
\frac{-\sqrt{\lad-2}}{\lad-1} \\ \\
\frac{1}{\lad-1}
\end{array}
\right),
$$
$$
\thinlines
\unitlength 0.5mm
\begin{picture}(30,20)(0,0)
\multiput(11,-2)(0,10){2}{\line(1,0){8}}
\multiput(10,7)(10,0){2}{\line(0,-1){8}}
\multiput(10,-2)(10,0){2}{\circle*{1}}
\multiput(10,8)(10,0){2}{\circle*{1}}
\put(15,3){\makebox(0,0){{\sm $\a \ab$}}}
\put(6,12){\makebox(0,0){$b$}}
\put(24,12){\makebox(0,0){$a$}}
\put(6,-6){\makebox(0,0){$b$}}
\put(24,-6){\makebox(0,0){$c$}}
\end{picture}
=
\left(
\begin{array}{c}
\thinlines
\unitlength 0.5mm
\begin{picture}(30,15)(0,0)
\multiput(11,-2)(0,10){2}{\line(1,0){8}}
\multiput(10,7)(10,0){2}{\line(0,-1){8}}
\multiput(10,-2)(10,0){2}{\circle*{1}}
\multiput(10,8)(10,0){2}{\circle*{1}}
\put(15,3){\makebox(0,0){$\a$}}
\put(6,12){\makebox(0,0){$b$}}
\put(24,12){\makebox(0,0){$a$}}
\put(6,-6){\makebox(0,0){$a$}}
\put(24,-6){\makebox(0,0){$2$}}
\end{picture}
\cdot
\thinlines
\unitlength 0.5mm
\begin{picture}(30,15)(0,0)
\multiput(11,-2)(0,10){2}{\line(1,0){8}}
\multiput(10,7)(10,0){2}{\line(0,-1){8}}
\multiput(10,-2)(10,0){2}{\circle*{1}}
\multiput(10,8)(10,0){2}{\circle*{1}}
\put(15,3){\makebox(0,0){$\ab$}}
\put(6,12){\makebox(0,0){$a$}}
\put(24,12){\makebox(0,0){$2$}}
\put(6,-6){\makebox(0,0){$b$}}
\put(24,-6){\makebox(0,0){$c$}}
\end{picture}
\\ \\
\thinlines
\unitlength 0.5mm
\begin{picture}(30,15)(0,0)
\multiput(11,-2)(0,10){2}{\line(1,0){8}}
\multiput(10,7)(10,0){2}{\line(0,-1){8}}
\multiput(10,-2)(10,0){2}{\circle*{1}}
\multiput(10,8)(10,0){2}{\circle*{1}}
\put(15,3){\makebox(0,0){$\a$}}
\put(6,12){\makebox(0,0){$b$}}
\put(24,12){\makebox(0,0){$a$}}
\put(6,-6){\makebox(0,0){$c$}}
\put(24,-6){\makebox(0,0){$2$}}
\end{picture}
\cdot
\thinlines
\unitlength 0.5mm
\begin{picture}(30,15)(0,0)
\multiput(11,-2)(0,10){2}{\line(1,0){8}}
\multiput(10,7)(10,0){2}{\line(0,-1){8}}
\multiput(10,-2)(10,0){2}{\circle*{1}}
\multiput(10,8)(10,0){2}{\circle*{1}}
\put(15,3){\makebox(0,0){$\ab$}}
\put(6,12){\makebox(0,0){$c$}}
\put(24,12){\makebox(0,0){$2$}}
\put(6,-6){\makebox(0,0){$b$}}
\put(24,-6){\makebox(0,0){$c$}}
\end{picture}
\end{array}
\right)
= 
\left(
\begin{array}{c}
\frac{-1}{\lad-1} \\ \\
\frac{1}{(\lad-1)\sqrt{\lad-2}}
\end{array}
\right).
$$
\pha{x} \par \noindent
Since these two vectors are proportional, they are transformed to 
two proportional 
vectors by a left vertical gauge transform for the double edges $b$-$b$,
 i.e., multiplication from 
the left by an element of $U(2)$. Since we should have 1's in the
$(bb^{1}, aa^{1})$-entry and the $(bb^{1}, cc^{1})$-entry, they can be 
transformed into the following pair,
$$
\left( \begin{array}{c}
0 \\ \frac{1}{\sqrt{\lad-1}} 
\end{array} \right), \quad {\rm and} \quad
\left( \begin{array}{c}
0 \\ \frac{\sqrt{\lad-4}}{\sqrt{\lad-2}} 
\end{array} \right),
$$
then we have 
$ (bb^{2},ca)=\frac{1}{\sqrt{\lad-1}}, (bb^{2},ac)=
\frac{\sqrt{\lad-4}}{\sqrt{\lad-2}}$ respectively, where we have 
omitted ``-entry''. 
The same procedure for the entries with vertical 
double edges \s{b}-\s{b} and \ss{b}-\ss{b} gives two pairs of vectors as
follows.
$$
\thinlines
\unitlength 0.5mm
\begin{picture}(30,20)(0,0)
\multiput(11,-2)(0,10){2}{\line(1,0){8}}
\multiput(10,7)(10,0){2}{\line(0,-1){8}}
\multiput(10,-2)(10,0){2}{\circle*{1}}
\multiput(10,8)(10,0){2}{\circle*{1}}
\put(15,3){\makebox(0,0){{\sm $\a \ab$}}}
\put(6,12){\makebox(0,0){\s{b}}}
\put(24,12){\makebox(0,0){\s{a}}}
\put(6,-6){\makebox(0,0){\s{b}}}
\put(24,-6){\makebox(0,0){$c$}}
\end{picture}
=
\left(
\begin{array}{c}
\thinlines
\unitlength 0.5mm
\begin{picture}(30,15)(0,0)
\multiput(11,-2)(0,10){2}{\line(1,0){8}}
\multiput(10,7)(10,0){2}{\line(0,-1){8}}
\multiput(10,-2)(10,0){2}{\circle*{1}}
\multiput(10,8)(10,0){2}{\circle*{1}}
\put(15,3){\makebox(0,0){$\a$}}
\put(6,12){\makebox(0,0){\s{b}}}
\put(24,12){\makebox(0,0){\s{a}}}
\put(6,-6){\makebox(0,0){$c$}}
\put(24,-6){\makebox(0,0){$4$}}
\end{picture}
\cdot
\thinlines
\unitlength 0.5mm
\begin{picture}(30,15)(0,0)
\multiput(11,-2)(0,10){2}{\line(1,0){8}}
\multiput(10,7)(10,0){2}{\line(0,-1){8}}
\multiput(10,-2)(10,0){2}{\circle*{1}}
\multiput(10,8)(10,0){2}{\circle*{1}}
\put(15,3){\makebox(0,0){$\ab$}}
\put(6,12){\makebox(0,0){$c$}}
\put(24,12){\makebox(0,0){$4$}}
\put(6,-6){\makebox(0,0){\s{b}}}
\put(24,-6){\makebox(0,0){$c$}}
\end{picture}
\\ \\
\thinlines
\unitlength 0.5mm
\begin{picture}(30,15)(0,0)
\multiput(11,-2)(0,10){2}{\line(1,0){8}}
\multiput(10,7)(10,0){2}{\line(0,-1){8}}
\multiput(10,-2)(10,0){2}{\circle*{1}}
\multiput(10,8)(10,0){2}{\circle*{1}}
\put(15,3){\makebox(0,0){$\a$}}
\put(6,12){\makebox(0,0){\s{b}}}
\put(24,12){\makebox(0,0){\s{a}}}
\put(6,-6){\makebox(0,0){\s{a}}}
\put(24,-6){\makebox(0,0){$4$}}
\end{picture}
\cdot
\thinlines
\unitlength 0.5mm
\begin{picture}(30,15)(0,0)
\multiput(11,-2)(0,10){2}{\line(1,0){8}}
\multiput(10,7)(10,0){2}{\line(0,-1){8}}
\multiput(10,-2)(10,0){2}{\circle*{1}}
\multiput(10,8)(10,0){2}{\circle*{1}}
\put(15,3){\makebox(0,0){$\ab$}}
\put(6,12){\makebox(0,0){\s{a}}}
\put(24,12){\makebox(0,0){$4$}}
\put(6,-6){\makebox(0,0){\s{b}}}
\put(24,-6){\makebox(0,0){$c$}}
\end{picture}
\end{array}
\right)
= 
\left(
\begin{array}{c}
\frac{1}{\sqrt{(\lad-1)(\lad-2)}} \\ \\
\frac{-1}{\sqrt{\lad-1}}
\end{array}
\right),
$$
$$
\thinlines
\unitlength 0.5mm
\begin{picture}(30,20)(0,0)
\multiput(11,-2)(0,10){2}{\line(1,0){8}}
\multiput(10,7)(10,0){2}{\line(0,-1){8}}
\multiput(10,-2)(10,0){2}{\circle*{1}}
\multiput(10,8)(10,0){2}{\circle*{1}}
\put(15,3){\makebox(0,0){{\sm $\a \ab$}}}
\put(6,12){\makebox(0,0){\s{b}}}
\put(24,12){\makebox(0,0){$c$}}
\put(6,-6){\makebox(0,0){\s{b}}}
\put(24,-6){\makebox(0,0){\s{a}}}
\end{picture}
=
\left(
\begin{array}{c}
\thinlines
\unitlength 0.5mm
\begin{picture}(30,15)(0,0)
\multiput(11,-2)(0,10){2}{\line(1,0){8}}
\multiput(10,7)(10,0){2}{\line(0,-1){8}}
\multiput(10,-2)(10,0){2}{\circle*{1}}
\multiput(10,8)(10,0){2}{\circle*{1}}
\put(15,3){\makebox(0,0){$\a$}}
\put(6,12){\makebox(0,0){\s{b}}}
\put(24,12){\makebox(0,0){$c$}}
\put(6,-6){\makebox(0,0){$c$}}
\put(24,-6){\makebox(0,0){$4$}}
\end{picture}
\cdot
\thinlines
\unitlength 0.5mm
\begin{picture}(30,15)(0,0)
\multiput(11,-2)(0,10){2}{\line(1,0){8}}
\multiput(10,7)(10,0){2}{\line(0,-1){8}}
\multiput(10,-2)(10,0){2}{\circle*{1}}
\multiput(10,8)(10,0){2}{\circle*{1}}
\put(15,3){\makebox(0,0){$\ab$}}
\put(6,12){\makebox(0,0){$c$}}
\put(24,12){\makebox(0,0){$4$}}
\put(6,-6){\makebox(0,0){\s{b}}}
\put(24,-6){\makebox(0,0){\s{a}}}
\end{picture}
\\ \\
\thinlines
\unitlength 0.5mm
\begin{picture}(30,15)(0,0)
\multiput(11,-2)(0,10){2}{\line(1,0){8}}
\multiput(10,7)(10,0){2}{\line(0,-1){8}}
\multiput(10,-2)(10,0){2}{\circle*{1}}
\multiput(10,8)(10,0){2}{\circle*{1}}
\put(15,3){\makebox(0,0){$\a$}}
\put(6,12){\makebox(0,0){\s{b}}}
\put(24,12){\makebox(0,0){$c$}}
\put(6,-6){\makebox(0,0){\s{a}}}
\put(24,-6){\makebox(0,0){$4$}}
\end{picture}
\cdot
\thinlines
\unitlength 0.5mm
\begin{picture}(30,15)(0,0)
\multiput(11,-2)(0,10){2}{\line(1,0){8}}
\multiput(10,7)(10,0){2}{\line(0,-1){8}}
\multiput(10,-2)(10,0){2}{\circle*{1}}
\multiput(10,8)(10,0){2}{\circle*{1}}
\put(15,3){\makebox(0,0){$\ab$}}
\put(6,12){\makebox(0,0){\s{a}}}
\put(24,12){\makebox(0,0){$4$}}
\put(6,-6){\makebox(0,0){\s{b}}}
\put(24,-6){\makebox(0,0){\s{a}}}
\end{picture}
\end{array}
\right)
= 
\left(
\begin{array}{c}
\frac{1}{\sqrt{\lad-1}} \\ \\
\frac{-\sqrt{\lad-2}}{\sqrt{\lad-1}}
\end{array}
\right),
$$
\pha{x} \par \noindent
$$
\thinlines
\unitlength 0.5mm
\begin{picture}(30,20)(0,0)
\multiput(11,-2)(0,10){2}{\line(1,0){8}}
\multiput(10,7)(10,0){2}{\line(0,-1){8}}
\multiput(10,-2)(10,0){2}{\circle*{1}}
\multiput(10,8)(10,0){2}{\circle*{1}}
\put(15,3){\makebox(0,0){{\sm $\a \ab$}}}
\put(6,12){\makebox(0,0){\ss{b}}}
\put(24,12){\makebox(0,0){\ss{a}}}
\put(6,-6){\makebox(0,0){\ss{b}}}
\put(24,-6){\makebox(0,0){$c$}}
\end{picture}
=
\left(
\begin{array}{c}
\thinlines
\unitlength 0.5mm
\begin{picture}(30,15)(0,0)
\multiput(11,-2)(0,10){2}{\line(1,0){8}}
\multiput(10,7)(10,0){2}{\line(0,-1){8}}
\multiput(10,-2)(10,0){2}{\circle*{1}}
\multiput(10,8)(10,0){2}{\circle*{1}}
\put(15,3){\makebox(0,0){$\a$}}
\put(6,12){\makebox(0,0){\ss{b}}}
\put(24,12){\makebox(0,0){\ss{a}}}
\put(6,-6){\makebox(0,0){$c$}}
\put(24,-6){\makebox(0,0){$4$}}
\end{picture}
\cdot
\thinlines
\unitlength 0.5mm
\begin{picture}(30,15)(0,0)
\multiput(11,-2)(0,10){2}{\line(1,0){8}}
\multiput(10,7)(10,0){2}{\line(0,-1){8}}
\multiput(10,-2)(10,0){2}{\circle*{1}}
\multiput(10,8)(10,0){2}{\circle*{1}}
\put(15,3){\makebox(0,0){$\ab$}}
\put(6,12){\makebox(0,0){$c$}}
\put(24,12){\makebox(0,0){$4$}}
\put(6,-6){\makebox(0,0){\ss{b}}}
\put(24,-6){\makebox(0,0){$c$}}
\end{picture}
\\ \\
\thinlines
\unitlength 0.5mm
\begin{picture}(30,15)(0,0)
\multiput(11,-2)(0,10){2}{\line(1,0){8}}
\multiput(10,7)(10,0){2}{\line(0,-1){8}}
\multiput(10,-2)(10,0){2}{\circle*{1}}
\multiput(10,8)(10,0){2}{\circle*{1}}
\put(15,3){\makebox(0,0){$\a$}}
\put(6,12){\makebox(0,0){\ss{b}}}
\put(24,12){\makebox(0,0){\ss{a}}}
\put(6,-6){\makebox(0,0){\ss{a}}}
\put(24,-6){\makebox(0,0){$4$}}
\end{picture}
\cdot
\thinlines
\unitlength 0.5mm
\begin{picture}(30,15)(0,0)
\multiput(11,-2)(0,10){2}{\line(1,0){8}}
\multiput(10,7)(10,0){2}{\line(0,-1){8}}
\multiput(10,-2)(10,0){2}{\circle*{1}}
\multiput(10,8)(10,0){2}{\circle*{1}}
\put(15,3){\makebox(0,0){$\ab$}}
\put(6,12){\makebox(0,0){\ss{a}}}
\put(24,12){\makebox(0,0){$4$}}
\put(6,-6){\makebox(0,0){\ss{b}}}
\put(24,-6){\makebox(0,0){$c$}}
\end{picture}
\end{array}
\right)
= 
\left(
\begin{array}{c}
\frac{1}{\sqrt{(\lad-1)(\lad-2)}} \\ \\
\frac{-1}{\sqrt{\lad-1}}
\end{array}
\right),
$$
\pha{x} \par \noindent
$$
\thinlines
\unitlength 0.5mm
\begin{picture}(30,20)(0,0)
\multiput(11,-2)(0,10){2}{\line(1,0){8}}
\multiput(10,7)(10,0){2}{\line(0,-1){8}}
\multiput(10,-2)(10,0){2}{\circle*{1}}
\multiput(10,8)(10,0){2}{\circle*{1}}
\put(15,3){\makebox(0,0){{\sm $\a \ab$}}}
\put(6,12){\makebox(0,0){\ss{b}}}
\put(24,12){\makebox(0,0){$c$}}
\put(6,-6){\makebox(0,0){\ss{b}}}
\put(24,-6){\makebox(0,0){\ss{a}}}
\end{picture}
=
\left(
\begin{array}{c}
\thinlines
\unitlength 0.5mm
\begin{picture}(30,15)(0,0)
\multiput(11,-2)(0,10){2}{\line(1,0){8}}
\multiput(10,7)(10,0){2}{\line(0,-1){8}}
\multiput(10,-2)(10,0){2}{\circle*{1}}
\multiput(10,8)(10,0){2}{\circle*{1}}
\put(15,3){\makebox(0,0){$\a$}}
\put(6,12){\makebox(0,0){\ss{b}}}
\put(24,12){\makebox(0,0){\ss{a}}}
\put(6,-6){\makebox(0,0){$c$}}
\put(24,-6){\makebox(0,0){$4$}}
\end{picture}
\cdot
\thinlines
\unitlength 0.5mm
\begin{picture}(30,15)(0,0)
\multiput(11,-2)(0,10){2}{\line(1,0){8}}
\multiput(10,7)(10,0){2}{\line(0,-1){8}}
\multiput(10,-2)(10,0){2}{\circle*{1}}
\multiput(10,8)(10,0){2}{\circle*{1}}
\put(15,3){\makebox(0,0){$\ab$}}
\put(6,12){\makebox(0,0){$c$}}
\put(24,12){\makebox(0,0){$4$}}
\put(6,-6){\makebox(0,0){\ss{b}}}
\put(24,-6){\makebox(0,0){\ss{a}}}
\end{picture}
\\ \\
\thinlines
\unitlength 0.5mm
\begin{picture}(30,15)(0,0)
\multiput(11,-2)(0,10){2}{\line(1,0){8}}
\multiput(10,7)(10,0){2}{\line(0,-1){8}}
\multiput(10,-2)(10,0){2}{\circle*{1}}
\multiput(10,8)(10,0){2}{\circle*{1}}
\put(15,3){\makebox(0,0){$\a$}}
\put(6,12){\makebox(0,0){\ss{b}}}
\put(24,12){\makebox(0,0){$c$}}
\put(6,-6){\makebox(0,0){\ss{a}}}
\put(24,-6){\makebox(0,0){$4$}}
\end{picture}
\cdot
\thinlines
\unitlength 0.5mm
\begin{picture}(30,15)(0,0)
\multiput(11,-2)(0,10){2}{\line(1,0){8}}
\multiput(10,7)(10,0){2}{\line(0,-1){8}}
\multiput(10,-2)(10,0){2}{\circle*{1}}
\multiput(10,8)(10,0){2}{\circle*{1}}
\put(15,3){\makebox(0,0){$\ab$}}
\put(6,12){\makebox(0,0){\ss{a}}}
\put(24,12){\makebox(0,0){$4$}}
\put(6,-6){\makebox(0,0){\ss{b}}}
\put(24,-6){\makebox(0,0){\ss{a}}}
\end{picture}
\end{array}
\right)
= 
\left(
\begin{array}{c}
\frac{1}{\sqrt{\lad-1}} \\ \\
\frac{-\sqrt{\lad-2}}{\sqrt{\lad-1}}
\end{array}
\right).
$$
\pha{x} \par \noindent
The first pair concerning to \s{b}-\s{b} double edges can be 
transformed by gauge unitary into the pair
$$
\left( \begin{array}{c}
0 \\  \frac{1}{\sqrt{\lad-2}}
\end{array} \right)
\quad {\rm and} \quad
\left( \begin{array}{c}
0 \\  1 
\end{array}  \right),
$$
where we note $(\si{b}\si{b}^{2}, cc^{1})=(\si{b}\si{b}^{2}, 
\si{a}\si{a}^{1})=1$ by this gauge.
Since second pair is equal to the first pair, it can be 
transformed the same pair of vectors by the left gauge unitary of the 
double edges \ss{b}-\ss{b}. Thus, we have
$$(\si{b}\si{b}^{2}, \si{a}c)=(\ssi{b}\ssi{b}^{2}, \ssi{a}c)
=\frac{1}{\sqrt{\lad-2}} $$
and
$$
(\si{b}\si{b}^{2}, c\si{a})=(\ssi{b}\ssi{b}^{2}, c\ssi{a})=1. $$
Along the same argument, we have
$$(*b,aa^{1})=1,\quad (b*,aa^{1})=\frac{1}{\sqrt{\lad-1}}$$
from the ``actual'' multiplications
\begeq
\thinlines
\unitlength 0.5mm
\begin{picture}(30,20)(0,0)
\multiput(11,-2)(0,10){2}{\line(1,0){8}}
\multiput(10,7)(10,0){2}{\line(0,-1){8}}
\multiput(10,-2)(10,0){2}{\circle*{1}}
\multiput(10,8)(10,0){2}{\circle*{1}}
\put(15,3){\makebox(0,0){{\sm $\a \ab$}}}
\put(6,12){\makebox(0,0){$*$}}
\put(24,12){\makebox(0,0){$a$}}
\put(6,-6){\makebox(0,0){$b$}}
\put(24,-6){\makebox(0,0){$a$}}
\end{picture}
&=&
\left( \frac{\sqrt{\lad-1}}{\la}, \frac{-1}{\la} \right), \\
\thinlines
\unitlength 0.5mm
\begin{picture}(30,20)(0,0)
\multiput(11,-2)(0,10){2}{\line(1,0){8}}
\multiput(10,7)(10,0){2}{\line(0,-1){8}}
\multiput(10,-2)(10,0){2}{\circle*{1}}
\multiput(10,8)(10,0){2}{\circle*{1}}
\put(15,3){\makebox(0,0){{\sm $\a \ab$}}}
\put(6,12){\makebox(0,0){$b$}}
\put(24,12){\makebox(0,0){$a$}}
\put(6,-6){\makebox(0,0){$*$}}
\put(24,-6){\makebox(0,0){$a$}}
\end{picture}
&=&
\left( \frac{1}{\la}, \frac{-1}{\la\sqrt{\lad-1}} \right),
\endeq
by the right gauge unitary of double edges $a$-$a$. \\
So far, we have the entries of $\a \ab$ as in \framebox{Table
\ref{aab13.tex}}, 
where the entries $g_{**}$'s mean that they have not been 
determined so far.


\begin{table}[h]
\begin{center}
\begin{tabular}[pos]{|@{}c@{}||@{}c@{}|@{}c@{}|@{}c@{}||@{} c@{}|@{}c@{}
|@{}c@{}| @{}c@{}|@{}c@{}|@{}c@{}|| @{}c@{}|@{}c@{}|@{}c@{}|| @{}c@{}
|@{}c@{}|@{}c@{}|c }
\hl
& $aa^1$ & $aa^2$ & $ca$ & $ac$ & $cc^1$ & $cc^2$ & $cc^3$ & 
\s{a}$c$ &
\ss{a}$c$ & $c$\s{a} & \s{a}\s{a} & \ss{a}\s{a} & $c$\ss{a} & 
\s{a}\ss{a} & \ss{a}\ss{a} \\ \hline\hline
$**$ & 1 & 0 &&&&&&&&&&&&&  \\ \hl
$*b$ & 0 & 1 && 1 &&&&&&&&&&&  \\ \hline\hline
$b*$ & 0 & $\frac{1}{\la_1}$ & $\frac{\la_2}{\la_1}$ &
&&&&&&&&&&& \\ \hl
$bb^1$ & 1 & 0 & 0 &0 &1&0&0&&&&&&&& \\ \hl
$bb^2$ & 0 & $\frac{-\la_2}{\la_1}$ & 
$\frac{1}{\la_1}$ & $\frac{\la_4}{\la_2}$ & 0 & 
\mcol{2}{@{}c@{}|}{$g_{bb}$} & &&&&&&& \\ \hl
$b$\s{b} & &&& $\frac{\bar{\rho}}{\la_2}$ & 0 & 
\mcol{2}{@{}c@{}|}{$g_{b b_{\sigma}}$} &&&1&&&&& \\ \hl
$b$\ss{b} &&&& $\frac{\rho}{\la_2}$ & 0 & 
\mcol{2}{@{}c@{}|}{$g_{b b_{\sigma^2}}$} & && &&& 1 && \\ \hline\hline
\s{b}$b$ & && $\rho$ & & 0 & \mcol{2}{@{}c@{}|}{$g_{b_{\sigma}b}$} & 
$\frac{1}{\la_2}$ & & &&&&& \\ \hl
\s{b}\s{b}$^1$ & &&&& 1 & 0 & 0 & 0 & & 0 & 1 &&&& \\ \hl
\s{b}\s{b}$^2$ & &&&& 0 & \mcol{2}{@{}c@{}|}{$g_{b_{\sigma} b_{\sigma}}$} & 
$\frac{1}{\la_2}$ & & 1 & 0 &&&& \\ \hl
\s{b}\ss{b} & &&&& 0 & \mcol{2}{@{}c@{}|}{$g_{b_{\sigma} b_{\sigma^2}}$} &
$\frac{\la_4}{\la_2}$ & &&&& $\frac{1}{\la_1}$ & 
$\frac{\la_2}{\la_1}$ & \\ \hl
\s{b}\ss{*} &&&&& &&&&& &&& $\frac{\la_2}{\la_1}$ & 
$\frac{-1}{\la_1}$ & \\ \hline\hline
\ss{b}$b$ &&& $\bar{\rho}$ & &0 & \mcol{2}{@{}c@{}|}{$g_{b_{\sigma^2} b}$} & & 
$\frac{1}{\la_2}$ & &&&&& \\ \hl
\ss{b}\s{b} &&&&& 0 & \mcol{2}{@{}c@{}|}{$g_{b_{\sigma^2} b_{\sigma}}$} &&
$\frac{\la_4}{\la_2}$ & $\frac{1}{\la_1}$ & 0 &
$\frac{\la_2}{\la_1}$ &&& \\ \hl
\ss{b}\ss{b}$^1$ &&&&& 1 & 0 & 0 & & 0 & & & & 0 & & 1 \\ \hl
\ss{b}\ss{b}$^2$ &&&&& 0 & \mcol{2}{@{}c@{}|}{$g_{b_{\sigma^2} b_{\sigma^2}}$} &&
$\frac{1}{\la_2}$ &&&& 1 &&0 \\ \hl
\ss{b}\s{*} & &&&& &&&&& $\frac{\la_2}{\la_1}$ & 0 & 
$\frac{-1}{\la_1}$ & && \\ \hline\hline
\s{*}\ss{b} &&&& &&&& 1 & & &&&& $-1$ & \\ \hl
\s{*}\s{*} &&&&& &&&&& & 1 &&&& \\ \hline\hline
\ss{*}\s{b} &&&& &&&& &1 &&& $-1$ &&& \\ \hl
\ss{*}\ss{*} &&&&& &&&&& &&&&& 1 \\ \hl
\end{tabular}
\end{center}
\caption{Connection $\alpha{\tilde \alpha}$ ($\la_n=\sqrt{\lad-n}$)}
\label{aab13.tex}
\end{table}
Now we will obtain the entries marked \dia \ in 
Table \ref{land13.tex}. Denote the vectors of entries of ``actual'' 
multiplication $\a \ab$ corresponding to $\left( 0, g_{??} \right)$
by $f_{??}$. We use the following data of $f_{??}$'s.
\begeq
 f_{b\si{b}}=
\thinlines
\unitlength 0.5mm
\begin{picture}(30,15)(0,0)
\multiput(11,-2)(0,10){2}{\line(1,0){8}}
\multiput(10,7)(10,0){2}{\line(0,-1){8}}
\multiput(10,-2)(10,0){2}{\circle*{1}}
\multiput(10,8)(10,0){2}{\circle*{1}}
\put(15,3){\makebox(0,0){\mbox{\sm $\a \ab$}}}
\put(6,12){\makebox(0,0){$b$}}
\put(24,12){\makebox(0,0){$c$}}
\put(6,-6){\makebox(0,0){\s{b}}}
\put(24,-6){\makebox(0,0){$c$}}
\end{picture}
&=&
\left(
\thinlines
\unitlength 0.5mm
\begin{picture}(30,15)(0,0)
\multiput(11,-2)(0,10){2}{\line(1,0){8}}
\multiput(10,7)(10,0){2}{\line(0,-1){8}}
\multiput(10,-2)(10,0){2}{\circle*{1}}
\multiput(10,8)(10,0){2}{\circle*{1}}
\put(15,3){\makebox(0,0){$\a$}}
\put(6,12){\makebox(0,0){$b$}}
\put(24,12){\makebox(0,0){$c$}}
\put(6,-6){\makebox(0,0){$c$}}
\put(24,-6){\makebox(0,0){$2$}}
\end{picture}
\cdot
\thinlines
\unitlength 0.5mm
\begin{picture}(30,15)(0,0)
\multiput(11,-2)(0,10){2}{\line(1,0){8}}
\multiput(10,7)(10,0){2}{\line(0,-1){8}}
\multiput(10,-2)(10,0){2}{\circle*{1}}
\multiput(10,8)(10,0){2}{\circle*{1}}
\put(15,3){\makebox(0,0){$\ab$}}
\put(6,12){\makebox(0,0){$c$}}
\put(24,12){\makebox(0,0){$2$}}
\put(6,-6){\makebox(0,0){\s{b}}}
\put(24,-6){\makebox(0,0){$c$}}
\end{picture},
\thinlines
\unitlength 0.5mm
\begin{picture}(30,15)(0,0)
\multiput(11,-2)(0,10){2}{\line(1,0){8}}
\multiput(10,7)(10,0){2}{\line(0,-1){8}}
\multiput(10,-2)(10,0){2}{\circle*{1}}
\multiput(10,8)(10,0){2}{\circle*{1}}
\put(15,3){\makebox(0,0){$\a$}}
\put(6,12){\makebox(0,0){$b$}}
\put(24,12){\makebox(0,0){$c$}}
\put(6,-6){\makebox(0,0){$c$}}
\put(24,-6){\makebox(0,0){$3$}}
\end{picture}
\cdot
\thinlines
\unitlength 0.5mm
\begin{picture}(30,15)(0,0)
\multiput(11,-2)(0,10){2}{\line(1,0){8}}
\multiput(10,7)(10,0){2}{\line(0,-1){8}}
\multiput(10,-2)(10,0){2}{\circle*{1}}
\multiput(10,8)(10,0){2}{\circle*{1}}
\put(15,3){\makebox(0,0){$\ab$}}
\put(6,12){\makebox(0,0){$c$}}
\put(24,12){\makebox(0,0){$3$}}
\put(6,-6){\makebox(0,0){\s{b}}}
\put(24,-6){\makebox(0,0){$c$}}
\end{picture},
\thinlines
\unitlength 0.5mm
\begin{picture}(30,15)(0,0)
\multiput(11,-2)(0,10){2}{\line(1,0){8}}
\multiput(10,7)(10,0){2}{\line(0,-1){8}}
\multiput(10,-2)(10,0){2}{\circle*{1}}
\multiput(10,8)(10,0){2}{\circle*{1}}
\put(15,3){\makebox(0,0){$\a$}}
\put(6,12){\makebox(0,0){$b$}}
\put(24,12){\makebox(0,0){$c$}}
\put(6,-6){\makebox(0,0){$c$}}
\put(24,-6){\makebox(0,0){$4$}}
\end{picture}
\cdot
\thinlines
\unitlength 0.5mm
\begin{picture}(30,15)(0,0)
\multiput(11,-2)(0,10){2}{\line(1,0){8}}
\multiput(10,7)(10,0){2}{\line(0,-1){8}}
\multiput(10,-2)(10,0){2}{\circle*{1}}
\multiput(10,8)(10,0){2}{\circle*{1}}
\put(15,3){\makebox(0,0){$\ab$}}
\put(6,12){\makebox(0,0){$c$}}
\put(24,12){\makebox(0,0){$4$}}
\put(6,-6){\makebox(0,0){\s{b}}}
\put(24,-6){\makebox(0,0){$c$}}
\end{picture}
\right) \\
&=&
\left(
\frac{\rho}{\la(\lad-2)},\frac{\tau}{\sqrt{3}}, \frac{1}{\lad-2}
\right),
\endeq
$$f_{b\ssi{b}}
=\left( \frac{\bar{\rho}}{\la(\lad-2)}, \frac{\bar \tau}{\sqrt{3}},
\frac{1}{\lad-2} \right), $$

$$f_{\si{b}b}
=\left( \frac{\bar{\rho}}{\la(\lad-2)}, \frac{\bar \tau}{\sqrt{3}},
\frac{1}{\lad-2} \right), $$

$$f_{\si{b}\ssi{b}}=
\left( \frac{(\lad-1){\bar \rho}^{2}}{\la(\lad-2)},
\frac{{\bar \tau}^{2}}{\sqrt{3}}, \frac{1}{(\lad-2)\sqrt{\lad-1}} 
\right),
$$
$$ f_{\ssi{b}b}=
\left(
\frac{\rho}{\la(\lad-2)},\frac{\tau}{\sqrt{3}}, \frac{1}{\lad-2}
\right), $$
$$f_{\ssi{b}\si{b}}=
\left( \frac{(\lad-1){\bar \rho}^{2}}{\la(\lad-2)},
\frac{{\bar \tau}^{2}}{\sqrt{3}}, \frac{1}{(\lad-2)\sqrt{\lad-1}} 
\right). $$
Note that $f_{b\ssi{b}}=f_{\si{b}b}$ and $f_{\ssi{b}b}=f_{b\si{b}}$, 
so they are transformed with keeping equality by the gauge transform 
of the triple edges $c$-$c$. Therefore, we see only $f_{b\si{b}}$, 
$f_{\si{b}b}$, $f_{\si{b}\ssi{b}}$ and $f_{\ssi{b}\si{b}}$. We have 
the following lemma.

\begin{lm}
The three vectors
$$u_{1}=\left( \frac{1}{\sqrt{3}}, \frac{1}{\la}, 
\frac{1}{\sqrt{\lad-2}} \right), $$
$$u_{2}=\left(\frac{\sqrt{\lad-2}}{3}, 
\frac{\sqrt{\lad-2}}{\la\sqrt{3}}, -\frac{\lad-3}{\sqrt{3}} \right), $$
and
$$\left( \frac{\lad-2}{\la\sqrt{3}}, -\frac{\lad-2}{3},0 \right)$$
form an orthonormal basis for ${\bf C}^{3}$ and 
$$
f_{\ssi{b}b}=f_{b\si{b}}=\frac{-1}{\sqrt{3}} u_{2}+
(\frac{\sqrt{\lad-2}}{2\la}+
\frac{\lad-3}{2\la}i) u_{3}, $$
$$
f_{b\ssi{b}}=f_{\si{b}b}=\frac{-1}{\sqrt{3}}u_{2}+ (\frac{\sqrt{\lad-2}}{2\la}-
\frac{\lad-3}{2\la}i) u_{3}, $$
$$f_{\si{b}\ssi{b}}= -\sqrt{\frac{\lad-4}{3}}u_{2}+
(-\frac{\lad-2}{2\sqrt{3}}+\frac{\sqrt{\lad-3}}{2}i)u_{3}, $$
$$f_{\ssi{b}\si{b}}= -\sqrt{\frac{\lad-4}{3}}u_{2}+
(-\frac{\lad-2}{2\sqrt{3}}-\frac{\sqrt{\lad-3}}{2}i)u_{3}. $$
\label{basis}
\end{lm}
{\it Proof} \\
Checked by elementary, but heavy computations, using 
$\la^{4}-5\lad+3=0$. \\
\qed \\

From Lemma \ref{basis}, we have $g_{??}$'s as the expression of 
$f_{??}$'s by the orthonormal basis $u_{2}$ and $u_{3}$ as follows,
\begeq
 g_{b \si{b}}&=&g_{\ssi{b}b}
=\left( \frac{-1}{\sqrt{3}}, \frac{\sqrt{\lad-2}}{2\la}+
\frac{\lad-3}{2\la}i \right), \\
g_{\si{b}b}&=&g_{b \ssi{b}}=
\left( \frac{-1}{\sqrt{3}}, \frac{\sqrt{\lad-2}}{2\la}-
\frac{\lad-3}{2\la}i \right), \\
g_{\si{b}\ssi{b}}&=& \left(-\sqrt{\frac{\lad-4}{3}}, 
-\frac{\lad-2}{2\sqrt{3}}+\frac{\sqrt{\lad-3}}{2}i \right), \\
g_{\ssi{b}\si{b}}&=& \left( -\sqrt{\frac{\lad-4}{3}}, 
-\frac{\lad-2}{2\sqrt{3}}-\frac{\sqrt{\lad-3}}{2}i \right).
\endeq
$g_{bb},g_{\si{b}\si{b}}$ and $g_{\ssi{b}\ssi{b}}$ are uniquely 
determined so that the matrices
$$
\thinlines
\unitlength 0.5mm
\begin{picture}(30,20)(0,0)
\multiput(11,-2)(0,10){2}{\line(1,0){8}}
\multiput(10,7)(10,0){2}{\line(0,-1){8}}
\multiput(10,-2)(10,0){2}{\circle*{1}}
\multiput(10,8)(10,0){2}{\circle*{1}}
\put(15,3){\makebox(0,0){{\sm $\a \ab$}}}
\put(6,12){\makebox(0,0){$b$}}
\put(24,12){\makebox(0,0){}}
\put(6,-6){\makebox(0,0){}}
\put(24,-6){\makebox(0,0){$c$}}
\end{picture}
=
\left( \begin{array}{cc}
\sqrt{\frac{\lad-4}{\lad-2}} & g_{bb} \\
\frac{\bar{\rho}}{\sqrt{\lad-2}} & g_{b\si{b}} \\
\frac{\rho}{\sqrt{\lad-2}}  & g_{b\ssi{b}} 
\end{array} \right),
$$ 
$$
\thinlines
\unitlength 0.5mm
\begin{picture}(30,20)(0,0)
\multiput(11,-2)(0,10){2}{\line(1,0){8}}
\multiput(10,7)(10,0){2}{\line(0,-1){8}}
\multiput(10,-2)(10,0){2}{\circle*{1}}
\multiput(10,8)(10,0){2}{\circle*{1}}
\put(15,3){\makebox(0,0){{\sm $\a \ab$}}}
\put(6,12){\makebox(0,0){\s{b}}}
\put(24,12){\makebox(0,0){}}
\put(6,-6){\makebox(0,0){}}
\put(24,-6){\makebox(0,0){$c$}}
\end{picture}
=
\left( \begin{array}{cc}
  g_{\si{b}b} & \frac{1}{\sqrt{\lad-2}} \\
g_{\si{b}\si{b}} & \frac{1}{\sqrt{\lad-2}}\\
   g_{\si{b}\ssi{b}} &\sqrt{\frac{\lad-4}{\lad-2}} 
   \end{array} \right),
$$
and
$$
\thinlines
\unitlength 0.5mm
\begin{picture}(30,20)(0,0)
\multiput(11,-2)(0,10){2}{\line(1,0){8}}
\multiput(10,7)(10,0){2}{\line(0,-1){8}}
\multiput(10,-2)(10,0){2}{\circle*{1}}
\multiput(10,8)(10,0){2}{\circle*{1}}
\put(15,3){\makebox(0,0){{\sm $\a \ab$}}}
\put(6,12){\makebox(0,0){\ss{b}}}
\put(24,12){\makebox(0,0){}}
\put(6,-6){\makebox(0,0){}}
\put(24,-6){\makebox(0,0){$c$}}
\end{picture}
=
\left( \begin{array}{cc}
  g_{\ssi{b}b} & \frac{1}{\sqrt{\lad-2}} \\
g_{\ssi{b}\si{b}} &  \sqrt{\frac{\lad-4}{\lad-2}}  \\
   g_{\ssi{b}\ssi{b}} & \frac{1}{\sqrt{\lad-2}}
   \end{array} \right)
$$
are unitaries, hence we have
\begeq
g_{bb}&=& \left( \frac{-1}{\sqrt{3}}, -\frac{\sqrt{\lad-2}}{\la} 
\right),  \\
g_{\si{b}\si{b}}&=&\left( \frac{\lad-3}{\sqrt{3}}, 0 \right), \\
g_{\ssi{b}\ssi{b}}&=&\left( \frac{\lad-3}{\sqrt{3}}, 0 \right).
\endeq
Now, the connection $(\a \ab-{\bf 1})$ is as in the \framebox{Table
\ref{aab-1.tex}}. 


\begin{table}[h]
\begin{center}
\begin{tabular}[pos]{|@{}c@{}||@{}c@{}|@{}c@{}|| @{}c@{}|@{}c@{}|@{}c@{}|
 @{}c@{}|@{}c@{}||@{}c@{}| @{}c@{}||@{}c@{}|@{}c@{}| }
\hl
& $aa^2$ & $ca$ & $ac$ &  $cc^2$ & $cc^3$ & \s{a}$c$ &\ss{a}$c$ &
 $c$\s{a} &  \ss{a}\s{a} & $c$\ss{a} & \s{a}\ss{a}  \\ \hline\hline
$*b$ & 1 &&&&&&&&&&  \\ \hline\hline
$b*$ & $\frac{1}{\la_1}$ & $\frac{\la_2}{\la_1}$ &
&&&&&&&& \\ \hl
$bb^2$ &  $\frac{-\la_2}{\la_1}$ & 
$\frac{1}{\la_1}$ & $\frac{\la_4}{\la_2}$ &  
\mcol{2}{c|}{$g_{bb}$} & &&&&& \\ \hl
$b$\s{b} & && $\frac{\bar{\rho}}{\la_2}$ &  
\mcol{2}{@{}c@{}|}{$g_{b \si{b}}$} &&&1&&& \\ \hl
$b$\ss{b} &&& $\frac{\rho}{\la_2}$ &  
\mcol{2}{@{}c@{}|}{$g_{b \ssi{b}}$} & & &&& 1 & \\ \hline\hline
\s{b}$b$ & & $\rho$ & &  \mcol{2}{@{}c@{}|}{$g_{\si{b}b}$} & 
$\frac{1}{\la_2}$ & &&&& \\ \hl
\s{b}\s{b}$^2$ & &&&  \mcol{2}{@{}c@{}|}{$g_{\si{b}\si{b} }$} & 
$\frac{1}{\la_2}$ & & 1  &&& \\ \hl
\s{b}\ss{b} & &&& \mcol{2}{@{}c@{}|}{$g_{\si{b} \ssi{b}}$} &
$\frac{\la_4}{\la_2}$ & &&& $\frac{1}{\la_1}$ & 
$\frac{\la_2}{\la_1}$  \\ \hl
\s{b}\ss{*} &&&&&&& &&& $\frac{\la_2}{\la_1}$ & 
$\frac{-1}{\la_1}$  \\ \hline\hline
\ss{b}$b$ && $\bar{\rho}$ & & \mcol{2}{@{}c@{}|}{$g_{b_{\sigma^2} b}$} & & 
$\frac{1}{\la_2}$ &&&& \\ \hl
\ss{b}\s{b} &&&& \mcol{2}{@{}c@{}|}{$g_{\ssi{b} \si{b}}$} &&
$\frac{\la_4}{\la_2}$ & $\frac{1}{\la_1}$ & 
$\frac{\la_2}{\la_1}$ && \\ \hl
\ss{b}\ss{b}$^2$ &&&& \mcol{2}{@{}c@{}|}{$g_{\ssi{b} \ssi{b}}$} &&
$\frac{1}{\la_2}$ &&& 1 & \\ \hl
\ss{b}\s{*} &&& &&&&& $\frac{\la_2}{\la_1}$ & 
$\frac{-1}{\la_1}$ & & \\ \hline\hline
\s{*}\ss{b} && &&&& 1 &  &&&& -1  \\ \hline\hline
\ss{*}\s{b} && &&&& &1 && -1 && \\ \hl
\end{tabular}
\end{center}
\caption{Connection $\alpha {\tilde \alpha -{\bf 1}}$ ($\la_n=\sqrt{\lad-n}$)}
\label{aab-1.tex}
\end{table}
Our aim is to show $(\a \ab-{\bf 1}) \cong \sigma(\a \ab-{\bf 1})
\sigma$. For this purpose, an expression of $\a \ab-{\bf 1}$ with
symmetry up to $\sigma$ is useful. We will re-choose another gauge
as in \framebox{Table \ref{aab-1ng.tex}},
\begin{table}[h]
\begin{center}
\begin{tabular}[pos]{|@{}c@{}|@{}c@{}||c||c||c@{} 
|@{}c@{}|@{}c@{}| c@{}|c||c@{}| c||c@{}|c@{}| }
\hl
&& 1& 1& 1& & & $s$ & ${\bar s}$ & ${\bar s}$ & $-{\bar s}$ &
$s$ & $-s$ \\ \hl
&& $aa^2$ & $ca$ & $ac$ &  $cc^2$ & $cc^3$ & \s{a}$c$ &\ss{a}$c$ &
  $c$\s{a} &  \ss{a}\s{a} & $c$\ss{a} & \s{a}\ss{a}  \\ \hline\hline
1 & $*b$ & 1 && 1 &&&&&&&&  \\ \hline\hline
1 & $b*$ & $\frac{1}{\la_1}$ & $\frac{\la_2}{\la_1}$ &
&&&&&&&& \\ \hl
1 & $bb^2$ &  $\frac{-\la_2}{\la_1}$ & 
$\frac{1}{\la_1}$ & $\frac{\la_4}{\la_2}$ &  
\mcol{2}{c|}{$g'_{bb}$} & &&&&& \\ \hl
$s^{2}$ & $b$\s{b} & && $\frac{\bar s}{\la_2}$ &  
\mcol{2}{@{}c@{}|}{$g'_{b \si{b}}$} &&& $s$ &&& \\ \hl
${\bar s}^{2}$ & $b$\ss{b} &&& $\frac{s}{\la_2}$ &  
\mcol{2}{@{}c@{}|}{$g'_{b \ssi{b}}$} & & &&& ${\bar s}$ & \\ 
\hline\hline
${\bar s}^{2}$ & \s{b}$b$ & & $s$ & &  \mcol{2}{@{}c@{}|}{$g'_{\si{b}b}$} & 
$\frac{\bar s}{\la_2}$ & &&&& \\ \hl
1 & \s{b}\s{b}$^2$ & &&&  \mcol{2}{@{}c@{}|}{$g'_{\si{b}\si{b} }$} & 
$\frac{s}{\la_2}$ & & ${\bar s}$  &&& \\ \hl
${\bar s}$ & \s{b}\ss{b} & &&& \mcol{2}{@{}c@{}|}{$g'_{\si{b} \ssi{b}}$}&
$\frac{\la_4}{\la_2}$ & &&& $\frac{1}{\la_1}$ & 
$\frac{-\la_2}{\la_1}$  \\ \hl
${\bar s}$ & \s{b}\ss{*} &&&&&&& &&& $\frac{\la_2}{\la_1}$ & 
$\frac{1}{\la_1}$  \\ \hline\hline
$s^{2}$ & \ss{b}$b$ && $\bar{s}$ & & \mcol{2}{@{}c@{}|}{$g'_{b_{\sigma^2} b}$} & & 
$\frac{s}{\la_2}$ &&&& \\ \hl
$s$ &\ss{b}\s{b} &&&& \mcol{2}{@{}c@{}|}{$g'_{\ssi{b} \si{b}}$} &&
$\frac{\la_4}{\la_2}$ & $\frac{1}{\la_1}$ & 
$\frac{-\la_2}{\la_1}$ && \\ \hl
1 & \ss{b}\ss{b}$^2$ &&&& \mcol{2}{@{}c@{}|}{$g'_{\ssi{b} \ssi{b}}$} &&
$\frac{\bar s}{\la_2}$ &&& $s$ & \\ \hl
$s$ & \ss{b}\s{*} &&& &&&&& $\frac{\la_2}{\la_1}$ & 
$\frac{1}{\la_1}$ & & \\ \hline\hline
1 & \s{*}\ss{b} && &&&& 1 &  &&&& 1  \\ \hline\hline
1 & \ss{*}\s{b} && &&&& &1 && 1 && \\ \hl
\end{tabular}
\end{center}
\caption{Connection $\alpha {\tilde \alpha} -{\bf 1}$ after taking symmetric 
gauge choice}
\label{aab-1ng.tex}
\end{table}
where $s={\bar \tau}$, $\la_n=\sqrt{\lad-n}$,  and numbers beside the name of edges denote
$1 \times 1$ unitaries corresponding to the edges, namely, we have
multiplied these numbers to the corresponding rows (resp. columns) in 
the previous table, and $g'_{**}$'s denote the vectors corresponding
to $g_{**}$'s after being multiplied by suitable gauge numbers
 respectively.
 By seeing this table, we easily see that 
entries {\it other than} $g'_{**}$'s are invariant to the transformation
of $$(\a \ab-{\bf 1}) \longrightarrow \sigma(\a \ab-{\bf 1})\sigma,$$
which acts on the table as the relabeling $xy \rightarrow 
\sigma(x)\sigma(y)$. The remaining problem is whether we have a
gauge unitary matrix $u{{c}\cho{c}}_2$ corresponding to double edges
$c$-$c$ such that
$$ g'_{**} \stackrel{\sm u{{c}\cho{c}}}{\longrightarrow}
 g'_{\sigma(*)\sigma(*)} $$
 or not.
 We can check by a simple computation that
 $$u{{c}\cho{c}}_2=
 \left( \begin{array}{cc}
 \frac{-1}{2}-\frac{(\lad-2)^{\frac{3}{2}}}{6} i & \frac{\lad-2}
 {\la\sqrt{3}} \\
 \frac{\lad-2}{\la\sqrt{3}} i & \frac{-1}{2}+
 \frac{(\lad-2)^{\frac{3}{2}}}{6}i
 \end{array} \right)
 $$
 gives rise to the transformation
 $$g'_{b\si{b}} \rightarrow g'_{\si{b}\ssi{b}} \rightarrow
 g'_{\ssi{b}b} \rightarrow g'_{b\si{b}},$$
 $$g'_{\si{b}b} \rightarrow g'_{\ssi{b}\si{b}} \rightarrow
 g'_{b\ssi{b}} \rightarrow g'_{\si{b}b},$$
 and
 $$ g'_{bb} \rightarrow g'_{\si{b}\si{b}}  \rightarrow
 g'_{\ssi{b}\ssi{b}}.$$
 Thus, we have proved the equivalence of connections
 $$\a \ab-{\bf 1} \cong \sigma(\a \ab-{\bf 1})\sigma.$$
 Finally we will check conditions 1) and 2). \label{last}
 Mutual inequivalence is obvious by seeing four graphs of 
 the connections appearing there. Namely, 
 connections producing the bimodules of different indices are
trivially mutually inequivalent. To prove inequivalence of connections 
which produce the bimodules of the same index, it is 
sufficient to show the existence of the unitary matrices of the connection 
of the form 
\thinlines
\unitlength 0.5mm
\begin{picture}(30,20)(0,-5)
\multiput(11,-2)(0,10){2}{\line(1,0){8}}
\multiput(10,7)(10,0){2}{\line(0,-1){8}}
\multiput(10,-2)(10,0){2}{\circle*{1}}
\multiput(10,8)(10,0){2}{\circle*{1}}
\put(6,12){\makebox(0,0){$x$}}
\put(24,-6){\makebox(0,0){$y$}}
\end{picture}
which have different sizes in each connection. We can check it only by 
seeing the four graphs.  About the indecomposability,  since it was 
irreducibility of the bimodules in our original lemma, all we must 
see is the irreducibility of bimodules made of connections here.
The bimodule $X^{\bf 1}={}_{N} N_{N}$ is trivially irreducible, and  
indecomposability of $\sigma$ and $\sigma^{2}$ follows.
To see the irreducibility of $X^{\a}$, consider the subfactor  $N 
\subset M$ 
constructed from the connection $\a$. Then $X^{\a} = {}_{N} M_{M}$.
By  Ocneanu's compactness argument, (see Section 3, Theorem 4)
\begeq 
\mbox{End}({}_{N} X^{\a}_{M}) &=&
\mbox{End}({}_{N} M_{M}) \\
&=& N' \cap M  \subset  \mbox{String}_{*}^{(1)} {\cal 
G}={\bf C}, 
\endeq
where ${\cal G}$ is the upper graph in Figure \ref{13graph},
thus irreducibility of $X^{\a}$, $X^{\sigma \a}$, and $X^{\sigma^{2} \a}$
follows. Similarly, we have
\begeq
& \mbox{End}({}_{N} X^{\a \ab}_{N}) = \mbox{End}({}_{N} M
 \otimes_{M} M_{N}) &= N' \cap M_{1} \subset \mbox{String}_{*}^{(2)} {\cal 
G} = {\bf C} \oplus {\bf C},
\endeq
where $N \subset M \subset M_{1} \cdots $ is Jones tower of $N 
\subset M$, 
thus irreducibility of $X^{\a \ab -{\bf 1}}$, $X^{\sigma (\a \ab -{\bf 1})}$, 
and $X^{\sigma^{2} (\a \ab -{\bf 1})}$ follows. Irreducibility of 
$X^{\a \ab \a -2\a}$ follows in the same way using
$${\rm String}_{*}^{(3)}{\cal G} = {\bf C} \oplus {\rm M}_{2}({\bf 
C}).$$
Now, the proposition holds and thus we have proved the theorem. \\
\qed

\section{Main theorem for the case of $(5+\sqrt{17})/2$}
In this section, we will give a proof for our main theorem 
for the case of index $(5+\sqrt{17})/2$ due to the first named author.

\begin{thm}
A subfactor with principal graph and dual principal graph as in Figure
\ref{17graph} exists.
\end{thm}
From the key lemma, we know that the above theorem follows from the 
next proposition. We define the connection $\sigma$ as
\begin{displaymath}
\sigma
\left(
\thinlines
\unitlength 0.5mm
\begin{picture}(30,30)
\multiput(11,-2)(0,10){2}{\line(1,0){8}}
\multiput(10,7)(10,0){2}{\line(0,-1){8}}
\multiput(10,-2)(10,0){2}{\circle*{1}}
\multiput(10,8)(10,0){2}{\circle*{1}}
\put(6,12){\makebox(0,0){$p$}}
\put(24,12){\makebox(0,0){$q$}}
\put(6,-6){\makebox(0,0){$r$}}
\put(24,-6){\makebox(0,0){$s$}}
\end{picture}
\right)
= \delta_{p,{\tilde{r}}}\delta_{q,{\tilde{s}}},
\end{displaymath}
here $p,q,r,s$ are the vertices on the upper graph in Figure 
\ref{17graph}
and we consider $\tilde{\tilde{x}}$ as $x$ and if $x$ is one of 
$e,f,g$, $\tilde{x}=x$. \par

\begin{prop}
Let $\alpha$ be the unique connection on the four graphs 
consisting of the pair of the graphs appearing in Figure 
\ref{17graph},
and $\sigma$ be the connection defined above.
Then, the following hold. \\
1) The eight connections
\begeq
{\bf 1}, \; \sigma, \;  \sigma^2, \;  (\a \ab -{\bf 1}), \; 
\sigma(\a \ab -{\bf 1}),
\; (\a \ab -{\bf 1})\sigma, \; \sigma(\a \ab -{\bf 1})\sigma,\; \\ 
(\a \ab)^{2}-3 \a \ab +{\bf 1}, \; 
\sigma((\a \ab)^{2}-3 \a \ab +{\bf 1})
\endeq
are indecomposable and mutually inequivalent. \\
2) The six connections
$$ \a, \; \sigma \a, \; \a \ab \a-2\a, \; \sigma(\a \ab \a-2\a),
\; (\a \ab)^{2}\a-4 \a \ab \a + 3\a, \; (\a \ab-{\bf 1})\sigma \a
$$
are irreducible and mutually inequivalent. \\
3) 
$$ \sigma(\a \ab -{\bf 1})\sigma \a \cong 
(\a \ab -{\bf 1})\sigma \a. $$
\end{prop}
{\it Proof} \\
The four graphs of the connection $\a$ and the Perron-Frobenius weights 
are as in  
\framebox{Figure \ref{17fourgraphs}}.
\begin{figure}[H]
\begin{center}
\thinlines
\unitlength 1.0mm
\begin{picture}(30,20)(0,-5)
\multiput(11,-2)(0,10){2}{\line(1,0){8}}
\multiput(10,7)(10,0){2}{\line(0,-1){8}}
\multiput(10,-2)(10,0){2}{\circle*{1}}
\multiput(10,8)(10,0){2}{\circle*{1}}
\put(6,12){\makebox(0,0){$V_0$}}
\put(24,12){\makebox(0,0){$V_1$}}
\put(6,-6){\makebox(0,0){$V_3$}}
\put(24,-6){\makebox(0,0){$V_2$}}
\put(15,12){\makebox(0,0){${\cal G}_0$}}
\put(15,-6){\makebox(0,0){${\cal G}_2$}}
\put(6,3){\makebox(0,0){${\cal G}_3$}}
\put(24,3){\makebox(0,0){${\cal G}_1$}}
\put(15,3){\makebox(0,0){$\a$}}
\end{picture} 
\end{center}

\begin{center}
	\setlength{\unitlength}{0.5mm}
	\begin{picture}(200,120)(-20,-60)
	\multiput(20,-25)(40,0){5}{\circle*{2}}
	\put(120,-25){\circle*{2}}
	\multiput(0,50)(40,0){6}{\circle*{2}}
	\multiput(100,50)(40,0){3}{\circle*{2}}
	\multiput(20,25)(40,0){5}{\circle*{2}}
	\put(120,25){\circle*{2}}
	\multiput(0,0)(40,0){5}{\circle*{2}}
	\put(100,0){\circle*{2}}
	\multiput(0,-50)(40,0){6}{\circle*{2}}
	\multiput(100,-50)(40,0){3}{\circle*{2}}
	\multiput(0,0)(40,0){5}{\line(4,-5){21}}
	\multiput(40,0)(40,0){4}{\line(-4,-5){21}}
	\multiput(100,0)(20,0){2}{\line(0,-1){25}}
	\multiput(20,-25)(40,0){5}{\line(4,-5){21}}
	\multiput(20,-25)(40,0){3}{\line(-4,-5){21}}
	\put(140,-25){\line(-8,-5){41}}
	\multiput(100,-25)(20,0){3}{\line(0,-1){25}}
	\put(180,-25){\line(0,-1){25}}
	\put(120,-25){\line(12,-5){60}}
	\multiput(0,50)(40,0){3}{\line(4,-5){21}}
	\multiput(40,50)(40,0){3}{\line(-4,-5){21}}
	\multiput(100,50)(20,0){3}{\line(0,-1){25}}
	\put(100,50){\line(8,-5){40}}
	\multiput(160,50)(20,0){2}{\line(-8,-5){40}}
	\put(160,50){\line(4,-5){21}}
	\put(200,50){\line(-4,-5){21}}
	\multiput(20,25)(40,0){4}{\line(4,-5){21}}
	\multiput(20,25)(40,0){5}{\line(-4,-5){21}}
	\multiput(100,25)(20,0){2}{\line(0,-1){25}}
	\put(25,-25){\makebox(1,1){$A$}}
	\put(65,-25){\makebox(1,1){$C$}}
	\put(106,-25){\makebox(1,1){$E$}}
	\put(124,-22){\makebox(1,1){$\tilde{C}$}}
	\put(145,-25){\makebox(1,1){$G$}}
	\put(185,-25){\makebox(1,1){$\tilde{A}$}}
	\put(4,50){\makebox(1,1){*}}
	\put(44,50){\makebox(1,1){$b$}}
	\put(84,50){\makebox(1,1){$d$}}
	\put(96,50){\makebox(1,1){$f$}}
	\put(124,50){\makebox(1,1){$\tilde{d}$}}
	\put(144,50){\makebox(1,1){$h$}}
	\put(164,50){\makebox(1,1){$\tilde{b}$}}
	\put(184,50){\makebox(1,1){$\tilde{h}$}}
	\put(204,50){\makebox(1,1){$\tilde{*}$}}
	\put(24,25){\makebox(1,1){$a$}}
	\put(64,25){\makebox(1,1){$c$}}
	\put(104,25){\makebox(1,1){$e$}}
	\put(126,25){\makebox(1,1){$\tilde{c}$}}
	\put(146,25){\makebox(1,1){$g$}}
	\put(184,25){\makebox(1,1){$\tilde{a}$}}
	\put(5,0){\makebox(1,1){1}}
	\put(45,0){\makebox(1,1){2}}
	\put(85,0){\makebox(1,1){3}}
	\put(105,0){\makebox(1,1){4}}
	\put(125,0){\makebox(1,1){5}}
	\put(165,0){\makebox(1,1){6}}
	\put(4,-50){\makebox(1,1){*}}
	\put(44,-50){\makebox(1,1){$b$}}
	\put(84,-50){\makebox(1,1){$d$}}
	\put(96,-50){\makebox(1,1){$f$}}
	\put(124,-50){\makebox(1,1){$\tilde{d}$}}
	\put(144,-50){\makebox(1,1){$h$}}
	\put(164,-50){\makebox(1,1){$\tilde{b}$}}
	\put(184,-50){\makebox(1,1){$\tilde{h}$}}
	\put(204,-50){\makebox(1,1){$\tilde{*}$}}
	\put(-10,37){\makebox(2,2){\large ${\cal G}_0$}}
	\put(-10,12){\makebox(2,2){\large ${\cal G}_1$}}
	\put(-10,-13){\makebox(2,2){\large ${\cal G}_2$}}
	\put(-10,-38){\makebox(2,2){\large ${\cal G}_3$}}
	\put(-15,50){\makebox(2,2){\large $V_0$}}
	\put(-15,25){\makebox(2,2){\large $V_1$}}
	\put(-15,0){\makebox(2,2){\large $V_2$}}
	\put(-15,-25){\makebox(2,2){\large $V_3$}}
	\put(-15,-50){\makebox(2,2){\large $V_0$}}
	\end{picture}
\end{center}
\mbox{The Perron-Frobenius weights:} \\
\mbox{$\mu(*)=\mu(\td{*})=1,\quad \mu(a)=\mu(\td{a})=\gbs,\quad
\mu(b)=\mu(\td{b})=\mu(h)=\mu(\td{h})=\gb-1,$ } \\
\mbox{$\mu(c)=\mu(\td{c})=\gbs^{3}-2\gbs, \quad 
\mu(d)=\mu(\td{d})=2\gb-1, $ }
\mbox{$\mu(e)=\gbs^{3}+\gbs, \quad \mu(f)=2\gb, $} \\
\mbox{$\mu(g)=\gbs^{3}-\gbs, \quad \mu(2)=\gb-1, \quad \mu(3)=2\gb-1, 
\quad \mu(4)=\gb+1,$ }    \\
\mbox{$\mu(5)=3\gb-2, \quad \mu(6)=\gb.$}
\caption{Four graphs of the connection $\alpha$}
\label{17fourgraphs}
\end{figure}
\noindent
Note that the Perron-Frobenius weights of the vertices in $V_{3}$ are 
the same as that of the vertices in $V_{1}$, and here we used 
$\gbf-5\gb+2=0$.
The biunitary connection $\alpha$ on these four graphs is 
determined uniquely as in \framebox{Table \ref{ab17.tex}}, as in 
$(5+\sqrt{13})/2$ case. We will also 
display the Table of the connection $\ab$ for use of later computations.  
\pagebreak
\begin{table}[H]
\hspace{-1.5cm}
\begin{center}
\begin{tabular}[pos]{|c||@{}c@{}||c|c||c|c||@{}c@{}||c|c|c||@{}c@{}|@{}c@{}|}
\hl
& $a1$ & $a2$ & $c2$ & $c3$ &  $e3$ & $e4$ & $e5$ & $\td{c}5$ & 
$g5$ & $g6$ & $\td{a}6$  \\ \hline\hline
$*A$ & 1 & 1 & &&&&&&&& \\ \hline\hline
$bA$ & 1 & $\frac{-1}{\gbs_1'}$ & $\frac{\gbs \gbs_2}{\gbs_1'}$
 & &&&&&&& \\ \hl
$bC$ & & $\frac{\gbs\gbs_2}{\gbs_1'}$ & $\frac{1}{\gbs_1'}$
 & 1 &  & &&&&& \\ \hline\hline
$dC$ & && 1 & $\frac{-1}{\gamma'}$ & $\frac{2\gbs \gbs_1}{\gamma'}$ 
 &&&&&& \\ \hl
$dE$ & &&& $\frac{2\gbs \gbs_1}{\gamma'}$ & 
$\frac{1}{\gamma'}$ & 1 & 1 &&&& \\ \hline\hline
$fE$ & && & & 1 & $-1$ &  $\frac{1}{\gbs_2'}$ & 
& $\frac{\gbs_{-1}}{\gbs_2'}$ && \\ \hl
$fG$ &&&&&& & $\frac{\gbs_{-1}}{\gbs_2'}$ & & $\frac{-1}{\gbs_2'}$ & 1 
& \\ \hline\hline
$hG$ &&&&&&&&& 1 &1 & \\ \hline\hline
$\td{h}\td{C}$ & &&&&&&&& 1 && \\ \hl
$\td{h}\td{A}$  & &&&&&&&&& 1  & \\ \hline\hline
$\td{d}E$ &&&&& 1 & 1 & $\frac{-\gbs_2}{\gamma}$ & 
$\frac{\gbs_{-1}}{\gamma}$ & && \\ \hl
$\td{d}\td{C}$   &&&&&&& $\frac{\gbs_{-1}}{\gamma}$ & 
$\frac{\gbs_2}{\gamma}$ &&& \\ \hline\hline
$\td{b}G$  &&&&&&&& 1 &&& 1 \\ \hline\hline
$\td{*}\td{A}$ &&&&&&&&&&& 1 \\ \hl
\end{tabular}
\end{center}
\hspace{-1.5cm}
\begin{center}
\begin{tabular}[pos]{|c||c|c||c|c||c|c|c||@{}c@{}||c|c||@{}c@{}|}
\hl 
& $1a$ & $2a$ & $2c$ & $3c$ &  $3e$ & $4e$ & $5e$ & $5\td{c}$ & 
$5g$ & $6g$ & $6\td{a}$  \\ 
\hline\hline
$A*$ & $\frac{1}{\gbs}$ & $ \frac{\gbs_1}{\gbs}$ & &&&&&&&& \\ 
\hl
$Ab$ & $\frac{\gbs_1}{\gbs}$ & $\frac{-1}{\gbs}$ & 1 & &&&&&&& 
\\
\hline\hline
$Cb$ & & 1 & $\frac{1}{\gbs \gbs_2'}$ & 
$\frac{\gbs_1'}{\sqrt{2} \gbs_2' }$ & &&&&&& \\ \hl
$Cd$ & && $\frac{\gbs_1'}{\sqrt{2} \gbs_2'}$ & 
$\frac{-1}{\gbs \gbs_2'}$ & 1 &&&&&& \\ \hline\hline
$Ed$ &&&& 1 & $\frac{1}{\gbs \gbs_{-1}'}$ &
$\frac{\gbs_1}{\sqrt{2} \gbs_{-1}}$ & $\frac{\gb}{\gbs_{-1}'}$ &&&& \\ \hl
$Ef$ & &&&& $\frac{\gbs \gbs_1}{\gbs_{-1}'}$ & 
$\frac{-\sqrt{2}}{\gbs_{-1}}$ & $\frac{\gbs_1}{\gbs_{-1}'}$ & & 
1 && \\ \hl
$E\td{d}$ &&&&& $\frac{\gamma'}{\gbs \gbs_{-1}'}$ & $\frac{\gbs_1}{\sqrt{2} \gbs_{-1}}$ 
& $\frac{-\gbs_2'}{\gbs_{-1}'}$ & 1 &&& \\ \hline\hline
$\td{C}\td{d}$ &&&&&&& 1 & 1 &&& \\ \hl
$\td{C}\td{h}$ &&&&&&&&& 1 && \\ \hline\hline
$Gf$ & &&&&&& 1 & & $\frac{-1}{\gbs_1}$ & 
$\frac{\gbs_2}{\gbs_1}$ & \\ \hl
$Gh$ &&&&&&&&& $\frac{\gbs_2}{\gbs_1}$ & $\frac{1}{\gbs_1}$ & \\ 
\hl
$G\td{b}$ &&&&&&&& 1 &&& 1 \\ \hline\hline
$\td{A}\td{h}$ &&&&&&&&&& 1 & \\ \hl
$\td{A}\td{*}$ &&&&&&&&&&& 1 \\ \hl
\end{tabular}
\end{center}
\mbox{where,} \\
\mbox{$\gbs_n = \sqrt{\gb-n}, \quad \gbs_n'=\gb-n$,} \\
\mbox{$\gamma = \sqrt{2\gb-1}, \quad \gamma' = 2\gb-1$.}
\caption{Connections $\a$ (upper) and $\ab$ (lower)}
\label{ab17.tex}
\end{table}
\newpage
First we check condition 3), namely we prove
$$ \sigma(\a \ab -{\bf 1})\sigma \a \cong 
(\a \ab -{\bf 1})\sigma \a. $$
up to vertical gauge choice. 
Now we compute the connection $\a \ab$. The four graphs on 
which the connection $\a \ab$ exists are as in 
\framebox{Figure \ref{aab17graph}}.
\begin{figure}[H]
\begin{center}
\thinlines
\unitlength 1.0mm
\begin{picture}(30,20)(0,-5)
\multiput(11,-2)(0,10){2}{\line(1,0){8}}
\multiput(10,7)(10,0){2}{\line(0,-1){8}}
\multiput(10,-2)(10,0){2}{\circle*{1}}
\multiput(10,8)(10,0){2}{\circle*{1}}
\put(6,12){\makebox(0,0){$V_2$}}
\put(24,12){\makebox(0,0){$V_3$}}
\put(6,-6){\makebox(0,0){$V_2$}}
\put(24,-6){\makebox(0,0){$V_3$}}
\put(15,12){\makebox(0,0){\large ${\cal G}_2$}}
\put(15,-6){\makebox(0,0){\large ${\cal G}_2$}}
\put(4,3){\makebox(0,0){\large ${\cal G}_3 {\cal G}_3^{t}$}}
\put(26,3){\makebox(0,0){\large ${\cal G}_1 {\cal G}_1^{t}$}}
\put(15,3){\makebox(0,0){$\a \ab$}}
\end{picture} 
\end{center}
\begin{center}
\setlength{\unitlength}{0.5mm}
	\begin{picture}(200,140)(-20,-10)
\multiput(0,100)(40,0){6}{\circle*{2}}
	\multiput(100,100)(40,0){3}{\circle*{2}}
	\multiput(20,75)(40,0){5}{\circle*{2}}
	\put(120,75){\circle*{2}}
	\multiput(20,50)(40,0){5}{\circle*{2}}
	\put(120,50){\circle*{2}}
	\multiput(0,25)(40,0){6}{\circle*{2}}
	\multiput(100,25)(40,0){3}{\circle*{2}}
	\multiput(0,0)(40,0){6}{\circle*{2}}
	\multiput(100,0)(40,0){3}{\circle*{2}}
		\put(4,100){\makebox(1,1){*}}
	\put(44,100){\makebox(1,1){$b$}}
	\put(84,100){\makebox(1,1){$d$}}
	\put(96,100){\makebox(1,1){$f$}}
	\put(124,100){\makebox(1,1){$\tilde{d}$}}
	\put(144,100){\makebox(1,1){$h$}}
	\put(164,100){\makebox(1,1){$\tilde{b}$}}
	\put(184,100){\makebox(1,1){$\tilde{h}$}}
	\put(204,100){\makebox(1,1){$\tilde{*}$}}
			\put(4,0){\makebox(1,1){*}}
	\put(44,0){\makebox(1,1){$b$}}
	\put(84,0){\makebox(1,1){$d$}}
	\put(96,0){\makebox(1,1){$f$}}
	\put(124,0){\makebox(1,1){$\tilde{d}$}}
	\put(144,0){\makebox(1,1){$h$}}
	\put(164,0){\makebox(1,1){$\tilde{b}$}}
	\put(184,0){\makebox(1,1){$\tilde{h}$}}
	\put(204,0){\makebox(1,1){$\tilde{*}$}}
	\put(24,75){\makebox(1,1){$a$}}
	\put(64,75){\makebox(1,1){$c$}}
	\put(104,75){\makebox(1,1){$e$}}
	\put(126,75){\makebox(1,1){$\tilde{c}$}}
	\put(146,75){\makebox(1,1){$g$}}
	\put(184,75){\makebox(1,1){$\tilde{a}$}}
	\put(24,50){\makebox(1,1){$a$}}
	\put(64,50){\makebox(1,1){$c$}}
	\put(104,50){\makebox(1,1){$e$}}
	\put(126,50){\makebox(1,1){$\tilde{c}$}}
	\put(146,50){\makebox(1,1){$g$}}
	\put(184,50){\makebox(1,1){$\tilde{a}$}}
	\put(4,25){\makebox(1,1){*}}
	\put(44,25){\makebox(1,1){$b$}}
	\put(84,25){\makebox(1,1){$d$}}
	\put(96,25){\makebox(1,1){$f$}}
	\put(124,25){\makebox(1,1){$\tilde{d}$}}
	\put(144,25){\makebox(1,1){$h$}}
	\put(164,25){\makebox(1,1){$\tilde{b}$}}
	\put(184,25){\makebox(1,1){$\tilde{h}$}}
	\put(204,25){\makebox(1,1){$\tilde{*}$}}
		\multiput(0,100)(40,0){3}{\line(4,-5){20}}
	\multiput(40,100)(40,0){3}{\line(-4,-5){21}}
	\multiput(100,100)(20,0){3}{\line(0,-1){25}}
	\put(100,100){\line(8,-5){40}}
	\multiput(160,100)(20,0){2}{\line(-8,-5){40}}
	\put(160,100){\line(4,-5){21}}
	\put(200,100){\line(-4,-5){21}}
	\dline(20,75) \dline(60,75) \dline(100,75) \dline(120,75)
	\dline(140,75) \dline(180,75)
	\put(18,75){\line(0,-1){25}}
		\multiput(20,75)(40,0){4}{\line(8,-5){40}}
		\multiput(60,75)(40,0){4}{\line(-8,-5){40}}
		\multiput(58,75)(40,0){3}{\line(0,-1){25}}
		\put(102,75){\line(0,-1){25}}
		\multiput(100,75)(20,0){2}{\line(4,-5){20}}
		\multiput(120,75)(20,0){2}{\line(-4,-5){20}}
		\multiput(20,50)(40,0){3}{\line(4,-5){21}}
		\put(180,50){\line(4,-5){21}}
	\multiput(20,50)(40,0){3}{\line(-4,-5){20}}
	\put(140,50){\line(-8,-5){40}}
	\multiput(100,50)(20,0){3}{\line(0,-1){25}}
	\put(180,50){\line(-4,-5){20}}
	\put(120,50){\line(8,-5){40}}
	\put(140,50){\line(8,-5){40}}
	\dline(0,25) \dline(40,25)  \dline(80,25) \dline(100,25)
	\dline(120,25) \dline(140,25) \dline(160,25) 
	\dline(180,25) \dline(200,25)
	\put(182,25){\line(0,-1){25}}
	\multiput(0,25)(40,0){3}{\line(8,-5){40}}
	\multiput(40,25)(40,0){3}{\line(-8,-5){40}}
	\multiput(42,25)(60,0){2}{\line(0,-1){25}}
	\multiput(82,25)(40,0){2}{\line(0,-1){25}}
	\multiput(80,25)(20,0){2}{\line(4,-5){20}}
	\multiput(100,25)(20,0){2}{\line(-4,-5){20}}
	\put(100,25){\line(12,-5){60}}
    \put(100,25){\line(8,-5){40}}
  \put(120,25){\line(12,-5){60}}
	\put(140,25){\line(-8,-5){40}}
	\put(140,25){\line(4,-5){20}}
	\put(160,25){\line(-12,-5){60}}
	\put(160,25){\line(-4,-5){20}}
	\put(180,25){\line(-12,-5){60}}
	\put(180,25){\line(4,-5){20}}	
	\put(200,25){\line(-4,-5){20}}
	\put(-10,87){\makebox(2,2){\large ${\cal G}_0$}}
	\put(-12,62){\makebox(2,2){\large ${\cal G}_1 {\cal G}_1^{t}$}}
	\put(-10,37){\makebox(2,2){\large ${\cal G}_0$}}
	\put(-12,12){\makebox(2,2){\large ${\cal G}_3 {\cal G}_3^{t}$}}
	\put(-15,100){\makebox(2,2){\large $V_0$}}
	\put(-15,75){\makebox(2,2){\large $V_1$}}
	\put(-15,50){\makebox(2,2){\large $V_1$}}
	\put(-15,25){\makebox(2,2){\large $V_0$}}
	\put(-15,0){\makebox(2,2){\large $V_0$}}
	\end{picture}
\caption{Four graphs of the connection $\a \ab$}
\label{aab17graph}
\end{center}
\end{figure}
The broken edges correspond to the trivial connection ${\bf 1}$.
We will now compute the connection $\a \ab-{\bf 1}$,which is determined only 
up to vertical gauges. As in the $(5 + \sqrt{13})/2$-case, we assume 
that $1 \times 1$ gauge transfrom unitaries corresponding to single 
vertical edges which connect different vertices in the graph ${\cal 
G}_{1}{\cal G}^{t}_{1} \cup {\cal G}_{3}{\cal G}^{t}_{3}$ 
to be $1$. Then we easily find 38 entries of $\a \ab-{\bf 1}$ by 
``actual'' multiplication of the connections $\a$ and $\ab$. Next, as 
in the $(5+\sqrt{13})/2$-case, we can compute all the entries of $\a 
\ab-{\bf 1}$ which involve the double edges in the graph of $\a \ab$ by a 
simple gauge transform, then we have $14$ entries listed in the 
 \framebox{Tables \ref{aab-117-1.ps}--\ref{aab-117-2.ps}} {\it other than}
 the entries marked ``(?)'', where $\gbs_n=\sqrt{\gb-n}$
and $\gamma=\sqrt{2\gb-1}$.

\newpage
\begin{table}[H]
\begin{center}
\begin{tabular}[pos]{|@{}c@{}||@{}c@{}|@{}c@{}|| @{}c@{}|
@{}c@{}|@{}c@{}|| @{}c@{}|c|c| @{}c@{}|c|| }
\hl
& $aa$ & $ca$ & $ac$ & $cc$ & $ec$ & $ce$ & $ee^1$ & $ee^2$ &
 $\td{c}e$ & $ge$  \\ \hline\hline
$*b$ & 1 && 1 &&& &&& &  \\ \hline\hline
$b*$ & $\frac{1}{\gbs_1}$ & $\frac{\gbs_2}{\gbs_1}$ &
&&&& &&& \\ \hl
$bb$ &  $-\frac{\gbs_2}{\gbs_1}$  (?) & 
$\frac{1}{\gbs_1}$ & $\frac{1}{\sqrt{2}\gbs}$ & 
$\frac{-\gamma}{\sqrt{2}\gbs}$  (?) &  &&&&&  \\ \hl
$bd$  &   & 
 & $\frac{1}{\sqrt{2}\gbs}$ & 
$\frac{-\gamma}{\sqrt{2}\gbs}$ &  & 1 &&&&  \\
 \hline\hline 
$db$ && 1 && $\frac{1}{\gb}$ & $\frac{\sqrt{\gbf-1}}{\gb}$ &&&&& \\ \hl
$dd$ &&&& $\frac{-\sqrt{\gbf-1}}{\gb}$  (?) & $\frac{1}{\gb}$ & 
$\frac{1}{\gb}\frac{\gbs_2}{\gbs_{-1}}$ &
$l_1$ & $l_2$  & &  \\ \hl
$df$ &&&& & & 
$\frac{\gb-2}{\gb-1}$ & $m_1$ & $m_2$
 & & \\ \hl
$d \td{d}$ &&&& & & 
$\frac{2 \gbs_1}{\gb+1}$ &
$n_1$ & $n_2$
& & \\
 \hline\hline
$fd$ &&&&& 1&& $p_{1}$ &$p_{2}$ && $\frac{\sqrt{3\gb-1}}
{2\gbs_{-1}}$ \\ \hl
$ff$ &&&&&& & $q_{1}$ & $q_{2}$ && $\frac{-1}{2}$  \\ \hl
$fh$ &&&&&& & & &&  \\ \hl
$f \td{d}$ &&&&&& & $r_{1}$ & $r_{2}$ && 
$\frac{-1}{\gbs_{-1}}$  \\ \hl
$f \td{b}$ &&&& &&& \mcol{2}{c|}{} && 
  \\ \hl
 \hl
 $hf$ &&&&&&&&&& 1 \\ \hl
 $h \td{b}$  &&&&&&&&&&  \\ \hl
 \hl
 $\td{h}\td{h}$ &&&&&&&&&&  \\ \hl
 $\td{h}\td{d}$ &&&&&&&&& &1  \\ \hl
 $\td{h}\td{*}$&&&&&&&&&&   \\ \hl
 \hl
 $\td{d}d$ &&&&& 1 && $s_{1}$ & $s_{2}$ & 
 $\frac{\gbs_2}{\gbs_{-1}}$ & \\ \hl
 $\td{d}f$ &&&&&&& $t_{1}$ & $t_{2}$ & 
 $\frac{1}{\gb-1}$ & \\ \hl
  $\td{d}\td{h}$ &&&&&&& \mcol{2}{c|}{} & 
 & \\ \hl
$\td{d}\td{d}$ &&&&&&& $u_{1}$ & $u_{2}$ & 
$\frac{\gbs_2}{\gbs_{-1}}$ &  \\ \hl
\hl
$\td{b}f$&&&&&&&&&1&  \\ \hl
$\td{b}h$&&&&&&&&&&  \\ \hline\hline
$\td{*}\td{h}$ &&&&&&&&&&  \\ \hl
\end{tabular}
\end{center}
 \caption{Connection $\alpha {\tilde \alpha} -{\bf 1}$ (left part of 
 diagram)}
\label{aab-117-1.ps}
\end{table}

\newpage

\begin{table}[H]
\begin{center}
\begin{tabular}[pos]{|@{}c@{}||
c|c||@{}c@{}|@{}c@{}|@{}c@{}|@{}c@{}||@{}c@{}| }
\hl
&  $e\td{c}$ & $g\td{c}$ & $eg$ & $\td{c}g$ &
 $gg$ & $\td{a}g$ & $g\td{a}$  \\ \hline\hline
$*b$  &&&& &&&  \\ \hline\hline
$b*$  &&& &&&& \\ \hl
$bb$  &&& &&&& \\ \hl
$bd$   &&& &&&& \\ \hline\hline 
$db$ 
&&&&&&& \\ \hl
$dd$  &&&& &&& \\ \hl
$df$ &&& 1 &&&& 
\\ \hl
$d \td{d}$  & 1 &&& &&& \\
 \hline\hline
$fd$   &&&& &&& \\ \hl
$ff$   &&& 
$\frac{-\gbs_1}{\sqrt{2}\gbs}$ && 
$\frac{\gbs_{-1}}{\sqrt{2}\gbs}$ (?)  && \\ \hl
$fh$   &&& 
$\frac{\gbs_{-1}}{\sqrt{2}\gbs}$ && 
$\frac{\gbs_1}{\sqrt{2}\gbs}$ && \\ \hl
$f \td{d}$   & $\frac{1}{\gb-2}$ & 
$\frac{\gbs_{-1}}{\gb-2}$ & &&  &&  \\ \hl
$f \td{b}$ 
 & $\frac{\gbs_{-1}}{\gb-2}$ & $\frac{-1}{\gb-2}$ &  &&&& 1 \\ \hl
 \hl
 $hf$  &&& && 1 && \\ \hl
 $h \td{b}$    && 1 &&& && 1  \\ \hl
 \hl
 $\td{h}\td{h}$  && & & & 1  (?) && \\ \hl
 $\td{h}\td{d}$    && 1 &&&&& \\ \hl
 $\td{h}\td{*}$  &&&& &&& 1  \\ \hl
 \hl
 $\td{d}d$  &&&& &&& \\ \hl
 $\td{d}f$   & && $-\sqrt{\frac{\gb+1}{2\gb-1}}$ & 
 $\sqrt{\frac{\gb+1}{2\gb-1}}$ &&& \\ \hl
  $\td{d}\td{h}$  &&&  $\sqrt{\frac{\gb+1}{2\gb-1}}$
& $-\sqrt{\frac{\gb+1}{2\gb-1}}$ &&& \\ \hl
$\td{d}\td{d}$  & 1 && &&&& \\ \hl \hl
$\td{b}f$ && && $\frac{-1}{\gbs_1}$ && 
$\frac{\gbs_2}{\gbs_1}$ & \\ \hl
$\td{b}h$ && && $\frac{\gbs_2}{\gbs_1}$ & &
$\frac{1}{\gbs_1}$ & \\ \hline\hline
$\td{*}\td{h}$   &&&& && 1 & \\ \hl
\end{tabular}
\end{center}
\caption{Connection $\alpha {\tilde \alpha} -{\bf 1}$ (right part of 
diagram)}
\label{aab-117-2.ps}
\end{table}



\newpage

The four entries marked (?) in Table \ref{aab-117-1.ps} can easily be 
computed by the unitarity of the $2 \times 2$ matrices which they 
are part of, and the entry marked (?) in Table \ref{aab-117-2.ps}
can be put equal to 1, because a gauge choice  corresponding ot the 
$\td{h}\td{h}$-edge in the vertical left graph will only be concerned 
with $(\td{h}\td{h}, gg)$-entry. \\
The only entries left to compute are the 18 entries $l_{1}, l_{2}, 
m_{1}, m_{2}, \ldots , u_{1}, u_{2}$ in Table \ref{aab-117-1.ps} and 
Table \ref{aab-117-2.ps}. They can also be obtained as in 
$(5+\sqrt{13})/2$ case, but here we will make a shortcut: Since the 
entries of the connection $\a \ab$ obtained by ``actual'' multiplication
are all real scalars, all the gauge choices involved in decomposing 
the connection $\a \ab$ into $(\a \ab -{\bf 1}) + {\bf 1}$ can be 
chosen to be matrices with real entries. Hence, $l_{1}, l_{2}, 
m_{1}, m_{2}, \ldots , u_{1}, u_{2}$ become real numbers. We still 
have a possibility of making a gauge choice of the double edges 
$e$-$e$ with an orthogonal matrix, i.e., we can make the following 
change:
$$ (l_{1}, l_{2}) \rightarrow (l_{1}, l_{2})v, \quad 
(m_{1}, m_{2}) \rightarrow (m_{1}, m_{2})v, \ldots , (u_{1}, u_{2})
 \rightarrow (u_{1}, u_{2})v $$
for some $v \in O(2)$. (A common orthogonal matrix for all the vectors 
in ${\bf R}^{2}$.)
Then, we can assume $l_{2}=0$ and $m_{2} \geq 0$, thus, we obtain

\begin{displaymath}
\thinlines
\unitlength 0.5mm
\begin{picture}(30,20)
\multiput(11,-2)(0,10){2}{\line(1,0){8}}
\multiput(10,7)(10,0){2}{\line(0,-1){8}}
\multiput(10,-2)(10,0){2}{\circle*{1}}
\multiput(10,8)(10,0){2}{\circle*{1}}
\put(6,12){\makebox(0,0){$d$}}
\put(24,-6){\makebox(0,0){$e$}}
\end{picture}
= \left(
\matrix{
\frac{1}{\gb}\sqrt{\frac{\gb-2}{\gb+1}} & l_{1} & l_{2} \cr
\frac{\gb-2}{\gb-1} & m_{1} & m_{2} \cr
\frac{2\sqrt{\gb-1}}{\gb+1} & n_{1} & n_{2} \cr
} \right)
= \left(
\matrix{
   \frac{1}{\gb}\sqrt{\frac{\gb-2}{\gb+1}} &
\sqrt{\frac{\gbf+4}{\gb(\gb+1)}} & 0 \cr
   \frac{\gb-2}{\gb-1} & \frac{-\sqrt{\gbf-4}}{\beta 
(\gb-1)\sqrt{\gbf+4}} & \sqrt{\frac{2\gbf}{(\gb-1)(\gbf+4)}} \cr
  \frac{2\sqrt{\gb-1}}{\gb+1} & 
-\sqrt{\frac{\gb(\gb-2)}{(\gb+1)(2\gb-1)(\gbf+4)}} & 
\frac{-4\sqrt{\gb(2\gb-1)}}{(\gb-1)\sqrt{(\gb-1)(\gbf+4)}} \cr
} \right)
\end{displaymath}
 by the orthogonality of the matrix. \\
 Now, all the gauge choices have been used up.  We know that there is 
 an orthogonal matrix $V \in O(3)$ such that
 $$
\thinlines
\unitlength 0.5mm
\begin{picture}(30,20)(0,0)
\multiput(11,-2)(0,10){2}{\line(1,0){8}}
\multiput(10,7)(10,0){2}{\line(0,-1){8}}
\multiput(10,-2)(10,0){2}{\circle*{1}}
\multiput(10,8)(10,0){2}{\circle*{1}}
\put(15,3){\makebox(0,0){{\sm $\a \ab$}}}
\put(6,12){\makebox(0,0){$d$}}
\put(24,12){\makebox(0,0){$e$}}
\put(6,-6){\makebox(0,0){$f$}}
\put(24,-6){\makebox(0,0){$e$}}
\end{picture}
=(0, m_{1}, m_{2})V, \quad
\thinlines
\unitlength 0.5mm
\begin{picture}(30,20)(0,0)
\multiput(11,-2)(0,10){2}{\line(1,0){8}}
\multiput(10,7)(10,0){2}{\line(0,-1){8}}
\multiput(10,-2)(10,0){2}{\circle*{1}}
\multiput(10,8)(10,0){2}{\circle*{1}}
\put(15,3){\makebox(0,0){{\sm $\a \ab$}}}
\put(6,12){\makebox(0,0){$f$}}
\put(24,12){\makebox(0,0){$e$}}
\put(6,-6){\makebox(0,0){$\td{d}$}}
\put(24,-6){\makebox(0,0){$e$}}
\end{picture}
=(0, r_{1}, r_{2})V, $$
$$
\thinlines
\unitlength 0.5mm
\begin{picture}(30,20)(0,0)
\multiput(11,-2)(0,10){2}{\line(1,0){8}}
\multiput(10,7)(10,0){2}{\line(0,-1){8}}
\multiput(10,-2)(10,0){2}{\circle*{1}}
\multiput(10,8)(10,0){2}{\circle*{1}}
\put(15,3){\makebox(0,0){{\sm $\a \ab$}}}
\put(6,12){\makebox(0,0){$d$}}
\put(24,12){\makebox(0,0){$e$}}
\put(6,-6){\makebox(0,0){$\td{d}$}}
\put(24,-6){\makebox(0,0){$e$}}
\end{picture}
=(0, n_{1}, n_{2})V, \quad
\thinlines
\unitlength 0.5mm
\begin{picture}(30,20)(0,0)
\multiput(11,-2)(0,10){2}{\line(1,0){8}}
\multiput(10,7)(10,0){2}{\line(0,-1){8}}
\multiput(10,-2)(10,0){2}{\circle*{1}}
\multiput(10,8)(10,0){2}{\circle*{1}}
\put(15,3){\makebox(0,0){{\sm $\a \ab$}}}
\put(6,12){\makebox(0,0){$\td{d}$}}
\put(24,12){\makebox(0,0){$e$}}
\put(6,-6){\makebox(0,0){$d$}}
\put(24,-6){\makebox(0,0){$e$}}
\end{picture}
=(0, s_{1}, s_{2})V, $$
 $$
 \thinlines
\unitlength 0.5mm
\begin{picture}(30,20)(0,0)
\multiput(11,-2)(0,10){2}{\line(1,0){8}}
\multiput(10,7)(10,0){2}{\line(0,-1){8}}
\multiput(10,-2)(10,0){2}{\circle*{1}}
\multiput(10,8)(10,0){2}{\circle*{1}}
\put(15,3){\makebox(0,0){{\sm $\a \ab$}}}
\put(6,12){\makebox(0,0){$f$}}
\put(24,12){\makebox(0,0){$e$}}
\put(6,-6){\makebox(0,0){$d$}}
\put(24,-6){\makebox(0,0){$e$}}
\end{picture}
=(0, p_{1}, p_{2})V, \quad
\thinlines
\unitlength 0.5mm
\begin{picture}(30,20)(0,0)
\multiput(11,-2)(0,10){2}{\line(1,0){8}}
\multiput(10,7)(10,0){2}{\line(0,-1){8}}
\multiput(10,-2)(10,0){2}{\circle*{1}}
\multiput(10,8)(10,0){2}{\circle*{1}}
\put(15,3){\makebox(0,0){{\sm $\a \ab$}}}
\put(6,12){\makebox(0,0){$\td{d}$}}
\put(24,12){\makebox(0,0){$e$}}
\put(6,-6){\makebox(0,0){$f$}}
\put(24,-6){\makebox(0,0){$e$}}
\end{picture}
=(0, s_{1}, s_{2})V $$
 where 
 $$\thinlines
\unitlength 0.5mm
\begin{picture}(30,20)(0,0)
\multiput(11,-2)(0,10){2}{\line(1,0){8}}
\multiput(10,7)(10,0){2}{\line(0,-1){8}}
\multiput(10,-2)(10,0){2}{\circle*{1}}
\multiput(10,8)(10,0){2}{\circle*{1}}
\put(15,3){\makebox(0,0){{\sm $\a \ab$}}}
\end{picture} $$
 denotes the $1 \times 3$ matrices obtained by ``actual'' multiplication 
 of $\a$ and $\ab$. It is clear from the definition of the 
 renormalization of connection, that $(\a \ab -{\bf 1})^{\sim}=\a 
 \ab -{\bf 1}$ without any gauge transformation. Together with the 
 fact that all the entries of the connection $\a \ab$ by ``actual'' 
 multiplication are real numbers, we have 
 
 $$
 \thinlines
\unitlength 0.5mm
\begin{picture}(30,20)(0,0)
\multiput(11,-2)(0,10){2}{\line(1,0){8}}
\multiput(10,7)(10,0){2}{\line(0,-1){8}}
\multiput(10,-2)(10,0){2}{\circle*{1}}
\multiput(10,8)(10,0){2}{\circle*{1}}
\put(15,3){\makebox(0,0){{\sm $\a \ab$}}}
\put(6,12){\makebox(0,0){$f$}}
\put(24,12){\makebox(0,0){$e$}}
\put(6,-6){\makebox(0,0){$d$}}
\put(24,-6){\makebox(0,0){$e$}}
\end{picture}
=\sqrt{\frac{\mu(d)}{\mu(f)}}
 \thinlines
\unitlength 0.5mm
\begin{picture}(30,20)(0,0)
\multiput(11,-2)(0,10){2}{\line(1,0){8}}
\multiput(10,7)(10,0){2}{\line(0,-1){8}}
\multiput(10,-2)(10,0){2}{\circle*{1}}
\multiput(10,8)(10,0){2}{\circle*{1}}
\put(15,3){\makebox(0,0){{\sm $\a \ab$}}}
\put(6,12){\makebox(0,0){$d$}}
\put(24,12){\makebox(0,0){$e$}}
\put(6,-6){\makebox(0,0){$f$}}
\put(24,-6){\makebox(0,0){$e$}}
\end{picture},
$$
 $$
 \thinlines
\unitlength 0.5mm
\begin{picture}(30,20)(0,0)
\multiput(11,-2)(0,10){2}{\line(1,0){8}}
\multiput(10,7)(10,0){2}{\line(0,-1){8}}
\multiput(10,-2)(10,0){2}{\circle*{1}}
\multiput(10,8)(10,0){2}{\circle*{1}}
\put(15,3){\makebox(0,0){{\sm $\a \ab$}}}
\put(6,12){\makebox(0,0){$\td{d}$}}
\put(24,12){\makebox(0,0){$e$}}
\put(6,-6){\makebox(0,0){$d$}}
\put(24,-6){\makebox(0,0){$e$}}
\end{picture}
=\sqrt{\frac{\mu(d)}{\mu(\td{d})}}
 \thinlines
\unitlength 0.5mm
\begin{picture}(30,20)(0,0)
\multiput(11,-2)(0,10){2}{\line(1,0){8}}
\multiput(10,7)(10,0){2}{\line(0,-1){8}}
\multiput(10,-2)(10,0){2}{\circle*{1}}
\multiput(10,8)(10,0){2}{\circle*{1}}
\put(15,3){\makebox(0,0){{\sm $\a \ab$}}}
\put(6,12){\makebox(0,0){$d$}}
\put(24,12){\makebox(0,0){$e$}}
\put(6,-6){\makebox(0,0){$\td{d}$}}
\put(24,-6){\makebox(0,0){$e$}}
\end{picture},
$$
 $$
 \thinlines
\unitlength 0.5mm
\begin{picture}(30,20)(0,0)
\multiput(11,-2)(0,10){2}{\line(1,0){8}}
\multiput(10,7)(10,0){2}{\line(0,-1){8}}
\multiput(10,-2)(10,0){2}{\circle*{1}}
\multiput(10,8)(10,0){2}{\circle*{1}}
\put(15,3){\makebox(0,0){{\sm $\a \ab$}}}
\put(6,12){\makebox(0,0){$\td{d}$}}
\put(24,12){\makebox(0,0){$e$}}
\put(6,-6){\makebox(0,0){$f$}}
\put(24,-6){\makebox(0,0){$e$}}
\end{picture}
=\sqrt{\frac{\mu(f)}{\mu(\td{d})}}
 \thinlines
\unitlength 0.5mm
\begin{picture}(30,20)(0,0)
\multiput(11,-2)(0,10){2}{\line(1,0){8}}
\multiput(10,7)(10,0){2}{\line(0,-1){8}}
\multiput(10,-2)(10,0){2}{\circle*{1}}
\multiput(10,8)(10,0){2}{\circle*{1}}
\put(15,3){\makebox(0,0){{\sm $\a \ab$}}}
\put(6,12){\makebox(0,0){$f$}}
\put(24,12){\makebox(0,0){$e$}}
\put(6,-6){\makebox(0,0){$\td{d}$}}
\put(24,-6){\makebox(0,0){$e$}}
\end{picture},
$$
hence,
\begeq
(p_{1}, p_{2}) &=& \sqrt{\frac{2\gb-1}{2\gb}} (m_{1}, m_{2}) \\
&=& \sqrt{\frac{2\gb-1}{2\gb}} \left(\frac{-\sqrt{\gbf-4}}{\beta 
(\gb-1)\sqrt{\gbf+4}}, \sqrt{\frac{2\gbf}{(\gb-1)(\gbf+4)}} \right), 
\\
(s_{1}, s_{2})&=&(n_{1}, n_{2}) \\
&=& \left( -\sqrt{\frac{\gb(\gb-2)}{(\gb+1)(2\gb-1)(\gbf+4)}},
\frac{-4\sqrt{\gb(2\gb-1)}}{(\gb-1)\sqrt{(\gb-1)(\gbf+4)}} \right),
\endeq
and
 $$ (t_{1}, t_{2}) = \sqrt{\frac{2\gb}{2\gb-1}} (r_{1}, r_{2}). \qquad 
 (\natural 1)$$
 We next determine $ (r_{1}, r_{2})$ and $(t_{1}, t_{2})$. 
 In the text, we denote the connection matrix,
e.g.,
\begin{center}
\thinlines
\unitlength 0.5mm
\begin{picture}(30,20)(0,-5)
\multiput(11,-2)(0,10){2}{\line(1,0){8}}
\multiput(10,7)(10,0){2}{\line(0,-1){8}}
\multiput(10,-2)(10,0){2}{\circle*{1}}
\multiput(10,8)(10,0){2}{\circle*{1}}
\put(6,12){\makebox(0,0){$b$}}
\put(24,-6){\makebox(0,0){$5$}}
\end{picture}
\end{center}
 by $M(b/5)$ for the convenience of space. By 
 orthogonality of the first and the last row in the $3 \times 3$ 
 matrix $M(f/e)$
in the Tables \ref{aab-117-1.ps}--\ref{aab-117-2.ps}, we have
$$ p_{1}r_{1} + p_{2} r_{2} = \frac{\sqrt{3\gb-1}}{2(\gb+1)}, \qquad 
(\natural 2) $$
and by orthogonality of the first two rows in the $3 \times 3$ matrix 
$M(\td{d}/e)$ in Table \ref{aab-117-2.ps}, we have 
$$ s_{1}t_{1} + s_{2}t_{2} = 
\frac{-1}{\gb-1}\sqrt{\frac{\gb-2}{\gb+1}}, $$
too, and together with $(\natural 1)$ we have
$$(s_{1}r_{1} + s_{2} r_{2} 
)=\frac{-1}{\gb-1}\sqrt{\frac{(\gb-2)(2\gb-1)}{(\gb+1)2\gb}}. \qquad 
(\natural 3)$$
Solving $(\natural 2)$ and $(\natural 3)$ with respect to $(r_{1}, 
r_{2})$ using the known values of $p_{1}$, $p_{2}$, $s_{1}$, and 
$s_{2}$ gives
$$(r_{1}, r_{2}) = \left( \frac{-\beta^{3}}{\sqrt{(\gb+1)(\gbf+4)}}, 
\sqrt{\frac{4\gb}{(\gb+1)(\gbf+4)}} \right), $$
and therefore,
$$(t_{1}, t_{2})= \sqrt{\frac{2\gb}{2\gb-1}} \left( 
\frac{-\beta^{3}}{\sqrt{(\gb+1)(\gbf+4)}},   
\sqrt{\frac{4\gb}{(\gb+1)(\gbf+4)}} \right). $$
 The four remaining entries $q_{1}$, $q_{2}$, $u_{1}$, $u_{2}$ can 
 now be computed using the orhtogonarity of the $3 \times 3$ 
 matrices $M(f/e)$ and $M(\td{d}/e)$.
We have 
$$(q_{1}, q_{2})= \left( \sqrt{\frac{2\gbf}{(\gb+1)(\gbf+4)}},
\sqrt{\frac{2(\gb+1)}{\gbf+4}} \right), $$
$$ (u_{1}, u_{2}) = \left( \sqrt{\frac{\gb(\gb+2)}{(\gb+1)(\gbf+4)}},
\frac{2\sqrt{2}}{\sqrt{\gbf+4}} \right). $$
Now, we have obtained all the entries of $\a \ab-{\bf 1}$. We can 
obtain $(\a \ab-{\bf 1})\sigma$ only by exchanging the vertices at the
bottom of the connection $\a \ab-{\bf 1}$ as below.
$$
 \thinlines
\unitlength 0.5mm
\begin{picture}(30,30)
\multiput(11,-2)(0,10){2}{\line(1,0){8}}
\multiput(10,7)(10,0){2}{\line(0,-1){8}}
\multiput(10,-2)(10,0){2}{\circle*{1}}
\multiput(10,8)(10,0){2}{\circle*{1}}
\put(6,12){\makebox(0,0){$p$}}
\put(24,12){\makebox(0,0){$q$}}
\put(6,-6){\makebox(0,0){$r$}}
\put(24,-6){\makebox(0,0){$s$}}
\end{picture} (\mbox{in $(\a \ab-{\bf 1})\sigma$})
:=
 \thinlines
\unitlength 0.5mm
\begin{picture}(30,30)
\multiput(11,-2)(0,10){2}{\line(1,0){8}}
\multiput(10,7)(10,0){2}{\line(0,-1){8}}
\multiput(10,-2)(10,0){2}{\circle*{1}}
\multiput(10,8)(10,0){2}{\circle*{1}}
\put(6,12){\makebox(0,0){$p$}}
\put(24,12){\makebox(0,0){$q$}}
\put(6,-6){\makebox(0,0){$\td{r}$}}
\put(24,-6){\makebox(0,0){$\td{s}$}}
\end{picture} (\mbox{in $(\a \ab-{\bf 1})$})
$$

\bigskip
Together
with the information of $\a$, we obtain all the entries of the 
connection $(\a \ab-{\bf 1})\sigma \a$. Now we show the landscape of 
them in \framebox{Table \ref{land17.tex}}.

The four graphs of this connection are as in 
\framebox{Figure \ref{alle17graph.tex}}.

\newpage

\begin{figure}[H]
\begin{center}
\thinlines
\unitlength 1.0mm
\begin{picture}(30,20)
\multiput(11,-2)(0,10){2}{\line(1,0){8}}
\multiput(10,7)(10,0){2}{\line(0,-1){8}}
\multiput(10,-2)(10,0){2}{\circle*{1}}
\multiput(10,8)(10,0){2}{\circle*{1}}
\put(6,12){\makebox(0,0){$V_0$}}
\put(24,12){\makebox(0,0){$V_1$}}
\put(6,-6){\makebox(0,0){$V_2$}}
\put(24,-6){\makebox(0,0){$V_3$}}
\put(15,12){\makebox(0,0){${\cal G}_0$}}
\put(15,-6){\makebox(0,0){${\cal G}_2$}}
\put(6,3){\makebox(0,0){${\cal H}_0$}}
\put(24,3){\makebox(0,0){${\cal H}_1$}}
\end{picture} 
\end{center}
\begin{center}
	\setlength{\unitlength}{0.5mm}
	\begin{picture}(200,180)(-20,-100)
	\multiput(20,-50)(40,0){5}{\circle*{2}}
	\put(120,-50){\circle*{2}}
	\multiput(0,50)(40,0){6}{\circle*{2}}
	\multiput(100,50)(40,0){3}{\circle*{2}}
	\multiput(20,25)(40,0){5}{\circle*{2}}
	\put(120,25){\circle*{2}}
	\multiput(0,-25)(40,0){5}{\circle*{2}}
	\put(100,-25){\circle*{2}}
	\multiput(0,-100)(40,0){6}{\circle*{2}}
	\multiput(100,-100)(40,0){3}{\circle*{2}}
	
	\multiput(20,25)(40,0){2}{\line(2,-1){100}}
	\put(20,25){\line(140,-50){140}}
	\put(60,25){\line(2,-5){20}}
	\put(60,25){\line(4,-5){40}}
	\multiput(59,25)(2,0){2}{\line(6,-5){60}}
	\put(100,25){\line(-6,-5){60}}
	\multiput(99,25)(1,0){3}{\line(-2,-5){20}}
	\multiput(99.5,25)(1,0){2}{\line(0,-1){50}}
	\multiput(99,25)(1,0){4}{\line(2,-5){20}}
	\put(100,25){\line(6,-5){60}}
	\put(120,25){\line(-4,-5){40}}
	\put(120,25){\line(-2,-5){20}}
	\multiput(119.5,25)(1,0){2}{\line(0,-1){50}}
	\put(120,25){\line(4,-5){40}}
	\put(140,25){\line(-140,-50){140}}
	\multiput(139,25)(2,0){2}{\line(-2,-1){100}}
	\multiput(139,25)(2,0){2}{\line(-6,-5){60}}
	\put(140,25){\line(-4,-5){40}}
	\multiput(139.5,25)(1,0){2}{\line(-2,-5){20}}
	\put(140,25){\line(2,-5){20}}
	\put(180,25){\line(-2,-5){20}}
	\put(180,25){\line(-6,-5){60}}
	
	\multiput(20,-50)(40,0){3}{\line(8,-5){80}}
	\multiput(20,-50)(40,0){2}{\line(12,-5){120}}
	\put(20,-50){\line(16,-5){160}}
	\multiput(60,-50)(40,0){3}{\line(6,-5){60}}
	\multiput(59,-50)(2,0){2}{\line(4,-5){40}}
	\put(60,-50){\line(2,-5){20}}
	\multiput(100,-50)(20,0){3}{\line(4,-5){40}}
	\multiput(99,-50)(1,0){3}{\line(2,-5){20}}
	\multiput(99,-50)(1,0){3}{\line(-2,-5){20}}
    \multiput(99,-50)(1,0){3}{\line(0,-1){50}}
    \put(100,-50){\line(-6,-5){60}}
	\put(120,-50){\line(0,-1){50}}
	\multiput(119.5,-50)(1,0){2}{\line(-2,-5){20}}
	\multiput(120,-50)(20,0){2}{\line(-4,-5){40}}
	\put(120,-50){\line(-8,-5){80}}
    \put(140,-50){\line(2,-5){20}}
	\multiput(139,-50)(2,0){2}{\line(-2,-5){20}}
	\put(140,-50){\line(-4,-5){40}}
	\multiput(139.5,-50)(2,0){2}{\line(-6,-5){60}}
	\put(140,-50){\line(-2,-1){100}}
	\multiput(140,-50)(40,0){2}{\line(-14,-5){140}}
	\put(180,-50){\line(-8,-5){80}}
	\put(180,-50){\line(-2,-5){20}}
	
	\multiput(0,-25)(40,0){5}{\line(4,-5){20}}
	\multiput(40,-25)(40,0){4}{\line(-4,-5){20}}
	\multiput(100,-25)(20,0){2}{\line(0,-1){25}}
	\multiput(0,50)(40,0){3}{\line(4,-5){20}}
	\multiput(40,50)(40,0){3}{\line(-4,-5){20}}
	\multiput(100,50)(20,0){3}{\line(0,-1){25}}
	\put(100,50){\line(8,-5){40}}
	\multiput(160,50)(20,0){2}{\line(-8,-5){40}}
	\put(160,50){\line(4,-5){21}}
	\put(200,50){\line(-4,-5){21}}
	\put(25,-48){\makebox(1,1){$A$}}
	\put(66,-48){\makebox(1,1){$C$}}
	\put(106,-49){\makebox(1,1){$E$}}
	\put(125,-46){\makebox(1,1){$\tilde{C}$}}
	\put(146,-49){\makebox(1,1){$G$}}
	\put(185,-50){\makebox(1,1){$\tilde{A}$}}
	\put(4,50){\makebox(1,1){*}}
	\put(44,50){\makebox(1,1){$b$}}
	\put(84,50){\makebox(1,1){$d$}}
	\put(96,50){\makebox(1,1){$f$}}
	\put(124,50){\makebox(1,1){$\tilde{d}$}}
	\put(144,50){\makebox(1,1){$h$}}
	\put(164,50){\makebox(1,1){$\tilde{b}$}}
	\put(184,50){\makebox(1,1){$\tilde{h}$}}
	\put(204,50){\makebox(1,1){$\tilde{*}$}}
	\put(25,26){\makebox(1,1){$a$}}
	\put(65,26){\makebox(1,1){$c$}}
	\put(105,26){\makebox(1,1){$e$}}
	\put(126,24){\makebox(1,1){$\tilde{c}$}}
	\put(146,24){\makebox(1,1){$g$}}
	\put(185,25){\makebox(1,1){$\tilde{a}$}}
	\put(6,-27){\makebox(1,1){1}}
	\put(47,-26){\makebox(1,1){2}}
	\put(86,-27){\makebox(1,1){3}}
	\put(105,-25){\makebox(1,1){4}}
	\put(126,-25){\makebox(1,1){5}}
	\put(165,-25){\makebox(1,1){6}}
	\put(0,-107){\makebox(1,1){*}}
	\put(40,-106){\makebox(1,1){$b$}}
	\put(80,-105){\makebox(1,1){$d$}}
	\put(100,-106){\makebox(1,1){$f$}}
	\put(120,-106){\makebox(1,1){$\tilde{d}$}}
	\put(140,-106){\makebox(1,1){$h$}}
	\put(160,-107){\makebox(1,1){$\tilde{b}$}}
	\put(180,-107){\makebox(1,1){$\tilde{h}$}}
	\put(200,-106){\makebox(1,1){$\tilde{*}$}}
	\put(-10,37){\makebox(2,2){${\cal G}_0$}}
	\put(-10,0){\makebox(2,2){${\cal H}_1$}}
	\put(-10,-38){\makebox(2,2){${\cal G}_2$}}
	\put(-10,-75){\makebox(2,2){${\cal H}_0$}}
	\put(-15,50){\makebox(2,2){$V_0$}}
	\put(-15,25){\makebox(2,2){$V_1$}}
	\put(-15,-25){\makebox(2,2){$V_2$}}
	\put(-15,-50){\makebox(2,2){$V_3$}}
	\put(-15,-100){\makebox(2,2){$V_0$}}
	\end{picture}
\end{center}
\caption{Four graphs of the connection $(\a \ab-{\bf 1})\sigma \a$}
\label{alle17graph.tex}
\end{figure}

\newpage

\thispagestyle{plain}

\begin{table}[H]
\vspace{-2cm}
\hspace*{-1.5cm}
\begin{tabular}[pos]{|@{}c@{}||@{}c@{}||@{}c@{}|@{}c@{}||@{}c@{}
|@{}c@{}|@{}c@{}|@{}c@{}||@{}c@{}|@{}c@{}|@{}c@{}|@{}c@{}||@{}c@{}|@{}c@{}
|@{}c@{}|@{}c@{}|@{}c@{}|@{}c@{}||
@{}c@{}|@{}c@{}|@{}c@{}|@{}c@{}|@{}c@{}|@{}c@{}|}
\hl
 & 1 & \mcol{2}{c|}{2} & \mcol{4}{c|}{3} & \mcol{4}{c|}{4} &
 \mcol{6}{c|}{5} & \mcol{6}{c|}{6} \\ \hl
& $g1$ & $e$ & $g$ & $c$ & $e$ & $\td{c}$ & $g$ & $c$ & $e$
& $\td{c}$ & $g$ & $a$ & $c$ & $e$ & $\td{c}$ & $g$ & $\td{a}$
& $a$ & $c$ & $e$ & $\td{c}$ & $g$ & $\td{a}$ \\ \hline\hline
$*G$ &&&&&&&&&&&& \bl &&&&&& \bl &&&&& \\ \hline\hline
$bE$ &&&& \bl &&&& \bl &&&& \bl& \bl \bl& &&&&&&&&& \\ \hl
$b\td{C}$ &&&& &&&& &&&& \bl& \bl \bl&  &&&&&&&&& \\ \hl
$bG$ &&&&& &&&& &&& \bl& \bl \bl & &&&& \bl& \bl&&&& \\ \hl
$b\td{A}$ &&&&& &&&& &&&&& &    &&& \bl& \bl&&&& \\ \hline\hline
$dC$ && \bl && \bl&\bl \bl\bl&&& &&&&& &   &&&&&& &&& \\ \hl
$dE$ && && \bl&\bl  \bl \bl&& & \bl&\bl\bl &&& & \bl \bl & \bl \bl
\bl \bl && &&&& &&& \\ 
{} && && \bl &\bl  \bl \bl&& & \bl&\bl\bl &&& & \bl \bl & \bl \bl
\bl \bl && &&&& &&& \\ 
{}&& && \bl &\bl  \bl \bl&& & \bl&\bl\bl &&& & \bl \bl & \bl \bl
\bl \bl && &&&& &&& \\ \hl
$d\td{C}$ && &&&&&& && &&&  \bl \bl & \bl \bl \bl \bl && &&&& 
&&& \\ 
\hl
$dG$ && &&&&&& && &&&  \bl \bl & \bl \bl \bl \bl && &&&\bl&\bl 
&&& \\ 
{} && &&&&&& && &&&  \bl \bl & \bl \bl \bl \bl && &&&\bl&\bl 
&&& \\ \hline\hline
$fA$ & \bl & \bl & \bl \bl &&&&& && &&& &&& &&&& 
&&& \\ \hl
$fC$ & & \bl & \bl \bl &&  \bl \bl \bl && \bl \bl & && &&& &&& &&&& 
&&& \\ 
{} & & \bl & \bl \bl &&  \bl \bl \bl && \bl \bl & && &&& &&& &&&& 
&&& \\ \hl
$fE$  & && &&  \bl \bl \bl && \bl \bl & & \bl \bl & & \bl && &
  \bl \bl \bl \bl  && \bl \bl &&&& &&& \\ 
{}  & && &&  \bl \bl \bl && \bl \bl & & \bl \bl & & \bl && &
  \bl \bl \bl \bl  && \bl \bl &&&& &&& \\ 
{}  & && &&  \bl \bl \bl && \bl \bl & & \bl \bl & & \bl && &
  \bl \bl \bl \bl  && \bl \bl &&&& &&& \\ \hl
$f \td{C}$ &&&&&&&&&&&&&& \bl \bl \bl \bl  && \bl \bl &&&& &&& \\ 
{} &&&&&&&&&&&&&& \bl \bl \bl \bl  && \bl \bl &&&& &&& \\ \hl
$fG$ &&&&&&&&&&&&&& \bl \bl \bl \bl  && \bl \bl &&&& \bl && \bl & 
\\ \hl
$f \td{A}$  &&&&&&&&&&&&&& && &&&& \bl && \bl & \\ \hline\hline
$hA$ & \bl & & \bl \bl &&&&& && &&& &&& &&&& &&& \\ \hl
$hC$  & & & \bl \bl & & & & \bl \bl & && &&& &&& &&&& &&& \\ \hl
$hE$   & & &&&&& \bl \bl & && & \bl && &&& \bl \bl &&&& &&& \\ \hl
$hG$   & & &&&&&  & && & && &&& \bl \bl &&&& && \bl & \\ \hline\hline
$\td{h}A$ & \bl & & \bl \bl &&&&& && &&& &&& &&&& &&& \\ \hl
$\td{h}C$  & & & \bl \bl & & & & \bl \bl & && &&& &&& &&&& &&& \\ \hl
$\td{h}E$   & & &&&&& \bl \bl & && & \bl && &&& \bl \bl &&&& &&& \\ \hl
$\td{h}G$   & & &&&&&  & && & && &&& \bl \bl &&&& && \bl & \\ 
\hline\hline
$\td{d}C$  & &\bl & &&\bl \bl \bl &\bl && && &&& &&& &&&& &&& \\ \hl
$\td{d}E$  & && &&\bl \bl \bl &\bl && & \bl \bl & \bl &&& & \bl \bl 
\bl \bl & \bl \bl & &&&& &&& \\ 
{} & && &&\bl \bl \bl &\bl && & \bl \bl & \bl &&& & \bl \bl 
\bl \bl & \bl \bl & &&&& &&& \\ 
{} & && &&\bl \bl \bl &\bl && & \bl \bl & \bl &&& & \bl \bl 
\bl \bl & \bl \bl & &&&& &&& \\ \hline
$\td{d}\td{C}$ &&&&&&&&&&&&&& \bl \bl \bl \bl & \bl \bl &&&& &&&& \\ \hl
$\td{d}G$ &&&&&&&&&&&&&& \bl \bl \bl \bl & \bl \bl &&&&& \bl & \bl && \\
{}&&&&&&&&&&&&&& \bl \bl \bl \bl & \bl \bl &&&&& \bl & \bl && \\ 
\hline\hline
$\td{b}E$ &&&&&& \bl &&&& \bl &&&&& \bl \bl && \bl && & &&& \\ \hl
$\td{b}\td{C}$ &&&&&& &&&& &&&&& \bl \bl && \bl && & &&& \\ \hl
$\td{b}G$ &&&&&& &&&& &&&&& \bl \bl && \bl && & &\bl && \bl \\ \hl
$\td{b}\td{A}$ &&&&&& &&&& &&&&& &&&& & &\bl && \bl \\ \hline\hline
$\td{*}G$  &&&&&& &&&& &&&&& && \bl &&& &&& \bl \\ \hl

\end{tabular}
\caption{Landscape of $(\a \ab-{\bf 1})\sigma \a$}
\label{land17.tex}
\end{table}
\newpage
Since the exact values of the connection take up too much room to be 
listed up in a Table, we will show them in the shape of unitary 
matrices. Table \ref{land17.tex} gives also an overview of the 
connection $\sigma(\a \ab-{\bf 1})\sigma \a$, because 
it is easy to check that $\sigma(\a \ab-{\bf 1})\sigma \a$ has 
exactly the same vertical edges as $(\a \ab-{\bf 1})\sigma \a$.

Below we list all the entries of $(\a \ab-{\bf 1})\sigma \a$. These 
entries can be obtained by direct multiplication of the connections $(\a 
\ab-{\bf 1})\sigma$ and $\a$ as explained in section 3. In the list 
we have labeled rows and columns of the unitary matrices according to 
those entries that have to be used in he direct multiplication, for 
instance, in the $2 \times 2$-matrix 
\thinlines
\unitlength 0.5mm
\begin{picture}(30,20)(0, -6)
\multiput(11,-2)(0,10){2}{\line(1,0){8}}
\multiput(10,7)(10,0){2}{\line(0,-1){8}}
\multiput(10,-2)(10,0){2}{\circle*{1}}
\multiput(10,8)(10,0){2}{\circle*{1}}
\put(6,12){\makebox(0,0){$b$}}
\put(24,-6){\makebox(0,0){$6$}}
\end{picture}
below, the entry with row-label $G^{{\td b}}$ and column-label $a_{\td 
a}$ is computed as follows:
\begeq
\thinlines
\unitlength 0.5mm
\begin{picture}(30,30)
\multiput(11,-2)(0,10){2}{\line(1,0){8}}
\multiput(10,7)(10,0){2}{\line(0,-1){8}}
\multiput(10,-2)(10,0){2}{\circle*{1}}
\multiput(10,8)(10,0){2}{\circle*{1}}
\put(6,12){\makebox(0,0){$b$}}
\put(24,12){\makebox(0,0){$a_{\tilde{a}}$}}
\put(5,-6){\makebox(0,0){$G^{\tilde{b}}$}}
\put(24,-6){\makebox(0,0){$6$}}
\end{picture}
&=& \quad
\thinlines
\unitlength 0.5mm
\begin{picture}(50,30)
\multiput(11,-2)(0,10){2}{\line(1,0){38}}
\multiput(10,7)(40,0){2}{\line(0,-1){8}}
\multiput(10,-2)(40,0){2}{\circle*{1}}
\multiput(10,8)(40,0){2}{\circle*{1}}
\put(6,12){\makebox(0,0){$b$}}
\put(54,12){\makebox(0,0){$a$}}
\put(5,-6){\makebox(0,0){${\tilde b}$}}
\put(54,-6){\makebox(0,0){${\td a}$}}
\put(30, 3){\makebox(0,0){{\sm $(\a \ab-{\bf 1})\sigma$}}}
\end{picture} \quad
\cdot
\thinlines
\unitlength 0.5mm
\begin{picture}(30,30)
\multiput(11,-2)(0,10){2}{\line(1,0){8}}
\multiput(10,7)(10,0){2}{\line(0,-1){8}}
\multiput(10,-2)(10,0){2}{\circle*{1}}
\multiput(10,8)(10,0){2}{\circle*{1}}
\put(6,12){\makebox(0,0){${\td b}$}}
\put(24,12){\makebox(0,0){${\tilde{a}}$}}
\put(5,-6){\makebox(0,0){$G$}}
\put(24,-6){\makebox(0,0){$6$}}
\put(15,3){\makebox(0,0){$\a$}}
\end{picture} \\
&=&
\thinlines
\unitlength 0.5mm
\begin{picture}(40,30)
\multiput(11,-2)(0,10){2}{\line(1,0){28}}
\multiput(10,7)(30,0){2}{\line(0,-1){8}}
\multiput(10,-2)(30,0){2}{\circle*{1}}
\multiput(10,8)(30,0){2}{\circle*{1}}
\put(6,12){\makebox(0,0){$b$}}
\put(44,12){\makebox(0,0){$a$}}
\put(5,-6){\makebox(0,0){$b$}}
\put(44,-6){\makebox(0,0){$a$}}
\put(25, 3){\makebox(0,0){{\sm $\a \ab-{\bf 1}$}}}
\end{picture} \quad
\cdot
\thinlines
\unitlength 0.5mm
\begin{picture}(30,30) 
\multiput(11,-2)(0,10){2}{\line(1,0){8}}
\multiput(10,7)(10,0){2}{\line(0,-1){8}}
\multiput(10,-2)(10,0){2}{\circle*{1}}
\multiput(10,8)(10,0){2}{\circle*{1}}
\put(6,12){\makebox(0,0){${\td b}$}}
\put(24,12){\makebox(0,0){${\tilde{a}}$}}
\put(5,-6){\makebox(0,0){$G$}}
\put(24,-6){\makebox(0,0){$6$}}
\put(15,3){\makebox(0,0){$\a$}}
\end{picture} \\
&& \\
&=& \quad ( -\frac{\sqrt {{\beta}^2-2}}{\sqrt{{\beta}^2-1}} ) \cdot 1,
\endeq
where the last equality is obtained from the tables 
\ref{aab-117-1.ps} and \ref{a13}. Sometimes the entries listed below 
appear at first glance to be different from the entries obtained by 
direct multiplication. However in all those cases, it is just a 
different representation of the same algebraic number.This can easily 
be checked using the following identities for 
$\gbs=\sqrt{\frac{5+\sqrt{17}}{2}}$:
$$
\begin{array}{cc}
\gb+1=\frac{2(\gb-1)^{2}}{\gb}, & \gb-4=\frac{2}{\gb-1}, \\
\gb+2=\frac{4\gbf}{(\gb-1)^{2}} & 5-\gb=\frac{2}{\gb}, \\
\gb+3=\frac{(\gb-1)^{2}}{2\gbf}(\gbf+4), & 
2\gb-1=\frac{\gb(\gb-1)}{2}, \\
\gb-2=\frac{2\gb}{\gb-1}, & 3\gb-1=(\gb-1)^{2}, \\
\gb-3=\frac{2(\gb-1)}{\gb}, & 3\gb-4=\frac{\gb-1}{2\gb}(\gbf+4).
\end{array}
$$
Here comes the list of entries of $(\a \ab-{\bf 1})\sigma \a$:


\begin{displaymath}
\thinlines
\unitlength 0.5mm
\begin{picture}(30,20)
\multiput(11,-2)(0,10){2}{\line(1,0){8}}
\multiput(10,7)(10,0){2}{\line(0,-1){8}}
\multiput(10,-2)(10,0){2}{\circle*{1}}
\multiput(10,8)(10,0){2}{\circle*{1}}
\put(6,12){\makebox(0,0){$*$}}
\put(24,12){\makebox(0,0){$a_{\tilde{a}}$}}
\put(5,-6){\makebox(0,0){$G^{\tilde{b}}$}}
\put(24,-6){\makebox(0,0){$6$}}
\end{picture}
= \thinlines
\unitlength 0.5mm
\begin{picture}(30,20)
\multiput(11,-2)(0,10){2}{\line(1,0){8}}
\multiput(10,7)(10,0){2}{\line(0,-1){8}}
\multiput(10,-2)(10,0){2}{\circle*{1}}
\multiput(10,8)(10,0){2}{\circle*{1}}
\put(6,12){\makebox(0,0){$*$}}
\put(24,12){\makebox(0,0){$a$}}
\put(6,-6){\makebox(0,0){${\tilde{b}}$}}
\put(24,-6){\makebox(0,0){$\tilde{a}$}}
\end{picture}
\cdot
 \thinlines
\unitlength 0.5mm
\begin{picture}(30,20)
\multiput(11,-2)(0,10){2}{\line(1,0){8}}
\multiput(10,7)(10,0){2}{\line(0,-1){8}}
\multiput(10,-2)(10,0){2}{\circle*{1}}
\multiput(10,8)(10,0){2}{\circle*{1}}
\put(6,12){\makebox(0,0){$\td{b}$}}
\put(24,12){\makebox(0,0){$\td{a}$}}
\put(6,-6){\makebox(0,0){$G$}}
\put(24,-6){\makebox(0,0){$6$}}
\end{picture}
= \quad 1,
\end{displaymath}

\begin{displaymath}
\thinlines
\unitlength 0.5mm
\begin{picture}(30,20)
\multiput(11,-2)(0,10){2}{\line(1,0){8}}
\multiput(10,7)(10,0){2}{\line(0,-1){8}}
\multiput(10,-2)(10,0){2}{\circle*{1}}
\multiput(10,8)(10,0){2}{\circle*{1}}
\put(6,12){\makebox(0,0){$\tilde{*}$}}
\put(24,12){\makebox(0,0){$\tilde{a}_g$}}
\put(6,-6){\makebox(0,0){$\tilde{G}^h$}}
\put(24,-6){\makebox(0,0){$6$}}
\end{picture}
=\quad 1,
\end{displaymath}

\begin{displaymath}
\thinlines
\unitlength 0.5mm
\begin{picture}(30,20)
\multiput(11,-2)(0,10){2}{\line(1,0){8}}
\multiput(10,7)(10,0){2}{\line(0,-1){8}}
\multiput(10,-2)(10,0){2}{\circle*{1}}
\multiput(10,8)(10,0){2}{\circle*{1}}
\put(6,12){\makebox(0,0){$*$}}
\put(24,12){\makebox(0,0){$a_{\tilde{c}}$}}
\put(5,-6){\makebox(0,0){$G^{\tilde{b}}$}}
\put(24,-6){\makebox(0,0){$5$}}
\end{picture}
=\quad1,
\end{displaymath}

\begin{displaymath}
\thinlines
\unitlength 0.5mm
\begin{picture}(30,20)
\multiput(11,-2)(0,10){2}{\line(1,0){8}}
\multiput(10,7)(10,0){2}{\line(0,-1){8}}
\multiput(10,-2)(10,0){2}{\circle*{1}}
\multiput(10,8)(10,0){2}{\circle*{1}}
\put(6,12){\makebox(0,0){$\tilde{*}$}}
\put(24,12){\makebox(0,0){$\tilde{a}_g$}}
\put(6,-6){\makebox(0,0){$G^{\tilde{h}}$}}
\put(24,-6){\makebox(0,0){$5$}}
\end{picture}
=\quad 1,
\end{displaymath}

\begin{displaymath}
\thinlines
\unitlength 0.5mm
\begin{picture}(30,20)
\multiput(11,-2)(0,10){2}{\line(1,0){8}}
\multiput(10,7)(10,0){2}{\line(0,-1){8}}
\multiput(10,-2)(10,0){2}{\circle*{1}}
\multiput(10,8)(10,0){2}{\circle*{1}}
\put(6,12){\makebox(0,0){$b$}}
\put(24,-6){\makebox(0,0){$6$}}
\end{picture}
=\quad
\bordermatrix{
&  {a}_{\tilde{a}} & c_{\tilde{a}}\cr
 {{\tilde{A}}^{\tilde{*}}} &\frac{1}{\sqrt{{\beta}^2-1}}  & 
\frac{\sqrt {{\beta}^2-2}}{\sqrt{{\beta}^2-1}}\cr
 {G^{\tilde{b}}} & -\frac{\sqrt {{\beta}^2-2}}{\sqrt{{\beta}^2-1}} 
&\frac{1}{\sqrt{{\beta}^2-1}} \cr},
\end{displaymath}


\begin{displaymath}
\thinlines
\unitlength 0.5mm
\begin{picture}(30,20)
\multiput(11,-2)(0,10){2}{\line(1,0){8}}
\multiput(10,7)(10,0){2}{\line(0,-1){8}}
\multiput(10,-2)(10,0){2}{\circle*{1}}
\multiput(10,8)(10,0){2}{\circle*{1}}
\put(6,12){\makebox(0,0){$\tilde{b}$}}
\put(24,-6){\makebox(0,0){6}}
\end{picture}
=\quad
\bordermatrix{
&  {\td{a}}_{\tilde{a}} & {\td c}_{\tilde{g}}\cr
 {{\tilde{A}}^{\tilde{h}}} &\frac{1}{\sqrt{{\beta}^2-1}}  & 
\frac{\sqrt {{\beta}^2-2}}{\sqrt{{\beta}^2-1}}\cr
 {G^{f}} & \frac{\sqrt {{\beta}^2-2}}{\sqrt{{\beta}^2-1}} & 
-\frac{1}{\sqrt{{\beta}^2-1}} \cr
},
\end{displaymath}


\begin{displaymath}
\thinlines
\unitlength 0.5mm
\begin{picture}(30,20)
\multiput(11,-2)(0,10){2}{\line(1,0){8}}
\multiput(10,7)(10,0){2}{\line(0,-1){8}}
\multiput(10,-2)(10,0){2}{\circle*{1}}
\multiput(10,8)(10,0){2}{\circle*{1}}
\put(6,12){\makebox(0,0){$b$}}
\put(24,-6){\makebox(0,0){$5$}}
\end{picture}
=\quad
\bordermatrix{
&  a_{\tilde{c}}&{c_2}_{\tilde{c}}  &c_e  \cr
 {G^{\tilde{b}}} & \frac{1}{\sqrt{2}\beta} & -\frac{\sqrt 
{{2\beta}^2-1}}{\sqrt{2}\beta}  & 0  \cr
 {E^{\tilde{d}}} & \frac{\sqrt{\beta^2+1}}{\sqrt{2}\beta} & 
\frac{1}{\sqrt{2}\beta}\sqrt{\frac{\beta^2+1}{2\beta^2-1}}  & 
-\sqrt{\frac{\beta^2-2}{2\beta^2-1}}  \cr
 {{\tilde{C}}^{\tilde{d}}} & \frac{\sqrt 
{{\beta}^2-2}}{\sqrt{2}\beta}  & 
\frac{1}{\sqrt{2}\beta}\sqrt{\frac{\beta^2-2}{2\beta^2-1}}  & 
\sqrt{\frac{\beta^2+1}{2\beta^2-1}} \cr
},
\end{displaymath}


\begin{displaymath}
\thinlines
\unitlength 0.5mm
\begin{picture}(30,20)
\multiput(11,-2)(0,10){2}{\line(1,0){8}}
\multiput(10,7)(10,0){2}{\line(0,-1){8}}
\multiput(10,-2)(10,0){2}{\circle*{1}}
\multiput(10,8)(10,0){2}{\circle*{1}}
\put(6,12){\makebox(0,0){$\tilde{b}$}}
\put(24,-6){\makebox(0,0){5}}
\end{picture}
=\quad
\bordermatrix{
& \tilde{a}_g  & \tilde{c}_g  &  \tilde{c}_e \cr
 {E^f} & \sqrt{\frac{\gb+1}{2\gb}}  & -\sqrt{\frac{\gb+1}{2\gb(\gb-2)}} & 
\sqrt{\frac{\gb-1}{2\gb(\gb-2)}} \cr
 {G^f} & -\frac{1}{\sqrt{2\gb}} & \frac{1}{\sqrt{2\gb(\gb-2)}} & 
\sqrt{\frac{\gbf-1}{2\gb(\gb-2)}}   \cr
 {\tilde{C}^{\tilde{h}}} &\frac{1}{\sqrt{\gb-1}} & 
\sqrt{\frac{\gb-2}{\gb-1}}  & 0  \cr
},
\end{displaymath}

\begin{displaymath}
\thinlines
\unitlength 0.5mm
\begin{picture}(30,20)
\multiput(11,-2)(0,10){2}{\line(1,0){8}}
\multiput(10,7)(10,0){2}{\line(0,-1){8}}
\multiput(10,-2)(10,0){2}{\circle*{1}}
\multiput(10,8)(10,0){2}{\circle*{1}}
\put(6,12){\makebox(0,0){$d$}}
\put(24,-6){\makebox(0,0){$6$}}
\end{picture}
=\quad
\bordermatrix{
&  e_g & c_{\tilde{a}} \cr
 {G^f} & 1 & 0\cr
 {G^{\tilde{b}}} & 0 & 1\cr},
\end{displaymath}


\begin{displaymath}
\thinlines
\unitlength 0.5mm
\begin{picture}(30,20)
\multiput(11,-2)(0,10){2}{\line(1,0){8}}
\multiput(10,7)(10,0){2}{\line(0,-1){8}}
\multiput(10,-2)(10,0){2}{\circle*{1}}
\multiput(10,8)(10,0){2}{\circle*{1}}
\put(6,12){\makebox(0,0){$\tilde{d}$}}
\put(24,-6){\makebox(0,0){6}}
\end{picture}
=\quad
\bordermatrix{
  & e_g & {\tilde{c}_g} \cr
 {G^f} & -\sqrt{\frac{\gb-2}{2\gb-1}}  & \sqrt{\frac{\gb+1}{2\gb-1}} 
\cr
 {G^h} & \sqrt{\frac{\gb+1}{2\gb-1}} & \sqrt{\frac{\gb-2}{2\gb-1}}  
\cr
},
\end{displaymath}


\begin{displaymath} 
\thinlines
\unitlength 0.5mm
\begin{picture}(30,20)
\multiput(11,-2)(0,10){2}{\line(1,0){8}}
\multiput(10,7)(10,0){2}{\line(0,-1){8}}
\multiput(10,-2)(10,0){2}{\circle*{1}}
\multiput(10,8)(10,0){2}{\circle*{1}}
\put(6,12){\makebox(0,0){$d$}}
\put(24,-6){\makebox(0,0){$5$}}
\end{picture} 
=
\end{displaymath}

\begin{flushleft}
$$\bordermatrix{
 &  {c_g} & c_e & {e_2}_e & {e_3}_e & e_{\tilde{c}} & e_g 
\cr
 {\tilde{C}^{\tilde{d}}} &\frac{-\sqrt{\gbf-1}}{2\gb-1} & 
\frac{1}{2\gb-1} & \sqrt{\frac{\gbf+4}{\gb(2\gb-1)}}  & 0  & 
\frac{1}{2\gb-1} & 0  \cr
 {E^{\tilde{d}}} & \frac{-\sqrt{2}(\gb+1)}{\beta^3} & 
\frac{-(\gb-2)}{\gb\sqrt{(\gb-1)^3}} & 
\frac{-\sqrt{2(\gbf+4)}}{(\gb-1)^2} & 0 & 
\frac{1}{\gb}\sqrt{\frac{\gb+1}{2\gb-1}} & 0 \cr
 {E^f} & 0 & \frac{1}{\gb-1} & 
\frac{-(\gb-2)}{\gb\sqrt{2(\gb-1)(\gbf+4)}} & 
\frac{\sqrt{\gb-1}}{\sqrt{2(\gbf+4)}} & 0  &  
\frac{\sqrt{\gb+1}}{\gb-2} \cr
 {E^d} & 0 & \frac{\gb}{\sqrt{(\gb-1)^3}} & 
\frac{-\sqrt{\gb(\gb-2)}}{\sqrt{(\gb-1)^3(\gbf+4)}} & 
\frac{-2(\gb-2)}{\sqrt{2(\gbf+4)}} & 0 & 0 \cr
 {G^f} & 0 & \sqrt{\frac{2}{\gb}} & 
\frac{-2}{\sqrt{\gb(\gb-1)(\gbf+4)}} & 
\sqrt{\frac{\gb(\gb+1)}{(\gb-2)(\gbf+4)}} & 0 & \frac{-1}{\gb-2}  \cr
 {G^{\tilde{b}}} & \frac{1}{\gb} & 0 & 0 & 0 & 
\frac{\sqrt{\gbf-1}}{\gb} & 0  \cr
}, $$ \end{flushleft} 


\begin{displaymath}
\thinlines
\unitlength 0.5mm
\begin{picture}(30,20)
\multiput(11,-2)(0,10){2}{\line(1,0){8}}
\multiput(10,7)(10,0){2}{\line(0,-1){8}}
\multiput(10,-2)(10,0){2}{\circle*{1}}
\multiput(10,8)(10,0){2}{\circle*{1}}
\put(6,12){\makebox(0,0){$\tilde{d}$}}
\put(24,-6){\makebox(0,0){5}}
\end{picture}  
= \end{displaymath}
 \begin{flushleft}
 $$\bordermatrix{
& {\tilde{c}}_g & {\tilde{c}}_e &{e_2}_e &  {e_3}_e & e_{\tilde{c}} & 
e_g \cr
 {{\tilde{C}}^{\tilde{d}}} & 0 & \frac{\gb-2}{\gb} & 
\frac{-1}{2\gb-1}\sqrt{\frac{\gb(\gb-2)}{\gbf+4}} & 
-\sqrt{\frac{4(\gb+1)}{(\gb+3)(2\gb-1)}} & \frac{\gb-2}{\gb} & 0  \cr
 {E^{\tilde{d}}} &  0 & \frac{-(\gb-2)}{\sqrt{(\gb+1)(2\gb-1)}} & 
\frac{\gb-2}{2\gb-1}\sqrt{\frac{\gb}{(\gbf+4)(\gb+1)}} & 
\sqrt{\frac{4(\gb-2)}{(\gb+3)(2\gb-1)}} & \sqrt{\frac{\gb+1}{2\gb-1}} 
& 0  \cr
{E^f} &  \frac{(\gb+1)\sqrt{\gb-1}}{\sqrt{2}\beta^3} & \frac{1}{2\gb} 
& \frac{-\sqrt{\gb-1}}{2\sqrt{\gb+3}} & \frac{\sqrt{\gb-1}}{\gb 
\sqrt{\gb+3}} & 0 & \frac{-\sqrt{\gb+1}}{\gb}  \cr
 {E^d} & 0 & \sqrt{\frac{\gb-2}{\gb+1}} & 
\sqrt{\frac{\gb(\gb+2)}{(\gb+1)(\gbf+4)}} & \frac{4}{\sqrt{2(\gbf+4)}} 
& 0 & 0  \cr
{G^f} & -\sqrt{\frac{\gbf-1}{2\beta^6}} & 
\frac{1}{\gb}\sqrt{\frac{\gb-1}{\gb-2}} & 
-\sqrt{\frac{\gb(\gb-1)}{\gbf+4}} & 
\frac{2\sqrt{\gb-1}}{\sqrt{\gb(\gbf+4)}} & 0 & \frac{1}{\gb} \cr
 {G^h} & \sqrt{\frac{\gb-2}{2\gb-1}} & 0 & 0& 0&0& 
\sqrt{\frac{\gb+1}{2\gb-1}} \cr
} ,$$
\end{flushleft}


\begin{displaymath}
\thinlines
\unitlength 0.5mm
\begin{picture}(30,20)
\multiput(11,-2)(0,10){2}{\line(1,0){8}}
\multiput(10,7)(10,0){2}{\line(0,-1){8}}
\multiput(10,-2)(10,0){2}{\circle*{1}}
\multiput(10,8)(10,0){2}{\circle*{1}}
\put(6,12){\makebox(0,0){$f$}}
\put(24,-6){\makebox(0,0){2}}
\end{picture}
=\quad
\bordermatrix{
 &  e_c & g_a & g_c  \cr
 {A^b} & \sqrt{\frac{\gb+1}{2(\gb-1)}}  
& -\frac{1}{\gb-1} & -\frac{1}{\sqrt{2(\gb-1)}}  \cr
 {C^d} &  \sqrt{\frac{\gb-1}{2\gb(\gb-2)}} & 
 0 &  \sqrt{\frac{\gbf-1}{2\gb(\gb-2)}} \cr
 {C^b} &\frac{\sqrt{\gb+1}}{2\gb} & 
\frac{\sqrt{\gb(\gb-2)}}{\gb-1} & -\frac{1}{2\gb}  \cr
},
\end{displaymath}

\begin{eqnarray*}
&&
\thinlines
\unitlength 0.5mm
\begin{picture}(30,20)
\multiput(11,-2)(0,10){2}{\line(1,0){8}}
\multiput(10,7)(10,0){2}{\line(0,-1){8}}
\multiput(10,-2)(10,0){2}{\circle*{1}}
\multiput(10,8)(10,0){2}{\circle*{1}}
\put(6,12){\makebox(0,0){$d$}}
\put(24,-6){\makebox(0,0){$3$}}
\end{picture}     \\
&=& \quad
\bordermatrix{
& {e_2}_e  &  {e_3}_e & e_c & c_{e}\cr
 {E^{\tilde{d}}} & \sqrt{\frac{\gbf+4}{\gb(\gb+1)}} & 0 & 0 & 
\frac{1}{\gb}\sqrt{\frac{\gb-2}{\gb+1}} \cr
 {E^f} & \frac{-1}{\gb-1}\sqrt{\frac{\gbf-4}{\gb(\gbf+4)}} & \gb 
\sqrt{\frac{2}{(\gb-1)(\gbf+4)}}  & 0 & \frac{\gb-2}{\gb-1}	\cr
 {E^d} &  -\frac{\gb-2}{(2\gb-1)(\gb-1)\sqrt{2(\gbf+4)}} & 
-\frac{\sqrt{2}(\gb-2)}{(2\gb-1)\sqrt{\gbf+4}} & \frac{4}{\sqrt{
\gb(\gb-1)}} 
& \frac{\gb}{\sqrt{(2\gb-1)^3(\gb+1)}}    \cr 
 {C^d} & -\frac{\gb(\gb-2)}{\sqrt{(2\gb-1)^3(\gbf+4)}} & 
-\frac{4(\gb-2)}{\sqrt{(\gbf+4)(2\gb-1)}} & -\frac{1}{2\gb-1} & 
\sqrt{\frac{\beta^6 (\gb-2)}{(2\gb-1)^3}}   \cr
}, 
\end{eqnarray*}


\begin{displaymath}
\thinlines
\unitlength 0.5mm
\begin{picture}(30,20)
\multiput(11,-2)(0,10){2}{\line(1,0){8}}
\multiput(10,7)(10,0){2}{\line(0,-1){8}}
\multiput(10,-2)(10,0){2}{\circle*{1}}
\multiput(10,8)(10,0){2}{\circle*{1}}
\put(6,12){\makebox(0,0){$\tilde{d}$}}
\put(24,-6){\makebox(0,0){3}}
\end{picture}
=\quad
\bordermatrix{
 &{e_2}_e &  {e_3}_e &  e_c &  {\tilde{c}}_e   \cr
{E^{\tilde{d}}} & -\sqrt{\frac{\gb(\gb-2)}{(\gb+1)(2\gb-1)(\gbf+4)}} 
& -\frac{2(\gb-2)}{\sqrt{2(\gbf+4)}} & 0 & 
\sqrt{\frac{\gb-2}{\gb+1}}   \cr
{E^f} &  
-\sqrt{\frac{2\gb}{2\gb-1}}\frac{\beta^3}{\sqrt{(\gb+1)(\gbf+4)}} & 
\sqrt{\frac{2\gb}{2\gb-1}}\sqrt{\frac{4\gb}{(\gb+1)(\gbf+4)}} & 0 & 
\sqrt{\frac{5-\gb}{\gb+1}}          \cr
{E^d} &  \frac{1}{2\gb-1}\sqrt{\frac{\gb(\gb+2)}{(\gb+1)(\gbf+4)}} & 
\frac{1}{2\gb-1}\frac{4}{\sqrt{2(\gbf+4)}} & 2\sqrt{\frac{2}{2\gb-1}} 
& \sqrt{\frac{\gb-2}{\gb+1}}\frac{1}{2\gb-1}       \cr
{C^d} &  
2\sqrt{\frac{2}{2\gb-1}}\sqrt{\frac{\gb(\gb+2)}{(\gb+1)(\gbf+4)}} & 
2\sqrt{\frac{2}{2\gb-1}}\frac{4}{\sqrt{2(\gbf+4)}} & -\frac{1}{2\gb-1} 
& \frac{\gbs(\gb-2)}{2\gb-1}         \cr
},
\end{displaymath}


\begin{displaymath}
\thinlines
\unitlength 0.5mm
\begin{picture}(30,20)
\multiput(11,-2)(0,10){2}{\line(1,0){8}}
\multiput(10,7)(10,0){2}{\line(0,-1){8}}
\multiput(10,-2)(10,0){2}{\circle*{1}}
\multiput(10,8)(10,0){2}{\circle*{1}}
\put(6,12){\makebox(0,0){$h$}}
\put(24,-6){\makebox(0,0){3}}
\end{picture}
=\quad
\bordermatrix{
&  g_c & g_e \cr
 {C^b} & 1 & 0\cr
 {E^f} & 0 & 1\cr},
\end{displaymath}


\begin{displaymath}
\thinlines
\unitlength 0.5mm
\begin{picture}(30,20)
\multiput(11,-2)(0,10){2}{\line(1,0){8}}
\multiput(10,7)(10,0){2}{\line(0,-1){8}}
\multiput(10,-2)(10,0){2}{\circle*{1}}
\multiput(10,8)(10,0){2}{\circle*{1}}
\put(6,12){\makebox(0,0){$\tilde{h}$}}
\put(24,-6){\makebox(0,0){3}}
\end{picture}
=\quad
\bordermatrix{
& g_c & g_e \cr
 {C^d} & -\frac{1}{2\gb-1}  & \frac{2\sqrt{2}}{\sqrt{2\gb-1}} \cr
 {E^d} & \frac{2\sqrt{2}}{\sqrt{2\gb-1}}  & \frac{1}{2\gb-1} \cr
},
\end{displaymath}


\begin{displaymath}
\thinlines
\unitlength 0.5mm
\begin{picture}(30,20)
\multiput(11,-2)(0,10){2}{\line(1,0){8}}
\multiput(10,7)(10,0){2}{\line(0,-1){8}}
\multiput(10,-2)(10,0){2}{\circle*{1}}
\multiput(10,8)(10,0){2}{\circle*{1}}
\put(6,12){\makebox(0,0){$f$}}
\put(24,-6){\makebox(0,0){3}}
\end{picture} 
= \end{displaymath}
\begin{flushleft}
$$ \bordermatrix{
&  g_c & g_e & {e_2}_e & {e_3}_e & e_c \cr
 {C^d} & \frac{-\sqrt{\gb+1}}{\gbf} & 
\frac{-2\sqrt{2}}{\sqrt{(\gb-1)^3}} & 
\frac{-2\sqrt{2}}{\sqrt{2\gb-1}}\frac{\beta^3}{\sqrt{(\gb+1)(\gbf+4)}} 
& \frac{2\sqrt{2}}{\sqrt{2\gb-1}}\sqrt{\frac{4\gb}{(\gb+1)(\gbf+4)}} & 
-\frac{1}{\gbf}        \cr
 {C^b} & \frac{-1}{\gb-2} & 0 & 0 & 0 & 
\frac{\sqrt{\gb+1}}{\gb-2}                                       \cr
 {E^{\tilde{d}}} & 0 & \sqrt{\frac{2\gb-1}{2\gb(\gb-3)}} & 
-\sqrt{\frac{(2\gb-1)(\gbf-4)}{2(\gbf+4)}}\frac{1}{\gb(\gb-1)} & 
\sqrt{\frac{\gb(2\gb-1)}{(\gb-1)(\gbf+4)}} & 0          \cr
 {E^f} & 0 & -\frac{1}{2} &  \sqrt{\frac{2\gbf}{(\gb+1)(\gbf+4)}} & 
\sqrt{\frac{2(\gb+1)}{\gbf+4}} & 0             \cr
 {E^d} & \frac{2\sqrt{\gbf-1}}{\beta^3} & 
\frac{-1}{(2\gb-1)\sqrt{\gb+1}} & 
-\frac{1}{2\gb-1}\frac{\beta^3}{\sqrt{(\gb+1)(\gbf+4)}} & 
\frac{1}{2\gb-1}\sqrt{\frac{4\gb}{(\gb+1)(\gbf+4)}} & 
2\frac{\sqrt{\gb-1}}{\beta^3}  \cr
},$$        
\end{flushleft}


\begin{displaymath}
\thinlines
\unitlength 0.5mm
\begin{picture}(30,20)
\multiput(11,-2)(0,10){2}{\line(1,0){8}}
\multiput(10,7)(10,0){2}{\line(0,-1){8}}
\multiput(10,-2)(10,0){2}{\circle*{1}}
\multiput(10,8)(10,0){2}{\circle*{1}}
\put(6,12){\makebox(0,0){$f$}}
\put(24,-6){\makebox(0,0){4}}
\end{picture}
=\quad
\bordermatrix{
&  g_e & e_{2e} &  e_{3e} \cr
 {E^{\tilde{d}}} & \sqrt{\frac{2\gb-1}{2\gb(\gb-3)}} & 
\frac{-(\gb-2)}{(\gb-1)\sqrt{2(\gbf+4)}} & 
\sqrt{\frac{\gb(2\gb-1)}{(\gb-1)(\gbf+4)}}  \cr
 {E^f} & \frac{1}{2} &  -\sqrt{\frac{2\gbf}{(\gb+1)(\gbf+4)}} & 
-\sqrt{\frac{2(\gb+1)}{\gbf+4}} \cr
 {E^d} & -\frac{1}{\sqrt{\gb+1}} & 
-\frac{\beta^3}{\sqrt{(\gb+1)(\gbf+4)}} & 
\sqrt{\frac{4\gb}{(\gb+1)(\gbf+4)}} \cr
},
\end{displaymath}


\begin{displaymath}
\thinlines
\unitlength 0.5mm
\begin{picture}(30,20)

\multiput(11,-2)(0,10){2}{\line(1,0){8}}
\multiput(10,7)(10,0){2}{\line(0,-1){8}}
\multiput(10,-2)(10,0){2}{\circle*{1}}
\multiput(10,8)(10,0){2}{\circle*{1}}
\put(6,12){\makebox(0,0){$h$}}
\put(24,-6){\makebox(0,0){5}}
\end{picture}
=\quad
\bordermatrix{
& g_g & g_e \cr
 {G^f} & -\sqrt{\frac{\gb-1}{2\gb(\gb-2)}}  & 
\sqrt{\frac{\gbf-1}{2\gb(\gb-2)}}\cr
 {E^f} & \sqrt{\frac{\gbf-1}{2\gb(\gb-2)}} & 
\sqrt{\frac{\gb-1}{2\gb(\gb-2)}} \cr
},
\end{displaymath}


\begin{displaymath}
\thinlines
\unitlength 0.5mm
\begin{picture}(30,20)
\multiput(11,-2)(0,10){2}{\line(1,0){8}}
\multiput(10,7)(10,0){2}{\line(0,-1){8}}
\multiput(10,-2)(10,0){2}{\circle*{1}}
\multiput(10,8)(10,0){2}{\circle*{1}}
\put(6,12){\makebox(0,0){$\tilde{h}$}}
\put(24,-6){\makebox(0,0){5}}
\end{picture}
=\quad
\bordermatrix{
&  g_g & g_e \cr
 {G^h} & 1 & 0\cr
 {E^d} & 0 & 1\cr},
\end{displaymath}


\begin{displaymath}
\thinlines
\unitlength 0.5mm
\begin{picture}(30,20)
\multiput(11,-2)(0,10){2}{\line(1,0){8}}
\multiput(10,7)(10,0){2}{\line(0,-1){8}}
\multiput(10,-2)(10,0){2}{\circle*{1}}
\multiput(10,8)(10,0){2}{\circle*{1}}
\put(6,12){\makebox(0,0){$f$}}
\put(24,-6){\makebox(0,0){5}}
\end{picture}
=\end{displaymath}
\begin{flushleft}
$$\bordermatrix{
& g_g & g_e & {e_2}_e & {e_3}_e & e_{\tilde{c}} & e_g \cr
 {E^{\tilde{d}}} & 0 & \frac{-(\gb-2)}{\sqrt{8\gb}} & 
\frac{\sqrt{(\gb-2)(\gbf-4)}}{\gb(\gb-1)\sqrt{2(\gbf+4)}} & 
-\sqrt{\frac{\gb(\gb-2)}{(\gb-1)(\gbf+4)}} & 
\sqrt{\frac{\gb+1}{2\gb-1}}  &  0 \cr
{E^f} &\frac{\gb+1}{2\gb}\sqrt{\frac{\gb-1}{\gb-2}} & 
\frac{-1}{2(\gb-2)} 
& \frac{\gbs}{2\sqrt{\gbf+4}} & 
\sqrt{\frac{(\gb+1)(\gb-1)}{\gb(\gb-2)(\gbf+4)}} & 0  
& \frac{-(\gb-1)}{2\gb}\sqrt{\frac{\gb+1}{\gb-2}}  \cr
{E^d} & 0 & \frac{-1}{\sqrt{\gb+1}} & 
\frac{-\beta^3}{\sqrt{(\gb+1)(\gbf+4)}} & 
\sqrt{\frac{4\gb}{(\gb+1)(\gbf+4)}} & 0 & 0  \cr
 \tilde{C}^{\tilde{d}} & 0 & \frac{1}{\sqrt{\gb-2}} & 
\frac{-\sqrt{(\gb+1)(\gbf-4)}}{\gb(\gb-1)\sqrt{2(\gbf+4)}} & 
\sqrt{\frac{2(\gb-1)}{\gbf+4}} & 
\sqrt{\frac{\gb-2}{2\gb-1}}    & 0          \cr
{{\tilde{C}^{\tilde{h}}}} & \sqrt{\frac{\gb-1}{2\gb}} & 0 & 0 & 0 & 
0 & \sqrt{\frac{\gb+1}{2\gb}}      \cr
{G^f} & \frac{-1}{2\gb}\sqrt{\frac{\gbf-1}{\gb-2}} & 
\frac{-1}{2}\sqrt{\frac{\gbf-1}{2\gb(\gb-2)}} & 
\sqrt{\frac{\gb(\gb-1)}{(\gbf+4)(\gb-2)}}  &  
\frac{(\gb+1)\sqrt{\gb-1}}{\sqrt{\gb(\gbf+4)(\gb-2)}} & 
0 & \frac{\gb-1}{2\gb}\frac{1}{\sqrt{\gb-2}}    \cr
},$$
\end{flushleft}
\begin{displaymath}
\thinlines
\unitlength 0.5mm
\begin{picture}(30,20)
\multiput(11,-2)(0,10){2}{\line(1,0){8}}
\multiput(10,7)(10,0){2}{\line(0,-1){8}}
\multiput(10,-2)(10,0){2}{\circle*{1}}
\multiput(10,8)(10,0){2}{\circle*{1}}
\put(6,12){\makebox(0,0){$b$}}
\put(24,12){\makebox(0,0){$c_e$}}
\put(6,-6){\makebox(0,0){$E^{\tilde{d}}$}}
\put(24,-6){\makebox(0,0){$3$}}
\end{picture}
=\quad1,
\end{displaymath}

\begin{displaymath}
\thinlines
\unitlength 0.5mm
\begin{picture}(30,20)
\multiput(11,-2)(0,10){2}{\line(1,0){8}}
\multiput(10,7)(10,0){2}{\line(0,-1){8}}
\multiput(10,-2)(10,0){2}{\circle*{1}}
\multiput(10,8)(10,0){2}{\circle*{1}}
\put(6,12){\makebox(0,0){$\tilde{b}$}}
\put(24,12){\makebox(0,0){$\tilde{c}_e$}}
\put(6,-6){\makebox(0,0){$E^f$}}
\put(24,-6){\makebox(0,0){$3$}}
\end{picture}
=\quad 1 ,
\end{displaymath}
\begin{displaymath}
\thinlines
\unitlength 0.5mm
\begin{picture}(30,20)
\multiput(11,-2)(0,10){2}{\line(1,0){8}}
\multiput(10,7)(10,0){2}{\line(0,-1){8}}
\multiput(10,-2)(10,0){2}{\circle*{1}}
\multiput(10,8)(10,0){2}{\circle*{1}}
\put(6,12){\makebox(0,0){$b$}}
\put(24,12){\makebox(0,0){$\tilde{c}_e$}}
\put(6,-6){\makebox(0,0){$E^{\tilde{d}}$}}
\put(24,-6){\makebox(0,0){$4$}}
\end{picture}
=\quad1,
\end{displaymath}

\begin{displaymath}
\thinlines
\unitlength 0.5mm
\begin{picture}(30,20)
\multiput(11,-2)(0,10){2}{\line(1,0){8}}
\multiput(10,7)(10,0){2}{\line(0,-1){8}}
\multiput(10,-2)(10,0){2}{\circle*{1}}
\multiput(10,8)(10,0){2}{\circle*{1}}
\put(6,12){\makebox(0,0){$\tilde{b}$}}
\put(24,12){\makebox(0,0){$\tilde{c}_e$}}
\put(6,-6){\makebox(0,0){$E^f$}}
\put(24,-6){\makebox(0,0){$4$}}
\end{picture}
=\quad -1 ,
\end{displaymath}


\begin{displaymath}
\thinlines
\unitlength 0.5mm
\begin{picture}(30,20)
\multiput(11,-2)(0,10){2}{\line(1,0){8}}
\multiput(10,7)(10,0){2}{\line(0,-1){8}}
\multiput(10,-2)(10,0){2}{\circle*{1}}
\multiput(10,8)(10,0){2}{\circle*{1}}
\put(6,12){\makebox(0,0){$d$}}
\put(24,-6){\makebox(0,0){$4$}}
\end{picture}
=\quad
 \bordermatrix{
&  c_e & {e_2}_e & {e_3}_e \cr
 {E^{\tilde{d}}} & \frac{1}{\gb}\sqrt{\frac{\gb-2}{\gb+1}} & 
\sqrt{\frac{\gbf+4}{\gb(\gb+1)}} & 0\cr
 {E^f} &-\frac{\sqrt{\gb+2}}{\gb-1} & 
\frac{1}{\gb-1}\sqrt{\frac{\gbf-4}{\gb(\gbf+4)}}& 
-\gb\sqrt{\frac{2}{(\gb-1)(\gbf+4)}}\cr
 {E^d} &\frac{\gb}{\sqrt{(\gb+1)(2\gb-1)}} & 
-\sqrt{\frac{\gb(\gb-2)}{(\gb+1)(2\gb-1)(\gbf+4)}} & 
-\frac{4}{\gb-1}\sqrt{\frac{2\gb-1}{2(3\gb-4)}} \cr
},
\end{displaymath}


\begin{displaymath}
\thinlines
\unitlength 0.5mm
\begin{picture}(30,20)
\multiput(11,-2)(0,10){2}{\line(1,0){8}}
\multiput(10,7)(10,0){2}{\line(0,-1){8}}
\multiput(10,-2)(10,0){2}{\circle*{1}}
\multiput(10,8)(10,0){2}{\circle*{1}}
\put(6,12){\makebox(0,0){$\tilde{d}$}}
\put(24,-6){\makebox(0,0){4}}
\end{picture}
=\quad
\bordermatrix{
& {\tilde{c}}_e &{e_2}_e &  {e_3}_e   \cr
{E^{\tilde{d}}} & \sqrt{\frac{\gb-2}{\gb+1}} & 
-\sqrt{\frac{\gb(\gb-2)}{(\gb+1)(2\gb-1)(\gbf+4)}} & 
-\frac{4}{\gb-1}\sqrt{\frac{2\gb-1}{2(3\gb-4)}}  \cr
{E^f} &  -\sqrt{\frac{5-\gb}{\gb+1}}& 
\sqrt{\frac{2\gb}{2\gb-1}}\frac{\beta^3}{\sqrt{(\gb+1)(\gbf+4)}} & 
-\sqrt{\frac{2\gb}{2\gb-1}}\sqrt{\frac{4\gb}{(\gb+1)(\gbf+4)}}                  
\cr
{E^d} &  \sqrt{\frac{\gb-2}{\gb+1}} & 
\sqrt{\frac{\gb(\gb+2)}{(\gb+1)(\gbf+4)}} & \frac{4}{\sqrt{2(\gbf+4)}}   \cr
},
\end{displaymath}

\begin{displaymath}
\thinlines
\unitlength 0.5mm
\begin{picture}(30,20)
\multiput(11,-2)(0,10){2}{\line(1,0){8}}
\multiput(10,7)(10,0){2}{\line(0,-1){8}}
\multiput(10,-2)(10,0){2}{\circle*{1}}
\multiput(10,8)(10,0){2}{\circle*{1}}
\put(6,12){\makebox(0,0){$h$}}
\put(24,12){\makebox(0,0){$g_g$}}
\put(6,-6){\makebox(0,0){$G^f$}}
\put(24,-6){\makebox(0,0){$6$}}
\end{picture}
=\quad 1 ,
\end{displaymath}

\begin{displaymath}
\thinlines
\unitlength 0.5mm
\begin{picture}(30,20)
\multiput(11,-2)(0,10){2}{\line(1,0){8}}
\multiput(10,7)(10,0){2}{\line(0,-1){8}}
\multiput(10,-2)(10,0){2}{\circle*{1}}
\multiput(10,8)(10,0){2}{\circle*{1}}
\put(6,12){\makebox(0,0){$h$}}
\put(24,12){\makebox(0,0){$g_g$}}
\put(6,-6){\makebox(0,0){$G^h$}}
\put(24,-6){\makebox(0,0){$6$}}
\end{picture}
=\quad 1 ,
\end{displaymath}


\begin{displaymath}
\thinlines
\unitlength 0.5mm
\begin{picture}(30,20)
\multiput(11,-2)(0,10){2}{\line(1,0){8}}
\multiput(10,7)(10,0){2}{\line(0,-1){8}}
\multiput(10,-2)(10,0){2}{\circle*{1}}
\multiput(10,8)(10,0){2}{\circle*{1}}
\put(6,12){\makebox(0,0){$f$}}
\put(24,12){\makebox(0,0){$g_a$}}
\put(6,-6){\makebox(0,0){$a^b$}}
\put(24,-6){\makebox(0,0){1}}
\end{picture}
= \quad *, 
\end{displaymath}

\begin{displaymath}
\thinlines
\unitlength 0.5mm
\begin{picture}(30,20)
\multiput(11,-2)(0,10){2}{\line(1,0){8}}
\multiput(10,7)(10,0){2}{\line(0,-1){8}}
\multiput(10,-2)(10,0){2}{\circle*{1}}
\multiput(10,8)(10,0){2}{\circle*{1}}
\put(6,12){\makebox(0,0){$f$}}
\put(24,-6){\makebox(0,0){$6$}}
\end{picture}
=\quad
\bordermatrix{
&  e_g & g_g  \cr
 G^f & &  \cr
 \tilde{A}^{\tilde{h}} &  &   \cr
 },
\end{displaymath}
(We do not use the values of these entries.)

\begin{displaymath}
\thinlines
\unitlength 0.5mm
\begin{picture}(30,20)
\multiput(11,-2)(0,10){2}{\line(1,0){8}}
\multiput(10,7)(10,0){2}{\line(0,-1){8}}
\multiput(10,-2)(10,0){2}{\circle*{1}}
\multiput(10,8)(10,0){2}{\circle*{1}}
\put(6,12){\makebox(0,0){$h$}}
\put(24,-6){\makebox(0,0){2}}
\end{picture}
=\quad
\bordermatrix{
& g_a & g_c \cr
 {A^b} & -\frac{1}{\gb-1}  & 
\frac{\sqrt{\gb(\gb-2)}}{{\gb-1}}\cr
 {C^b} & \frac{\sqrt{\gb(\gb-2)}}{{\gb-1}} & \frac{1}{\gb-1} 
\cr
},
\end{displaymath}


\begin{displaymath}
\thinlines
\unitlength 0.5mm
\begin{picture}(30,20)
\multiput(11,-2)(0,10){2}{\line(1,0){8}}
\multiput(10,7)(10,0){2}{\line(0,-1){8}}
\multiput(10,-2)(10,0){2}{\circle*{1}}
\multiput(10,8)(10,0){2}{\circle*{1}}
\put(6,12){\makebox(0,0){$\tilde{h}$}}
\put(24,-6){\makebox(0,0){2}}
\end{picture}
=\quad
\bordermatrix{
&  g_a & g_c \cr
 {A^*} & 1 & 0\cr
 {C^d} & 0 & 1\cr},
\end{displaymath}

\begin{displaymath}
\thinlines
\unitlength 0.5mm
\begin{picture}(30,20)
\multiput(11,-2)(0,10){2}{\line(1,0){8}}
\multiput(10,7)(10,0){2}{\line(0,-1){8}}
\multiput(10,-2)(10,0){2}{\circle*{1}}
\multiput(10,8)(10,0){2}{\circle*{1}}
\put(6,12){\makebox(0,0){$h$}}
\put(24,12){\makebox(0,0){$g_a$}}
\put(6,-6){\makebox(0,0){$A^b$}}
\put(24,-6){\makebox(0,0){$1$}}
\end{picture}
=\quad 1 ,
\end{displaymath}

\begin{displaymath}
\thinlines
\unitlength 0.5mm
\begin{picture}(30,20)
\multiput(11,-2)(0,10){2}{\line(1,0){8}}
\multiput(10,7)(10,0){2}{\line(0,-1){8}}
\multiput(10,-2)(10,0){2}{\circle*{1}}
\multiput(10,8)(10,0){2}{\circle*{1}}
\put(6,12){\makebox(0,0){$\tilde{h}$}}
\put(24,12){\makebox(0,0){$g_a$}}
\put(6,-6){\makebox(0,0){$A^b$}}
\put(24,-6){\makebox(0,0){$1$}}
\end{picture}
=\quad 1 ,
\end{displaymath}

\begin{displaymath}
\thinlines
\unitlength 0.5mm
\begin{picture}(30,20)
\multiput(11,-2)(0,10){2}{\line(1,0){8}}
\multiput(10,7)(10,0){2}{\line(0,-1){8}}
\multiput(10,-2)(10,0){2}{\circle*{1}}
\multiput(10,8)(10,0){2}{\circle*{1}}
\put(6,12){\makebox(0,0){$h$}}
\put(24,12){\makebox(0,0){$g_e$}}
\put(6,-6){\makebox(0,0){$E^f$}}
\put(24,-6){\makebox(0,0){$4$}}
\end{picture}
=\quad -1 ,
\end{displaymath}

\begin{displaymath}
\thinlines
\unitlength 0.5mm
\begin{picture}(30,20)
\multiput(11,-2)(0,10){2}{\line(1,0){8}}
\multiput(10,7)(10,0){2}{\line(0,-1){8}}
\multiput(10,-2)(10,0){2}{\circle*{1}}
\multiput(10,8)(10,0){2}{\circle*{1}}
\put(6,12){\makebox(0,0){$\tilde{h}$}}
\put(24,12){\makebox(0,0){$g_e$}}
\put(6,-6){\makebox(0,0){$E^d$}}
\put(24,-6){\makebox(0,0){$4$}}
\end{picture}
=\quad 1 ,
\end{displaymath}


\begin{displaymath}
\thinlines
\unitlength 0.5mm
\begin{picture}(30,20)
\multiput(11,-2)(0,10){2}{\line(1,0){8}}
\multiput(10,7)(10,0){2}{\line(0,-1){8}}
\multiput(10,-2)(10,0){2}{\circle*{1}}
\multiput(10,8)(10,0){2}{\circle*{1}}
\put(6,12){\makebox(0,0){$d$}}
\put(24,12){\makebox(0,0){$e$}}
\put(6,-6){\makebox(0,0){$C$}}
\put(24,-6){\makebox(0,0){$2$}}
\end{picture}
=\quad1,
\end{displaymath}


\begin{displaymath}
\thinlines
\unitlength 0.5mm
\begin{picture}(30,20)
\multiput(11,-2)(0,10){2}{\line(1,0){8}}
\multiput(10,7)(10,0){2}{\line(0,-1){8}}
\multiput(10,-2)(10,0){2}{\circle*{1}}
\multiput(10,8)(10,0){2}{\circle*{1}}
\put(6,12){\makebox(0,0){$\tilde{d}$}}
\put(24,12){\makebox(0,0){$e_c$}}
\put(6,-6){\makebox(0,0){$C^d$}}
\put(24,-6){\makebox(0,0){$2$}}
\end{picture}
=\quad 1.
\end{displaymath}
\pha{x} \par 
Here we will display three matrices of the connection $((\a \ab-{\bf 
1})\sigma \a)^{\sim} = \a \sigma(\a \ab-{\bf 1})$ for easiness of later 
procedure. These matrices are computed from the entries of $(\a \ab-{\bf 
1})\sigma \a $ and the Perron-Frobenius weights of the horizontal 
graphs by renormalization rule as in section 4.


\begin{eqnarray*}
 (1) && \phantom{XXXX}
\thinlines
\unitlength 0.5mm
\begin{picture}(30,20)
\multiput(11,-2)(0,10){2}{\line(1,0){8}}
\multiput(10,7)(10,0){2}{\line(0,-1){8}}
\multiput(10,-2)(10,0){2}{\circle*{1}}
\multiput(10,8)(10,0){2}{\circle*{1}}
\put(6,12){\makebox(0,0){$G$}}
\put(24,-6){\makebox(0,0){$e$}}
\end{picture}         \\
&=&\quad
\bordermatrix{
& 5_{e_{2}} & 5_{e_{3}} & 5_{\td{c}} & 5_{g} & 6_{g}  \cr
d^{f} & 
-\frac{1}{\gb(\gb-1)}\sqrt{\frac{(2\gb-1)(\gbf-4)}{2(\gbf+4)}} & 
\sqrt{\frac{\gb(2\gb-1)}{(\gb-1)(\gbf+4)}} &  0 &  
-\frac{\sqrt{2\gb-1}}{2(\gb-1)} &  \sqrt{\frac{2\gb-1}{\gbf-1}} \cr
d^{\td{b}} & 0 & 0 & 1 & 0 & 0 \cr
\td{d}_{f} & -\frac{\beta^3}{\sqrt{(\gb+1)(\gbf+4)}} & 
\sqrt{\frac{4\gb}{(\gb+1)(\gbf+4)}} &  0 & \frac{1}{\sqrt{\gbf-1}} & 
-\sqrt{\frac{\gb-2}{\gbf-1}} \cr
\td{d}^{h} & 0 & 0 & 0 & \sqrt{\frac{\gb-2}{\gb-1}} & 
\frac{1}{\sqrt{\gb-1}} \cr
f_{1}^{f} & \sqrt{\frac{2\gb}{(\gb+1)(\gbf+4)}} & 
\sqrt{\frac{2(\gb+1)}{\gbf+4}} & \frac{1}{2\sqrt{\gb-1}} & 
\frac{-1}{\sqrt{\gb+1}}  \cr
}, 
\end{eqnarray*}


\begin{displaymath}
\thinlines
\unitlength 0.5mm
\begin{picture}(30,20)
\multiput(11,-2)(0,10){2}{\line(1,0){8}}
\multiput(10,7)(10,0){2}{\line(0,-1){8}}
\multiput(10,-2)(10,0){2}{\circle*{1}}
\multiput(10,8)(10,0){2}{\circle*{1}}
\put(6,12){\makebox(0,0){$A$}}
\put(24,-6){\makebox(0,0){$g$}}
\end{picture}
=\quad
\bordermatrix{
& *_{a} & 2_{a} & 2_{c} \cr
f^{b} & \sqrt{\frac{\gb-2}{\gb}} & -\sqrt{\frac{\gb-2}{(\gb-1)\gb}} &  
-\frac{1}{\sqrt{\gb-1}}   \cr
h^{b} & \frac{1}{\gbs} & -\frac{1}{\sqrt{\gb(\gb-1)}} & 
\sqrt{\frac{\gb-2}{\gb-1}} \cr
{\tilde{h}}^{*} & \frac{1}{\gbs} & \sqrt{\frac{\gb-1}{\gb}} & 0 \cr
},
\end{displaymath}


\begin{eqnarray*}
&&\thinlines
\unitlength 0.5mm
\begin{picture}(30,20)
\multiput(11,-2)(0,10){2}{\line(1,0){8}}
\multiput(10,7)(10,0){2}{\line(0,-1){8}}
\multiput(10,-2)(10,0){2}{\circle*{1}}
\multiput(10,8)(10,0){2}{\circle*{1}}
\put(6,12){\makebox(0,0){$C$}}
\put(24,-6){\makebox(0,0){$e$}}
\end{picture} \\
&=&\quad
\bordermatrix{
& 3_{e_2} & 3_{e_3} & 3_{c} & 2_{c} \cr
d_{d} &  -\sqrt{\frac{\gb(\gb-2)}{(\gb+1)(2\gb-1)(\gbf+4)}} & 
\frac{-2(\gb-2)}{\sqrt{2(\gbf+4)}} &  \frac{-1}{2\beta \sqrt{\gb-1}} & 
\frac{\beta}{2\sqrt{\gb-2}} \cr
f_{d} & 
-\sqrt{\frac{2\gb}{2\gb-1}}\frac{\beta^3}{\sqrt{(\gb+1)(\gbf+4)}} & 
\sqrt{\frac{2\gb}{2\gb-1}}\sqrt{\frac{4\gb}{(\gb+1)(\gbf+4)}} &    
\frac{-1}{2\beta^3} &  \frac{1}{\sqrt{2}(\gb-2)}  \cr
\tilde{d_{1}}^{d} & \sqrt{\frac{\gb(\gb+2)}{(\gb+1)(\gbf+4)}} & 
\frac{4}{\sqrt{2(\gbf+4)}} & \frac{-1}{2\beta \sqrt{\gb-1}} & 
\frac{\beta}{2\sqrt{\gb-2}} \cr
f^{b} & 0 & 0 & \frac{\gb-1}{\sqrt{2}(\gb-2)} & \frac{1}{\gbs(\gb-2)}  \cr
},
\end{eqnarray*}
Now we will prove that $(\alpha \tilde{\alpha}-\bf{1})\sigma
\alpha$ and $ \sigma(\alpha \tilde{\alpha}-\bf{1})\sigma \alpha$ are
equivalent up to {\it vertical} gauge choice. What we should do is to 
construct gauge transformation matrices for each vertical edges.
We write $u{p\cho{q}}_{m,l}$ 
for the $m \times m$ unitary gauge 
transformation matrix coming from the edges $p$-$q$ of multiplicity $m$ 
in the left vertical graphs ${\cal H}_{0}$ and $u{r\cho{s}}_{n,r}$ for 
the $n \times n$ unitary gauge transformation matrix coming from the 
edges $r$-$s$ of multiplicity $n$ in the right vertical graph ${\cal 
H}_{1}$. Let
\begin{displaymath}
\left( \thinlines
\unitlength 0.5mm
\begin{picture}(30,20)
\multiput(11,-2)(0,10){2}{\line(1,0){8}}
\multiput(10,7)(10,0){2}{\line(0,-1){8}}
\multiput(10,-2)(10,0){2}{\circle*{1}}
\multiput(10,8)(10,0){2}{\circle*{1}}
\put(6,12){\makebox(0,0){$x$}}
\put(24,12){\makebox(0,0){$z$}}
\put(6,-6){\makebox(0,0){$y$}}
\put(24,-6){\makebox(0,0){$w$}}
\put(6,3){\makebox(0,0){$\xi$}}
\put(24,3){\makebox(0,0){$\eta$}}
\end{picture} \right)_{\xi,\eta}
\end{displaymath}
to be a $n \times m$ matrix of the connection $(\alpha 
\tilde{\alpha}-\bf{1})\sigma \alpha$, where $n$ and $m$ is the 
multiplicities of the edges $x$-$y$ and $z$-$w$ respectively, and 
\begin{displaymath}
\left( \thinlines
\unitlength 0.5mm
\begin{picture}(30,20)
\multiput(11,-2)(0,10){2}{\line(1,0){8}}
\multiput(10,7)(10,0){2}{\line(0,-1){8}}
\multiput(10,-2)(10,0){2}{\circle*{1}}
\multiput(10,8)(10,0){2}{\circle*{1}}
\put(12,1){$\sim$}
\put(6,12){\makebox(0,0){$x$}}
\put(24,12){\makebox(0,0){$z$}}
\put(6,-6){\makebox(0,0){$y$}}
\put(24,-6){\makebox(0,0){$w$}}
\put(6,3){\makebox(0,0){$\xi$}}
\put(24,3){\makebox(0,0){$\eta$}}
\end{picture} \right)_{\xi,\eta}
\end{displaymath}
to be a $n \times m$ matrix of the connection $\sigma(\alpha 
\tilde{\alpha}-\bf{1})\sigma \alpha$. Then, the gauge matrices which 
we are going to construct should satisfy the equality
$$
\label{eqn}
\left( \thinlines
\unitlength 0.5mm
\begin{picture}(30,20)
\multiput(11,-2)(0,10){2}{\line(1,0){8}}
\multiput(10,7)(10,0){2}{\line(0,-1){8}}
\multiput(10,-2)(10,0){2}{\circle*{1}}
\multiput(10,8)(10,0){2}{\circle*{1}}
\put(12,1){$\sim$}
\put(6,12){\makebox(0,0){$x$}}
\put(24,12){\makebox(0,0){$x$}}
\put(6,-6){\makebox(0,0){$y$}}
\put(24,-6){\makebox(0,0){$w$}}
\put(6,3){\makebox(0,0){$\xi$}}
\put(24,3){\makebox(0,0){$\eta$}}
\end{picture} \right)_{\xi,\eta}
=u{x\cho{y}}_{n,l}
\left( \thinlines
\unitlength 0.5mm
\begin{picture}(30,20)
\multiput(11,-2)(0,10){2}{\line(1,0){8}}
\multiput(10,7)(10,0){2}{\line(0,-1){8}}
\multiput(10,-2)(10,0){2}{\circle*{1}}
\multiput(10,8)(10,0){2}{\circle*{1}}
\put(6,12){\makebox(0,0){$x$}}
\put(24,12){\makebox(0,0){$z$}}
\put(6,-6){\makebox(0,0){$y$}}
\put(24,-6){\makebox(0,0){$w$}}
\put(6,3){\makebox(0,0){$\xi$}}
\put(24,3){\makebox(0,0){$\eta$}}
\end{picture} \right)_{\xi,\eta}
u{z\cho{w}}_{m,r}
$$
for all pair of vertical edges $(xy, zw)$. Notice that multiplying the connection
 $\sigma$ from the left means simply  
 changing the upper vertices of the connection as $ x \leftrightarrow 
\tilde{x}$, and then the above equality is equivalent to 
$$
\left( \thinlines
\unitlength 0.5mm
\begin{picture}(30,20)
\multiput(11,-2)(0,10){2}{\line(1,0){8}}
\multiput(10,7)(10,0){2}{\line(0,-1){8}}
\multiput(10,-2)(10,0){2}{\circle*{1}}
\multiput(10,8)(10,0){2}{\circle*{1}}
\put(6,12){\makebox(0,0){$\td{x}$}}
\put(24,12){\makebox(0,0){$\td{z}$}}
\put(6,-6){\makebox(0,0){$y$}}
\put(24,-6){\makebox(0,0){$w$}}
\put(6,3){\makebox(0,0){$\xi$}}
\put(24,3){\makebox(0,0){$\eta$}}
\end{picture} \right)_{\xi,\eta}
=u{x\cho{y}}_{n,l}
\left( \thinlines
\unitlength 0.5mm
\begin{picture}(30,20)
\multiput(11,-2)(0,10){2}{\line(1,0){8}}
\multiput(10,7)(10,0){2}{\line(0,-1){8}}
\multiput(10,-2)(10,0){2}{\circle*{1}}
\multiput(10,8)(10,0){2}{\circle*{1}}
\put(6,12){\makebox(0,0){$x$}}
\put(24,12){\makebox(0,0){$z$}}
\put(6,-6){\makebox(0,0){$y$}}
\put(24,-6){\makebox(0,0){$w$}}
\put(6,3){\makebox(0,0){$\xi$}}
\put(24,3){\makebox(0,0){$\eta$}}
\end{picture} \right)_{\xi,\eta}
u{y\cho{w}}_{m,r}.
$$



Note that the vertices $e$, $f$ and $g$ are fixed by taking $\sim$.
We easily know that 

\begin{displaymath}
 u{z\cho{w}}_{n,l} 
= u{\tilde{z}\cho{w}}_{n,l}^t
= u{z\cho{w}}_{n,r},
\end{displaymath}
\phantom{x} \\
Note that the multiplicity $n$ of the edges $z$-$w$ is equal to 
that of the edges $\tilde{z}$-$w$.
 Now we begin to construct a candidate for the list of gauge 
 transformation matrices. 
First, for the connections
$M(*/6)$ and $M(\tilde{*}/6)$, we fix the gauges for the simple edges 
as 

$$ 
u{*\cho{G}}_{1,l}=u{a\cho{6}}_{1,r}=u{\td{*}\cho{G}}_{1,l}=
u{\tilde{a}\cho{6}}_{1,r}=(1)_{1,1}, 
$$ 
here the matrices are all $1 \times 1$.  
From the next matrices, we always fix the gauges for simple edges to 
$1 \times 1$ matrices $(1)_{1,1}$,
unless otherwise specified.  We denote it simply by 1. \par
Next we fix the gauges for the connections $M(*/5)$,
$M(\tilde{*}/5)$ and $M(b/6)$, $M(\tilde{b}/6)$. We put
$u{b\cho{G}}_{1,l}=u{\tilde{b}\cho{G}}_{1,l}=-1$.\par 
For $M(b/5)$ and $M(\tilde{b}/5)$, we fix gauges as follows:

\begin{displaymath}
\thinlines
\unitlength 0.5mm
\begin{picture}(30,20)
\multiput(11,-2)(0,10){2}{\line(1,0){8}}
\multiput(10,7)(10,0){2}{\line(0,-1){8}}
\multiput(10,-2)(10,0){2}{\circle*{1}}
\multiput(10,8)(10,0){2}{\circle*{1}}
\put(6,12){\makebox(0,0){$\tilde{b}$}}
\put(24,-6){\makebox(0,0){$5$}}
\end{picture}
=
\left( \begin{array}{ccc}
1 & &  \\
   & -1 & \\
   &  & 1  
\end{array} \right)
\thinlines
\unitlength 0.5mm
\begin{picture}(30,20)
\multiput(11,-2)(0,10){2}{\line(1,0){8}}
\multiput(10,7)(10,0){2}{\line(0,-1){8}}
\multiput(10,-2)(10,0){2}{\circle*{1}}
\multiput(10,8)(10,0){2}{\circle*{1}}
\put(6,12){\makebox(0,0){$b$}}
\put(24,-6){\makebox(0,0){$5$}}
\end{picture}
\left( \begin{array}{cc}
1 &   \\
   &   u{\tilde{c}\cho{5}}_{2,r}
 \end{array}  \right),
\end{displaymath}

\begin{displaymath}
\left( \begin{array}{cc}
1 &   \\
   &   u{\tilde{c}\cho{5}}_{2,r}
 \end{array}  \right)
=
\left( \matrix{
1 & 0 & 0  \cr
0 & \frac{1}{\gb} & \frac{\sqrt{\gbf-1}}{\gb}  \cr
0 & \frac{\sqrt{\gbf-1}}{\gb} & \frac{-1}{\gb}  \cr
} \right).
\end{displaymath}

In the same way we get
$$u{c\cho{5}}_{2,r}
= \left( \matrix{
\frac{1}{\gb} & \frac{\sqrt{\gbf-1}}{\gb}  \cr
\frac{\sqrt{\gbf-1}}{\gb} & \frac{-1}{\gb}  \cr
} \right).
$$
For $M(d/6)$ and $M(\tilde{d}/6)$,

\begin{displaymath}
\thinlines
\unitlength 0.5mm
\begin{picture}(30,20)
\multiput(11,-2)(0,10){2}{\line(1,0){8}}
\multiput(10,7)(10,0){2}{\line(0,-1){8}}
\multiput(10,-2)(10,0){2}{\circle*{1}}
\multiput(10,8)(10,0){2}{\circle*{1}}
\put(6,12){\makebox(0,0){$\tilde{d}$}}
\put(24,-6){\makebox(0,0){$6$}}
\end{picture}
=
u{d\cho{G}}_{2,l}
\thinlines
\unitlength 0.5mm
\begin{picture}(30,20)
\multiput(11,-2)(0,10){2}{\line(1,0){8}}
\multiput(10,7)(10,0){2}{\line(0,-1){8}}
\multiput(10,-2)(10,0){2}{\circle*{1}}
\multiput(10,8)(10,0){2}{\circle*{1}}
\put(6,12){\makebox(0,0){$d$}}
\put(24,-6){\makebox(0,0){$6$}}
\end{picture},
\end{displaymath}

\begin{displaymath}
u{d\cho{G}}_{2,l}
=
\left(\matrix{
   -\sqrt{\frac{\gb-2}{2\gb-1}}  & \sqrt{\frac{\gb+1}{2\gb-1}} \cr
   \sqrt{\frac{\gb+1}{2\gb-1}} & \sqrt{\frac{\gb-2}{2\gb-1}}  \cr
}  \right)
=
u{\tilde{d}\cho{G}}_{2,l},
\end{displaymath}
by symmetricity of this matrix.
To check $M(d/5)$ and  $M(\tilde{d}/5)$, we use $M(G/e)$. See the 
matrix 
(1).
Note that the multiplication by $\sigma$ from the left to $(\a \ab 
-{\bf 1})\sigma \a$ corresponds to that by $\sigma$ from the right on 
$\ab \sigma (\a \ab-{\bf 1})$, which causes the permutation of the entries of
connection matrix $M(G/e)$ as follows.

\begin{center}
\thinlines
\unitlength 0.5mm
\begin{picture}(30,20)
\multiput(11,-2)(0,10){2}{\line(1,0){8}}
\multiput(10,7)(10,0){2}{\line(0,-1){8}}
\multiput(10,-2)(10,0){2}{\circle*{1}}
\multiput(10,8)(10,0){2}{\circle*{1}}
\put(6,12){\makebox(0,0){$G$}}
\put(24,12){\makebox(0,0){$e$}}
\put(24,-6){\makebox(0,0){$e$}}
\put(6,-6){\makebox(0,0){$d$}}
\end{picture}
{$\longrightarrow$}
\thinlines
\unitlength 0.5mm
\begin{picture}(30,20)
\multiput(11,-2)(0,10){2}{\line(1,0){8}}
\multiput(10,7)(10,0){2}{\line(0,-1){8}}
\multiput(10,-2)(10,0){2}{\circle*{1}}
\multiput(10,8)(10,0){2}{\circle*{1}}
\put(6,12){\makebox(0,0){$G$}}
\put(24,12){\makebox(0,0){$e$}}
\put(24,-6){\makebox(0,0){$e$}}
\put(6,-6){\makebox(0,0){$\tilde{d}$}}
\end{picture}
\end{center}

We denote the connection matrix made from $M(G/e)$ by multiplying 
$\sigma$ 
by $M(G/e)_{\sim}$.  Since the vertex $e$ is fixed by multiplying
$\sigma$, we should fix the gauge matrix so that $M(G/e)$ and 
$M(G/e)_{\sim}$
are transferred each other.  By the effect of multiplying $\sigma$, we
see that $M(G/e)_{\sim}$ is made from $M(G/e)$ by exchanging $d^f$
(resp. $d^{\tilde{b}}$)-row and $\tilde{d}^f$ (resp. 
$\tilde{d}^h$)-row,
i.e., we have the following relation:

$$
\thinlines
\unitlength 0.5mm
\begin{picture}(30,20)
\multiput(11,-2)(0,10){2}{\line(1,0){8}}
\multiput(10,7)(10,0){2}{\line(0,-1){8}}
\multiput(10,-2)(10,0){2}{\circle*{1}}
\multiput(10,8)(10,0){2}{\circle*{1}}
\put(6,12){\makebox(0,0){$G$}}
\put(24,-6){\makebox(0,0){$e$}}
\put(6,-6){\makebox(0,0){$\sim$}}
\end{picture}
=
\left( \matrix{
0 & 0 & 1 & 0 & 0 \cr
0 & 0 & 0 & 1 & 0 \cr
1 & 0 & 0 & 0 & 0 \cr
0 & 1 & 0 & 0 & 0 \cr
0 & 0 & 0 & 0 & 1
 }\right)
\thinlines
\unitlength 0.5mm
\begin{picture}(30,20)
\multiput(11,-2)(0,10){2}{\line(1,0){8}}
\multiput(10,7)(10,0){2}{\line(0,-1){8}}
\multiput(10,-2)(10,0){2}{\circle*{1}}
\multiput(10,8)(10,0){2}{\circle*{1}}
\put(6,12){\makebox(0,0){$G$}}
\put(24,-6){\makebox(0,0){$e$}}
\end{picture},
$$
here $\sim$ at the lower left corner in the left hand side square 
means 
the changing of the labels and replacing of columns according to the
labels.  Now we get the gauge as usual. Note that gauge matrices for 
upside down edges are the same to those of normal position.

$$
\thinlines
\unitlength 0.5mm
\begin{picture}(30,20)
\multiput(11,-2)(0,10){2}{\line(1,0){8}}
\multiput(10,7)(10,0){2}{\line(0,-1){8}}
\multiput(10,-2)(10,0){2}{\circle*{1}}
\multiput(10,8)(10,0){2}{\circle*{1}}
\put(6,12){\makebox(0,0){$G$}}
\put(24,-6){\makebox(0,0){$e$}}
\put(6,-6){\makebox(0,0){$\sim$}}
\end{picture}
=
\left( \begin{array}{ccc}
  u{\tilde{d}\cho{G}}_{2,l}  & &  \\
  &  u{d\cho{G}}_{2,l}  &  \\
& & 1    
 \end{array}  \right)   
\thinlines
\unitlength 0.5mm
\begin{picture}(30,20)
\multiput(11,-2)(0,10){2}{\line(1,0){8}}
\multiput(10,7)(10,0){2}{\line(0,-1){8}}
\multiput(10,-2)(10,0){2}{\circle*{1}}
\multiput(10,8)(10,0){2}{\circle*{1}}
\put(6,12){\makebox(0,0){$G$}}
\put(24,-6){\makebox(0,0){$e$}}
\end{picture}
\left( \begin{array}{cc}
  u{e\cho{5}}_{4,r} &  \\
  &1
 \end{array}  \right) $$
\begin{flushleft} 
$$ {} =
\left( \matrix{
\frac{(\gb-2)^3}{2\gbf \sqrt{2(\gbf+4)}}& \frac{-(\gb-2)}{\sqrt{2(\gbf+4)}} &
\sqrt{\frac{\gb+1}{2\gb-1}} &\frac{\sqrt{\gb-2}}{2(\gb-1)} &
\frac{-(\gb-2)}{2(\gb-1)} \cr
\frac{-\sqrt{2}(\gb-2)}{\gb\sqrt{(\gb-1)(\gbf+4)}} & 
\frac{\sqrt{2(\gb-1)}}{\sqrt{\gbf+4}} & \sqrt{\frac{\gb-2}{2\gb-1}} & 
\frac{-1}{\sqrt{2\gb}} & \frac{1}{\sqrt{\gb-1}} \cr
\frac{\gb(\gb-2)}{(\gb-1)\sqrt{2(\gbf+4)}}& 
\frac{-\sqrt{2}(\gb-2)}{(\gb-1)\sqrt{\gbf+4}} & 0 &
 \frac{\gb-2}{\sqrt{\gbf-1}} & 
 \frac{\sqrt{2\gb-1}}{\sqrt{\gbf-1}} \cr
  \frac{\sqrt{\gb(\gb-2)}}{\sqrt{\gbf+4}}& \frac{-2\sqrt{\gb-2}}
{\sqrt{\gb(\gbf+4)}} &  0 & \frac{\sqrt{2}}{\gbs} & 0 \cr
\frac{\sqrt{2}\gb}{\sqrt{(\gb+1)(\gbf+4)}}  & 
\frac{\sqrt{2(\gb+1)}}{\sqrt{\gbf+4}} & 0 & \frac{1}{2\sqrt{\gb-1}}  & 
\frac{-1}{\sqrt{\gb+1}} 
} \right)
 \left( \begin{array}{cc}
 u{e\cho{5}}_{4,r}  &  \\
  &1
  \end{array}  \right),$$
\end{flushleft}

\begin{eqnarray*}
&& \left( \begin{array}{cc}
  u{e\cho{5}}_{4,r}   &  \\
  &1
  \end{array}  \right)   \\
&=&
\left( \begin{array}{@{}c@{}|@{}c@{}|@{}c@{}|@{}c@{}|c@{}}
\frac{2(\gb-2)}{\gbf+4} & \frac{2(\gb+4)}{\gbf+4} & 
\frac{-\sqrt{\gb(\gb-2)}}{\sqrt{\gbf+4}} & 
\frac{-2}{\sqrt{\gb(\gb-1)(\gbf+4)}} & 0 \\ 
\frac{2(\gb+4)}{\gbf+4} & \frac{\gb}{\gbf+4} & \frac{2\gbs}
{\sqrt{(2\gb-1)(\gbf+4)}} & \frac{(\gb-1)\sqrt{\gb-1}}
{\sqrt{\gb(\gbf+4)}}  & 0 \\
\frac{-\sqrt{\gb(\gb-2)}}{\sqrt{\gbf+4}} & \frac{2\gbs}
{\sqrt{(2\gb-1)(\gbf+4)}} & 0 & \frac{\sqrt{2}}{\gbs} & 0 \\
\frac{-2}{\sqrt{\gb(\gb-1)(\gbf+4)}} & \frac{(\gb-1)\sqrt{\gb-1}}
{\sqrt{\gb(\gbf+4)}} & \frac{\sqrt{2}}{\gbs}  & \frac{-1}{\gb-2} &
0 \\
0 & 0 & 0 & 0 & 1 \end{array}  \right).
\end{eqnarray*} 

It is too hard to obtain above gauge matrix by calculating all the
elements by the multiplication of matrices.  Note that,

\begin{eqnarray*}
\thinlines
\unitlength 0.5mm
\begin{picture}(30,20)
\multiput(11,-2)(0,10){2}{\line(1,0){8}}
\multiput(10,7)(10,0){2}{\line(0,-1){8}}
\multiput(10,-2)(10,0){2}{\circle*{1}}
\multiput(10,8)(10,0){2}{\circle*{1}}
\put(6,12){\makebox(0,0){$G$}}
\put(24,-6){\makebox(0,0){$e$}}
\put(6,-6){\makebox(0,0){$\sim$}}
\end{picture}
&=&
\left( \matrix{
0 & 0 & 1 & 0 & 0 \cr
0 & 0 & 0 & 1 & 0 \cr
1 & 0 & 0 & 0 & 0 \cr
0 & 1 & 0 & 0 & 0 \cr
0 & 0 & 0 & 0 & 1 \cr
 }\right) 
\thinlines
\unitlength 0.5mm
\begin{picture}(30,20)
\multiput(11,-2)(0,10){2}{\line(1,0){8}}
\multiput(10,7)(10,0){2}{\line(0,-1){8}}
\multiput(10,-2)(10,0){2}{\circle*{1}}
\multiput(10,8)(10,0){2}{\circle*{1}}
\put(6,12){\makebox(0,0){$G$}}
\put(24,-6){\makebox(0,0){$e$}}
\end{picture}
 \\
 &=&
\left( \begin{array}{ccc}
  u{\tilde{d}\cho{G}}_{2,l}  & &  \\
  &  u{d\cho{G}}_{2,l}  &  \\
& & 1    
 \end{array}  \right)
\thinlines
\unitlength 0.5mm
\begin{picture}(30,20)
\multiput(11,-2)(0,10){2}{\line(1,0){8}}
\multiput(10,7)(10,0){2}{\line(0,-1){8}}
\multiput(10,-2)(10,0){2}{\circle*{1}}
\multiput(10,8)(10,0){2}{\circle*{1}}
\put(6,12){\makebox(0,0){$G$}}
\put(24,-6){\makebox(0,0){$e$}}
\end{picture}
\left( \begin{array}{cc}
  u{e\cho{5}}_{4,r} &  \\
  &1
 \end{array}  \right),
\end{eqnarray*}

\begin{eqnarray*}
&&
\left( \matrix{
0 & 0 & 1 & 0 & 0 \cr
0 & 0 & 0 & 1 & 0 \cr
1 & 0 & 0 & 0 & 0 \cr
0 & 1 & 0 & 0 & 0 \cr
0 & 0 & 0 & 0 & 1 \cr
 } \right) 
\thinlines
\unitlength 0.5mm
\begin{picture}(30,20)
\multiput(11,-2)(0,10){2}{\line(1,0){8}}
\multiput(10,7)(10,0){2}{\line(0,-1){8}}
\multiput(10,-2)(10,0){2}{\circle*{1}}
\multiput(10,8)(10,0){2}{\circle*{1}}
\put(6,12){\makebox(0,0){$G$}}
\put(24,-6){\makebox(0,0){$e$}}
\end{picture}
 \\
&=&
\left( \matrix{
0 & 0 & 1 & 0 & 0 \cr
0 & 0 & 0 & 1 & 0 \cr
1 & 0 & 0 & 0 & 0 \cr
0 & 1 & 0 & 0 & 0 \cr
0 & 0 & 0 & 0 & 1 \cr
 }\right)  
 \left( \begin{array}{ccc}
  u{\tilde{d}\cho{G}}^t_{2,l}  & &  \\
  &  u{d\cho{G}}^t_{2,l}  &  \\
& & 1    
 \end{array}  \right)
 \left( \matrix{
0 & 0 & 1 & 0 & 0 \cr
0 & 0 & 0 & 1 & 0 \cr
1 & 0 & 0 & 0 & 0 \cr
0 & 1 & 0 & 0 & 0 \cr
0 & 0 & 0 & 0 & 1 \cr
 }\right) 
 \thinlines
\unitlength 0.5mm
\begin{picture}(30,20)
\multiput(11,-2)(0,10){2}{\line(1,0){8}}
\multiput(10,7)(10,0){2}{\line(0,-1){8}}
\multiput(10,-2)(10,0){2}{\circle*{1}}
\multiput(10,8)(10,0){2}{\circle*{1}}
\put(6,12){\makebox(0,0){$G$}}
\put(24,-6){\makebox(0,0){$e$}}
\end{picture}
\left( \begin{array}{cc}
  u{e\cho{5}}^t_{4,r} &  \\
  &1
 \end{array}  \right)  \\
 &=&
\left( \begin{array}{ccc}
  u{d\cho{G}}^{t}_{2,l}  & &  \\
  &  u{\tilde{d}\cho{G}}^{t}_{2,l}  &  \\
& & 1    
 \end{array}  \right)
\thinlines
\unitlength 0.5mm
\begin{picture}(30,20)
\multiput(11,-2)(0,10){2}{\line(1,0){8}}
\multiput(10,7)(10,0){2}{\line(0,-1){8}}
\multiput(10,-2)(10,0){2}{\circle*{1}}
\multiput(10,8)(10,0){2}{\circle*{1}}
\put(6,12){\makebox(0,0){$G$}}
\put(24,-6){\makebox(0,0){$e$}}
\end{picture}
\left( \begin{array}{cc}
  u{e\cho{5}}^t_{4,r} &  \\
  &1
 \end{array}  \right),
\end{eqnarray*}
and $u{d\cho{G}}_{2,l}=u{\tilde{d}\cho{G}}_{2,l}=u{d\cho{G}}_{2,l}^t$, 
by comparing the
above two equations we know that the matrix for
 the gauge $u{e\cho{5}}_{4,r}$ is symmetric.  Thus it is
enough first to check the (5,5)-entry is 1, and then to calculate 
the 
$(1,1)$, $(2,1)$, $(3,1)$, $(4,1)$, $(2,2)$, $(3,2)$, $(3,3)$, 
$(2,4)$, $(3,4)$ and $(4,4)$-entries. \par
We continue to fix gauge transformation matrices.

\begin{eqnarray*}
&& 
\thinlines
\unitlength 0.5mm
\begin{picture}(30,20)
\multiput(11,-2)(0,10){2}{\line(1,0){8}}
\multiput(10,7)(10,0){2}{\line(0,-1){8}}
\multiput(10,-2)(10,0){2}{\circle*{1}}
\multiput(10,8)(10,0){2}{\circle*{1}}
\put(6,12){\makebox(0,0){$\tilde{d}$}}
\put(24,-6){\makebox(0,0){$5$}}
\end{picture}
=
\left( \begin{array}{ccc}
  -1 & &  \\
&  u{d\cho{E}}_{3,l} &  \\
 &  & u{d\cho{G}}_{2,l}
  \end{array}  \right)
\thinlines
\unitlength 0.5mm
\begin{picture}(30,20)
\multiput(11,-2)(0,10){2}{\line(1,0){8}}
\multiput(10,7)(10,0){2}{\line(0,-1){8}}
\multiput(10,-2)(10,0){2}{\circle*{1}}
\multiput(10,8)(10,0){2}{\circle*{1}}
\put(6,12){\makebox(0,0){$d$}}
\put(24,-6){\makebox(0,0){$5$}}
\end{picture}
\left( \begin{array}{cc}
 u{c\cho{5}}_{2,r}  &  \\
 &  u{e\cho{5}}_{4,r}    
  \end{array}  \right)  \\
&=&
\left( \begin{array}{ccc}
  -1 & &  \\
&  u{d\cho{e}}_{3,l} &  \\
 &  & u{d\cho{G}}_{2,l}
  \end{array}  \right)  \\  
&&\left( \matrix{
0 & \gb-4 & \frac{\gb-2}{\sqrt{\gb(2\gb-1)(\gbf+4)}} & \frac{3\gb+4}
{\sqrt{\gb(2\gb-1)(\gbf+4)}} & \frac{-\gb}{2\gb-1} & 0 \cr
 * & \frac{-(4\gb-5)}{2\gb \sqrt{\gb-1}} & \frac{3\sqrt{2}}{(\gb-1)
\sqrt{\gbf+4}} & \frac{-\sqrt{2}(3\gb+1)}{(\gb-1)^2 \sqrt{\gbf+4}} & 
\frac{\sqrt{2\gb(\gb-2)}}{(\gb-1)^2} & * \cr
* & \frac{-(\gb-2)}{2\gbf} & \frac{\sqrt{(\gb-2)^3}}{\beta^3 \sqrt
{\gbf+4}} & \frac{4 \sqrt{\gb-2}}{\sqrt{\gb(\gbf+4)}} &
\frac{\gb-2}{\gb-1} & * \cr
* & \frac{-1}{\sqrt{(\gb-1)^3}} & \frac{\sqrt{2}(\gb-2)^2}{(\gb-1)
\sqrt{\gbf+4}} & \frac{\sqrt{2}(\gb-2)}{(\gb-1)\sqrt{\gbf+4}} & 
\frac{-1}{\sqrt{(\gb-1)^3}} & * \cr
\frac{\sqrt{2(\gbf-1)}}{\beta^3} & * & * & * & 0 & \frac{2}{\gb-1} \cr
\frac{1}{\gbf}  & * & * & * & 0 & \frac{\gb+1}{\gb\sqrt{\gb-1}} 
} \right),
\end{eqnarray*}

\begin{eqnarray*}
&& \left( \begin{array}{ccc}
  -1 & &  \\
&   u{d\cho{E}}_{3,l}  &  \\
 &  &  u{d\cho{G}}_{2,l}
  \end{array}  \right)    \\
&=&
\left( \matrix{
-1 & 0 & 0 & 0 & 0 & 0 \cr
0 & \frac{\gb-2}{\gb+1} & \frac{\gbs(\gb-2)}{\sqrt{2(\gbf-1)}} & 
\frac{-1}{\gb+1} & 0 & 0 \cr
0 & \frac{-\sqrt{\gb-2}}{2\sqrt{\gb+1}} & \sqrt{\frac{\gb-2}{\gbf-1}} 
&
 \frac{\gbs(\gb-2)}{\sqrt{2(\gbf-1)}} & 0 & 0 \cr
0 & \frac{-\gb}{\gb+1} & \frac{\sqrt{\gb-2}}{2\sqrt{\gb+1}} & 
\frac{-(\gb-2)}{\gb+1} & 0 & 0 \cr
0 & 0 & 0 & 0 & \frac{-(\gb-2)}{\gb} & \frac{2\sqrt{\gb-1}}{\gb} \cr
0 & 0 & 0 & 0 & \frac{2\sqrt{\gb-1}}{\gb} & \frac{\gb-2}{\gb}  \cr
} \right).
\end{eqnarray*}

Here, to obtain the above matrix for the gauges $u{d\cho{E}}$ and
$u{d\cho{G}}$, we use {\it Mathematica} to see the signs of 
non-zero entries.  First we calculate the $(6,5)$, $(6,6)$-entries
and $(1,1)$-entry (equal to $-1$), and also check 
$ \{(6,5) $-entry$ \}^2+\{(6,6) $-entry$ \}^2=1 $, then we know
 the entries of first and last rows and of first column are equal to 0
except for 
(1,1),(6,5) and (6,6). Next we calculate $(2,2),\dots,(2,4)$ and 
$(4,2),
\dots, (4,4)$ and check that square sums are equal to 1 respectively, 
and also
calculate (5,2)-entry is equal to 0.  Then (3,2)-entry is determined 
by using 
the fact that the square sum of the entries in the second column is
equal to 1, and by sign of the entry obtained by the numerical
calculation of the product of matrices
 by {\it Mathematica}.  Then we have (3,3) and (3,4)-entries by the
orthogonal relation of the
second column and the third, and forth. Now we know that the matrix is
block diagonal, and we have the rest (5,5) and (5,6)-entries by
unitarity and signs obtained by {\it Mathematica}. To execute the above
calculation, we do not use the 
entries $*$ in the multiplication matrix of $M(d/5)$ and the right
gauges.  
Note $u{d\cho{\tilde{C}}}=u{\tilde{d}\cho{\tilde{C}}}=-1$. \par
We see that $u{d\cho{G}}$ is the same matrix as we have already 
gotten.  
From now, 
we consider that we are always checking it when the matrices for the
gauges whose matrices have been already obtained appear with the 
matrices 
of new gauges.

We check for the rest connections as before.

\begin{displaymath}
\thinlines
\unitlength 0.5mm
\begin{picture}(30,20)
\multiput(11,-2)(0,10){2}{\line(1,0){8}}
\multiput(10,7)(10,0){2}{\line(0,-1){8}}
\multiput(10,-2)(10,0){2}{\circle*{1}}
\multiput(10,8)(10,0){2}{\circle*{1}}
\put(6,12){\makebox(0,0){$A$}}
\put(24,-6){\makebox(0,0){$g$}}
\put(6,-6){\makebox(0,0){$\sim$}}
\end{picture}
= \quad
\thinlines
\unitlength 0.5mm
\begin{picture}(30,20)
\multiput(11,-2)(0,10){2}{\line(1,0){8}}
\multiput(10,7)(10,0){2}{\line(0,-1){8}}
\multiput(10,-2)(10,0){2}{\circle*{1}}
\multiput(10,8)(10,0){2}{\circle*{1}}
\put(6,12){\makebox(0,0){$A$}}
\put(24,-6){\makebox(0,0){$g$}}
\end{picture}
\left( \begin{array}{cc}
  1  &  \\
 &  u{g\cho{2}}_{2,r}  
  \end{array}  \right),
\end{displaymath}

\begin{displaymath}
\left( \begin{array}{cc}
  1  &  \\
 &   u{g\cho{2}}_{2,r}   
  \end{array}  \right)
=
\left( \matrix{
1 &  0 & 0 \cr
0 & \frac{-1}{\gb-1} & \frac{\sqrt{\gb(\gb-2)}}{\gb-1} \cr
0 & \frac{\sqrt{\gb(\gb-2)}}{\gb-1} & \frac{1}{\gb-1}  \cr
} \right),
\end{displaymath}

\begin{eqnarray*}
\thinlines
\unitlength 0.5mm
\begin{picture}(30,20)
\multiput(11,-2)(0,10){2}{\line(1,0){8}}
\multiput(10,7)(10,0){2}{\line(0,-1){8}}
\multiput(10,-2)(10,0){2}{\circle*{1}}
\multiput(10,8)(10,0){2}{\circle*{1}}
\put(6,12){\makebox(0,0){$f$}}
\put(24,-6){\makebox(0,0){$2$}}
\end{picture}
&=&
\left( \begin{array}{cc}
  1  &  \\
 &   u{f\cho{C}}_{2,l}  
  \end{array}  \right)
\thinlines
\unitlength 0.5mm
\begin{picture}(30,20)
\multiput(11,-2)(0,10){2}{\line(1,0){8}}
\multiput(10,7)(10,0){2}{\line(0,-1){8}}
\multiput(10,-2)(10,0){2}{\circle*{1}}
\multiput(10,8)(10,0){2}{\circle*{1}}
\put(6,12){\makebox(0,0){$f$}}
\put(24,-6){\makebox(0,0){$2$}}
\end{picture}
\left( \begin{array}{cc}
  1  &  \\
 &   u{g\cho{2}}_{2,r}   
  \end{array}  \right)  \\
&=&
\left( \begin{array}{cc}
  1  &  \\
 &   u{f\cho{C}}_{2,l}   
  \end{array}  \right)
\left( \matrix{
\frac{\sqrt{2}}{\sqrt{\gb-2}} & \frac{-1}{\gb-1} & \frac{\sqrt{\gb-2}}
{2\gbs} \cr
\frac{\gb+1}{2\gb} & \frac{-\sqrt{\gb-2}}{2\gbs} & 
\frac{\gb}{2(\gb-2)}
\cr 
* & * & * \cr } \right),
\end{eqnarray*}

\begin{displaymath}
\left( \begin{array}{cc}
  1  &  \\
 &   u{f\cho{C}}_{2,l}  
  \end{array}  \right)
=
\left( \matrix{
1 & 0 & 0 \cr
0 & \frac{1}{\gb-2} & \frac{\sqrt{\gb+1}}{\gb-2} \cr
0 & \frac{\sqrt{\gb+1}}{\gb-2} & \frac{-1}{\gb-2}  \cr } \right),
\end{displaymath}

\begin{eqnarray*}
\thinlines
\unitlength 0.5mm
\begin{picture}(30,20)
\multiput(11,-2)(0,10){2}{\line(1,0){8}}
\multiput(10,7)(10,0){2}{\line(0,-1){8}}
\multiput(10,-2)(10,0){2}{\circle*{1}}
\multiput(10,8)(10,0){2}{\circle*{1}}
\put(6,12){\makebox(0,0){$C$}}
\put(24,-6){\makebox(0,0){$e$}}
\put(6,-6){\makebox(0,0){$\sim$}}
\end{picture}
&=&
\left( \begin{array}{ccc}
1 & & \\
 & 1 & \\
 & & u{f\cho{C}}_{2,l}  
\end{array}  \right) 
\thinlines
\unitlength 0.5mm
\begin{picture}(30,20)
\multiput(11,-2)(0,10){2}{\line(1,0){8}}
\multiput(10,7)(10,0){2}{\line(0,-1){8}}
\multiput(10,-2)(10,0){2}{\circle*{1}}
\multiput(10,8)(10,0){2}{\circle*{1}}
\put(6,12){\makebox(0,0){$e$}}
\put(24,-6){\makebox(0,0){$c$}}
\end{picture}
\left( \begin{array}{cc}
   u{e\cho{3}}_{3,r}     & \\
 & 1  
  \end{array}  \right)
   \\
&=&
\left( \matrix{
\frac{-(\gb-2)}{(\gb-1)\sqrt{2(\gbf+4)}}   
&\sqrt{\frac{\gb(\gb+2)}{(\gb+1)(\gbf+4)}} & 
\frac{-\gbs\sqrt{\gb-1}}{\sqrt{\gbf+4}} & * \cr
-\frac{2(\gb-2)}{\sqrt{2(\gbf+4)}}  & \frac{4}{\sqrt{2(\gbf+4)}} & 
\frac{2\sqrt{\gb-1}}{\sqrt{\gb(\gbf+4)}} & * \cr
\frac{-1}{2\gbs\sqrt{\gb-1}} & \frac{-1}{2\gbs\sqrt{\gb-1}} &
\frac{-2\sqrt{2}}{(\gb-2)^2} & \frac{3\gb-4}{2\gbs(\gb-2)} \cr
* & * & \frac{1}{(\gb-2)\gbs} & \frac{1}{\sqrt{2}(\gb-2)} \cr
} \right)
\left( \begin{array}{cc}
  u{e\cho{3}}_{3,r}    & \\
 & 1  
  \end{array}  \right),
\end{eqnarray*}

\begin{displaymath}
\left( \begin{array}{cc}
 u{e\cho{3}}_{3,r}       & \\
 & 1  
  \end{array}  \right)
=
\left( \matrix{
\frac{2(\gb-2)}{\gbf+4} & \frac{-2(\gb-2)^2}{\gbf+4} & 
\frac{-\sqrt{(\gb-2)^3}}{\sqrt{\gbf+4}} & 0 \cr
\frac{-2(\gb-2)^2}{\gbf+4} & \frac{-2(7\gb+2)}{(\gb-1)(\gbf+4)} & 
\frac{2\sqrt{(\gb-2)^3}}{\gb\sqrt{\gbf+4}} & 0 \cr
\frac{-\sqrt{(\gb-2)^3}}{\sqrt{\gbf+4}} & \frac{2\sqrt{(\gb-2)^3}}
{\gb\sqrt{\gbf+4}} & \frac{-2}{\gb} & 0 \cr
0 & 0 & 0 & 1 \cr
} \right),
\end{displaymath}

\begin{eqnarray*}
\thinlines
\unitlength 0.5mm
\begin{picture}(30,20)
\multiput(11,-2)(0,10){2}{\line(1,0){8}}
\multiput(10,7)(10,0){2}{\line(0,-1){8}}
\multiput(10,-2)(10,0){2}{\circle*{1}}
\multiput(10,8)(10,0){2}{\circle*{1}}
\put(6,12){\makebox(0,0){$\tilde{d}$}}
\put(24,-6){\makebox(0,0){$3$}}
\end{picture}
&=&
\left( \begin{array}{cc}
  u{\tilde{d}\cho{E}}_{3,l}    & \\
 & 1  
  \end{array}  \right)
\thinlines
\unitlength 0.5mm
\begin{picture}(30,20)
\multiput(11,-2)(0,10){2}{\line(1,0){8}}
\multiput(10,7)(10,0){2}{\line(0,-1){8}}
\multiput(10,-2)(10,0){2}{\circle*{1}}
\multiput(10,8)(10,0){2}{\circle*{1}}
\put(6,12){\makebox(0,0){$d$}}
\put(24,-6){\makebox(0,0){$3$}}
\end{picture}
\left( \begin{array}{cc}
  u{e\cho{3}}_{3,r}    & \\
 & 1  
  \end{array}  \right)  \\
&=&
\left( \begin{array}{cc}
  u{\tilde{d}\cho{E}}_{3,l}    & \\
 & 1  
  \end{array}  \right)
\left( \matrix{
\frac{2(\gb-2)}{\sqrt{2(\gbf+4)}(\gb-1)} & \frac{-8\sqrt{2}\gb}
{\sqrt{\gbf+4}(\gbf-1)} & * & \frac{\sqrt{\gb-2}}{\gb\sqrt{\gb+1}} \cr
\frac{-\sqrt{(\gb-2)^3}}{\sqrt{\gb(\gbf+4)}} & 
\frac{-4\sqrt{2(\gb-1)}}
{\gb\sqrt{\gbf+4}} & * & \frac{\gb-2}{\gb-1} \cr
* & * & *  & * \cr } \right),
\end{eqnarray*}

\begin{displaymath}
\left( \begin{array}{cc}
   u{\tilde{d}\cho{E}}_{3,l}     & \\
  & 1  
  \end{array}  \right)
=
\left( \matrix{
\frac{\gb-2}{\gb+1} & \frac{-(\gb-2)}{4\sqrt{\gb-1}} & \frac{-\gb}{\gb+1} &
 0 \cr
\frac{\gb(\gb-2)}{2\sqrt{(\gb-1)^3}} & \frac{\gb}{(\gb-1)^2} & 
\frac{\gb-2}{4\sqrt{\gb-1}} &  0 \cr
\frac{-1}{\gb+1} & \frac{\gb(\gb-2)}{2\sqrt{(\gb-1)^3}} &
\frac{-(\gb-2)}{\gb+1} & 0 \cr
0 & 0 & 0 & 1 \cr} \right),
\end{displaymath}

\begin{displaymath}
\thinlines
\unitlength 0.5mm
\begin{picture}(30,20)
\multiput(11,-2)(0,10){2}{\line(1,0){8}}
\multiput(10,7)(10,0){2}{\line(0,-1){8}}
\multiput(10,-2)(10,0){2}{\circle*{1}}
\multiput(10,8)(10,0){2}{\circle*{1}}
\put(6,12){\makebox(0,0){$\tilde{h}$}}
\put(24,-6){\makebox(0,0){$3$}}
\end{picture}
=
\thinlines
\unitlength 0.5mm
\begin{picture}(30,20)
\multiput(11,-2)(0,10){2}{\line(1,0){8}}
\multiput(10,7)(10,0){2}{\line(0,-1){8}}
\multiput(10,-2)(10,0){2}{\circle*{1}}
\multiput(10,8)(10,0){2}{\circle*{1}}
\put(6,12){\makebox(0,0){$h$}}
\put(24,-6){\makebox(0,0){$3$}}
\end{picture}
u{g\cho{3}}_{2,r},
\end{displaymath}

\begin{displaymath}
 u{g\cho{3}}_{2,r}
=
\left( \matrix{
\frac{-1}{2\gb-1} & \frac{2\sqrt{2}}{\sqrt{2\gb-1}} \cr
\frac{2\sqrt{2}}{\sqrt{2\gb-1}} & \frac{1}{2\gb-1} 
} \right),
\end{displaymath}

\begin{eqnarray*}
&& \thinlines
\unitlength 0.5mm
\begin{picture}(30,20)
\multiput(11,-2)(0,10){2}{\line(1,0){8}}
\multiput(10,7)(10,0){2}{\line(0,-1){8}}
\multiput(10,-2)(10,0){2}{\circle*{1}}
\multiput(10,8)(10,0){2}{\circle*{1}}
\put(6,12){\makebox(0,0){$f$}}
\put(24,-6){\makebox(0,0){$3$}}
\end{picture}
=
\left( \begin{array}{cc}
  u{f\cho{C}}_{2,l}    & \\
 &  u{f\cho{E}}_{3,l}    
  \end{array}  \right)
\thinlines
\unitlength 0.5mm
\begin{picture}(30,20)
\multiput(11,-2)(0,10){2}{\line(1,0){8}}
\multiput(10,7)(10,0){2}{\line(0,-1){8}}
\multiput(10,-2)(10,0){2}{\circle*{1}}
\multiput(10,8)(10,0){2}{\circle*{1}}
\put(6,12){\makebox(0,0){$f$}}
\put(24,-6){\makebox(0,0){$3$}}
\end{picture}
\left( \begin{array}{cc}
 u{g\cho{3}}_{2,r}      & \\
 &   u{e\cho{3}}_{3,r}   
  \end{array}  \right)  \\
&=&
\left( \begin{array}{cc}
   u{f\cho{C}}_{2,l}    & \\
 &    u{f\cho{E}}_{3,l}     
  \end{array}  \right)
\left( \matrix{
\frac{1}{\gbf} & * & * & * & \frac{-2\sqrt{\gb+1}}{\gb(\gb-2)} \cr
\frac{-(\gb-1)\sqrt{\gb+1}}{\gbf} & * & * & * & \frac{4\gb-1}{\gbf} 
\cr
* & \frac{1}{\sqrt{2}\gbs(\gb-1)} & \frac{1}{\sqrt{\gb+3}} & 
\frac{-2\sqrt{2}(\gb-1)}{\gb\sqrt{\gbf+4}} & * \cr
* & \frac{-1}{\gb(\gb-1)} & \frac{-\sqrt{2}(\gb-2)}{2\gb\sqrt{\gbf+4}}
& \frac{-(\gb+4)}{\sqrt{\gb(\gbf+4)}} & * \cr
* & * & * & * & * \cr } \right),
\end{eqnarray*}

\begin{eqnarray*}
\left( \begin{array}{cc}
   u{f\cho{C}}_{2,l}    & \\
 &  u{f\cho{E}}_{3,l} 
  \end{array}  \right)  
=
\left( \matrix{
\frac{1}{\gb-2} & \frac{\sqrt{\gb+1}}{\gb-2} & 0 & 0 & 0 \cr
\frac{\sqrt{\gb+1}}{\gb-2} & \frac{-1}{\gb-2}  & 0 & 0 & 0 \cr
0 & 0 & \frac{-1}{\gb+1} & \frac{-\sqrt{2}\gbs}{\gb+1} & 
\frac{\gb}{\gb+1} \cr
0 & 0 & \frac{-\sqrt{2}\gbs}{\gb+1} & \frac{-\gb}{2(\gb-1)} & 
\frac{-\sqrt{2}\gbs}{\gb+1} \cr
0 & 0 & \frac{\gb}{\gb+1} & \frac{-\sqrt{2}\gbs}{\gb+1} & \frac{-1}
{\gb+1} \cr } \right),
\end{eqnarray*}

\begin{eqnarray*}
\thinlines
\unitlength 0.5mm
\begin{picture}(30,20)
\multiput(11,-2)(0,10){2}{\line(1,0){8}}
\multiput(10,7)(10,0){2}{\line(0,-1){8}}
\multiput(10,-2)(10,0){2}{\circle*{1}}
\multiput(10,8)(10,0){2}{\circle*{1}}
\put(6,12){\makebox(0,0){$f$}}
\put(24,-6){\makebox(0,0){$4$}}
\end{picture}
&=&
 u{f\cho{E}}_{3,l} 
\thinlines
\unitlength 0.5mm
\begin{picture}(30,20)
\multiput(11,-2)(0,10){2}{\line(1,0){8}}
\multiput(10,7)(10,0){2}{\line(0,-1){8}}
\multiput(10,-2)(10,0){2}{\circle*{1}}
\multiput(10,8)(10,0){2}{\circle*{1}}
\put(6,12){\makebox(0,0){$f$}}
\put(24,-6){\makebox(0,0){$4$}}
\end{picture}
\left( \begin{array}{cc}
-1 &  \\
& u{e\cho{4}}_{2,r}
\end{array} \right)   \\
&=& 
\left( \matrix{
\frac{\gbs}{2\sqrt{2}} & \frac{-(\gb-2)}{\sqrt{2(\gbf+4)}} & * \cr
\frac{-1}{2} & \frac{2\gbs}{\sqrt{\gbf+4}} & * \cr
\frac{1}{\sqrt{\gb+1}} & \frac{\gb-2}{\sqrt{2(\gbf+4)}} & * \cr
} \right)
\left( \begin{array}{cc}
-1 &  \\
& u{e\cho{4}}_{2,r} 
\end{array} \right), 
\end{eqnarray*}

\begin{displaymath}
\left( \begin{array}{cc}
-1 &  \\
& u{e\cho{4}}_{2,r} 
\end{array} \right)
=
\left( \matrix{
-1 & 0 & 0 \cr
0 & \frac{-4\gb}{\gbf+4} & \frac{-(\gbf-4)}{\gbf+4} \cr
0 & \frac{-(\gbf-4)}{\gbf+4} & \frac{4\gb}{\gbf+4} \cr } \right),
\end{displaymath}

here note $u(g\mbox{-}4)=-1$

\begin{displaymath}
\thinlines
\unitlength 0.5mm
\begin{picture}(30,20)
\multiput(11,-2)(0,10){2}{\line(1,0){8}}
\multiput(10,7)(10,0){2}{\line(0,-1){8}}
\multiput(10,-2)(10,0){2}{\circle*{1}}
\multiput(10,8)(10,0){2}{\circle*{1}}
\put(6,12){\makebox(0,0){$\tilde{h}$}}
\put(24,-6){\makebox(0,0){$5$}}
\end{picture}
=
\thinlines
\unitlength 0.5mm
\begin{picture}(30,20)
\multiput(11,-2)(0,10){2}{\line(1,0){8}}
\multiput(10,7)(10,0){2}{\line(0,-1){8}}
\multiput(10,-2)(10,0){2}{\circle*{1}}
\multiput(10,8)(10,0){2}{\circle*{1}}
\put(6,12){\makebox(0,0){$h$}}
\put(24,-6){\makebox(0,0){$5$}}
\end{picture}
u{g\cho{5}}_{2,r},
\end{displaymath}

\begin{displaymath}
u{g\cho{5}}_{2,r}
=
\left( \matrix{
\frac{-1}{\gb-2} & \frac{\sqrt{\gb+1}}{\gb-2} \cr
\frac{\sqrt{\gb+1}}{\gb-2} & \frac{1}{\gb-2}   \cr } \right),
\end{displaymath}

\begin{displaymath}
\thinlines
\unitlength 0.5mm
\begin{picture}(30,20)
\multiput(11,-2)(0,10){2}{\line(1,0){8}}
\multiput(10,7)(10,0){2}{\line(0,-1){8}}
\multiput(10,-2)(10,0){2}{\circle*{1}}
\multiput(10,8)(10,0){2}{\circle*{1}}
\put(6,12){\makebox(0,0){$f$}}
\put(24,-6){\makebox(0,0){$5$}}
\end{picture}
=
\left( \begin{array}{ccc}
u{f\cho{E}}_{3,l} & & \\
& u{f\cho{\tilde{C}}}_{2,l} & \\
& & 1 
\end{array} \right)
\thinlines
\unitlength 0.5mm
\begin{picture}(30,20)
\multiput(11,-2)(0,10){2}{\line(1,0){8}}
\multiput(10,7)(10,0){2}{\line(0,-1){8}}
\multiput(10,-2)(10,0){2}{\circle*{1}}
\multiput(10,8)(10,0){2}{\circle*{1}}
\put(6,12){\makebox(0,0){$f$}}
\put(24,-6){\makebox(0,0){$5$}}
\end{picture}
\left( \begin{array}{cc}
u{g\cho{5}}_{2,r}   & \\
& u{e\cho{5}}_{4,r} 
\end{array} \right),
\end{displaymath}

\begin{eqnarray*}
&&
\thinlines
\unitlength 0.5mm
\begin{picture}(30,20)
\multiput(11,-2)(0,10){2}{\line(1,0){8}}
\multiput(10,7)(10,0){2}{\line(0,-1){8}}
\multiput(10,-2)(10,0){2}{\circle*{1}}
\multiput(10,8)(10,0){2}{\circle*{1}}
\put(6,12){\makebox(0,0){$f$}}
\put(24,-6){\makebox(0,0){$5$}}
\end{picture}
\left( \begin{array}{cc}
 u{g\cho{5}}_{2,r}   & \\
& u{e\cho{5}}_{4,r} 
\end{array} \right)   \\
&=&
\left( \matrix{
* & * & * & * & * & *  \cr
* & \frac{3}{2(\gb-2)} & \frac{3}{\gbs\sqrt{\gbf+4}} & \frac{-\gbs}
{(\gb-2)\sqrt{\gbf+4}} & \frac{-(\gb-2)}{2\sqrt{2\gb-1}} & * \cr
* & \frac{-1}{2\sqrt{2}\gbs} & \frac{\sqrt{2}}{(\gb-1)\sqrt{\gbf+4}} &
\frac{-(\gb-1)}{\sqrt{2(\gbf+4)}} & \frac{\gbs}{\sqrt{(\gb-1)^3}} & * 
\cr
\frac{2\sqrt{\gb-1}}{(\gb-2)^2} & * & * & * & * & \frac{\gb+1}
{\sqrt{2}\gbs(\gb-2)} \cr
\frac{-1}{\sqrt{(\gb-2)^3}} & * & * & * & * & \frac{-\sqrt{\gb+1}}
{\sqrt{2}\gbs(\gb-2)} \cr
\frac{-2}{(\gb-2)^2} & \frac{-\sqrt{\gb+1}}{2(\gb-2)} & \frac{\gb-1}
{\sqrt{2(\gbf+4)}} & \frac{\sqrt{\gb-1}(\gb+1)}
{\gbs\sqrt{(\gbf+4)(\gb-2)}} & * & \frac{1}{\sqrt{(\gb-2)^3}} \cr
} \right),
\end{eqnarray*}

\begin{eqnarray*}
\left(  \begin{array}{ccc}
u{f\cho{E}}_{3,l} & & \\
& u{f\cho{\tilde{C}}}_{2,l} & \\
& & 1 
\end{array} \right)   
=
\left( \matrix{
\frac{-1}{\gb+1} & \frac{-(\gb-2)^2}{4\sqrt{2}\gbs} & 
\frac{\gb}{\gb+1}
& 0 & 0 & 0 \cr
\frac{-(\gb-2)^2}{4\sqrt{2}\gbs} & \frac{-(\gb-1)}{\gb+1} & 
\frac{-(\gb-2)}{2\sqrt{\gb+1}} & 0 & 0 & 0 \cr
\frac{\gb}{\gb+1} & \frac{-(\gb-2)}{2\sqrt{\gb+1}} & \frac{-1}{\gb+1} 
& 
0 & 0 & 0 \cr
0 & 0 & 0 & \frac{1}{\gb-2} & \frac{\gb+1}{2\sqrt{2}\gbs} & 0 \cr
0 & 0 & 0 & \frac{\gb+1}{2\sqrt{2}\gbs} & \frac{-1}{\gb-2} & 0 \cr
0 & 0 & 0 & 0 & 0 & 1 \cr
} \right),
\end{eqnarray*}

\begin{displaymath}
\thinlines
\unitlength 0.5mm
\begin{picture}(30,20)
\multiput(11,-2)(0,10){2}{\line(1,0){8}}
\multiput(10,7)(10,0){2}{\line(0,-1){8}}
\multiput(10,-2)(10,0){2}{\circle*{1}}
\multiput(10,8)(10,0){2}{\circle*{1}}
\put(6,12){\makebox(0,0){$b$}}
\put(24,12){\makebox(0,0){$c_e$}}
\put(6,-6){\makebox(0,0){$E^f$}}
\put(24,-6){\makebox(0,0){$3$}}
\end{picture}
=
\thinlines
\unitlength 0.5mm
\begin{picture}(30,20)
\multiput(11,-2)(0,10){2}{\line(1,0){8}}
\multiput(10,7)(10,0){2}{\line(0,-1){8}}
\multiput(10,-2)(10,0){2}{\circle*{1}}
\multiput(10,8)(10,0){2}{\circle*{1}}
\put(6,12){\makebox(0,0){$\tilde{b}$}}
\put(24,12){\makebox(0,0){$\tilde{c}_e$}}
\put(6,-6){\makebox(0,0){$E^{\tilde{d}}$}}
\put(24,-6){\makebox(0,0){$3$}}
\end{picture}
\times (-1),
\end{displaymath}

$$u(c\mbox{-}4)=u(\tilde{c}\mbox{-}4)=-1.$$

From the above computations, we can extract a list of 32 unitary 
matrices $u{x\cho{y}}_{m,l}$ labeled by the edges $x$-$y$ of ${\cal 
H}_{0}$ and 23 unitary matrices $u{z\cho{w}}_{n,r}$ labeled by the 
edges $z$-$w$ of ${\cal H}_{1}$. We now have to check the equality 

$$
\left( \thinlines
\unitlength 0.5mm
\begin{picture}(30,20)
\multiput(11,-2)(0,10){2}{\line(1,0){8}}
\multiput(10,7)(10,0){2}{\line(0,-1){8}}
\multiput(10,-2)(10,0){2}{\circle*{1}}
\multiput(10,8)(10,0){2}{\circle*{1}}
\put(6,12){\makebox(0,0){$\td{x}$}}
\put(24,12){\makebox(0,0){$\td{z}$}}
\put(6,-6){\makebox(0,0){$y$}}
\put(24,-6){\makebox(0,0){$w$}}
\end{picture} \right)_{xy,zw}
=u{x\cho{y}}_{n,l}
\left( \thinlines
\unitlength 0.5mm
\begin{picture}(30,20)
\multiput(11,-2)(0,10){2}{\line(1,0){8}}
\multiput(10,7)(10,0){2}{\line(0,-1){8}}
\multiput(10,-2)(10,0){2}{\circle*{1}}
\multiput(10,8)(10,0){2}{\circle*{1}}
\put(6,12){\makebox(0,0){$x$}}
\put(24,12){\makebox(0,0){$z$}}
\put(6,-6){\makebox(0,0){$y$}}
\put(24,-6){\makebox(0,0){$w$}}
\end{picture} \right)_{xy,zw}
u{z\cho{w}}_{m,r}
$$
for {\it all} pairs of edges $(xy, 
zw)$ from the vertical graphs. This amounts to check 280 equalities of 
real numbers. First we go back and check that all the matrix 
identities listed above for the construction of the gauge matrices 
holds for {\it all} entries and not just for the subset of entries 
needed to produce the candidates for $u{x\cho{y}}_{n,l}$'s and 
$u{z\cho{w}}_{m,r}$'s. Then, we see that a few equalities are left to 
be checked, namely, the equality of gauge transform between $M(d/4)$ and 
$M(\td{d}/4)$, and that of 
some other matrices. The checking of the former is done as follows.


\begin{eqnarray*}
\thinlines
\unitlength 0.5mm
\begin{picture}(30,20)
\multiput(11,-2)(0,10){2}{\line(1,0){8}}
\multiput(10,7)(10,0){2}{\line(0,-1){8}}
\multiput(10,-2)(10,0){2}{\circle*{1}}
\multiput(10,8)(10,0){2}{\circle*{1}}
\put(6,12){\makebox(0,0){$\tilde{d}$}}
\put(24,-6){\makebox(0,0){4}}
\end{picture}
&=&
 u{\tilde{d}\cho{E}}_{3,l}
\thinlines
\unitlength 0.5mm
\begin{picture}(30,20)
\multiput(11,-2)(0,10){2}{\line(1,0){8}}
\multiput(10,7)(10,0){2}{\line(0,-1){8}}
\multiput(10,-2)(10,0){2}{\circle*{1}}
\multiput(10,8)(10,0){2}{\circle*{1}}
\put(6,12){\makebox(0,0){$d$}}
\put(24,-6){\makebox(0,0){$4$}}
\end{picture}
\left( \begin{array}{cc}
 -1 &  \\
 & u{e\cho{4}}_{2,r}
\end{array}
\right)  \\ 
&=&
\left( \matrix{
\frac{-(\gb-2)}{\sqrt{\gb+1}} & * & 
\frac{\gbs(\gb-2)}{\sqrt{(\gb+1)(\gbf+4)}}
  \cr
\frac{\sqrt{5-\gb}}{\sqrt{\gb+1}} & * & \frac{-2\beta^3}
{\sqrt{(\gbf-1)(\gbf+4)}} \cr
\frac{-(\gb-2)}{\sqrt{\gb+1}} & * & 
\frac{2\sqrt{2}}{(\gb+1)\sqrt{\gbf+4}}
 \cr
}  \right)
\left( \begin{array}{cc}
 -1 &  \\
 & u{e\cho{4}}_{2,r}
\end{array}
\right),
\end{eqnarray*}

\begin{displaymath}
\left( \begin{array}{cc}
 -1 &  \\
 & u{e\cho{4}}_{2,r}
\end{array}
\right)
=
\left(
\matrix{
-1 & 0 & 0 \cr
0 & \frac{-4\gb}{\gbf+4} & \frac{-(\gbf-4)}{\gbf+4} \cr
0 & \frac{-(\gbf-4)}{\gbf+4} & \frac{4\gb}{\gbf+4} \cr } \right).
\end{displaymath}
We can see that $u{e\cho{4}}_{2,r}$ here is the same matrix as 
when it appeared first.
The latter equalities are of scalar matrices or $2 \times 2$
matrices which we can check at a glance that the gauge matrix is 
common with 
which we already used, in either case it is easy enough not to write
down. \par

All the above identities, we have checked using {\it Mathematica} and 
of course we have made repeatably use of the identity $\gbf -5 \gb + 
2 =0 $.
At last, we have obtained the equivalence of the connections
$$ (\alpha \tilde{\alpha}-\bf{1})\sigma \alpha \cong
 \sigma(\alpha \tilde{\alpha}-\bf{1})\sigma \alpha $$
 up to the vertical gauge choice. \\
 
Finally we will check conditions 1) and 2).
Along the same argument of the proof of the previous theorem for 
$(5+\sqrt{13})/2$ case,
we see the indecomposability {\it other than}
 for $(\a \ab-{\bf 1})\sigma \a$
and mutually inequivalence of all.
 In Figure \ref{alle17graph.tex},
$*$ in $V_{0}$ is vertex of only one edge in ${\cal K}$. Thus, using 
Cororally \ref{irred}, we have 
indecomposability of the connection $(\a \ab-{\bf 1})\sigma \a$.
 Now, the proposition holds and thus
 we have proved the theorem. \\
 \qed

\end{document}